\documentclass[letter,12pt,openany]{book}

\usepackage[utf8]{inputenc}
\usepackage[T1]{fontenc}
\usepackage[english]{babel}
\usepackage{times}
\usepackage{amsfonts}
\usepackage{fourier}
\usepackage{amsmath,amsfonts,amsthm}
\usepackage{stmaryrd}
\usepackage{graphicx}
\usepackage{enumerate}
\usepackage{amsthm}
\usepackage{amssymb}
\usepackage{array}
\usepackage{epstopdf}
\usepackage{ragged2e}
\usepackage{sectsty}
\usepackage{float} 			
\usepackage{newclude}
\usepackage{relsize}		
\usepackage{hyperref}

\usepackage[all]{xy}
\usepackage{ar}

\usepackage[nottoc]{tocbibind}

\usepackage{tikz}
\usepackage{verbatim}
\usetikzlibrary{shapes.geometric}
\usetikzlibrary{shapes.arrows}
\usepackage{array}

\let\oldsection\section
\renewcommand{\section}{
  \renewcommand{\theequation}{\thesection.\arabic{equation}}
  \oldsection}
\let\oldsubsection\subsection
\renewcommand{\subsection}{
  \renewcommand{\theequation}{\thesubsection.\arabic{equation}}
  \oldsubsection}

\newtheorem{thm}{Theorem}[section]
\newtheorem*{thm*}{Theorem}
\newtheorem{lem}[thm]{Lemma}
\newtheorem{cor}[thm]{Corolary}
\newtheorem{prop}[thm]{Proposition}
\newtheorem{defn}[thm]{Definition}
\newtheorem{obs}[thm]{Observation}
\newtheorem{ejem}[thm]{Example}

\newcommand{\nin}{\not\in}
\newcommand{\barra}{\,\middle|\,}
\newcommand{\set}[1]{\left\{#1\right\}}
\newcommand{\SET}[2]{\left\{#1\,\middle|\,#2\right\}}
\newcommand{\parentesis}[1]{\left(#1\right)}
\newcommand{\corchetes}[1]{\left[#1\right]}
\newcommand{\valorabs}[1]{\left|#1\right|}

\newcommand{\prodint}[1]{\left\langle #1 \right\rangle}

\newcommand{\dual}[1]{\parentesis{ #1 }^{\ast}}
\newcommand{\floor}[1]{\lfloor #1 \rfloor}
 
\newcommand{\R}{\mathbb{R}}
\newcommand{\C}{\mathbb{C}}
\newcommand{\Cc}{\mathbb{C}\setminus\set{0}}

\newcommand{\hC}{\hat{\mathbb{C}}}

\newcommand{\Z}{\mathbb{Z}}
\newcommand{\N}{\mathbb{N}}
\newcommand{\Q}{\mathbb{Q}}

\newcommand{\Ss}{\mathbb{S}}
\newcommand{\D}{\mathbb{D}}

\newcommand{\Hh}{\mathbb{H}}
\newcommand{\id}{\mathit{id}}

\newcommand{\CP}{\mathbb{CP}}

\newcommand{\HhC}{\Hh_{\C}}

\newcommand{\dih}{\text{Dih}_\infty}
\newcommand{\rot}{\text{Rot}_\infty}
\newcommand{\so}{\text{SO}(3)}
\newcommand{\epa}{\text{Epa}\parentesis{\C}}
\newcommand{\autc}{\text{Aut}\parentesis{\C^\ast}}
\newcommand{\pslr}{\text{PSL}\parentesis{2,\R}}
\newcommand{\PSL}{\text{PSL}\parentesis{3,\C}}
\newcommand{\psl}{\text{PSL}\parentesis{2,\C}}
\newcommand{\PSLN}{\text{PSL}\parentesis{n+1,\C}}
\newcommand{\SLN}{\text{SL}\parentesis{n+1,\C}}
\newcommand{\SL}{\text{SL}\parentesis{3,\C}}
\newcommand{\sL}{\text{SL}\parentesis{2,\C}}
\newcommand{\GLN}{\text{GL}\parentesis{n+1,\C}}
\newcommand{\GL}{\text{GL}\parentesis{3,\C}}

\newcommand{\MC}{\mathcal{M}_{3\times 3}\parentesis{\C}}
\newcommand{\qp}{\text{QP}\parentesis{2,\C}}
\newcommand{\QP}{\text{QP}\parentesis{3,\C}}
\newcommand{\QPN}{\text{QP}\parentesis{n+1,\C}}

\newcommand{\pu}{\text{PU}\parentesis{2,1}}
\newcommand{\fix}{\text{Fix}}
\newcommand{\kernel}{\text{Ker}}
\newcommand{\core}{\text{Core}}
\newcommand{\rank}{\text{rank}}

\newcommand{\KulL}{\Lambda_{\text{Kul}}}
\newcommand{\KulD}{\Omega_{\text{Kul}}}
\newcommand{\Eq}{\text{Eq}}
\newcommand{\CG}{\Lambda_{\text{CG}}}

\newcommand{\FraL}{\Lambda_{\text{F}}}
\newcommand{\FraD}{\Omega_{\text{F}}}
\newcommand{\grL}{\Lambda_{\text{Gr}}}

\newcommand*{\linproy}[1]{\overleftrightarrow{#1}}

\newcommand{\bola}[3]{\text{B}^{#1}_{#2}\parentesis{#3}}

\newcommand{\titem}[1]{{\scriptsize \textbf{#1}}}
\newcommand{\titemT}[1]{{\small \textbf{#1}}}

\begin{document}
 
\title{Three Generalizations Regarding Limit Sets for Complex Kleinian Groups}
\author{Gerardo Mauricio Toledo Acosta}
\date{\today}

\maketitle

\tableofcontents

\chapter*{Acknowledgements}

Undertaking and completing this PhD has been truly a life-changing journey and a dream come true. This would not be possible without the support and guidance I received from many people and institutions.\\

First, I would like to express my sincere gratitude to my advisor Dra. M\'onica Moreno Rocha for the continuous support during this five year period. Her guidance helped me in the time of research and writing. I'm also very grateful for helping me getting support to assist to several great and very useful workshops and events in these years.\\

I also want to say a big thank you to my co-advisor Dr. \'Angel Cano Cordero. His patience, dedication and ideas motivated me and helped me greatly, especially in the hardest parts of the work. \\

My sincere thanks also goes to Dr. Manuel Cruz L\'opez for his guidance during the first three years. The work in Chapter \ref{chapter_measures} started based on his insightful ideas and suggestions.\\

I would like to thank the rest of the thesis committee: Dr. Xavier Gomez-Mont, Dr. Carlos Cabrera and Dr. Waldemar Barrera for the helpful advice, insightful comments and their crucial remarks that shaped this final dissertation.\\

I want to thank CONACYT for the PhD scholarship during the four years of the program. I also want to thank Dr. Jos\'e Seade for the scholarship during this last year through his project 282937 \emph{Din\'amica y Geometr\'ia Real y Compleja}.\\

The work of this thesis would not have come to a successful completion without the support of the community and staff of CIMAT and IMATE, Cuernavaca. Both places were great environments for working, learning and developing the ideas of this work.\\

Finally, but not least, I would like to express my gratitude to my wife Mar\'ia and my family: my parents, my sister and her husband, my brother and my little nephew and niece. I owe them a lot, they kept me going through all these years.

\chapter*{Introduction}

Kleinian groups are discrete subgroups of $\psl$, the group of biholomorphic automorphisms of the complex projective line $\CP^1$, acting properly and discontinuously on a non-empty region of $\CP^1$. Equivalently, these groups can be viewed as groups of orientation-preserving isometries of the hyperbolic 3-space $\Hh^3_\R$ or as groups of conformal autormorphisms of the sphere $\Ss^2$. Kleinian groups have been thoroughly studied since the end of the 19th century. They were first studied by Lazarus Fuchs, when he wanted to understand whether the solutions of certain ordinary differential equations were algebraic or not. Felix Klein then improved on Fuchs' solution, thus starting the field of Kleinian Groups.\\

Kleinian groups have important applications to Riemann surfaces, Teichm\"uller theory, holomorphic dynamics and automorphic forms, among others. For further reading, see \cite{maskit}, \cite{taniguchi}, \cite{hubbard2006teichmuller}, \cite{hubbard2016teichmuller}.\\

There are two ways to generalize Kleinian groups to higher dimensions: One can either look at conformal automophisms of the $n$-sphere $\Ss^n$ or one can look at holomorphic automorphisms of the complex projective space $\CP^n$. This distinction occurs because in higher dimensions, conformality is not equivalent to holomorphicity in general. In the first case one studies groups of isometries of real hyperbolic spaces, there is a rich body of knowledge in this subject thanks to the work of people like Ahlfors, Sullivan, Kapovich, Mostow, Thurston and many others. In the second case one is dealing with an area of mathematics that is still in its childhood. It has been studied by Jos\'e Seade, Angel Cano, Juan Pablo Navarrete, Waldemar Barrera and others (see \cite{ckg_libro} for a detailed introduction).\\

When dealing with complex Kleinian groups, there are several differences with the complex 1-dimensional case. One of these differences is the concept of limit set. In the classic case of Kleinian groups there are several equivalent characterizations for the limit set: 
	\begin{itemize}
	\item It is the closure of fixed points of loxodromic elements.
	\item It is the set of accumulation points of orbits.
	\item It is the complement of the maximal region where the action of the group is proper and discontinuous.
	\item It is the maximal region where the action is equicontinuous.
	\end{itemize}
In complex dimension 2, these notions no longer coincide in general. Each of these statements might define a different limit set for the action of a complex Kleinian group. Because of this distinction there are several possible definitions of limit set. One of these definitions is the Kulkarni limit set (see definition \ref{defn_kulkarni}), it is the most appropriate definition for our needs and it's one of the most used notions of limit set, at least in complex dimension 2. It was first defined on \cite{kulkarni}, as a limit set of a group acting on very general topological spaces and then it was adapted for complex Kleinian groups on \cite{ckg_libro}.\\

Some other definitions of limit set are: the Conze-Guivarc'h limit set (see \cite{Conze}) and the Myrberg limit set (see \cite{bcn16}). There are relations between them under certain hypotheses. For example, the Kulkarni limit set and the Conze-Guivarc'h limit set are dual to each other in complex dimension 2 (see \cite{bgn18} and \cite{tesisadriana}).\\

In chapter \ref{chapter_frances} we propose a definition for a new limit set for the action of a discrete subgroup of $\PSLN$ acting on $\CP^n$ that will be called the \emph{Frances limit set}. This limit set has the advantage of being purely dimensional, this will mean that the limit set is made up entirely of subspaces of the same dimension. The Kulkarni limit set does not have this property, as it can be the union of a line and a point (see Example \ref{ej_falso_hopf}). In this chapter we will prove some additional properties and give a few examples of this new limit set.\\

The Sullivan Dictionary gives relations between classical Kleinian groups and iteration of rational maps on $\CP^1$. It was first described by Dennis Sullivan in his 1985 paper \cite{sulldict}, where he introduced quasiconformal methods to the setting of rational maps and thus translated Ahlfors' finiteness theorem into the problem of wandering domains. The dictionary provides a conceptual framework for understanding the connections between dynamics of rational maps and Kleinian groups. In most cases, this dictionary suggests analogies which motivate research in both areas; in other cases, similar proofs can be given of related results in the two subjects. Some examples of this dictionary are the following (see \cite{qc_surgery}):

\begin{center}
\begin{tabular}{p{6cm}p{6cm}}
\textbf{Complex Dynamics} & \textbf{Kleinian Groups}\\
\hline
Julia Set & Limit Set\\
Fatou Set & Region of discontinuity\\
Rational or Entire function & Non-elemental group\\
Rational function with degree $\geq 2$ & Non-elemental, finitely generated group\\
No wandering domains theorem & Ahlfors finiteness theorem\\
Blaschke product & Fuchsian group\\
Hyperbolic rational map with two completely invariant components & Schottky group\\
Non-attracting cycles & Non-elliptic fixed points\\
\end{tabular}
\end{center}

This dictionary has been enriched by Curtis McMullen and others, for further reading see \cite{qc_surgery}. When one considers the generalization of Kleinian groups to complex dimension 2, it is natural to look for a similar dictionary between complex Kleinian groups and the dynamics of iteration of holomorphic functions of $\CP^2$. This latter field has been extensively studied by John Erik Fornaess and Nessim Sibony, among others (for further reading see \cite{forsib1}, \cite{forsib2}).\\ 

In chapter \ref{chapter_solvable} we study one of these possible correspondences. In the classic theory of Kleinian groups, \emph{elemental groups} are discrete subgroups of $\psl$ such that the limit set is a finite set (in other words, it is empty or it consists of 1 or 2 points). In complex dimension 2, the Kulkarni limit set is either a finite union of complex lines (1, 2 or 3) or it contains an infinite number of complex lines. Analogously, the limit set either contains a finite number of lines in general position (1, 2, 3 or 4) or it contains infinitely many lines in general position (see \cite{bcn16}). Therefore, one could define elemental groups in complex dimension 2 as discrete subgroups of $\PSL$ such that its Kulkarni limit set contains a finite number of lines or discrete subgroups of $\PSL$ such that its Kulkarni limit set contains a finite number of lines in general position. In bigger dimensions the situation is unknown. Another way one could define elemental subgroups of $\PSL$ is to consider groups with reducible action. In this chapter, we propose solvable subgroups of $\PSL$ as elemental groups.\\

We will describe the dynamics of discrete solvable subgroups of $\PSL$ by describing the Kulkarni limit set. We will also describe the groups themselves, by giving their representations.\\

One of the main goals in the study of complex Kleinian groups is the description of their dynamics through the description of the Kulkarni limit set. This study has been developing through the last 10 years in the following way:

\begin{itemize}
\item In \cite{navarrete2006limit}, it was proven that, for a discrete subgroup $\Gamma\subset\pu$ the Kulkarni limit set is made up of the tangent lines to the unitary open ball in the points of the Chen-Greenberg limit set. 
\item In \cite{bcn2011}, it was studied the case of subgroups of $\PSL$ leaving no proper subspace invariant, it was found that the Kulkarni limite set agreed with the equicontinuity region.
\item In \cite{bcn11cuatrolineas}, it was described the Kulkarni limit sets consisting of 4 lines in general position, furthermore, it was proven that this limit set coincides with the equicontinuity region in this case.
\item In \cite{bcn14unalinea} was described the Kulkarni limit set of subgroups of $\PSL$ such that the limit set contains exactly one line.
\item In \cite{ppar}, the purely parabolic subgroups of $\PSL$ are studied.
\end{itemize}

Consider the following classes of groups.

\begin{center}
\begin{tabular}{ll}
A: & Solvable Groups \\
B: & Groups with 1 line in the limit set \\
C: & Groups with 1 line and 1 point in the limit set \\
D: & Groups with 2 lines in general position in the limit set \\
E: & Groups with 3 lines in general position in the limit set \\
F: & Groups with 4 lines in general position in the limit set \\
G: & Purely parabolic groups
\end{tabular}
\end{center}

\begin{figure}[H]
\begin{center}	
\includegraphics[height=35mm]{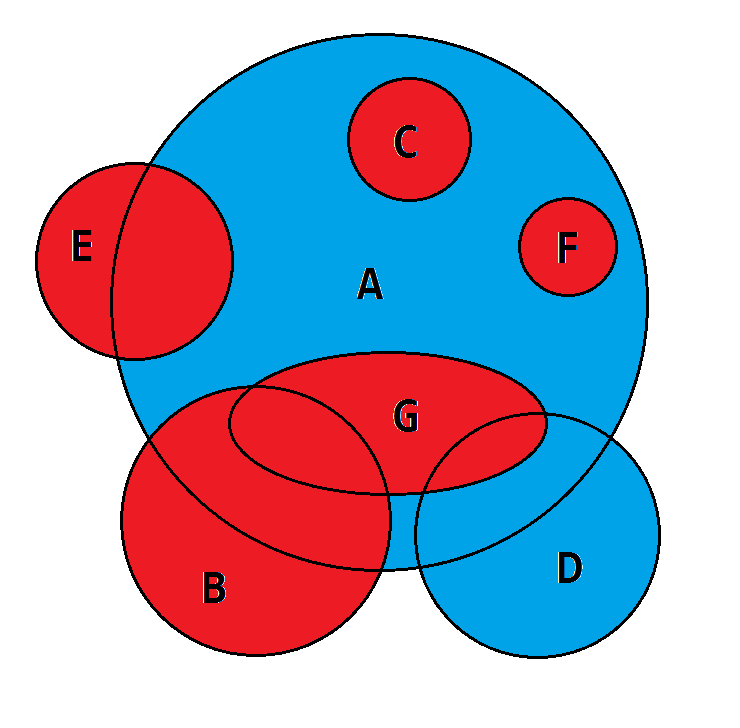}
\caption{Different classes of groups.}
\label{fig_introduccion_clases_de_grupos}
\end{center}	
\end{figure}

Summarizing, the description of the Kulkarni limit set has been provided for groups B, C, E, F and G (colored in red in figure \ref{fig_introduccion_clases_de_grupos}).

In chapter \ref{chapter_solvable} we give a full description of solvable groups (colored in blue in figure \ref{fig_introduccion_clases_de_grupos}). With this description, the full description of the Kulkarni limit set for general discrete subgroups of $\PSL$ will be almost complete. The main result we prove in this Chapter \ref{chapter_solvable} is the following theorem (see Theorem \ref{thm_main_solvable}).

\begin{thm*}
Let $\Gamma\subset\PSL$ be a solvable complex Kleinian group such that its Kulkarni limit set does not contain exactly four lines in general position. Then, there exists a non-empty open region $\Omega_\Gamma\subset \CP^2$ such that
	\begin{enumerate}[(i)]
	\item $\Omega_\Gamma$ is the maximal open set where the action is proper and discontinuous.
	\item $\Omega_\Gamma$ is homeomorphic to one of the following regions: $\C^2$, $\C^2\setminus\set{0}$, $\C\times\parentesis{\Hh^+\cup\Hh^-}$ or $\C\times\C^\ast$.
	\item $\Gamma$ is finitely generated and $\rank(\Gamma)\leq 4$.
	\item The group $\Gamma$ can be written as
		$$\Gamma=\Gamma_p\rtimes \Gamma_{LP}\rtimes\Gamma_{L},$$
	where $\Gamma_p$ is the subgroup of $\Gamma$ consisting of all the parabolic elements of $\Gamma$, $\Gamma_{LP}$ is the subgroup of $\Gamma$ consisting of all the loxo-parabolic elements of $\Gamma$ and $\Gamma_{L}$ is the subgroup of $\Gamma$ consisting of all the strongly loxodromic and  complex homotheties of $\Gamma$.
	\item The group $\Gamma$, up to a finite index subgroup, leaves a full flag invariant. 
	\end{enumerate}
\end{thm*}

In chapter \ref{chapter_measures} we study the existence of quasi-invariant measures of probability supported on the Kulkarni limit set (or the Conze-Guivarc'h limit set) of a discrete subgroup of $\PSL$. In complex dimension 1, these measures are the Patterson-Sullivan measures. They were first hinted in \cite{beardon1971inequalities}, when Alan Beardon was studying the Hausdorff dimension of the limit sets of Fuchsian groups. In \cite{patterson}, Patterson introduced a measure of probability whose exponent of convergence was the desired Hausdorff dimension. In \cite{sullivan79}, Sullivan generalized these measures for Kleinian groups. In this chapter we study the generalization of these measures to complex dimension 2 and present some advances towards said generalization. As we will prove, the existence of these measures depends on the finiteness of the entropy volume of the Kobayashi metric. We prove the following:\\

\begin{thm*}
Let $\Gamma\subset\PSL$ be a strongly irreducible complex Kleinian group acting on $\Omega=\KulD(\Gamma)$ and let $z\in\Omega$ such that $e(\Omega,z)<\infty$. Then, there exists a family of quasi-invariant measures supported on the Kulkarni limit set of $\Gamma$. 
\end{thm*}

We also give a strategy to estimate the entropy volume of the Kobayashi metric and present some partial results in this direction.\\

In chapter \ref{chapter_preliminaries} we introduce the definitions and properties needed to understand complex Kleinian groups and we cover some additional topics, like the Kobayashi metric needed for Chapter \ref{chapter_measures} or the Plucker embedding needed for Chapter \ref{chapter_frances}.\\
\chapter{Preliminaries}
\label{chapter_preliminaries}

In this chapter we review briefly the concepts and tools needed for the rest of this work. In general, we present the results without proof.

\section{The complex projective space $\CP^n$}

The complex projective space $\CP^n$ is defined by $\parentesis{\C^{n+1}\setminus\set{0}}\bigg/\sim$, where $\sim$ is the equivalence relation on $\C^3\setminus\set{0}$ given by
	$$Z\sim W\;\Longleftrightarrow\; \exists \lambda\in\C^\ast,\text{ such that } Z=\lambda W,$$

where $\C^\ast=\Cc$. $\CP^n$ is a compact complex $n$-manifold, it is the space of complex projective lines through the origin. We will denote by $\corchetes{z}\in\CP^n$ the projectivization of $z\in\C^{n+1}$, that is, the equivalence class of $z$.\\ 

In the same way, if $A\subset\C^{n+1}$ is a subset, we denote by $\corchetes{A}\subset\CP^n$ to the projectivization of $A$. If $A$ is a vector subspace of $\C^{n+1}$ then we say that $\corchetes{A}$ is a projective subspace of $\CP^n$.\\

If $x\in\CP^n$ is a point (resp. a subset $A\subset\CP^n$) we denote by $\mathbf{x}\in\C^{n+1}$ (resp. $\mathbf{A}\subset\C^{n+1}$) to any lift of $x$ (resp. $A$).\\

In the particular case when $n=2$ we call $\CP^2$ the \emph{complex projective plane}.

\begin{defn}\label{defn_linea}
Let $p,q\in\CP^n$, we define the complex projective line $\overleftrightarrow{p,q}$ as the intersection of all the projective subspaces $A\subset\CP^n$ such that $\set{p,q}\subset A$.
\end{defn}

\begin{prop}
If $p,q\in\CP^n$, let $\mathbf{p},\mathbf{q}\in\C^{n+1}$ are lifts and if $\langle\mathbf{p},\mathbf{q}\rangle\subset\C^{n+1}$ is the complex vector subspace generated by $\mathbf{p}$ and $\mathbf{q}$, then
	$$\overleftrightarrow{p,q}=\corchetes{\prodint{\mathbf{p},\mathbf{q}}}.$$
\end{prop}

From the previous description of a complex projective line it follows that, for any two distinct points in $\CP^n$, there exists a unique complex projective line containing both points. Analogously, if $\ell_1$ and $\ell_2$ are two distinct complex projective lines in $\CP^2$, then $\ell_1\cap \ell_2$ contains exactly one point. 

\begin{defn}
Consider an array of $k\geq 2$ complex projective lines in $\CP^n$. If exactly two lines of the array pass through a point, we call this point \emph{regular}. If more than two lines pass through a point, we call this point \emph{singular}. We say that the array is in \emph{general position} if it only contains regular points. 
\end{defn}

\begin{figure}[H]
\begin{center}	
\includegraphics[height=35mm]{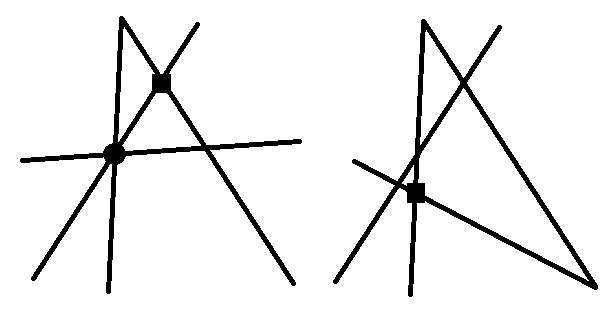}	
\label{fig_general_position}
\caption{An array of lines not in general position (left) and in general position (right). The \emph{round} point is a singular point and the \emph{square} points are regular.}
\end{center}	
\end{figure}

In Figure \ref{fig_general_position} we have an array of lines in general position (right) and an array not in general position (left).

In the previous figure and every time we draw complex projective lines, we have to keep in mind that the drawings are just illustrative and don't exactly represent the lines.\\

The following classical \emph{decomposition} of $\CP^n$ in terms of $\C^n$ will be very useful through the rest of this work.

\begin{prop}\label{prop_CPn_Cn}
$\CP^n$ is homeomorphic to $\C^n\cup\CP^{n-1}$.
\end{prop}

This proposition says that we can realize $\CP^n$ as $\C^n$ together with a hyperspace at infinity in the same way as we can see $\CP^1$ as $\C$ together with a point at infinity.

\section{The group of automorphisms of $\CP^2$}

Let $\MC$ be the group of all $3\times 3$ matrices with complex coefficients and let $\SL\subset\MC$ be the subgroup of matrices with determinant equal to $1$. The elements of $SL(3,\C)$ act linearly on $\C^3$, therefore they induce a transformation on $\CP^2$ in the following way: if $A\in SL(3,\C)$, we define $\gamma_A:\CP^2\rightarrow\CP^2$ as 
	$$\gamma_A(z)=[A\mathbf{z}]$$
where $\mathbf{z}\in\C^3$ is a lift of $z$ as before. It is clear that $\gamma_A$ is well defined and is an automorphism of $\CP^2$. Furthermore, it is well known that every biholomorphism of $\CP^2$ arises in this way. Then we have a group homomorphism between $\SL$ and $\text{Aut}(\CP^2)$, the group of biholomorphic automorphisms of $\CP^2$.\\

On the other hand, two matrices $A,B\in SL(3,\C)$ induce the same automorphism on $\CP^2$ if and only if $A=\omega B$ where $\omega$ is a cubic root of the unity. If we denote by $C_3=\set{\omega,\omega^2,1}$ the group of cubic roots of the unity, then the kernel of the group homomorphism is $C_3$, and therefore we have the following group isomorphism
	$$\PSL:=\SL\bigg/ C_3 \cong \text{Aut}(\CP^2).$$	
From now on, we will simply say \emph{automorphisms} instead of \emph{biholomorphic automorphisms}. We denote by $\fix(g)\subset\CP^2$ the set of fixed points of an automorphism $g\in\PSL$.

\begin{obs}\label{obs_eigenvect_fixpoints}
If $\gamma\in\SL$ and $v\in\C^3$ is an eigenvector of $\gamma$ then $[v]\in\fix{\corchetes{\gamma}}$. Furthermore, linearly independent eigenvectors of $\gamma$ define distinct fixed points of $\corchetes{\gamma}$.
\end{obs}

In the same way we define $\PSLN$, the group of automorphisms of $\CP^n$, as
	$$\PSLN:=\SLN\bigg/ C_n \cong \text{Aut}(\CP^n).$$

As before, if $A=(a_{ij})\in\SLN$ we denote by $[A]=\corchetes{a_{ij}}\in\PSLN$ the projectivization of $A$. If $g\in\PSLN$, we denote by $\mathbf{g}\in\SLN$ to any lift of $g$.\\

For $a_1,...,a_{n+1}\in\C^\ast$, we will denote by $A=\text{Diag}(a_1,...,a_{n+1})$ to the diagonal element
	$$A=\corchetes{\begin{array}{ccc}
	a_1 & ... & 0 \\
	... & \ddots & ...\\
	0 & ... & a_{n+1}\\ 
	\end{array}
	}\in\PSL,$$ 
or to any of its lifts $\mathbf{A}\in\SL$.\\

\subsection{Quasi-projective maps}

Quasi-projective maps will be a very useful tool to describe the Kulkarni limit set and the equicontinuity region of a discrete subgroup of $\PSL$ (see Propositions \ref{prop_descripcion_Eq} and \ref{prop_convergencia_qp}). For a more detailed overview, see Section 3 of \cite{cano2010equicontinuity} or Section 7.4 of \cite{ckg_libro}.\\

Let $M\in\MC$, consider the kernel of $M$, that is
	$$\kernel(M)=\SET{z\in\C^3}{Mz=0\in\C^3}.$$ 
Consider the projectivization of this set, $\corchetes{\kernel(M)\setminus\set{0}}$. Then $M$ induces a well defined map 
	$$\corchetes{M}:\CP^2\setminus \corchetes{\kernel(M)\setminus\set{0}}\rightarrow \CP^2,$$
given by
	$$\corchetes{M}(z)=\corchetes{M\mathbf{z}}.$$
We define in this way the quasi-projective maps $\QP$ as
	$$\QP=\parentesis{\MC\setminus\set{0}}\bigg/ \C^\ast.$$
This space is the closure of $\PSL$ in the space $\MC$ and therefore every sequence of elements in $\PSL$ converge to an element of $\QP$. For an element $T\in\QP$ we define its kernel as
	$$\kernel(T)=\corchetes{\kernel(\mathbf{T})\setminus\set{0}}.$$ 
We define $\QPN$ in a similar way. 

\section{Complex Kleinian groups}

In this section we define complex Kleinian groups and give some of their pro-perties. We also classify the automorphism of $\CP^2$ in three main types: elliptic, parabolic and loxodromic. For a more complete background see \cite{ckg_libro}. We start by giving some general definitions.

\begin{defn}
An element $g$ of a group $G$ is called a torsion element if it has finite order, i.e., if there is a positive integer $m$ such that $g^m=\id$. A group is called \emph{torsion free group} if the only torsion element it contains is the identity.  
\end{defn}

\begin{defn}
Let $G$ be a group acting on a space $X$ and let $x\in X$, we define the \emph{isotropy group} of $x$ by
	$$\text{Isot}(x,G)=\SET{g\in G}{gx=x}.$$
\end{defn}

\begin{defn}
Let $\Gamma$ be a subgroup of $\PSL$ and let $\Omega\subset\CP^2$ be an open $\Gamma$-invariant set. We say that the action of $\Gamma$ on $\Omega$ is \emph{proper and discontinuous} if, for every compact set $K\subset\Omega$, 
	$$\valorabs{\SET{\gamma\in\Gamma}{K\cap\gamma(K)\neq\emptyset}}<\infty.$$
\end{defn}

\begin{defn}
A discrete subgroup $\Gamma$ of $\PSL$ is said to be a \emph{complex Kleinian group} if it acts properly and discontinuously on some non-empty open subset of $\CP^2$.
\end{defn}

In the following example we introduce a family of complex Kleinian groups that will play an important role in Chapter \ref{chapter_solvable}.

\begin{ejem}\label{defn_inuoe}
Let us define:
	\begin{align*}
	\text{Sol}_0^4&=\SET{\corchetes{\begin{array}{ccc}
	\lambda & 0 & a \\
	0 & \valorabs{\lambda}^{-2} & b \\
	0 & 0 & 1\\
	\end{array}}}{(\lambda,a,b)\in\C^{\ast}\times\C\times\R}\\
	\text{Sol}_1^4&=\SET{\corchetes{\begin{array}{ccc}
	\varepsilon & a & b \\
	0 & \alpha & c \\
	0 & 0 & 1\\
	\end{array}}}{\alpha,a,b,c\in\R\text{, }\alpha>0\text{, }\varepsilon=\pm 1}\\
	\text{Sol'}_1^4&=\SET{\corchetes{\begin{array}{ccc}
	1 & a & b+i\text{log}\alpha \\
	0 & \alpha & c \\
	0 & 0 & 1\\
	\end{array}}}{\alpha,a,b,c\in\R\text{, }\alpha>0}.
	\end{align*}
We identify $\C^2$ with the affine chart $\SET{[z_1:z_2:1]\in\CP^2}{z_1,z_2\in\C^2}$. Consider a compact complex projective 2-manifold 
	$$M=\parentesis{\C\times\Hh}/\Gamma,$$
where $\Gamma=\pi_1(M)$ is a torsion free group contained in one of the Lie groups $\text{Sol}_0^4$, $\text{Sol}_1^4$ or $\text{Sol'}_1^4$. The manifold $M$ is called an \emph{Inuoue surface} and $\Gamma$, a \emph{fundamental Groups of Inuoe surfaces}. The group $\Gamma$ is an example of complex Kleinian groups (see \cite{brunella_inuoe} and Section 8.3 of \cite{ckg_libro} for further details).	
\end{ejem}   

There is no standard definition of limit set for the action of subgroups of $\PSL$, as in the case of subgroups of $\psl$. We first define the Kulkarni limit set (see \cite{kulkarni}), this definition has the advantage that the group acts properly and discontinuously on its complement, however these region is not maximal in general.

\begin{defn}\label{defn_kulkarni}
Let $\Gamma$ be a discrete subgroup of $\PSLN$ acting on $\CP^n$. Let us define the following sets:
	\begin{itemize}
	\item Let $L_0(\Gamma)$ be the closure of the set of points in $\CP^n$ with infinite group of isotropy.
	\item Let $L_1(\Gamma)$ be the closure of the set of accumulation points of orbits of points in $\CP^n\setminus L_0(\Gamma)$.
	\item Let $L_2(\Gamma)$ be the closure of the set of accumulation points of orbits of compact subsets of $\CP^n\setminus\parentesis{L_0(\Gamma)\cup L_1(\Gamma)}$. 
	\end{itemize}
We define the \emph{Kulkarni limit set} of $\Gamma$ as
	$$ \KulL(\Gamma):=L_0(\Gamma)\cup L_1(\Gamma)\cup L_1(\Gamma).$$
The \emph{Kulkarni region of discontinuity} of $\Gamma$ is defined as
	$$\KulD(\Gamma):=\CP^n\setminus\KulL(\Gamma).$$
\end{defn}

From the previous definition is clear that $\KulL(\Gamma)$ is a closed subset. For simplicity we will write $\KulL(g)$ instead of $\KulL\parentesis{\prodint{g}}$ for any $g\in\PSLN$, analogously for $\KulD(g)$.\\

In the following proposition we summarize some properties of the Kulkarni limit set.

\begin{prop}\label{prop_kulk_tiene_una_linea}
\mbox{}
\begin{enumerate}[(i)]
\item $\KulL(\Gamma)$ is a $\Gamma$-invariant set and it might be empty. $\KulD(\Gamma)$ is a $\Gamma$-invariant open set and might also be empty, besides $\Gamma$ acts properly and discontinuously on $\KulD(\Gamma)$. 
\item If $\Gamma\subset\PSL$ is a discrete infinite subgroup then $\KulL(\Gamma)$ contains always a complex projective subspace of dimension at least 1.
\end{enumerate}
\end{prop}

The proof of (i) is given in \cite{kulkarni} and the proof of (ii) is given in \cite{ckg_libro}.

\begin{defn}\label{defn_elemental}
Let $\Gamma\subset\PSL$ be a discrete subgroup. We say that $\Gamma$ is \emph{elemental} if $\KulL(\Gamma)$ is a finite union of complex projective subspaces.
\end{defn}

The following theorem gives a description of $\KulL(\Gamma)$ when $\Gamma$ is elemental (see \cite{bcn2011}). In Chapter \ref{chapter_solvable} we will give a more detailed description of this situation.

\begin{thm}\label{thm_opciones_kul_elemental}
Let $\Gamma\subset\PSL$ be a discrete elemental subgroup, then one of the following conclusions hold:
	\begin{enumerate}[(i)]
	\item $\KulL(\Gamma)$ is a complex projective line.
	\item $\KulL(\Gamma)$ is the union of two complex projective lines.
	\item $\KulL(\Gamma)$ is the union of three complex projective lines in general position.
	\item $\KulL(\Gamma)$ is the union of a complex projective line and a point not belonging to this line. 
	\end{enumerate}
\end{thm}

\begin{defn}
The \emph{equicontinuity region} for a familiy $\Gamma$ of automorphisms of $\CP^n$, denoted $\Eq(\Gamma)$, is defined to be the set of points $z\in\CP^n$ for which there is an open neighborhood $U$ of $z$ such that $\Gamma$ restricted to $U$ is a normal family. 
\end{defn}

The following propositions, will be very useful to compute the Kulkarni limit set in terms of the quasi-projective limits (see Proposition 7.4.1 of \cite{ckg_libro}, Proposition 2.5 and Corollary 2.6 of \cite{clu2017} respectively). The next proposition is known as the $\lambda$-lemma.

\begin{prop}\label{prop_convergencia_qp}
Let $\set{\gamma_k}\subset\PSLN$ a sequence of distinct elements, then there is a subsequence of $\set{\gamma_k}$, still denoted by $\set{\gamma_k}$, and a quasi-projective map $\gamma\in\QPN$ such that $\gamma_k\overset{k\rightarrow\infty}{\longrightarrow}\gamma$ uniformly on compact sets of $\CP^n\setminus\corchetes{\text{Ker}(\gamma)}$.
\end{prop}

\begin{prop}\label{prop_descripcion_Eq}
Let $\Gamma\subset\PSLN$ be a group, we say $\gamma\in\QPN$ is a limit of $\Gamma$, in symbols $\gamma\in\text{Lim}(\Gamma)$, if there is a sequence $\set{\gamma_m}\subset\Gamma$ of distinct elements satisfying $\gamma_m\rightarrow\gamma$. Thus we have
	$$\Eq(\Gamma)=\CP^n\setminus\overline{\underset{\gamma\in\text{Lim}(\Gamma)}{\bigcup}\kernel(\gamma)}.$$
\end{prop}

\begin{prop}\label{prop_eq_in_kuld}
Let $\Gamma\subset\PSLN$ be a discrete group, then $\Gamma$ acts pro-perly discontinuously on $\Eq(\Gamma)$. Moreover 
	$$\Eq(\Gamma)\subset\KulD(\Gamma).$$
\end{prop}

The following theorem (see Theorem 6.3.3 of \cite{ckg_libro}) describes when the equicontinuity region and $\KulD(\Gamma)$ coincide. This suggest that the Kulkarni limit set is the \emph{right} notion of limit set in complex dimension 2.

\begin{thm}
Let $\Gamma\subset\PSL$ be a complex Kleinian group such that $\KulL(\Gamma)$ contain four lines in general position. Then:
	\begin{enumerate}[(i)]
	\item $\KulD(\Gamma)=\Eq(\Gamma)$.
	\item If the group acts without fixed points or fixed lines then $\KulD(\Gamma)$ is the largest open set where $\Gamma$ acts properly and discontinuously. 
	\end{enumerate}	
\end{thm}

We give now some more definitions and general results of complex Kleinian groups which will be used in the following chapters. These can be found in \cite{ckg_libro}, \cite{bcn2011}, \cite{bcn16} and \cite{cs2014}.

\begin{prop}
Let $\Gamma\subset\PSLN$ be a discrete group acting properly and discontinuously on some open subset $\Omega\subset\CP^n$. Then, for any compact set $K\subset \Omega$, the set of accumulation points of $\Gamma$-orbits is a subset of $\CP^n\setminus\Omega$.
\end{prop}

\begin{defn}\label{defn_irreducible}
Let $\Gamma\subset\PSL$ be a discrete subgroup.
\begin{enumerate}[(i)]
\item We say that the action of $\Gamma$ is \emph{irreducible} if there is no proper non-empty subspace, invariant under the action of $\Gamma$. In any other case, we say that the action is \emph{reducible}.
\item We say that the action of $\Gamma$ is \emph{strongly irreducible} if there is no proper non-empty subspace, invariant under the action of some finite index subgroup of $\Gamma$.
\end{enumerate}
\end{defn}

\begin{defn}\label{defn_suspension}
If $\Gamma\subset\psl$ is a discrete subgroup, the \emph{suspension} of $\Gamma$ in $\PSL$ is the subgroup given by
	$$\SET{\corchetes{\begin{array}{cc}
	\mathbf{\gamma} & 0 \\ 
	0 & 1
	\end{array}}}{\gamma\in\Gamma}$$
where $\mathbf{\gamma}\in\sL$ is a lift of $\gamma\in\Gamma$. The suspension can also be defined as 
	$$\SET{\corchetes{\begin{array}{cc}
	1 & 0 \\ 
	0 & \mathbf{\gamma}
	\end{array}}}{\gamma\in\Gamma}$$
\end{defn}

\begin{ejem}\label{ejem_irreducible}
\mbox{}\begin{enumerate}[(i)]
\item Suspensions and upper triangular groups are examples of complex Kleinian groups with reducible action.
\item Some examples of groups with strongly irreducible action are: Schottky-like groups (see Definition \ref{defn_schottky_like}) and Veronesse groups (see Chapter 2 of \cite{tesisalex}).
\end{enumerate}

\end{ejem}

\subsection{Classification of elements of $\PSL$}
\label{sec_clasificacionPSL3C}

In this subsection we define the classes and subclasses of elements of $\PSL$. This classification will be needed in Chapter \ref{chapter_solvable}. For examples of the different classes of elements, their Kulkarni limit sets and more details, see Chapter 3 of \cite{ckg_libro}.\\

Before defining the different types of elements of $\PSL$, consider the following description of unitary balls in $\CP^2$.

\begin{obs}
We can identify the set 
	$$T=\set{[z_1:z_2:z_3]\in\CP^2\barra \valorabs{z_1}^2+\valorabs{z_2}^2=\valorabs{z_3}^2}$$
with the unitary $3$-sphere in $\C^2$. We consider this 3-sphere the unitary 3-sphere in $\CP^2$, observe that $\partial\HhC^2 = T$. In the same way, we can identify 
	$$B=\set{[z_1:z_2:z_3]\in\CP^2\barra \valorabs{z_1}^2+\valorabs{z_2}^2<\valorabs{z_3}^2}$$ 
with the unitary open $4$-ball in $\C^2$ and 
	$$D=\set{[z_1:z_2:z_3]\in\CP^2\barra \valorabs{z_1}^2+\valorabs{z_2}^2\leq\valorabs{z_3}^2}$$ 
with the unitary closed $4$-ball.
\end{obs}

\begin{defn}\label{defn_3esferas}
We call the images of $T$ (resp. $B$) under elements of $\PSL$ the \emph{$3$-spheres} (resp. open $4$-balls) of $\CP^2$.
\end{defn}

\subsubsection{Elliptic Elements}

\begin{figure}[H]
\begin{center}	
\includegraphics[height=40mm]{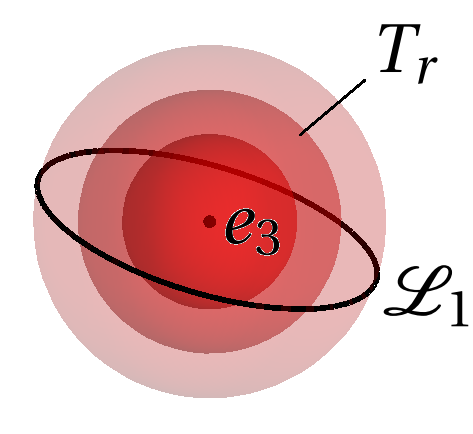}
\caption{Foliation of $\CP^2\setminus\parentesis{\mathcal{L}_1\cup\set{e_3}}$ by the family of 3-spheres $\set{T_r}_{r>0}$.}	
\label{fig_foliacion_CP2}
\end{center}	
\end{figure}

We define the family of 3-spheres of $\CP^2$
	$$T_r=\set{[z_1:z_2:z_3]\in\CP^2\barra \valorabs{z_1}^2+\valorabs{z_2}^2=r\valorabs{z_3}^2}$$ 
for $r>0$. Then $\set{T_r}_{r>0}$ is a foliation of $\CP^2\setminus\parentesis{\linproy{e_1,e_2}\cup\set{e_3}}$. We think of $e_3$ as the origin of $\CP^2$ and $\mathcal{L}_1=\linproy{e_1,e_2}$ as the line at infinity (see figure \ref{fig_foliacion_CP2}).

\begin{obs}
The elements of the family $\set{T_r}_{r>0}$ are 3-spheres centered in $e_3$, in fact, $T_r=g_r(T)$ where

	$$g_r=\corchetes{
	\begin{array}{ccc}
 	\lambda_1 & 0 & 0 \\
 	0 & \lambda_2 & 0 \\
 	0 & 0 & \lambda_3
	\end{array}}\in\PSL,$$
	
where $\valorabs{\lambda_1}=\valorabs{\lambda_2}$ and $\lambda_3=\frac{\lambda_1}{\sqrt{r}}$. On the other hand, the images of $T$ under arbitrary elements of $\PSL$ are not centered in $e_3$ in general.
\end{obs}

\begin{defn}
We say that the element $g\in \PSL$ is \emph{elliptic} if it preserves each leaf of the foliation $\set{T_r}_{r>0}$. In other words, there exists $h\in \PSL$ such that $h^{-1}gh\parentesis{T_r}=T_r$ for every $r>0$.
\end{defn}

We have the following characterizations of elliptic elements of $\PSL$ (see Section 4.2.1 of \cite{ckg_libro}).

\begin{prop}
An element $g\in \PSL$ is elliptic if and only if it has a lift $\mathbf{g}\in SL(3,\C)$ such that $\mathbf{g}$ is diagonalizable and every eigenvalue $\lambda_i$ satisfies $\valorabs{\lambda_i}=1$. In other words, it has a canonical Jordan form $D$ such that

	$$D=\parentesis{
	\begin{array}{ccc}
 	e^{2\pi i \theta_1} & 0 & 0 \\
 	0 & e^{2\pi i \theta_2} & 0 \\
 	0 & 0 & e^{2\pi i \theta_3}
	\end{array}},$$
	
with $\theta_1,\theta_2,\theta_3\in\R$.
\end{prop}

\begin{prop}
An element $g\in \PSL$ is elliptic if and only if $\KulL(g)=\emptyset$ or $\KulL(g)=\CP^2$ (according to whether the order of $g$ is finite or infinite).
\end{prop}

Subgroups of $\PSL$ containing elliptic elements of infinite order cannot be discrete.

\subsubsection{Parabolic Elements}

\begin{defn}
An element $g\in \PSL$ is \emph{parabolic} if there exist a family of $g$-invariant 3-spheres $\mathcal{F}$ and a point $z_0\in\fix(g)$ such that
	\begin{enumerate}
	\item For any different $T_1$ and $T_2\in\mathcal{F}$, it holds $T_1\cap T_2=\set{z_0}$.
	\item $\overline{\bigcup \mathcal{F}}$ is a closed 4-ball in $\CP^2$.
	\end{enumerate}
\end{defn}

\begin{obs}
Unlike the case of parabolic elements of $\psl$, a parabolic element of $\PSL$ might have two fixed points or a whole complex projective line made of fixed points.
\end{obs} 

There are two subclasses of parabolic elements in $\PSL$.

\begin{defn}
A parabolic transformation is:
	\begin{itemize}
	\item \emph{Unipotent}, if it has a lift in $SL(3,\C)$ such that every eigenvalue is $1$.
	\item \emph{Ellipto-parabolic}, if it is not unipotent.
	\end{itemize}
\end{defn}

\begin{obs}\label{obs_forma_jordan_parab}
From the previous definiton it follows that an element $g\in \PSL$ is parabolic if and only if it has a non-diagonalizable lift $\mathbf{g}\in\SL$ such that every eigenvalue has norm $1$. If $g$ is unipotent, the canonical Jordan form of $\mathbf{g}$ is
	$$
	\begin{array}{ccc}
	\parentesis{
	\begin{array}{ccc}
 	1 & 1 & 0 \\
 	0 & 1 & 0 \\
 	0 & 0 & 1
	\end{array}} & \text{or} & \parentesis{
	\begin{array}{ccc}
 	1 & 1 & 0 \\
 	0 & 1 & 1 \\
 	0 & 0 & 1
	\end{array}}
	\end{array}.		
	$$
	
If $g$ is ellipto-parabolic, the canonical Jordan form of $\mathbf{g}$ is
	$$
	\parentesis{
	\begin{array}{ccc}
 	\lambda & 1 & 0 \\
 	0 & \lambda & 0 \\
 	0 & 0 & \lambda^{-2}
	\end{array}},$$
with $\valorabs{\lambda}=1$ and $\lambda\neq 1$.
\end{obs}

\begin{defn}
Let $g\in\PSL$ be an ellipto-parabolic element with a lift conjugate to 
	$$\parentesis{
	\begin{array}{ccc}
 	e^{2\pi i \theta} & 1 & 0 \\
 	0 & e^{2\pi i \theta} & 0 \\
 	0 & 0 & e^{-4\pi i \theta}
	\end{array}},$$
we say that $g$ is \emph{rational (resp. irrational) ellipto-parabolic} if $\theta\in\Q$ (resp. $\theta\in\R\setminus\Q$). 
\end{defn}

\subsubsection{Loxodromic Elements}

\begin{defn}
An element $g\in \PSL$ is \emph{loxodromic} if there is an open set $W\subset \CP^2$ such that 
	$$g\parentesis{W\cup \mathbb{X}_i}\subset \mathbb{X}_i$$
for $i=1$ or $i=2$, where $\mathbb{X}_1$ and $\mathbb{X}_2$ are the connected components of $\CP^2\setminus W$.  
\end{defn}


Now we define the four subclasses of loxodromic elements of $\PSL$.

\begin{defn} 
A loxodromic element is called:
\begin{enumerate}
\item \emph{Loxo-parabolic} if it is conjugated to an element $h\in\PSL$ such that
	$$\mathbf{h}=\parentesis{
	\begin{array}{ccc}
 	\lambda & 1 & 0 \\
 	0 & \lambda & 0 \\
 	0 & 0 & \lambda^{-2}
	\end{array}},$$
with $\valorabs{\lambda}\neq 1$.
\item \emph{Complex homothety} if it is conjugated to an element $h\in\PSL$ such that
	$$\mathbf{h}=\parentesis{
	\begin{array}{ccc}
 	\lambda & 0 & 0 \\
 	0 & \lambda & 0 \\
 	0 & 0 & \lambda^{-2}
	\end{array}},$$
with $\valorabs{\lambda}\neq 1$.
\item \emph{Rational screw} (resp. \emph{irrational screw}) if it is conjugated to an element $h\in\PSL$ such that
	$$\mathbf{h}=\parentesis{
	\begin{array}{ccc}
 	\lambda_1 & 0 & 0 \\
 	0 & \lambda_2 & 0 \\
 	0 & 0 & \lambda_3
	\end{array}},$$
with $\valorabs{\lambda_1}=\valorabs{\lambda_2}\neq \valorabs{\lambda_3}$ and $\frac{\lambda_1}{\lambda_2}=e^{2\pi i x}$ with $x\in \Q$ (resp. $x\in\R\setminus\Q$).
\item \emph{Strongly loxodromic} if it is conjugated to an element $h\in\PSL$ such that
	$$\mathbf{h}=\parentesis{
	\begin{array}{ccc}
 	\lambda_1 & 0 & 0 \\
 	0 & \lambda_2 & 0 \\
 	0 & 0 & \lambda_3
	\end{array}},$$
where the elements $\set{\valorabs{\lambda_1},\valorabs{\lambda_2},\valorabs{\lambda_3}}$ are pairwise different.
\end{enumerate}
\end{defn}

We summarize the classification of elements of $\PSL$ in the following table. 

\begin{center}
\begin{tabular}{|l|l|l|c|}
  \hline
   & R & \small Pairwise distinct unitary eigenvalues. &  \\ \cline{2-3}
   
  E & CR & \small Two different eigenvalues (all eigenvalues are unitary). & Diag \\ \hline
  \hline
   &  U & \small All eigenvuales equal to 1. &  \\ \cline{2-3}
  P & EP & \small Not unipotent. & N-Diag \\ \hline
  \hline
   & LP & \small Two different eigenvalues, Jordan block of size $2\times 2$, $\valorabs{\lambda_1}\neq\valorabs{\lambda_2}$. & \\ \cline{2-3}
   & S & \small Three different eigenvalues, $\valorabs{\lambda_1}=\valorabs{\lambda_2}\neq \valorabs{\lambda_3}$ & Diag\\ \cline{2-3} 
  L & CH & \small Two different eigenvalues, $\valorabs{\lambda_1}\neq\valorabs{\lambda_2}$. & \\ \cline{2-3}
   & SL & \small Three different eigenvalues, $\valorabs{\lambda_1}\neq\valorabs{\lambda_2}\neq\valorabs{\lambda_3}\neq\valorabs{\lambda_1}$. & \\ 
  \hline
\end{tabular}
\end{center}

\begin{center}
\begin{tabular}{|r|r|r|r|}
\hline
E & Elliptic & R & Regular\\
P & Parabolic & CR & Complex Reflection\\
L & Loxodromic & U & Unipotent\\
 & & EP & Ellipto-Parabolic\\
Diag & Diagonalizable & LP & Loxo-Parabolic \\
N-Diag & Non-diagonalizable & S & Screw \\
& & CH & Complex Homothety \\
& & SL & Strongly-Loxodromic \\ 
\hline
\end{tabular}
\end{center}

\section{Groups with a control group}

In this section we describe a construction useful for reducing the action of a group $\Gamma\subset\PSL$ on $\CP^2$ to the action of a subgroup of $\psl$ on a complex line in $\CP^2$, thus, simplifying the study of the dynamics of $\Gamma$. The details and proofs of this construction are given in chapter 5 of \cite{ckg_libro}.\\

Consider a subgroup $\Gamma\subset\PSL$ acting on $\CP^2$ with a global fixed point $p\in\CP^2$. Let $\ell\subset\CP^2\setminus\set{p}$ be a projective complex line, we define the projection
	$$\pi=\pi_{p,\ell}:\CP^2\rightarrow\ell$$
given by $\pi(x)=\ell\cap\linproy{p,x}$. This function is holomorphic, and it allows us to define the group homomorphism
	$$\Pi=\Pi_{p,\ell}:\Gamma\rightarrow\text{Bihol}(\ell)\cong\psl$$
given by $\Pi(g)(x)=\pi(g(x))$ for $g\in\Gamma$ (see Lemma 6.11 of \cite{cs2014}). 

\begin{center}
\begin{figure}[h]
\label{fig_defn_grupo_control}
\includegraphics[height=32mm]{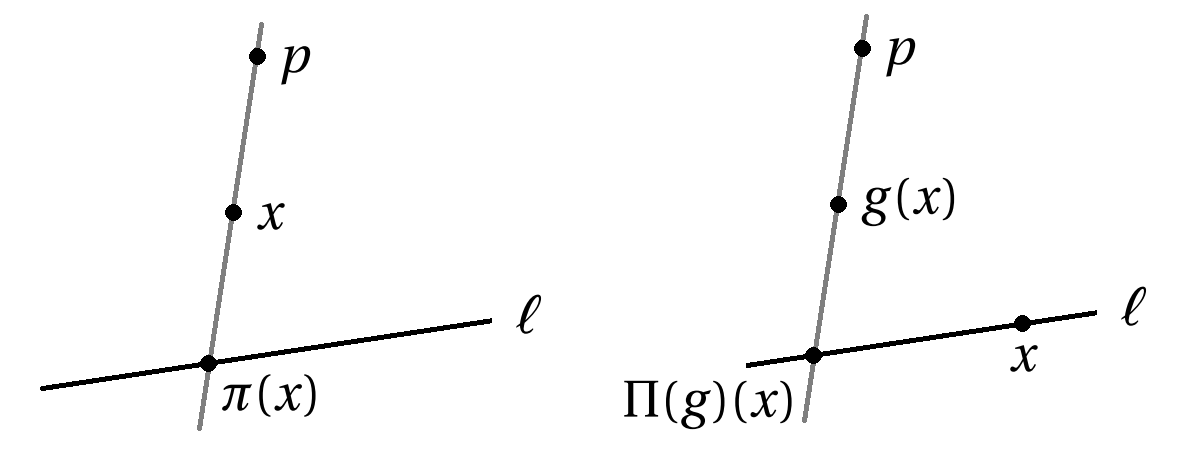}
\caption{Definition of $\pi_{p,\ell}$ and $\Pi_{p,\ell}$.}
\end{figure}
\end{center}

If we choose another line, $\ell'\subset\CP^2\setminus\set{p}$, we obtain a projection $\pi'=\pi_{p,\ell'}$ and a group homomorphisms $\Pi'=\Pi_{p,\ell'}$. The homomorphisms $\Pi$ and $\Pi'$ are equivalent in the sense that there exists a biholomorphism $h:\ell\rightarrow\ell'$ inducing an automorphism $H$ of $\psl$ such that $H\circ\Pi=\Pi'$. The line $\ell$ is called \emph{the horizon}.\\

To simplify the notation we will write $\text{Ker}(\Gamma)$ instead of $\text{Ker}(\Pi)\cap \Gamma$.

\begin{defn}
Let $\Gamma\subset\PSL$ be a discrete subgroup. We say that the group $\Gamma$ is \emph{weakly semi-controllable} if it acts with a fixed point in $\CP^2$. In this case, a choice of a horizon $\ell$ determines a \emph{control group} $\Pi(\Gamma)\subset\psl$, which is well-defined and independent of $\ell$ up to and automorphism of $\psl$.
\end{defn}

\subsection{Non-discrete subgroups of $\psl$}\label{subsec_groups_no_discretos}

In this subsection we give a quick overview of non-discrete subgroups of $\psl$. This will be needed in Chapter \ref{chapter_solvable} as the control group of discrete subgroups of $\PSL$ may be conjugated to a non-discrete subgroup of $\psl$. For more details, see Chapter 5 of \cite{ckg_libro} and Section 2 of \cite{cs2014}.

\begin{defn}
Let $\Sigma$ be a non-discrete subgroup of $\psl$. The \emph{limit set of $\Sigma$ in the sense of Greenberg}, denoted by $\Lambda_{Gr}(\Sigma)$, is the intersection with $\mathbb{S}^2_{\infty}:=\partial\Hh^3_\R$ of the set of accumulation points of de orbits of points in $\Hh^3_\R$.
\end{defn}

\begin{defn}
A non-discrete group $\Sigma\subset\psl$ is elemental if its equicontinuity set omits at most 2 points in $\CP^1$.
\end{defn}

\begin{thm}
Let $\Sigma\subset\psl$ be a non-discrete subgroup then 
	$$\Lambda_{Gr}(\Sigma)=\CP^1\setminus\text{Eq}(\Sigma).$$
\end{thm}

\begin{thm}\label{thm_no_discreto_no_elem_cerradura}
Let $\Sigma\subset\psl$ be non-discrete subgroup such that $\valorabs{\Lambda_{Gr}(\Sigma)}\geq 2$. Then $\Lambda_{Gr}(\Sigma)$ is the closure of fixed points of loxodromic elements of $\Sigma$.
\end{thm}

\begin{ejem}\label{ej_defns_dih_rot_so3}
The following are important examples of non-discrete subgroups of $\PSL$ that will appear in the sequel (see Section 2 and 13 of \cite{cs2014}):
\begin{itemize}
\item We denote by $\rot$ the group generated by all rotations around the origin.
\item The infinite dihedral group, is given by
	$$\dih=\prodint{\rot,\;z\mapsto -z}.$$ 
It satisfies
	$$\Eq\parentesis{\rot}=\Eq\parentesis{\dih}=\Ss^2.$$
\item The special orthogonal group, denoted by $\so$, can be embedded in $\psl$ as follows:
	$$\so=\SET{\corchetes{\begin{array}{cc}
	a & -c\\
	c & \overline{a}	
	\end{array}}}{\valorabs{a}^2+\valorabs{c}^2=1}.$$
This is a purely elliptical group satisfying 
	$$\Eq\parentesis{\so}=\CP^1.$$
\item We denote by $\epa$ to the group of all affine M\"obius transformations which are either parabolic or elliptic,
	$$\epa=\SET{\corchetes{\begin{array}{cc}
	a & b\\
	0 & a^{-1}	
	\end{array}}}{\valorabs{a}=1\text{, }b\in \C}.$$
	It satisfies $\Eq\parentesis{\epa}=\C$.
\item The group $\autc$ is the group of all M\"obius transformations that leave $\C^\ast$ invariant. It can be embedded in $\psl$ as follows
	$$\autc=\prodint{\SET{\corchetes{\begin{array}{cc}
	a & 0\\
	0 & a^{-1}	
	\end{array}}}{a\in\C^\ast},\corchetes{\begin{array}{cc}
	0 & 1\\
	1 & 0	
	\end{array}}}.$$
	It satisfies $\Eq\parentesis{\autc}=\C^\ast$. 
\end{itemize}
\end{ejem}

In the following proposition we summarize Theorem 2.14, Corollary 13.4 and Proposition 13.15 of \cite{cs2014}. This proposition lists all the possibilities for the Greenberg limit set of a non-discrete subgroup of $\psl$ and describe each possibility.

\begin{prop}\label{prop_opciones_conj_lim_greenberg}
Let $\Sigma\subset\psl$ be a non-discrete group, then one of the following cases happen:
	\begin{itemize}
	\item $\grL(\Sigma)=\emptyset$ if and only if $\Sigma$ is finite or $\Sigma$ is conjugated to a dense subgroup of $\so$, $\rot$, or $\dih$.
	\item $\valorabs{\grL(\Sigma)}=1$ if and only if $\Sigma$ is conjugated to a subgroup of $\epa$, whose closure contains parabolic elements.
	\item $\valorabs{\grL(\Sigma)}=2$ if and only if $\Sigma$ is conjugated to a subgroup of $\autc$, whose closure contains loxodromic elements.
	\item $\grL(\Sigma)=\Ss^1$ if and only if $\Sigma$ is conjugated to a subgroup of $\pslr$.	
	\item $\grL(\Sigma)=\CP^1$ if and only if $\Sigma$ is conjugated to a subgroup of $\psl$.
	\end{itemize}
\end{prop}

\section{The Conze-Guivarc'h limit set}

In this section we define another limit set for the action of discrete subgroups of $\PSL$, this definition is based on \cite{Conze}. This limit set is smaller than the Kulkarni limit set. 

\begin{defn}\label{defn_proximal}
We say that an element $\gamma\in GL(3,\C)$ is \emph{proximal} if it has an unique eigenvalue $\lambda_0$ such that $\valorabs{\lambda_0}>\valorabs{\lambda}$ for any other eigenvalue $\lambda$ of $\gamma$. An eigenvector $v_0$ corresponding to $\lambda_0$ is called a \emph{dominant eigenvector} for $\lambda_0$.\\
We say that $v_0\in\C^3$ is a \emph{dominant vector} for $\Gamma$ if there exists an element $\gamma\in\Gamma$ such that $v_0$ is a dominant eigenvector for $\gamma_0$. 
\end{defn}

\begin{defn}\label{defn_CG1}
Let $\Gamma\subset\PSL$ be a discrete group. Denote by
	$$\CG(\Gamma)= \overline{\SET{\corchetes{v_0}}{v_0\text{ is a dominant vector for }\Gamma}}.$$
We call $\CG(\Gamma)$ the \emph{Conze-Guivarc'h limit set} of $\Gamma$.
\end{defn}

\begin{defn}\label{defn_minimal}
Let $F\subset\CP^n$ be a $\Gamma$-invariant closed subset. We say that $F$ is \emph{minimal} if it contains no proper closed $\Gamma$-invariant subspace.
\end{defn}



\begin{obs}
Observe that $\CG(\Gamma)\subset\KulL(\Gamma)$. This follows from the fact that $\CG(\Gamma)$ has two types of points: attractive fixed points of loxodromic elements and accumulation points of these attractive fixed points, the former are in $L_0(\Gamma)$ and the latter in $L_1(\Gamma)$.
\end{obs}

\begin{obs}\label{obs_cg_noelemental}
\mbox{}\begin{enumerate}
\item From Definition \ref{defn_CG1}, it follows that $\CG(\Gamma)=\emptyset$ whenever $\Gamma$ contains no loxodromic elements. 
\item From Theorem \ref{thm_opciones_kul_elemental} it follows that if $\Gamma$ is elemental, $\CG(\Gamma)$ contains 1, 2 or 3 points. When $\Gamma$ is non-elemental, $\CG(\Gamma)$ is a perfect set (see \cite{bcn2011}). 
\end{enumerate}
\end{obs} 

\subsection{The dual space $\dual{\CP^2}$}
\label{subsec_dualization}

In this subsection we give a quick review of the dual projective complex plane and the action of a discrete subgroup of $\PSL$ in this dual space, for more details see \cite{bgn18} and \cite{tesisadriana}. Dualization will be useful in Section \ref{cor_HC3_no_hay_en_no_conmutativos}, when we need to determine the set of acumulation points of orbits of compact sets. Dualization also provide a relation between the Kulkarni limit set and the Conze-Guivarc'h limit set (see \cite{bgn18} and \cite{tesisadriana}).\\

Denote by $\dual{\CP^2}$ to the space of all complex projective lines in $\CP^2$. Each line in $\CP^2$ is the projectivization of a complex plane in $\C^3$ passing through the origin, this plane is the set of points $(x,y,z)\in\C^3$ such that
	$$Ax + By + Cz = 0$$
where $A,B,C\in\C^3$. So, we associate with every plane in $\C^3$ passing through the origin the point $(A,B,C)\in\C^3$, this is well defined up to multiplication by elements of $\C^{\ast}$. Then, projecting we have
	$$\dual{\CP^2}\cong \CP^2,$$
the equivalence is given by the map
	$$[A:B:C]\mapsto[Ax+By+Cz=0].$$
	
Under this identification, the natural action of $g\in\PSL$ on $\dual{\CP^2}$ is given by
	$$g\cdot\ell \longleftrightarrow \parentesis{\mathbf{g}^{-1}}^T\parentesis{\begin{array}{c}
	A\\
	B\\
	C\end{array}},$$
where the line $\ell\in\dual{\CP^2}$ is identified with $[A:B:C]\in\CP^2$.
	
\section{Grassmanians and flags}

In this subsection we define the Grassmanian, which will be needed in Chapter \ref{chapter_frances}. We also state some of their properties (ver \cite{harris}). 

\begin{defn}
Let $0\leq k<n$, we define the \emph{Grassmanian} $Gr(k,n)$ as the space of all the $k$-dimensional projective subspaces of $\CP^n$.
\end{defn}

$Gr(k,n)$ is a compact connected complex manifold of dimension $k(n-k)$. We define the function
	$$\psi:Gr(n,k)\mapsto \mathbb{P}\parentesis{\bigwedge^{k+1} \C^{n+1}}$$
given by 
	$$V\overset{\psi}{\mapsto} \corchetes{v_1\wedge...\wedge v_{k+1}}$$
where $\set{v_1,...,v_{k+1}}$ is a base of $V$. This function is well defined since, if $\set{w_1,...,w_{k+1}}$ is another base for $V$ and if $W\in\mathcal{M}_{k+1}\parentesis{\C}$ is the matrix of change of basis, then
	$$v_1\wedge...\wedge v_{k+1}=\det(W)\parentesis{w_1\wedge...\wedge w_{k+1}}.$$
The function $\psi$ is an embedding, called the \emph{Pl\"ucker embedding}. The homogeneous coordinates in $\mathbb{P}\parentesis{\bigwedge^{k+1} \C^{n+1}}$ are called the \emph{Pl\"ucker coordinates} in $Gr(k,n)$.

\begin{obs}
If $\gamma\in\SLN$, then $\gamma$ acts on $\bigwedge^i \C^{n+1}$ in the following way
	$$\gamma\parentesis{v_1\wedge...\wedge v_i}=\gamma v_1\wedge...\wedge \gamma v_i,$$
where $v_1,...,v_i\in\C^{n+1}$ and $i<n+1$. 
\end{obs}

\begin{defn}
A \emph{flag} in $\CP^n$ is a sequence of $k$-dimensional projective subspaces $V_k\subset\CP^n$ such that
	$$V_1\varsubsetneq V_2 \varsubsetneq ... \varsubsetneq V_r $$
where $r\leq n$. We denote this flag as $\set{V_1,...,V_r}$. If $\text{dim}(V_k)=k$ for all $k=1,...,r=n$ then the flag is called \emph{complete} or \emph{full}. 
\end{defn}

\section{Intrinsic metrics}\label{sec_kobayashi}

In this section we review two instrinsic metrics, the Kobayashi metric and the Eisenman metric. These metrics, along with a few others, are defined on complex manifolds and depend only on the complex structure of the manifold and nothing else.\\

In subsection \ref{subsec_kobayashi} we study the main facts about the Kobayashi metric, these definitions and properties will be used in Chapter \ref{chapter_measures}. For details and proofs, see \cite{kobayashi} and \cite{kobayashi05}. In Subsection \ref{subsec_eisenman}, we will define the Eisenman metric, which shares similarities and properties with the Kobayashi metric, for more details, see \cite{graham1985some}. The main difference between them is that the Kobayashi volume uses complex poly-discs and the Eisenman volume uses complex balls embedded in the manifold.\\

\subsection{The Kobayashi metric}\label{subsec_kobayashi}

Let $M$ be a complex manifold. Let us define a pseudo-distance $d_M$ on $M$ in the following way: Given two points $p,q\in M$, we choose points 
	$$p=p_0,p_1,...,p_{k-1},p_k=q\in M,$$ 
points 
	$$a_1,...,a_k,b_1,...,b_k\in\D$$ 
and holomorphic functions $f_1,...,f_k:\D\rightarrow M$ such that $f_i(a_i)=p_{i-1}$ and $f_i(b_i)=p_i$ for $i=1,...,k$. We call a \emph{chain of holomorphic disks} to every such collection of points and functions. The length of the chain of holomorphic disks $\alpha$ is defined as
	$$\ell(\alpha)=\rho(a_1,b_1)+...+\rho(a_k,b_k).$$
We define then
	$$d_M(p,q)=\inf_\alpha \ell(\alpha),$$
where the infimum is taken over all the chains of holomorphic disks $\alpha$. It can be shown that $d_M$ is a pseudo-distance, that is, there can be two distinct points $p,q\in M$ such that $d_M(p,q)=0$.\\

\begin{figure}[H]
\begin{center}	
\includegraphics[height=38mm]{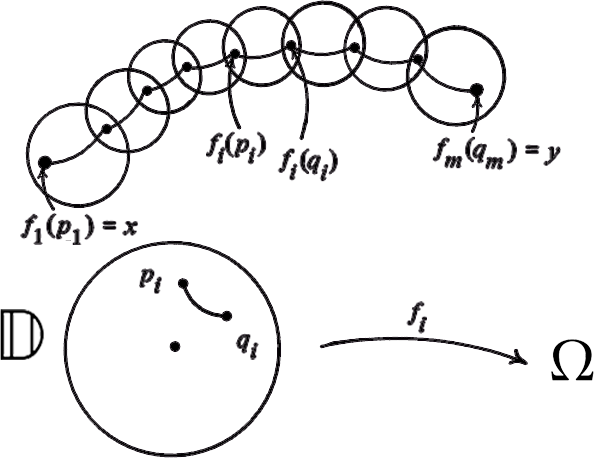}
\end{center}	
\end{figure}

\begin{defn}
If $M$ is a complex manifold such that $d_M$ is a metric, we say that $M$ is \emph{Kobayashi hyperbolic}. A Kobayashi hyperbolic space $M$ is \emph{complete} if it is Cauchy complete respect to $d_M$.  
\end{defn}

If $M$ is a Kobayashi hyperbolic complex manifold, we will always denote by $d_M$ the Kobayashi metric on $M$. 
Finally, for $z\in M$ and $r>0$, we denote by $\bola{r}{M}{z}$ to the open ball respect to the Kobayashi metric, 
	$$\bola{r}{M}{z}=\SET{w\in M}{d_M(z,w)<r}.$$

\begin{prop}\label{prop_dist_contractante}
Let $M$ and $N$ be two complex manifolds and let $f:M\rightarrow N$ be a holomorphic function. Then
	$$d_M(p,q)\geq d_N\parentesis{f(p),f(q)}$$
for any $p,q\in M$.
\end{prop}

As a immediate consequence of the previous result we have the following proposition.

\begin{prop}\label{prop_kob_isometrias}
Let $\Gamma\subset \PSL$ be a subgroup and let $\Omega\subset\CP^2$ be a Kobayashi hyperbolic subspace such that $\Omega$ is $\Gamma$-invariant. Then $\Gamma$ acts by isometries on $\Omega$, with respect to the Kobayashi metric $d_\Omega$.
\end{prop}
	
We will need the following theorem (see Theorem 1.3 of \cite{bcn2011}).

\begin{thm}\label{thm_complemento_kob_hyp}
Let $\Gamma\subset\PSL$ be a discrete subgroup without fixed points and without fixed lines, suppose that $\KulL(\Gamma)$ has more than 4 projective complex lines in general position, then each connected component of $\KulD(\Gamma)$ is complete Kobayashi hyperbolic.
\end{thm}

\subsection{The Eisenman metric}\label{subsec_eisenman}

The construction of these metrics was modelled on the construction of the Kobayashi metric. As we will see in Chapter \ref{chapter_measures}, there are some advantages of this metric over the Kobayashi metric. There are not many sources studying these metrics, unlike the other intrinsic metrics (see \cite{graham1985some}).\\

Again, let $M$ be an $n$-dimensional complex manifold. We denote by $B_n$ the unitary open ball in $\C^n$ and $\text{B}(0,R)\subset\C^n$ the usual open ball with radius $R>0$ centered in the origin. For $p\in M$, we denote by $T_p M$ the holomorphic tangent space to $M$ at $p$. Let $\bigwedge^k T_p M$ be the $k$-th exterior power of $T_p M$. The decomposable elements of $\bigwedge^k T_p M$ will be denoted by $D_p^k M$ (i.e. the elements of $\bigwedge^k T_p M$ that can be expressed as the exterior product of elements of $k$ elements of $T_p M$).

\begin{defn}
Let $k\in\Z$, $1\leq k\leq n$. Let $\alpha\in D_p^k M$ and $p\in M$. The \emph{intrinsic Eisenman norm} of $\alpha$ is
	$$\text{E}_k(p,\alpha)=\inf\SET{R^{-2k}}{\begin{array}{c}\exists\;\text{holomorphic } f:B_k(R)\rightarrow M\text{ such that }\\
	f(0)=p\text{ and }f_\ast \parentesis{\frac{\partial}{\partial z_1}\wedge ... \wedge \frac{\partial}{\partial z_k}(0)}=\alpha
	\end{array}}.$$
\end{defn}

Let $p\in M$ and let $\set{w_1,...,w_n}$ be complex coordinates near $p$. When $k=n$, we can re-formulate $E_k$ in the following way
	
	$$\text{E}_n\parentesis{p,\frac{\partial}{\partial w_1}\wedge ... \wedge \frac{\partial}{\partial w_n}(p)}=\inf\SET{\valorabs{Jf(0)}^{-2}}{\begin{array}{c}\exists\;\text{holomorphic } f:B_n\rightarrow M\\
	\text{ such that }f(0)=p
	\end{array}}$$

where $Jf$ denotes the complex Jacobian determinant of $f$. Now we introduce the \emph{Eisenman volume form} $\tau_M$ on $M$ defined by
	 $$\tau_M(p)=\text{E}_n\parentesis{p,\frac{\partial}{\partial w_1}\wedge ... \wedge \frac{\partial}{\partial w_n}(p)}\parentesis{\frac{i}{2}}^n dw_1\wedge d\overline{w}_1\wedge \cdots \wedge dw_n\wedge d\overline{w}_n.$$
Let $U\subset M$ be an open subset, then the \emph{intrinsic Eisenman volume} of $U$ is given by
	$$\text{Vol}_M(U)=\int_U \tau_M.$$ 
If there is no ambiguity, we will write $\text{Vol}(U)$ instead of $\text{Vol}_M(U)$. The Eisenman metric and volume satisfies the same contraction property given by Proposition \ref{prop_dist_contractante} (see Lemma 2.16 of \cite{graham1985some}).
\chapter{Dynamics of solvable groups}\label{chapter_solvable}
\label{chapter_solvable}

In this chapter we study the dynamics of discrete solvable subgroups of $\PSL$. These groups present \emph{simple} dynamics contrary to the \emph{rich} dynamics of strongly irreducible discrete subgroups of $\PSL$, which have been studied extensively (see \cite{bcn2011} or \cite{ckg_libro}). The study of solvable groups will help to complete the classification and understanding of the dynamics of all complex Kleinian groups of $\PSL$.\\

In the first section we give the necessary general definitions and results for the next sections. In Section \ref{sec_solvable_invariant_flag} we give a dynamical characterization of solvable subgroups of $\PSL$, in Section \ref{sec_solvables_triangularizables} we prove that solvable groups always contain a finite index triangularizable group. Therefore, if we want to describe solvable groups, a first step is to describe the discrete triangular subgroups of $\PSL$. In Section \ref{sec_commutative_triangular}, we first give a full description of commutative discrete triangular subgroups of $\PSL$ and, in Section \ref{sec_non_commutative_triangular}, we describe the non-commutative case. The description of the dynamics of solvable groups is summarized in Theorem \ref{thm_main_solvable}.

\section{Notation and properties of solvable groups}

In this section we give some definitions and algebraic results necessary for the next sections. We will use some standard concepts from algebraic geometry, see for example Chapter 2 of \cite{harris}, Chapter 2 of \cite{oni} or Section 1.2 of \cite{hartshorne2013algebraic}.\\

If $g$ and $h$ are two elements of a group $G$, we denote its commutator by
	$$\corchetes{g,h}=g^{-1}h^{-1}gh\in G.$$
If $g$ and $h$ commute then $\corchetes{g,h}=\id$. In the same way, we define the commutator subgroup 
	$$\corchetes{G,G}=\SET{\corchetes{g,h}}{g,h\in G}\subset G.$$ 
If $G$ is commutative then $\corchetes{G,G}$ is the trivial group.

\begin{defn}
Let $G$ be a group. The derived series $G^{(i)}$ of $G$ is defined inductively in the following way
	$$\begin{array}{cc}
	G^{(0)}=G, & G^{(i+1)}=\left[G^{(i)},G^{(i)}\right].
	\end{array}$$
We say that $G$ is \emph{solvable} if, for some $n\geq 0$, we have $G^{(n)}=\set{ \id}$. The least integer $n$ such that $G^{(n)}=\{ \id\}$ is called the \emph{solvability length} of $G$.
\end{defn}

\begin{ejem}\label{ej_solvable}\mbox{}
	\begin{enumerate}
	\item The infinite dihedral group $\dih$ is solvable (see definition \ref{ej_defns_dih_rot_so3}).
	\item The group $\so$ is not solvable.
	\item Any triangular group is solvable, with solvability length at most 3. In particular, the fundamental groups of Inoue surfaces are solvable groups (see example \ref{defn_inuoe}).
	\item Cyclic groups are solvable.
	\item Suspensions of solvable groups are solvable.

	\end{enumerate}
\end{ejem}

\begin{obs}
It follows directly from the definition above that a subgroup $G\subset\PSL$ is solvable if and only if any of its lifts $\mathbf{G}\subset\GL$ is solvable. In the same way, $G$ is triangular if and only if any of its lifts $\mathbf{G}\subset\GL$ is triangular.
\end{obs}

The previous observation allows us to treat indistinctly a complex Kleinian group and any of its lifts in the matrix group $\SL$ when we prove that solvable groups are virtually triangularizable.\\ 

The next result is stated in Section 1.2.1 of \cite{borel}.

\begin{lem}\label{lem_producto_zariski_denso}
Let $\Gamma$ be an algebraic group, then $\Gamma\times\Gamma$ is Zariski dense in $\overline{\Gamma}\times\overline{\Gamma}$.
\end{lem}

\begin{prop}\label{prop_cerradura_solvable}
Let $\Gamma\subset\GLN$ be a subgroup group, if $\overline{\Gamma}$ is the Zariski closure of $\Gamma$, then $\Gamma$ is solvable if and only if $\overline{\Gamma}$ is solvable.
\end{prop}

\begin{proof}
Consider the continuous function (respect to the Zariski topology) $c:\GLN\times\GLN\rightarrow\GLN$ given by $c(g,h)=g^{-1}h^{-1}gh$. Using lemma \ref{lem_producto_zariski_denso} we know that $\Gamma\times\Gamma$ is Zariski dense in $\overline{\Gamma}\times\overline{\Gamma}$. Then $c\parentesis{\Gamma\times\Gamma}$ is Zariski dense in $c\parentesis{\overline{\Gamma}\times\overline{\Gamma}}$, that is, 
	\begin{equation}\label{eq_dem_cerradura_solvable_1}
	\overline{\corchetes{\Gamma,\Gamma}}=\corchetes{\overline{\Gamma},\overline{\Gamma}}.
	\end{equation}
Now, using induction on $k\in\N$, we will verify that
	\begin{equation}\label{eq_dem_cerradura_solvable_2}
	\overline{\Gamma^{(k)}}=\overline{\Gamma}^{(k)}.
	\end{equation}
Using (\ref{eq_dem_cerradura_solvable_1}) we have:
	\begin{align*}
	\overline{\Gamma}^{(1)}&=\corchetes{\overline{\Gamma},\overline{\Gamma}}=\overline{\corchetes{\Gamma,\Gamma}}=\overline{\Gamma^{(1)}}\\
	\overline{\Gamma}^{(2)}&=\corchetes{\overline{\Gamma}^{(1)},\overline{\Gamma}^{(1)}}=\corchetes{\overline{\Gamma^{(1)}},\overline{\Gamma^{(1)}}}=\overline{\corchetes{\Gamma^{(1)},\Gamma^{(1)}}}=\overline{\Gamma^{(2)}}\\
	 &\cdots\\
	\overline{\Gamma}^{(k)}&=\corchetes{\overline{\Gamma}^{(k-1)},\overline{\Gamma}^{(k-1)}}=\corchetes{\overline{\Gamma^{(k-1)}},\overline{\Gamma^{(k-1)}}}=\overline{\corchetes{\Gamma^{(k-1)},\Gamma^{(k-1)}}}=\overline{\Gamma^{(k)}}.
	\end{align*}		
This proves (\ref{eq_dem_cerradura_solvable_2}).\\

Now, assume $\Gamma$ is solvable and that the solvability length of $\Gamma$ is $k$, then $\Gamma^{(k)} =\set{ \id }$ and therefore, 
	$$\overline{\Gamma}^{(k)}=\overline{\Gamma^{(k)}}=\overline{\set{\id}}.$$
Since $\set{\id}$ is closed in the Zariski topology, then $\overline{\set{\id}}=\set{\id}$ and therefore $\overline{\Gamma}^{(k)}=\set{\id}$, which proves that $\overline{\Gamma}$ is a solvable group.\\

Finally, if $\overline{\Gamma}$ is solvable, with solvability length $k$, then $\overline{\Gamma}^{(k)}=\set{\id}$. From this and (\ref{eq_dem_cerradura_solvable_2}) it follows that $\Gamma^{(k)}=\set{\id}$, which proves that $\Gamma$ is solvable, with solvability length $k$.
\end{proof}

Now we list some results and definitions necessary for the rest of this chapter. The following theorem is proven in Chapter 1 of \cite{series}. 

\begin{thm}\label{teo_elem_solvable}
Let $G$ be a discrete group of M\"obius transformations, then $G$ is elemental if and only if it is solvable.
\end{thm}

The next theorem is known as the \emph{topological Tits alternative} (see Theorem 1 of \cite{tits} or Theorem 1.3 of \cite{tits2}).

\begin{thm}\label{teo_tits2}
Let $K$ be a local field and let $\Gamma\subset\text{GL}(n,K)$ be a subgroup. Then, either, $\Gamma$ contains an open solvable group or $\Gamma$ contains a dense free subgroup.  
\end{thm}

Before stating Theorem \ref{teo_punto_fijo_Borel}, we need the following definition.

\begin{defn}
Let $G$ be an algebraic group, $V$ a variety, and let $\alpha:G\times V\rightarrow V$ be an action of the group $G$ in $V$, $(g,x)\mapsto gx=\alpha(g,x)$. One says that $G$ \emph{acts morphically} on $V$ if the action $\alpha$ satisfies the following axioms:
	\begin{enumerate}[(i)]
	\item $\alpha(e,x)=x$, for any $x\in V$, where $e\in G$ is the identity element.
	\item $\alpha(g,hx)=\alpha(gh,x)$ for any $g,h\in G$ and $x\in V$.	
	\end{enumerate}	  
\end{defn}

\begin{ejem}
Clearly $\PSLN$ acts morphically on $\CP^n$.
\end{ejem}

The next theorem is called the \emph{Borel fixed point theorem} (see Theorem 10.4 of \cite{borel}).

\begin{thm}\label{teo_punto_fijo_Borel}
Let $G$ be a connected solvable group acting morphically on a non-empty complete variety $V$. Then $G$ has a fixed point in $V$. 
\end{thm}

The next theorem is proved on \cite{auslander}.

\begin{thm}\label{thm_solubles_fg}
Let $G\subset\GLN$ be a solvable matrix group, then $G$ is finitely generated.
\end{thm}

\section{Solvable groups have an invariant full flag}
\label{sec_solvable_invariant_flag}
In this section we are going to show that solvable subgroups of $\PSL$ have an invariant full flag in $\CP^2$.\\

Recall that Schottky groups are the ``simplest'' examples of Kleinian groups and enjoy very interesting properties (see Chapter 4 of \cite{mumford2002indra}). Unfortunately they are not always realizable in the higher dimensional setting (see \cite{ckg_libro}). For this reason a weaker form of Schottky groups are introduced.

\begin{figure}[H]
\begin{center}
\label{dinamica_schottky_like}
\includegraphics[height=30mm]{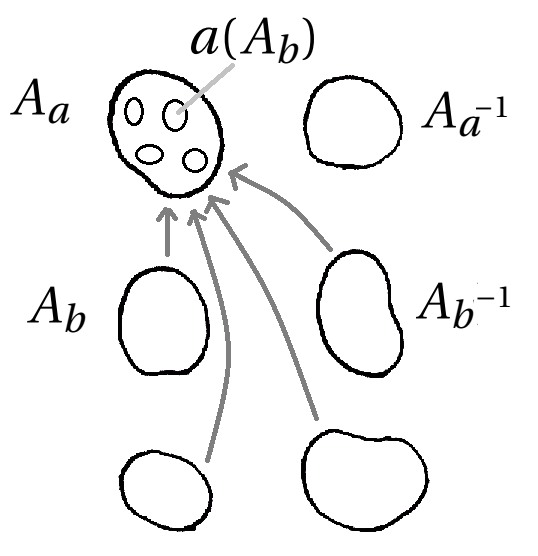}
\caption{Ping-pong dynamics in a Schottky-like group.}
\end{center}
\end{figure}

\begin{defn}\label{defn_schottky_like}
Let $\Sigma\subset \PSLN$ be a finite set which is symmetric (i.e., $a^{-1} \in\Sigma $ for all $a\in \Sigma $) and $A_\sigma= \{A_a\}_{a \in \Sigma }$ a family of compact non-empty pairwise disjoint subsets of $\CP^n$ such that for each $a \in \Sigma$ we have  
	\begin{equation}\label{eq_defn_schottky_like}
	\bigcup_{b\in \Sigma\setminus\{a^{-1}\}}a(A_b) \subset  A_a.
	\end{equation}
The group $\Gamma$ generated by $\Sigma$ is called a \emph{Schottky-like group}. The dynamics determined by (\ref{eq_defn_schottky_like}) are known as \emph{ping-pong dynamics} (see Figure \ref{dinamica_schottky_like}).
\end{defn}

\begin{ejem}
Every Schottky group of $\psl$ acting on $\CP^1$ is a Schottky-like group.
\end{ejem}

As a consequence of the ping-pong dynamics it is possible to deduce that a Schottky-like group are free and discrete.

\begin{lem}\label{l:exfolox}
Let $\gamma_1,\gamma_2\in \PSL$ be loxodromic elements and let 
	$$\tau_1^+,\tau_2^+,\tau_1^-,\tau_2^-\in \QP\setminus \PSL$$
be elements such that all the following conditions hold:
	\begin{itemize}
	\item $\gamma_i^{\pm n} \underset{n \longrightarrow \infty}{\rightarrow} \tau_i^\pm$,
	\item $\text{Im}(\tau_1^+)$ and $\text{Im}(\tau_2^+)$ are distinct points,
	\item $\text{Im}(\tau_2^-)\nsubseteq \kernel(\tau_1^+)$  and $\text{Im}(\tau_1^-)\nsubseteq	\kernel(\tau_2^+)$	
	\end{itemize}
then $\langle \gamma_1,\gamma_2\rangle$ contains a strongly loxodromic element.
\end{lem}

\begin{proof}
	Let us define  $\sigma_n=\gamma_1^n\gamma_2^{-n}$, then it is clear that  
	\begin{align*}
	\sigma_n & \underset{n \rightarrow \infty}{\longrightarrow} \tau_1^+\tau_2^-\\
	\sigma_n^{-1} & \underset{n \rightarrow \infty}{\longrightarrow} \tau_2^+\tau_1^-
	\end{align*}
	to conclude the proof observe that $\text{Im}(\tau_1^+\tau_2^-)=\text{Im}(\tau_1^+)$ and $ \text{Im}(\tau_2^+\tau_1^-)=\text{Im}(\tau_2^+)$.
\end{proof}

\begin{lem} \label{l:hsc}
Let $a,x,w,z\in \C^\ast$, $\theta\in\R$ and $y,u\in \C$ such that $\valorabs{y}+\valorabs{u}\neq 0$ and $\valorabs{a}<1$, then the group
	$$\left 	\langle 
\gamma_1=
	\begin{bmatrix}
	ae^{2\pi i \theta} & 0 &0\\
	0 & a & 0\\
	0 & 0 & a^{-2}e^{-2\pi i \theta} 
	\end{bmatrix},\,
	\gamma_2= 	\begin{bmatrix}
	x & 0 &y\\
	0 & z & u\\
	0 & 0 & w
	\end{bmatrix}
	\right 
	\rangle$$
is not discrete.
	\end{lem}
	
\begin{proof}
A straightforward calculation shows that
$$\gamma_1^n\gamma_2\gamma_1^{-n}\gamma_2=
	\begin{bmatrix}
1 & 0 & -\frac{y}{w} \parentesis{1-a^{3 n}e^{4\pi i \theta}} \\
0 & 1 & -\frac{u}{w} \parentesis{1-a^{3 n} e^{2\pi i \theta}} \\
0 & 0 & 1 \\
\end{bmatrix}
	\xymatrix{
	\ar[r]_{n \rightarrow \infty}&} 
\begin{bmatrix}
1 & 0 & \frac{-y }{w} \\
0 & 1 & \frac{-u}{w} \\
0 & 0 & 1 \\
\end{bmatrix}$$
which shows the assertion.
	\end{proof}

\begin{lem}  \label{l:5pg}
Let $\Gamma\subset \PSL$ be a complex Kleinian group with strongly irreducible action on $\CP^2$ and $\ell\subset\CP^2$ a complex line. Then $\mathcal{L}=\SET{\gamma(l)}{\gamma\in \Gamma}$ contains infinite lines in general position.
\end{lem}

\begin{proof}
Since the action of $\Gamma$ on $\CP^2$ is strongly irreducible, we deduce that $\mathcal{L}$ cannot have exactly two lines in general position (otherwise $\bigcap\mathcal{L}$ is a fixed point of $\Gamma$ and therefore, a non-empty invariant space for the action of $\Gamma$, see Figure \ref{fig_dos_lineas_pos_gral}). 

\begin{figure}[H]
\begin{center}
\label{fig_dos_lineas_pos_gral}
\includegraphics[height=30mm]{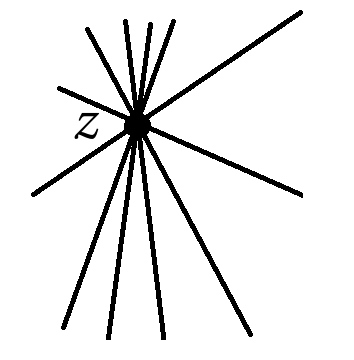}
\caption{Two lines in general position.}
\end{center}
\end{figure}

So $\mathcal{L}$ must contain at least three lines in general position. On the other hand, $\CP^2\setminus \cup \mathcal{L}$ is $\Gamma$-invariant, then, $\CP^2\setminus \cup \mathcal{L}\subset \Eq(\Gamma)$ (see Theorem 3.5 in \cite{bcn2011}). As a consequence of the main Theorem in \cite{barrera2010estimates}, the number of lines contained in $\CP^2\setminus \Eq(\Gamma)$ lying in general position is 1,2,3,4 or infinite, where the case of finite lines implies that the action of $\Gamma$ on $\CP^2$ is reducible, therefore $\mathcal{L}$ contains infinite lines in general position.
\end{proof}

\begin{lem}
Let $\Gamma\subset \PSL$ be an irreducible complex Kleinian group, then $\Gamma$ contains a strongly loxodromic element.
\end{lem}

\begin{proof}
By Lemma 4.1 in \cite{ppar}, we know $\Gamma$ contains a loxodromic element $\gamma$. Let $\tau_+,\tau_-\in \QP\setminus \PSL$  be such that 
	$$\gamma^{\pm n} \xymatrix{\ar[r]_{n \rightarrow \infty}&} \tau_\pm $$ 
as pseudo projective transformations. Since $\gamma$ is loxodromic we can assume $\text{Im}(\tau_+)$ is a point and $\text{Im}(\tau_+)\nin \kernel(\tau_+)$. Let us assume that $\gamma$ is not strongly loxodromic. Then, we have two possibilities:\\

\begin{itemize}
\item $\gamma$ is a complex homothety. In this case $\text{Im}(\tau_+)=\kernel(\tau_-)$ and $\text{Im}(\tau_-)=\kernel(\tau_+) $. Since $\Gamma$ is irreducible there is $\tau\in\Gamma$ such that  $\tau(\text{Im}(\tau_+))\neq  \text{Im}(\tau_+)$. Applying Lemma \ref{l:hsc} to $\prodint{ \gamma,  \tau\gamma \tau^{-1}}$ we deduce $\tau(\kernel(\tau_+))\neq \kernel(\tau_1)$. To conclude the proof in this case just apply Lemma \ref{l:exfolox} to $\gamma$ and $\tau\gamma \tau$.

\item $\gamma$ is loxo-parabolic. In this case $\fix(\gamma)=\{p_+,p_-\}$ and $\gamma$ leaves invariant exactly 2 complex lines, say $\ell_+$ and $\ell_-$. Let $\rho_+,\rho_- \in \QP\setminus \PSL$ such that $\gamma^{\pm n} 	\xymatrix{
	\ar[r]_{n \rightarrow \infty}&} \rho_\pm$, thus we can assume $\text{Im}(\rho_\pm)=p_\pm$, $\kernel(\rho_\pm )=\ell_\pm$ and $\ell_-=\langle p_+,p_-\rangle$. Applying Lemma \ref{l:5pg} we conclude that there is $\tau\in \Gamma$ such that $\tau \ell_-\cap \ell_-\cap\{p_+,p_-\}=\emptyset $. To conclude, apply Lemma \ref{l:exfolox} to 
	\begin{itemize}
	\item $\prodint{\tau \gamma^{-1}\tau^{-1}, \gamma}$, if $\tau p_+\in \ell_-$,
	\item $\prodint{\tau \gamma\tau^{-1}, \gamma^{-1}}$, if $\tau p_+\in \ell_+$,
	\item $\prodint{\tau \gamma\tau^{-1}, \gamma}$ in any other case.
	\end{itemize}
\end{itemize}
\end{proof}

Now the proof of the following corollary is straightforward.

\begin{cor}
Let $\Gamma\subset \PSL$ be an irreducible complex Kleinian group, then $\CG(\Gamma)$ is the closure of repelling fixed point of strongly loxodromic elements.
\end{cor}

\begin{lem}\label{lem_irred_contiene_schottky}
Let $\Gamma\subset\PSL$ be an irreducible complex Kleinian group, then $\Gamma$ contains a Schottky-like group. 
\end{lem}

\begin{proof}
Let  $\gamma_1\in \Gamma $ be a strongly loxodromic element and  $x_a,x_s,x_r$ be respectively its attracting, saddle and repulsive point. Applying Lemma \ref{l:5pg} to $\ell=\linproy{x_a,x_r}$ we deduce there is $\tau\in \Gamma$ such that 
 $$\tau\parentesis{\ell}\cap \ell\cap  \set{x_a,x_r}=\emptyset$$ 

\begin{figure}[H]
\begin{center}
\label{fig_moviendo_linea_FL}
\includegraphics[height=30mm]{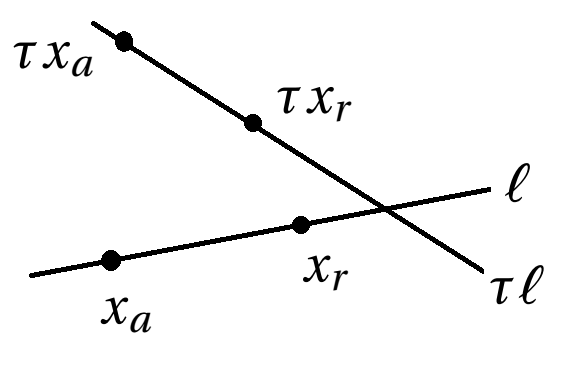}
\caption{The element $\tau$ moving the line $\ell$.}
\end{center}
\end{figure} 
 
Thus $\gamma_2=\tau \gamma_1 \tau^{-1}$ is a strongly loxodromic element and $\tau x_a,\tau x_s,\tau x_r$ are its attracting, saddle and repulsive points respectively. Clearly, $\valorabs{\set{x_a,x_r, \tau x_a,\tau x_r}}=4$. 
Now, it is trivial that we can find $N\in \N$ and four mutually disjoint closed balls $B_1^+,B_1^-, B_2^+$ and $B_2^-$ such that $\text{Int}(B_i^+)\cap \fix(\gamma_i)$ (resp. $\text{Int}(B_i^-)\cap \fix(\gamma_i)$) is the attracting (resp. repelling) fixed point of $\gamma_i$ (see Figure \ref{fig_bolas_schottky_type}). Furthermore, 
\begin{align*}
\gamma_i^{ \pm N}(B_j^{\pm}) &\subset \text{Int}(B_i^\pm ),\;\;\; i\neq j\\
\gamma_i^{ \pm N}(B_i^{\pm}) &\subset \text{Int}(B_i^\pm ) \\
\end{align*}
thus $\langle \gamma_1^N,\gamma_2^N \rangle$ is a Schottky-like group of rank two, which concludes the proof.

\begin{figure}[H]
\begin{center}
\label{fig_bolas_schottky_type}
\includegraphics[height=30mm]{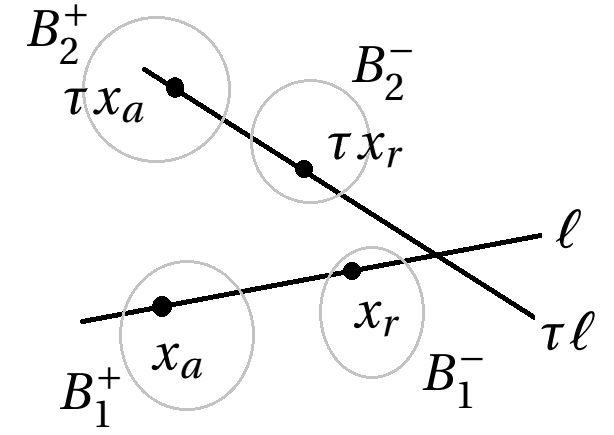}
\caption{The balls $B_1^+,B_1^-, B_2^+,B_2^-$ around the points $x_a,x_r, \tau x_a,\tau x_r$.}
\end{center}
\end{figure} 

\end{proof}

This last lemma shows that complex Kleinian groups with irreducible action on $\CP^2$ contain a free subgroup of rank at least 2. As a consequence of the topological Tits alternative (Theorem \ref{teo_tits2}), they cannot contain solvable groups. Thus, we have the following corollary.

\begin{cor}
Let $\Gamma\subset\PSL$ be a solvable complex Kleinian group, then $\Gamma$ contains a finite index subgroup leaving a full flag invariant.
\end{cor}

The converse is also true, if $\Gamma\subset\PSL$ is a group leaving a full flag inva-riant, up to conjugation, we can assume that the invariant flag is $\set{e_1,\linproy{e_1,e_2}}$. Then $\Gamma$ is upper triangular and therefore, solvable. This proof the following dynamic characterization of solvable groups.

\begin{cor}\label{cor_solubles_bandera_invariante}
Let $\Gamma\subset\PSL$ be a discrete group. $\Gamma$ is solvable if and only if $\Gamma$ contains a finite index subgroup leaving a full flag invariant.
\end{cor}

\section{Solvable groups are virtually triangularizable}\label{sec_solvables_triangularizables}

In this section we will prove one of the main theorems of this chapter (Theorem \ref{teo_solvable_triangularizables}). We will give two proofs, one is algebraic and the other one is dynamical.\\ 

The algebraic proof is the most appropriate because it is valid for solvable subgroups of $\PSLN$ acting on $\CP^n$ for any $n\geq 2$. On the other hand, the dynamical proof is only valid for the case $n=2$, since it strongly depends on what is known for the classic case of subgroups of $\psl$. However, we give this proof to make a case about the importance of Schottky groups in this dimension.\\

Before stating the main theorem of this section, we prove several auxiliary lemmas that will be needed.

\begin{lem}\label{lem_existen_lox_ajenos}
Let $\Sigma\subset\psl$ be a non elemental subgroup (discrete or non-discrete) acting on $\CP^1$ without fixed points. Then there are two loxodromic elements $\gamma_1,\gamma_2\in\Sigma$ such that 
	\begin{equation}\label{eq_lem_existen_lox_ajenos_1}
	\text{Fix}(\gamma_1)\cap \text{Fix}(\gamma_2)=\emptyset.
	\end{equation}
\end{lem}

\begin{proof}
Let's assume that $\Sigma$ is discrete, we denote by $\Lambda(\Sigma)$ the limit set of $\Sigma$. If $\Sigma$ is not discrete, the proof is similar.\\

Since $\Sigma$ is a non-elemental group then $\valorabs{\Lambda(\Sigma)}=\infty$ and, $\Lambda(\Sigma)$ is the closure of all the fixed points of loxodromic elements. Therefore, $\Sigma$ has an infinite amount of loxodromic elements.\\

Suppose that there are not two loxodromic elements in $\Sigma$ such that (\ref{eq_lem_existen_lox_ajenos_1}) holds, then we have two cases:
	\begin{enumerate}[(i)]
	\item All loxodromic elements of $\Sigma$ share the same fixed points, that is
		$$\text{Fix}(\gamma)=\set{a,b}$$
	for any loxodromic element $\gamma\in\Sigma$, and for some points $a,b\in\CP^1$. Then
		$$\Lambda(\Sigma)=\overline{\set{a,b}}=\set{a,b},$$
	contradicting that $\Sigma$ is non-elemental.
	\item All loxodromic elements of $\Sigma$ share one fixed point, that is,  
		\begin{equation}\label{eq_lem_existen_lox_ajenos_2}
		\text{Fix}(\gamma)=\set{p,q_\gamma}
		\end{equation}			
	for any loxodromic element $\gamma\in\Sigma$, and for some point $p\in\CP^1$.\\
	
	Since $\Sigma$ acts on $\CP^1$ without fixed points then $\exists$ $\tau\in\Sigma$ such that 
		\begin{equation}\label{eq_lem_existen_lox_ajenos_3}
		\tau(p)\neq p.
		\end{equation}				
	Let $\gamma\in\Sigma$ be a loxodromic element such that
		\begin{equation}\label{eq_lem_existen_lox_ajenos_4}
		\tau(q_\gamma)\neq p.
		\end{equation} 
	If there is no loxodromic element satisfying (\ref{eq_lem_existen_lox_ajenos_4}), then $\tau(q_\gamma)=p$ for any loxodromic $\gamma\in\Sigma$, contradicting that $\tau\in\psl$ is injective.\\
	
	Consider the conjugation $\tau\gamma\tau^{-1}\in\Sigma$, it is a loxodromic element of $\Sigma$ with fixed points
		$$\text{Fix}(\tau\gamma\tau^{-1})=\set{\tau(p),\tau(q_\gamma)}.$$
	But, as a consequence of (\ref{eq_lem_existen_lox_ajenos_3}) and (\ref{eq_lem_existen_lox_ajenos_4}) we have that $p\nin\text{Fix}(\tau\gamma\tau^{-1})$. This contradicts the assumption (\ref{eq_lem_existen_lox_ajenos_2}).
	\end{enumerate}
Therefore, none of the previous two cases can occur. This proves the lemma.
\end{proof}

\begin{lem}\label{lem_suf_grande_schottky}
Let $\gamma_1,\gamma_2$ be two loxodromic elements of a subgroup $\Sigma\subset\psl$. If 
	\begin{equation}\label{eq_lem_suf_grande_schottky_1}
	\text{Fix}(\gamma_1)\cap \text{Fix}(\gamma_2)=\emptyset
	\end{equation}			
then there exists a large enough integer $n>0$ such that the subgroup $\prodint{\gamma_1^n,\gamma_2^n}\subset\Sigma$ is a Schottky group.
\end{lem}

\begin{proof}
Let's denote the fixed points of $\gamma_i$ by $\text{Fix}(\gamma_i)=\set{a_i,b_i}$ for $i=1,2$. We can take a large enough $n>0$ such that there are open balls $U_1,V_1,U_2,V_2\subset\CP^1$ around the points $a_1$, $b_1$, $a_2$, $b_2$ respectively such that $\prodint{\gamma_1^n,\gamma_2^n}\subset\Sigma$ is a Schottky group. We can take this balls as a consequence of the hypothesis (\ref{eq_lem_suf_grande_schottky_1}) and the fact that $\gamma_1,\gamma_2$ are loxodromic elements. 
\end{proof}

Now we state and prove the main theorem of this section.

\begin{thm}\label{teo_solvable_triangularizables}
Let $G\subset\GL$ be a discrete solvable subgroup then $G$ is virtually triangularizable. 
\end{thm}

We will give the two proofs in the following two subsections.

\subsection*{Dynamical proof}

First we give the dynamical proof.

\begin{proof}
Let $G\subset\GL$ be a solvable discrete subgroup. Since $G$ is solvable, it is not strongly irreducible; if it were strongly irreducible it would contain a free subgroup of rank at least 2, contradicting that its solvability (see Lemma \ref{lem_irred_contiene_schottky}).\\

Since $G$ is not strongly irreducible, there exists a non-empty proper subspace $p\subset\CP^2$ (i.e. $p$ is a point or a line) such that the orbit $G(p)$ is finite. Assume first that $p$ is a point.\\

Let $G_1=\text{Isot}(p,G)\subset G$, then $G_1$ is finite index subgroup of $G$, which acts on $\CP^2$ fixing $p$. Then $G_1$ is weakly semi-controllable, let $\ell\subset\CP^2$ be a line such that $p\nin\ell$ and let 
	$$\Pi=\Pi_{p,\ell}:\PSL\rightarrow\text{Aut}(\ell)\cong\psl$$
the group homomorphism defined in Subsection \ref{subsec_groups_no_discretos}. Denote $\Sigma = \Pi(G_1)\subset\psl$, the subgroup $\Sigma$ can be either discrete or non-discrete.\\

If $\Sigma$ is discrete then, by Theorem \ref{teo_elem_solvable}, $\Sigma$ is elemental and therefore one of the following cases hold:
	\begin{enumerate}\setcounter{enumi}{0}
	\item $\Lambda\parentesis{\Sigma}=\emptyset$, or
	\item $\valorabs{\Lambda\parentesis{\Sigma}}\in\set{1,2}$.
	\end{enumerate}
If $\Sigma$ is non-discrete, it can be elemental or non-elemental. If it is elemental, there are two cases:
	\begin{enumerate}\setcounter{enumi}{2}
	\item $\Lambda_{Gr}\parentesis{\Sigma}=\emptyset$, or
	\item $\valorabs{\Lambda_{Gr}\parentesis{\Sigma}}\in\set{1,2}$.
	\end{enumerate}
If $\Sigma$ is non-elemental, there are two cases (see Proposition 13.15 of \cite{cs2014}):
	\begin{enumerate}\setcounter{enumi}{4}
	\item $\Lambda_{Gr}\parentesis{\Sigma}=\mathbb{S}^1$, or
	\item $\Lambda_{Gr}\parentesis{\Sigma}=\CP^1$.
	\end{enumerate}	

Now we analyze each of the previous 6 cases.
	
	\begin{enumerate}
	\item Observe that all points in $\CP^1$ have finite $G$-orbit. If there would be a point with infinite $G$-orbit, and since $\CP^1$ is compact, there would be a subsequence of $\set{\gamma q}_{\gamma\in\Sigma}$ converging to a point in $\CP^1$, and this point would have to be in $\Lambda(\Sigma)$, contradicting that $\Lambda(\Sigma)=\emptyset$.\\
	
	Let $x\in\CP^1$ be an arbitrary point and let us define
		$$G_0=\Pi^{-1}\parentesis{\text{Isot}(x,\Sigma)}.$$
	Clearly $G_0$ is a subgroup of $G$, furthermore, it has finite index in G, since the cosets of $\text{Isot}(x,\Sigma)$ in $\Sigma$ correspond to the $\Sigma$-orbit of $x$, which is finite, as we already stated.\\
	
	Observe that the elements of $\text{Isot}(x,\Sigma)$ are elements $\gamma\in\Sigma$ such that $\gamma(x)=x$, but, by the definition of $\Pi$, the elements of $\Pi^{-1}\parentesis{\text{Isot}(x,\Sigma)}$ are elements leaving invariant the line $\overleftrightarrow{p,x}$, besides fixing the point $p$. This means that $G_0$ leaves invariant the full flag $\set{p,\overleftrightarrow{p,x}}$.\\
	
	We have proved that $G$ has a finite index subgroup which is triangularizable.  
	\item Let $x\in\Lambda(\Sigma)$, define
		$$G_0 = \Pi^{-1}\parentesis{\text{Isot}(x,\Sigma)}.$$
		Observe that $G_0$ is a subgroup of $G$. Furthermore, $\text{Isot}(x,\Sigma)$ has finite index in $\Sigma$, since the cosets of $\text{Isot}(x,\Sigma)$ are the points of $\Lambda(\Sigma)$, which are finite.\\
		
		By the same argument used in case 1, $G_0$ leaves invariant the full flag $\set{p,\overleftrightarrow{p,x}}$. Which means that $G$ has a triangularizable subgroup of finite index.
	\item Since $\Sigma$ is not discrete, it is an infinite group. Since $\Lambda_{Gr}(\Sigma)=\emptyset$, then $\Sigma$ is a  subgroup of $\dih$ or $\so$ (see Theorem 1.14 of \cite{cs2014}). We have the following two subcases:
		\begin{enumerate}
		\item If $\Sigma$ is a subgroup of $\dih$, then the point $\set{0}$ has a finite $\Sigma$-orbit, analogously to the previous cases, we define the subgroup of $G$,
			$$G_0=\Pi^{-1}\parentesis{\text{Isot}(0,\Sigma)},$$
		which has finite index in $G$ and leaves a full flag invariant.
		\item If $\Sigma$ is a subgroup of $\so$, it can be dense or $\Sigma$ is finite (see Proposition \ref{prop_opciones_conj_lim_greenberg}). We consider both cases:
     \begin{enumerate}
      \item A similar proof of Lemma 4.3 of \cite{ppar} dismisses that $\Sigma$ could contain a dense subgroup of $\so$. 
      \item If $\Sigma=\Pi(\Gamma)$ is finite then $\Gamma_0:=\kernel(\Gamma)$ has finite index in $\Gamma$ and, by the definition of $\Pi$, leaves the full flag $\set{e_1,\linproy{e_1,e_2}}$ invariant.
      \end{enumerate}

		\end{enumerate}
	
	\item We use the same argument, used in case 2, on the limit set $\Lambda_{Gr}(\Sigma)$.			
	
	\item We have two cases:
		\begin{enumerate}
		\item If $\Sigma$ has a global fixed point, we denote it by $x\in\Lambda_{Gr}(\Sigma)$ and define 
			$$G_0 = \Pi^{-1}\parentesis{\Sigma} = G_1.$$
		Analogously to the other cases, this finite index subgroup of $G$ leaves invariant the full flag $\set{p,\overleftrightarrow{p,x}}$. Therefore, $G$ is virtually triangularizable.
		\item If $\Sigma$ doesn't have a global fixed point, and since the limit set $\Lambda_{Gr}(\Sigma)=\Ss^1$ is the closure of fixed points of $\Sigma$, these fixed points are dense in $\Ss^1$. Therefore as consequence of Lemma \ref{lem_existen_lox_ajenos}, we can choose two loxodromic elements $\gamma_1,\gamma_2\in\Sigma$ such that
			$$\text{Fix}(\gamma_1)\cap \text{Fix}(\gamma_2)=\emptyset.$$
		Therefore using Lemma \ref{lem_suf_grande_schottky}, we can guarantee that, for a sufficiently large $n$, the subgroup $\prodint{\gamma_1^n,\gamma_2^n}\subset\Sigma$ is a Schottky group. This means that $\Sigma$ contains a free subgroup of rank 2, contradicting that $\Sigma$ is solvable (see Theorem \ref{teo_tits2}).
		\end{enumerate}
	
	\item As in the last case, we have the same two cases:
		\begin{enumerate}
		\item If $\Sigma$ has a global fixed point, this point determines the invariant full flag, in the same way as in case (a) of the previous case.
		\item If $\Sigma$ doesn't have a global fixed point, Lemmas \ref{lem_existen_lox_ajenos} and \ref{lem_suf_grande_schottky} guarantee that $\Sigma$ contains a Schottky subgroup and therefore that $\Sigma$ is solvable.
		\end{enumerate}
\end{enumerate}
\end{proof}

\subsection*{Algebraic proof}

Now we give the algebraic proof of Theorem \ref{teo_solvable_triangularizables}. The proof is given for discrete solvable subgroups $G\subset\GL$ but the same arguments are valid for the general case of discrete solvable subgroups $G\subset\GLN$.\\ 

\begin{proof}
Let $G\subset\GL$ be a discrete and solvable subgroup. Since $G$ is discrete, $G$ is not connected, however the Zariski closure $\overline{G}$ is a solvable subgroup (Proposition \ref{prop_cerradura_solvable}), which acts morphically on $\CP^2$. Denote by $G_0$ the connected component of $\overline{G}$ which contains the identity, then $G_0$ is a solvable connected group, furthermore it is a finite index subgroup of $\overline{G}$. By Theorem \ref{teo_punto_fijo_Borel}, $G_0$ has a global fixed point in $\CP^2$.\\

Up to conjugation by an element of $\GL$, we can assume that this fixed point is $e_1=\corchetes{1:0:0}$ and therefore every element of $G_0$ has the form

	$$\parentesis{\begin{array}{ccc}
	a_{11} & a_{12} & a_{13} \\ 
	0 & a_{22} & a_{23} \\ 
	0 & a_{32} & a_{33}
	\end{array}}.$$
Let $a: G_0\rightarrow a(G_0)\subset GL(2,C)$ be the group morphism given by
	$$\parentesis{\begin{array}{ccc}
	a_{11} & a_{12} & a_{13} \\ 
	0 & a_{22} & a_{23} \\ 
	0 & a_{32} & a_{33}
	\end{array}} \overset{a}{\mapsto} \parentesis{\begin{array}{cc} 
	a_{22} & a_{23} \\ 
    a_{32} & a_{33}
	\end{array}}.$$	
Since $G_0$ is solvable and $a$ is a suprajective group morphism, then $a(G_0)$ is solvable. Repeating the same argument we applied before to $G\subset\GL$, we now have a connected solvable subgroup with finite index $H_0\subset H=a(G_0)$ acting with a fixed point on $\CP^1$. Therefore, up to conjugation, every element of $H_0$ has the form
	
	$$\parentesis{\begin{array}{cc} 
	a_{22} & a_{23} \\ 
	0 & a_{33}
	\end{array}}.$$

Taking the inverse image $a^{-1}(H_0)$ we get a finite index subgroup of $G$ with the form
	 
	$$\parentesis{\begin{array}{ccc}
	a_{11} & a_{12} & a_{13} \\ 
	0 & a_{22} & a_{23} \\ 
	0 & 0 & a_{33}
	\end{array}}.$$

This proofs that $\overline{G}$ is virtually triangularizable and, since $G$ is a finite index subgroup of $\overline{G}$, then $G$ is virtually triangularizable.
\end{proof}

Theorem \ref{teo_solvable_triangularizables} states that every solvable subgroup $\Gamma\subset\PSL$ has a finite index subgroup $\Gamma_0\subset\Gamma$ such that, up to conjugation, is upper triangular. Since $\KulL(\Gamma_0)=\KulL(\Gamma)$ (see Proposition 3.6 of \cite{bcn16}), we can restrict our attention to the upper triangular subgroups of $\PSL$. We will study first the commutative triangular groups in Section \ref{sec_commutative_triangular} and then the non-commutative triangular groups in Section \ref{sec_non_commutative_triangular}.

\section{Commutative triangular groups}\label{sec_commutative_triangular}

In this section we describe the commutative triangular groups. The following notation and definitions will be used through the rest of the chapter.\\

We denote the upper triangular elements of $\PSL$ by

	$$U_{+}=\SET{\corchetes{\begin{array}{ccc}
	a_{11} & a_{12} & a_{13} \\
	0 & a_{22} & a_{23} \\
	0 & 0 & a_{33}
	\end{array}}}{a_{11}a_{22}a_{33}=1,\text{ }a_{ij}\in\C}.$$

Now we define the group morphisms $\lambda_{12},\lambda_{23},\lambda_{13}:(U_{+},\cdot)\rightarrow(\C^{\ast},\cdot)$ which are given by
	\begin{align*}
	\lambda_{12}\parentesis{\corchetes{a_{ij}}}&= a_{11} a^{-1}_{22}\\
	\lambda_{23}\parentesis{\corchetes{a_{ij}}}&= a_{22} a^{-1}_{33}\\
	\lambda_{13}\parentesis{\corchetes{a_{ij}}}&= a_{11} a^{-1}_{33}.	
	\end{align*}

To simplify the notation we will write $\text{Ker}\parentesis{\lambda_{ij}}$ instead of $\text{Ker}\parentesis{\lambda_{ij}}\cap \Gamma$ for subgroups $\Gamma\subset U_{+}$. We also define the projections
	$$\pi_{kl}\parentesis{\corchetes{ a_{ij}}}= a_{kl}.$$
	
Whenever we have a discrete subgroup $\Gamma\subset U_+$, we have a finite index torsion free subgroup $\Gamma'\subset\Gamma$ such that $\lambda_{12}(\Gamma')$ and $\lambda_{23}(\Gamma')$ are torsion free groups as well (see Lemma 5.8 of \cite{ppar}). This subgroup and the original group satisfy $\KulL(\Gamma)=\KulL(\Gamma')$ (see Proposition 3.6 of \cite{bcn16}). Therefore we can assume for the rest of the chapter that all discrete subgroups $\Gamma\subset U_+$ are torsion free.\\

The following immediate result will be used often.

\begin{prop}\label{prop_grupos_conm_sintorsion_Zk}
If $\Gamma\subset\PSL$ is a torsion free, commutative subgroup then
	$$\Gamma\cong\Z^r$$
where $r=\text{rank}(\Gamma)$.
\end{prop}

The following result (see Lemma 5.11 of \cite{ppar}) describe the form of the upper triangular commutative subgroups of $\PSL$.

\begin{lem}\label{lem_7casos_conmutativos}
Let $\Gamma\subset U_+$ be a commutative group, then there is a matrix $\tau\in\SL$ such that one of the following cases occurs:
	\begin{enumerate}
	\item Each element of $\tau\Gamma\tau^{-1}$ has the form
		$$\parentesis{\begin{array}{ccc}
		\alpha^{-2} & 0 & 0 \\
		0 & \alpha & \beta \\
		0 & 0 & \alpha
		\end{array}}.
		$$		
	\item Each element of $\tau\Gamma\tau^{-1}$ has the form
		$$\parentesis{\begin{array}{ccc}
		\alpha & 0 & \beta \\
		0 & \alpha^{-2} & 0 \\
		0 & 0 & \alpha
		\end{array}}.
		$$
	\item The group $\tau\Gamma\tau^{-1}$ is diagonal.
	\item Each element of $\tau\Gamma\tau^{-1}$ has the form
		$$\parentesis{\begin{array}{ccc}
		\alpha & 0 & \beta \\
		0 & \alpha & \gamma \\
		0 & 0 & \alpha
		\end{array}},
		$$
	where $\alpha$ is a cubic root of the unity.
	\item Each element of $\tau\Gamma\tau^{-1}$ has the form
		$$\parentesis{\begin{array}{ccc}
		\alpha & \beta & \gamma \\
		0 & \alpha & 0 \\
		0 & 0 & \alpha
		\end{array}},
		$$
	where $\alpha$ is a cubic root of the unity.
	\item Each element of $\tau\Gamma\tau^{-1}$ has the form
		$$\parentesis{\begin{array}{ccc}
		\alpha & \beta & \gamma \\
		0 & \alpha & \mu \\
		0 & 0 & \alpha
		\end{array}},
		$$
	where $\alpha$ is a cubic root of the unity.
	\end{enumerate}
In all cases, $\beta,\gamma,\mu\in\C$ and $\alpha\in\C^\ast$.
\end{lem} 

\begin{obs}
The previous lemma states that if we have an upper triangular commutative group, the group has one of the 6 described forms. However, the opposite doesn't necessarily holds: if we have a group of the sixth form, it is not commutative in general. 
\end{obs}

Observe that cases 4, 5 and 6 stated in Lemma \ref{lem_7casos_conmutativos} are purely parabolic, they have been already studied in \cite{ppar}. Cases 1, 2 and 3 can be purely parabolic or they can have loxodromic elements, therefore we will only study cases 1, 2 and 3, assuming that they contain loxodromic elements.\\ 

So, now we describe, in each of the 3 cases, how the groups should be in order to be commutative and discrete. We will also describe the Kulkarni limit set in each case. 

\subsection{Case 1}

This case will be treated in greater detail since the other cases will have simila-rities with this one. First we describe the form of the groups of this case.\\

First, we describe the explicit form of these groups.

\begin{prop}\label{prop_conmutativo_c1_descripcion}
Let $\Gamma\subset U_+$ be a commutative subgroup such that each element of $\Gamma$ has the form
	$$\corchetes{\begin{array}{ccc}
		\alpha^{-2} & 0 & 0 \\
		0 & \alpha & \beta \\
		0 & 0 & \alpha
		\end{array}},
		$$
for some $\alpha\in\C^\ast$ and $\beta\in\C$. Then there exists an additive subgroup $W\subset(\C,+)$ and a group morphism $\mu:(W,+)\rightarrow(\C^{\ast},\cdot)$ such that
	$$
	\Gamma=\Gamma_{W,\mu}=\SET{\corchetes{\begin{array}{ccc}
	\mu(w)^{-2} & 0 & 0 \\
	0 & \mu(w) & w\mu(w) \\
	0 & 0 & \mu(w)
	\end{array}}}{w\in W}.	
	$$
\end{prop}

\begin{proof}
Let $\zeta:(\Gamma,\cdot)\rightarrow(\C,+)$ be the group homomorphism given by 
	$$\corchetes{\alpha_{ij}}\overset{\zeta}{\mapsto}\alpha_{23}\alpha^{-1}_{33}.$$
Clearly, we have $\text{Ker}(\zeta)=\set{\id}$. Thus we can define the group homomorphism $\mu:(\zeta(\Gamma),+)\rightarrow(\C^{\ast},\cdot)$ as
	$$x\overset{\mu}{\mapsto}\pi_{22}\parentesis{\zeta^{-1}(x)}.$$
Define the additive group $W=\zeta\parentesis{\Gamma}$. It is straight forward to verify that 
	$$
	\Gamma=\SET{\corchetes{\begin{array}{ccc}
	\mu(w)^{-2} & 0 & 0 \\
	0 & \mu(w) & w\mu(w) \\
	0 & 0 & \mu(w)
	\end{array}}}{w\in W}.	
	$$
\end{proof}

For a commutative group $\Gamma_{W,\mu}\subset U_+$ with the form given by the previous proposition, and for $w\in W$, we denote 
	$$\gamma_w=\corchetes{\begin{array}{ccc}
	\mu(w)^{-3} & 0 & 0 \\
	0 & 1 & w \\
	0 & 0 & 1
	\end{array}}\in \Gamma_{W,\mu}.$$

\begin{prop}\label{prop_conmutativo_caso1_rangoG_rangoW}
Let $\Gamma=\Gamma_{W,\mu}\subset U_+$ be a commutative subgroup with the form given by Proposition \ref{prop_conmutativo_c1_descripcion}. If $\text{rank}(\Gamma)=r$, then $\text{rank}(W)=r$.
\end{prop}

\begin{proof}
Let $r=\text{rank}(\Gamma)$ and $\gamma_1,...,\gamma_r\in\Gamma$ such that $\Gamma=\prodint{\gamma_1,...,\gamma_r}$. Let $w_1,...,w_r\in W$ such that $\gamma_j=\gamma_{w_j}$, for $j=1,...,r$. Let $w\in W$ and consider $\gamma_w\in \Gamma$, then there exist $n_1,...,n_r\in\Z$ such that $\gamma=\gamma_1^{n_1}\cdots \gamma_r^{n_r}$, that is
	\begin{equation}\label{eq_dem_prop_conmutativo_caso1_rangoG_rangoW_1}
	\gamma=\corchetes{\begin{array}{ccc}
	\mu(w_1)^{-3n_1}...\mu(w_r)^{-3n_r} & 0 & 0 \\
	0 & 1 & n_1 w_1+...+n_r w_r \\
	0 & 0 & 1
	\end{array}}.
	\end{equation}
On the other hand,
	\begin{equation}\label{eq_dem_prop_conmutativo_caso1_rangoG_rangoW_2}
	\gamma=\corchetes{\begin{array}{ccc}
	\mu(w) & 0 & 0 \\
	0 & 1 & w \\
	0 & 0 & 1
	\end{array}}.
	\end{equation}
From (\ref{eq_dem_prop_conmutativo_caso1_rangoG_rangoW_1}) and (\ref{eq_dem_prop_conmutativo_caso1_rangoG_rangoW_2}), it follows
	$$w=n_1 w_1+...+n_r w_r.$$
This means that $W=\prodint{w_1,...,w_r}$ and therefore, $\text{rank}(W)\leq r$. To prove that $\text{rank}(W)= r$, assume without loss of generality that $w_r=n_1 w_1 +...+n_{r-1} w_{r-1}$, then
	$$\gamma_r=\corchetes{\begin{array}{ccc}
	\mu(w_1)^{-3n_1}...\mu(w_{r-1})^{-3n_{r-1}} & 0 & 0 \\
	0 & 1 & n_1 w_1+...+n_{r-1} w_{r-1} \\
	0 & 0 & 1
	\end{array}}=\gamma_1^{n_1}\cdots \gamma_r^{n_{r-1}}$$
contradicting that $\text{rank}(\Gamma)=r$. This completes the proof.
\end{proof}

\begin{lem}\label{lem_caso1_disc_disc}
Let $W\subset\C$ be a discrete additive subgroup and $\mu:(W,+)\rightarrow(\C^{\ast},\cdot)$ a group morphism. Then $\Gamma_{W,\mu}$ is a discrete subgroup of $\PSL$.
\end{lem}

\begin{proof}
Suppose that $\Gamma=\Gamma_{W,\mu}$ is not discrete. Then there is a sequence of distinct elements $\set{g_n}\subset \Gamma$ converging to the identity in $\PSL$. Denote this sequence by
	$$g_n=\corchetes{\begin{array}{ccc}
	\mu(w_n)^{-2} & 0 & 0 \\
	0 & \mu(w_n) & w_n\mu(w_n) \\
	0 & 0 & \mu(w_n)
	\end{array}}=\corchetes{\begin{array}{ccc}
	\mu(w)^{-3} & 0 & 0 \\
	0 & 1 & w_n \\
	0 & 0 & 1
	\end{array}},$$
for some sequence $\set{w_n}\subset W$. Then $w_n\rightarrow 0\in W$, but this cannot occur because $W$ is discrete, unless $w_n=0$ for all $n$, which implies that $\mu(w_n)=1$ for all $n$, since $\mu$ is a group morphism. Then $\set{g_n}=\set{id}$ is a constant sequence, contradicting that it is a sequence of distinct elements. This proves that $\Gamma$ is discrete.
\end{proof}

We will need the following lemma to produce examples of non-discrete a-dditive subgroups $W$ such that $\Gamma_{W,\mu}$ is discrete. 

\begin{lem}\label{lem_span_no_discreto}
If $\alpha$ and $\beta$ are two rationally independent real numbers then 
	$$W=\text{Span}_{\Z}(\alpha,\beta)$$
is a non-discrete additive subgroup of $\C$.
\end{lem}

\begin{proof}
Let $h:(W,+)\rightarrow(\Ss^1,\cdot)$ be a group homomorphism given by
	$$x\mapsto e^{2\pi i \frac{x}{\alpha}}.$$
Since $\set{\alpha,\beta}$ are rationally independent, so are $\set{1,\frac{\beta}{\alpha}}$. As a consequence of this, $h(W)$ is a sequence of distinct elements in $\Ss^1$, which is compact and therefore there a subsequence, denoted by $\set{g_n}\subset h(W)$, such that
	$$g_n=e^{2\pi i q_n \frac{\beta}{\alpha}}\rightarrow \xi \in \Ss^1,$$
for some $\set{q_n}\subset\Z$ and some $\xi \in \Ss^1$. Define the sequence $\set{h_n}\subset \Ss^1$ by $h_n=g_n g_{n+1}^{-1}$, then $h_n\rightarrow 1$. Denoting 
	$$h_n=e^{2\pi i r_n \frac{\beta}{\alpha}}$$
for some $\set{r_n}\subset\Z$, and taking the logarithm of the sequence we have
	\begin{equation}\label{eq_lem_span_no_discreto_1}
	2\pi i r_n \frac{\beta}{\alpha} + 2\pi i s_n \rightarrow 0
	\end{equation}		
for some logarithm branches defined by $\set{s_n}\subset\Z$. As a consequence of (\ref{eq_lem_span_no_discreto_1}) we have a sequence $\set{r_n \beta + s_n \alpha}\subset{W}$ converging to $0$ and therefore, $W$ is not discrete.
\end{proof}

In the following two examples we have discrete groups $\Gamma_{W,\mu}$ such that $W$ is not discrete.

\begin{ejem}\label{ej_caso1_rango2}
Let $W=\text{Span}_{\Z}(1,\sqrt{2})$. This is an abelian and free group of rank 2, besides $W$ is not discrete (see Lemma \ref{lem_span_no_discreto}).\\

Since $W$ is a free group, to define any group morphism in $W$, it is enough to define it in a generator set. Let $\mu:(W,+)\rightarrow(\C^{\ast},\cdot)$ be the group morphism given by
	$$\mu(1)=1,\;\;\;\;\valorabs{\mu(\sqrt{2})}\neq 1.$$
Suppose that $\Gamma=\Gamma_{W,\mu}$ is non-discrete. Then there exists a non-constant sequence $\set{g_n}\subset\Gamma$ converging to the identity in $W$. As we did in the proof of lemma \ref{lem_caso1_disc_disc}, we denote
	$$g_n=\corchetes{\begin{array}{ccc}
	\mu(w)^{-3} & 0 & 0 \\
	0 & 1 & w_n \\
	0 & 0 & 1
	\end{array}},$$	
then $w_n\rightarrow 0$. If we write $w_n=p_n+q_n\sqrt{2}$, then 
	$$\valorabs{p_n},\;\valorabs{q_n}\rightarrow \infty.$$
On the other hand, 
	\begin{align*}
	\mu(w_n) &= \mu\parentesis{p_n+q_n\sqrt{2}} = \mu\parentesis{p_n}\mu\parentesis{q_n\sqrt{2}}\\
		&= \mu\parentesis{1}^{p_n}\mu\parentesis{\sqrt{2}}^{q_n} = \mu\parentesis{\sqrt{2}}^{q_n}.
	\end{align*}
Since $\valorabs{\mu\parentesis{\sqrt{2}}}\neq 1$, either $\valorabs{\mu\parentesis{w_n}}\rightarrow 0$ or 
$\valorabs{\mu\parentesis{w_n}}\rightarrow \infty$ and therefore, $\valorabs{\mu\parentesis{w_n}}^{-3}\rightarrow \infty$ or 
$\valorabs{\mu\parentesis{w_n}}^{-3}\rightarrow 0$ respectively. In both cases is not possible to have $g_n\rightarrow\text{id}$. This implies that $\Gamma$ is discrete.
\end{ejem} 

\begin{ejem}\label{ej_caso1_rango3}
Let $W=\text{Span}_{\Z}(1,\sqrt{2},i)$. This is again a non-discrete, abelian and free group of rank 3. Let $\mu:(W,+)\rightarrow(\C^{\ast},\cdot)$ be a group morphism given by
	$$\mu(1)=1,\;\;\;\;\valorabs{\mu(\sqrt{2})}\neq 1.$$
and $\mu(i)$ arbitrary. Analogously to the previous example, assume that $\Gamma=\Gamma_{W,\mu}$ is not discrete, then there is a non-constant sequence $\set{g_n}\subset\Gamma$ converging to the identity in $W$ and this sequence determines a sequence $\set{w_n}\subset W$ such that $w_n\rightarrow 0$.\\
Denote $w_n=p_n+q_n\sqrt{2}+r_n i$ and, since $i$ is $\R$-linearly independent of $\set{1,\sqrt{2}}$, we have that $r_n\rightarrow 0$ and $\set{p_n}$, $\set{q_n}$ are as in the previous example. Then it holds again that 
	$$\valorabs{\mu(w_n)} = \mu\parentesis{\sqrt{2}}^{q_n}.$$
Using the same argument as in the last example, $\set{g_n}$ cannot converge to the identity, contradicting that $\Gamma$ is not discrete. Therefore, $\Gamma$ is discrete.
\end{ejem}

Examples \ref{ej_caso1_rango2} and \ref{ej_caso1_rango3} show non-discrete additive groups $W$ with range 2 and 3 such that $\Gamma_{W,\mu}$ is discrete. Proposition \ref{prop_conmutativo_c1_rangoW} says that we cannot found such a group with range 4 or higher. Besides, this proposition gives the full description of discrete commutative subgroups of $U_+$ belonging to the case 1. In order to prove that proposition, we first need to determine the equicontinuity region for case 1 groups.\\

Propositions \ref{prop_convergencia_qp}, \ref{prop_descripcion_Eq} and \ref{prop_eq_in_kuld} remark the importance of knowing all the possible quasi-projective maps that the sequences of distinct elements of $\Gamma$ can converge to. In the following table we list all of these possible quasi-projective limits.

\begin{table}[H]\label{fig_c1_casos_qp}
\begin{center}
  \begin{tabular}{ | l | c | c | c | c | }
    \hline
    Case & $\tau$ & Conditions & Ker($\tau$) & Im($\tau$) \\ \hline
    (i) & $\corchetes{\begin{array}{ccc}
	1 & 0 & 0\\
	0 & 0 & 0 \\
	0 & 0 & 0\\
	\end{array}}$ & \begin{tabular}{c}
	\small $w_n\rightarrow b\in\C$ and $\mu(w_n)\rightarrow 0$ \\
	\small or \\
	\small $w_n\rightarrow \infty$, $\mu(w_n)\rightarrow 0$ and $w_n \mu(w_n)^3\rightarrow 0$\\
	\end{tabular}	
	 & $\linproy{e_2,e_3}$ & $\set{e_1}$  \\ \hline
    (ii) & $\corchetes{\begin{array}{ccc}
	0 & 0 & 0\\
	0 & 1 & b \\
	0 & 0 & 1\\
	\end{array}}$ & \small $w_n\rightarrow b\in\C$ and $\mu(w_n)\rightarrow \infty$ & $\set{e_1}$ & $\linproy{e_2,e_3}$  \\ \hline
    (iii) & $\corchetes{\begin{array}{ccc}
	0 & 0 & 0\\
	0 & 0 & 1 \\
	0 & 0 & 0\\
	\end{array}}$ & \begin{tabular}{c}
\small $w_n\rightarrow \infty$ and $\mu(w_n)\rightarrow \infty$ \\
\small or \\
\small $w_n\rightarrow \infty$ and $\mu(w_n)\rightarrow a\in\C^{\ast}$\\
\small or \\
\small $w_n\rightarrow \infty$, $\mu(w_n)\rightarrow 0$ and $w_n \mu(w_n)^3\rightarrow\infty$\\
	\end{tabular}
   & $\linproy{e_1,e_2}$ & $\set{e_2}$  \\ \hline
    (iv) & $\corchetes{\begin{array}{ccc}
	1 & 0 & 0\\
	0 & 0 & b \\
	0 & 0 & 0\\
	\end{array}}$ & \small $w_n\rightarrow \infty$, $\mu(w_n)\rightarrow 0$ and $w_n \mu(w_n)^3\rightarrow b\in\C^\ast$ & $\set{e_2}$ & $\linproy{e_1,e_2}$ \\
    \hline
  \end{tabular}
\caption{Quasi-projective limits of sequences of distinct elements in $\Gamma$.}
\end{center}
\end{table}

\begin{prop}\label{prop_eq_c1}
Let $\Gamma\subset\PSL$ be a commutative discrete group with the form given in Proposition \ref{prop_conmutativo_c1_descripcion}. If $\Gamma$ contains loxodromic elements then
	$$\Eq(\Gamma)=\CP^2\setminus\parentesis{\linproy{e_1,e_2}\cup\linproy{e_2,e_3}}.$$
\end{prop}

\begin{proof}
We use Proposition \ref{prop_descripcion_Eq} to determine $\Eq(\Gamma)$. Let 
	$$\gamma_w=\corchetes{\begin{array}{ccc}
	\mu(w)^{-3} & 0 & 0\\
	0 & 1 & w \\
	0 & 0 & 1\\
	\end{array}}\in\Gamma$$ 
be a loxodromic element, then $\valorabs{\mu(w)}\neq 1$. Let us suppose, without loss of generality, that $\valorabs{\mu(w)}>1$. Consider the sequence $\set{\gamma_w^n}_{n\in\N}\subset\Gamma$, then
	\begin{equation}\label{eq_prop_eq_c1_1}
	\gamma_w^n\rightarrow\tau_1=\corchetes{\begin{array}{ccc}
	0 & 0 & 0 \\
	0 & 0 & 1 \\
	0 & 0 & 0 \\
	\end{array}},\;\;\;\text{ with Ker}(\tau_1)=\linproy{e_1,e_2}.
	\end{equation}
Considering the sequence $\set{\gamma_w^{-n}}_{n\in\N}\subset\Gamma$ instead, we have 
	\begin{equation}\label{eq_prop_eq_c1_2}
	\gamma_w^{-n}\rightarrow\tau_2=\corchetes{\begin{array}{ccc}
	1 & 0 & 0 \\
	0 & 0 & 0 \\
	0 & 0 & 0 \\
	\end{array}},\;\;\;\text{ with Ker}(\tau_2)=\linproy{e_2,e_3}.
	\end{equation}			
Proposition \ref{prop_descripcion_Eq} together with (\ref{eq_prop_eq_c1_1}) and (\ref{eq_prop_eq_c1_2}) imply that 
	\begin{equation}\label{eq_prop_eq_c1_3}
	\CP^{2}\setminus\parentesis{\linproy{e_1,e_2}\cup\linproy{e_2,e_3}}\subset\Eq(\Gamma).
	\end{equation} 
Proposition \ref{prop_descripcion_Eq} and the table in Figure \ref{fig_c1_casos_qp} imply that 
	\begin{equation}\label{eq_prop_eq_c1_4}
	\Eq(\Gamma)\subset \CP^{2}\setminus\parentesis{\linproy{e_1,e_2}\cup\linproy{e_2,e_3}}.
	\end{equation}	
(\ref{eq_prop_eq_c1_3}) and (\ref{eq_prop_eq_c1_4}) prove the proposition.
\end{proof}

\begin{obs}\label{obs_caso1_llenando_las_2_lineas_ya}
The previous proposition implies that, for a group $\Gamma$ of this first case, if we determine that 
	$$\linproy{e_1,e_2}\cup\linproy{e_2,e_3}\subset \KulL(\Gamma)$$
then, using Proposition \ref{prop_eq_in_kuld}, we have in fact, 
	$$\KulL(\Gamma)=\linproy{e_1,e_2}\cup\linproy{e_2,e_3}.$$
This observation will be very useful when we determine the Kulkarni limit set in the different subcases in the proof of  Theorem \ref{thm_case1_kulkarni}.
\end{obs}

Now we give the description of discrete subgroups of case 1.

\begin{prop}\label{prop_conmutativo_c1_rangoW}
Let $\Gamma=\Gamma_{W,\mu}\subset U_+$ be a group as described in Proposition \ref{prop_conmutativo_c1_descripcion}. The group $\Gamma$ is discrete if and only if $\text{rank}(W)\leq 3$ and the morphism $\mu$ satisfies the following condition:
\begin{enumerate}[(C)]
\item Whenever we have a sequence $\set{w_k}\in W$ of distinct elements such that $w_k\rightarrow 0$, either $\mu(w_k)\rightarrow 0$ or $\mu(w_k)\rightarrow \infty$. 
\end{enumerate}
\end{prop} 

\begin{proof}
First, assume that $\text{rank}(W)\leq 3$ and that the group morphism $\mu$ satisfies the condition (C). First, let us suppose that $W$ is discrete, then by Lemma \ref{lem_caso1_disc_disc}, $\Gamma$ is discrete.\\

Now, assume that $W$ is not discrete. Suppose that $\Gamma$ is not discrete, then let $\set{\gamma_k}\subset\Gamma$ be a sequence of distinct elements such that $\gamma_k\rightarrow\id$, denote
$$\gamma_k=\corchetes{\begin{array}{ccc}
	\mu(w_k)^{-3} & 0 & 0 \\
	0 & 1 & w_k \\
	0 & 0 & 1
	\end{array}}.$$
Then $\mu(w_k)$ converges to some cubic root of the unity and $w_k\rightarrow 0$. Since $\mu$ satisfies condition (C) then $\mu(w_k)\rightarrow 0$ or $\mu(w_k)\rightarrow \infty$ contradicting that $\mu(w_k)$ converges to some cubic root of the unity. This contradiction proves that $\Gamma$ is discrete.\\

Now assume that $\Gamma$ is discrete, by Propositions \ref{prop_eq_in_kuld} and \ref{prop_eq_c1}, $\Gamma$ acts properly and discontinuously on 
	$$\Eq(\Gamma)=\CP^2\setminus\parentesis{\linproy{e_1,e_2}\cup \linproy{e_2,e_3}}\cong \C\times\C^{\ast}.$$  
Consider the universal covering
	$$\pi=(\id,\text{exp}):\C\times\C\rightarrow\C\times\C^\ast,$$
where $\text{exp}(z)=e^z$ for $z\in\C$. The group $\Gamma$ can be written as $\Gamma\cong \Gamma_1\times\Gamma_2$ with
	\begin{align*}
	\Gamma_1 &= \SET{\corchetes{\begin{array}{cc}
	\mu(w)^{-3} & 0\\
	0 & 1 \end{array}}}{w\in W}\\
	\Gamma_2 &= \SET{\corchetes{\begin{array}{cc}
	1 & w\\
	0 & 1 \end{array}}}{w\in W}
	\end{align*}
with the multiplicative group $\Gamma_1$ acting on $\C^\ast$ and the additive group $\Gamma_2$ acting on $\C$. Let $\tilde{\Gamma}$ and $\tilde{\Gamma}_2$ be covering groups of $\Gamma$ and $\Gamma_2$ respectively (see Theorem 9.1 of \cite{bredon}), let $\tilde{\Gamma}_1=\Gamma_1$. As a consequence of this, there is a group morphism, induced by $\pi$ and still denoted by $\pi$, given by
	\begin{align*}
	\pi=(\id,\text{exp}):\tilde{\Gamma}\cong\tilde{\Gamma}_1\times\tilde{\Gamma}_2 &\rightarrow \Gamma_1\times\Gamma_2\\
	(\alpha,\beta) &\mapsto (\alpha,e^\beta) 
	\end{align*}

Since $\Gamma_1$ and $\Gamma_2$ are both commutative and torsion free then, by Proposition \ref{prop_grupos_conm_sintorsion_Zk}
	$$\Gamma_1\cong\Z^{k_1}\text{ and }\Gamma_2\cong\Z^{k_2}$$
with $k_1=\text{rank}(\Gamma_1)$ and $k_2=\text{rank}(\Gamma_2)$. Then $\Gamma \cong \Gamma_1\times\Gamma_2\cong \Z^{k_1}\times\Z^{k_2}=\Z^k$, with $k=k_1+k_2=\text{rank}(\Gamma)$. Observe that the kernel of the group morphism $\pi$ is given by
	$$\kernel(\pi)=\kernel(\id)\times\ker(\text{exp})\cong \Z$$
then,
	$$\tilde{\Gamma}\cong \kernel(\pi)\times \Gamma\cong \Z\times \Z^k,$$
and therefore,
	\begin{equation}\label{eq_dem_prop_conmutativo_c1_rangoW_1}
	\text{rank}\parentesis{\Gamma}=k+1.	
	\end{equation}	 
On the other hand, since $\Gamma$ acts properly and discontinuously on $\C\times\C^\ast$, then $\tilde{\Gamma}$ acts properly and discontinuously on $\C\times\C$, which is simply connected. Then, by Theorem \ref{thm_obdim_2}, $\text{rank}\parentesis{\tilde{\Gamma}}\leq 4$. This, together with (\ref{eq_dem_prop_conmutativo_c1_rangoW_1}) yields $\text{rank}(\Gamma)\leq 3$. Using Proposition \ref{prop_conmutativo_caso1_rangoG_rangoW} we conclude that $\text{rank}(W)\leq 3$.\\

Now we will verify that $\mu$ satisfies the condition (C). Let $\set{w_k}\subset W$ be a sequence of distinct elements such that $w_k\rightarrow 0$ (in the case that $W$ is non-discrete). Consider the sequence $\set{\mu(w_k)}\subset\C^\ast$, assume that it does not converge to $0$ or $\infty$, then there are open neighborhoods $U_0$ and $U_\infty$ of $0$ and $\infty$ respectively such that $\set{\mu(w_k)}\subset \CP^1\setminus\parentesis{U_0\cup U_\infty}$. Since $\CP^1$ is compact and $\CP^1\setminus\parentesis{U_0\cup U_\infty}$ is a closed subset of $\CP^1$, then $\CP^1\setminus\parentesis{U_0\cup U_\infty}$ is compact and therefore, there is a converging subsequence of $\set{\mu(w_k)}$, still denoted the same way. Let $z\in\C^\ast$ such that $\mu(w_k)\rightarrow z$, then
	$$\gamma_k=\corchetes{\begin{array}{ccc}
	\mu(w_k)^{-3} & 0 & 0 \\
	0 & 1 & w_k \\
	0 & 0 & 1
	\end{array}}\rightarrow \corchetes{\begin{array}{ccc}
	z^{-3} & 0 & 0 \\
	0 & 1 & 0 \\
	0 & 0 & 1
	\end{array}}.$$
Then $\Gamma$ cannot be discrete. This proves that $\mu$ satisfies the condition (C). 
\end{proof}

We say that $z\in\C$ is a rational rotation (resp. irrational rotation) if $z=e^{2\pi i \theta}$ for some $\theta\in\Q$ (resp. $\theta\in\R\setminus\Q$).\\ 

Now we describe the Kulkarni limit set of these groups. In order to do this, we will divide all groups of case 1 into subcases according to the following diagram:


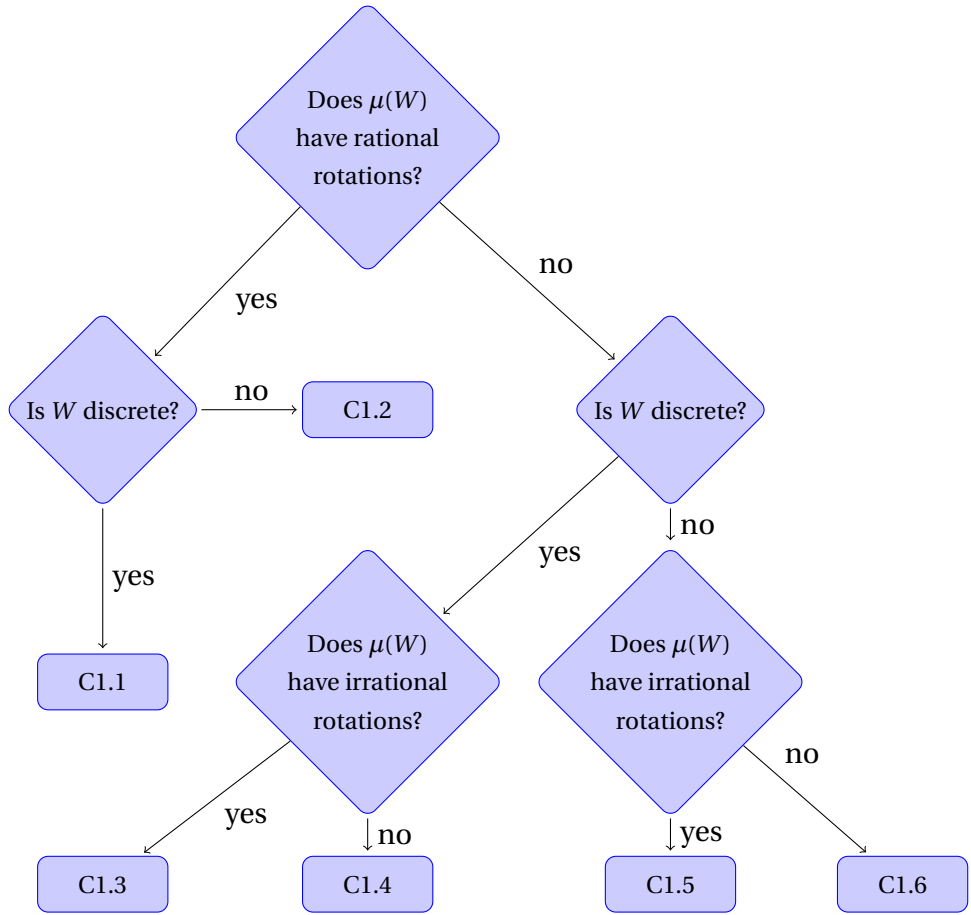
\begin{figure}[H]
\begin{tikzpicture} [
    auto,
    decision/.style = { diamond, draw=blue, thin, fill=blue!20,
                        text width=6em, text badly centered,
                        inner sep=1pt, rounded corners },
    block/.style    = { rectangle, draw=blue, thin, 
                        fill=blue!20, text width=4em, text centered,
                        rounded corners, minimum height=2em },
    line/.style     = { draw, thin , ->, shorten >=2pt },
  ]   
  \matrix [column sep=4mm, row sep=5mm] {
     &  \node [decision] (rotrac)
                        {\footnotesize Does $\mu(W)$ have rational rotations?};
                              &    &  & \\
     \node [decision] (disc1) {\footnotesize Is $W$ discrete?};  &  \node [block] (c2) {\footnotesize C1.2}; & \node [decision] (disc2) {\footnotesize Is $W$ discrete?}; & & \\
    \node [block] (c1) {\footnotesize C1.1};   &    \node [decision] (irrrot1) {\footnotesize Does $\mu(W)$ have irrational rotations?};   &  \node [decision] (irrrot2) {\footnotesize Does $\mu(W)$ have irrational rotations?}; & & \\
   \node [block] (c3) {\footnotesize C1.3}; & \node [block] (c4) {\footnotesize C1.4}; & \node [block] (c5) {\footnotesize C1.5}; & \node [block] (c6) {\footnotesize C1.6}; & \\ 
  };
  \begin{scope} [every path/.style=line]
    \path (rotrac)   --    node {yes} (disc1);
    \path (rotrac)   --    node {no} (disc2);
    \path (disc1)    --    node {yes} (c1);
	\path (disc1)    --    node {no} (c2);    
    \path (disc2)    --    node {yes} (irrrot1);
    \path (disc2)    --    node {no} (irrrot2);
    \path (irrrot1)  --    node {yes} (c3);
    \path (irrrot1)  --    node {no} (c4);
    \path (irrrot2)  --    node {yes} (c5);
    \path (irrrot2)  --    node [near start] {no} (c6);
  \end{scope}
\end{tikzpicture}
\caption{Diagram of the six subcases for case 1.}
\label{fig_casos}
\end{figure}

Before describing the Kulkarni limit set of case 1 groups, we will need the following proposition (see chapter 2 of \cite{TPR}).

\begin{prop}\label{prop_forma_subgrupos_aditivos}
The closed subgroup $H\subset\R^n$ is additive if and only if
	$$H\cong \R^p\oplus\Z^q$$
with non-negative integers $p,q$ such that $p+q\leq n$. 
\end{prop}

\begin{obs}
If $\Gamma$ is cyclic then $W$ is cyclic (see Proposition \ref{prop_conmutativo_caso1_rangoG_rangoW}), denote $W=\prodint{w}$. The group $\Gamma$ is either generated by a loxo-parabolic or an ellipto-parabolic element depending on whether $\valorabs{\mu(w)}\neq 1$ or $\valorabs{\mu(w)}=1$ respectively. According to Proposition 4.2.10 and 4.2.19 of \cite{ckg_libro}, 
	$$\KulL(\Gamma)=\begin{cases}
	\linproy{e_1,e_2}\cup\linproy{e_2,e_3}, & \valorabs{\mu(w)}\neq 1\\
	\linproy{e_1,e_2}, & \valorabs{\mu(w)}= 1 
	\end{cases}.$$
Therefore, we can assume that the group $\Gamma$ is not cyclic.
\end{obs}

\begin{lem}\label{lem_casoc1_L0}
Let $\Gamma\subset U_+$ be commutative discrete group with the form given in Proposition\ref{prop_conmutativo_c1_rangoW}, then
	$$L_0(\Gamma)=\begin{cases}
	\set{e_1,e_2},& \text{If }\mu(W)\text{ does not contain rational rotations}\\	
	\linproy{e_1,e_2},& \text{If }\mu(W)\text{ contains rational rotations}	
	\end{cases}	
	$$
\end{lem}

\begin{proof}
In both cases $\set{e_1,e_2}\subset L_0(\Gamma)$ since both are global fixed points of $\Gamma$. A direct computation shows that if $X\in L_0(\Gamma)$ then $X\in\linproy{e_1,e_2}$, that is, $L_0(\Gamma)\subset \linproy{e_1,e_2}$.\\

It is straight-forward to verify that, if $\mu(W)$ contains rational rotations, then every point in $\linproy{e_1,e_2}$ is in $L_0(\Gamma)$ and therefore $L_0(\Gamma)=\linproy{e_1,e_2}$. Also, if $\mu(W)$ doesn't contain rational rotations, then $L_0(\Gamma)=\set{e_1,e_2}$. 
\end{proof}

\begin{lem}\label{lem_casoc11}
Let $\Gamma=\Gamma_{W,\mu}\subset U_+$ be a commutative discrete group with the form given in Proposition\ref{prop_conmutativo_c1_rangoW}, if $W$ is discrete and $\mu(W)$ contains rational rotations then
	$$\KulL(\Gamma)=\linproy{e_1,e_2}.$$	
\end{lem}

\begin{proof}
Since $\mu(W)$ contains rational rotations, $L_0(\Gamma)=\linproy{e_1,e_2}$ (see Lemma \ref{lem_casoc1_L0}). Since $W$ is discrete, then $\Gamma$ has no sequences with quasi-projective limit with the form (ii) in Table \ref{fig_c1_casos_qp}, this implies $\KulL(\Gamma)\subset\linproy{e_1,e_2}$ (see Proposition \ref{prop_eq_in_kuld}). All of this yields
	$$\linproy{e_1,e_2}=L_0(\Gamma)\subset \KulL(\Gamma)\subset\linproy{e_1,e_2}.$$
This proves the proposition.
\end{proof}

\begin{lem}\label{lem_casoc1_KulEq}
Let $\Gamma=\Gamma_{W,\mu}\subset U_+$ be a commutative discrete group with the form given in Proposition\ref{prop_conmutativo_c1_rangoW}, if $\mu(W)$ contains rational rotations and $W$ is not discrete, then
	$$\KulL(\Gamma)=\linproy{e_1,e_2}\cup\linproy{e_2,e_3}.$$
\end{lem}

\begin{proof}
Since $W$ is not discrete, there is a sequence $\set{w_n}\subset W$ such that $w_n\rightarrow 0$ and therefore, either $\mu(w_n)\rightarrow \infty$ or $\mu(w_n)\rightarrow 0$ (see Proposition \ref{prop_conmutativo_c1_rangoW}). We can assume without loss of generality that the former happens (otherwise, consider the sequence $\set{-w_n}$ instead).\\

Then, the quasi-projective limit of the sequence $\set{\gamma_{w_n}}$ is 
	$$\tau=\corchetes{\begin{array}{ccc}
	0 & 0 & 0\\
	0 & 1 & 0 \\
	0 & 0 & 1\\
	\end{array}}.$$
Let $z\in \CP^2\setminus L_0(\Gamma)$, then $z\nin\kernel(\tau)=\set{e_1}$ (see Lemma \ref{lem_casoc1_L0}). Therefore, the set of accumulation points of points in $\CP^2\setminus L_0(\Gamma)$ is $\text{Im}(\tau)=\linproy{e_2,e_3}$. Then, 
 $$\linproy{e_1,e_2}\cup\linproy{e_2,e_3}\subset L_0(\Gamma)\cup L_1(\Gamma).$$
Using Observation \ref{obs_caso1_llenando_las_2_lineas_ya} we conclude that 
	$$\KulL(\Gamma)=\linproy{e_1,e_2}\cup\linproy{e_2,e_3},$$
whenever $W$ contains rational rotations and is not discrete.
\end{proof}

As an immediate consequence of the last lemma we conclude that groups $\Gamma$ belonging to the case C1.2 satisfy $\KulL(\Gamma)=\linproy{e_1,e_2}\cup\linproy{e_2,e_3}$.\\

We say that the commutative group $\Gamma=\Gamma_{W,\mu}$ satisfy the condition \titemT{(F)} if there is a sequence $\set{w_k}\subset W$ such that $w_k\rightarrow\infty$, $\mu(w_k)\rightarrow 0$ and $w_k\mu(w_k)^3\rightarrow b\in\C^\ast$.\\

If $W$ is discrete, using Table \ref{fig_c1_casos_qp} we conclude that $\linproy{e_2,e_3}$ cannot be contained in $\KulL(\Gamma)$. Then, if $\mu(W)$ contains no rational rotations, either $\KulL(\Gamma)=\set{e_1,e_2}$ or $\KulL(\Gamma)=\linproy{e_1,e_2}$. This argument, together with Table \ref{fig_c1_casos_qp} proves the following lemma.

\begin{lem}\label{lem_casoc13c14}
Let $\Gamma=\Gamma_{W,\mu}\subset U_+$ be a commutative discrete group with the form given in Proposition\ref{prop_conmutativo_c1_rangoW}, if $\mu(W)$ contains no rational rotations and $W$ is discrete, then
	$$\KulL(\Gamma)=\begin{cases}
	\linproy{e_1,e_2},& \text{ if }\Gamma\text{ satisfy condition }\text{\titemT{(F)}}\\
	\set{e_1,e_2},& \text{ any other case}
	\end{cases}.$$
\end{lem}

This previous lemma describes cases C1.3 and C1.4.
 
\begin{lem}\label{lem_casoc15c16}
Let $\Gamma=\Gamma_{W,\mu}\subset U_+$ be a commutative discrete group with the form given in Proposition\ref{prop_conmutativo_c1_rangoW}, if $\mu(W)$ contains no rational rotations and $W$ is not discrete, then
	$$\KulL(\Gamma)=\begin{cases}
	\linproy{e_1,e_2}\cup\linproy{e_2,e_3},& \text{ if }\Gamma\text{ satisfy condition }\text{\titemT{(F)}}\\
	\set{e_1}\cup\linproy{e_2,e_3},& \text{ any other case}
	\end{cases}.$$
\end{lem}

\begin{proof}
Since $\mu(W)$ contains no rational rotations, $L_0(\Gamma)=\set{e_1,e_2}$ (see Lemma \ref{lem_casoc1_L0}). Since $W$ is not discrete then there is a sequence $\set{w_k}\subset W$ such that $w_k\rightarrow 0$, and then $\mu(w_k)\rightarrow \infty$ or $\mu(w_k)\rightarrow 0$ (see Proposition\ref{prop_conmutativo_c1_rangoW}). In the former case, we can conclude, using Table \ref{fig_c1_casos_qp}, that $\linproy{e_2,e_3}\subset\KulL(\Gamma)$. In the latter case we can consider the sequence $\set{-w_k}$ which satisfies $\mu(-w_k)\rightarrow \infty$ and we conclude again that $\linproy{e_2,e_3}\subset\KulL(\Gamma)$.\\

If the group $\Gamma$ satisfies condition \titemT{(F)}, using Table \ref{fig_c1_casos_qp}, it follows that $\linproy{e_1,e_2}\subset\KulL(\Gamma)$ and using Observation \ref{obs_caso1_llenando_las_2_lineas_ya}, we conclude that $\KulL(\Gamma)= \linproy{e_1,e_2}\cup\linproy{e_2,e_3}$. Analogously, if $\Gamma$ doesn't satisfy condition \titemT{(F)}, $\KulL(\Gamma)= \set{e_1}\cup\linproy{e_2,e_3}$.
\end{proof}

This previous lemma describes cases C1.5 and C1.6. Finally, Lemmas \ref{lem_casoc11}, \ref{lem_casoc1_KulEq}, \ref{lem_casoc13c14} and \ref{lem_casoc15c16} together prove the following theorem. 
 
\begin{thm}\label{thm_case1_kulkarni}
Let $\Gamma\subset\PSL$ be commutative discrete group with the form given in Proposition\ref{prop_conmutativo_c1_rangoW}, then
	$$
	\KulL(\Gamma) = \begin{cases}
	\set{e_1,e_2}, & \text{Cases C1.3 or C1.4 with condition \titemT{(F)} not holding.} \\
	\linproy{e_1,e_2}, & \begin{cases}
						 \text{Cases C1.3 or C1.4, satisfying condition \titemT{(F)}} &\\
						 \text{Case C1.1} & \\	
						 \end{cases}	  \\
	\set{e_1}\cup\linproy{e_2,e_3}, & \text{Cases C1.5 or C1.6 with condition \titemT{(F)} not holding}.\\
	\linproy{e_1,e_2}\cup\linproy{e_2,e_3}, & \begin{cases}
						 \text{Cases C1.5 or C1.6, satisfying condition \titemT{(F)}} &\\
						 \text{Case C1.2} & \\	
						 \end{cases}	  \\
	\end{cases}
	$$
\end{thm}

In the following example we give a group belonging to the case C1.6 and determine its Kulkarni limit set. This example is important because it was believed that only fundamental groups of Hopf surfaces had a Kulkarni limit set consisting of a line and a point. This is an example of a group with this limit set, which is not a fundamental group of a Hopf surface (since these groups are cyclic, see \cite{kato1975topology}).  

\begin{ejem}\label{ej_falso_hopf}
Let $W=\text{Span}_{\Z}\set{1,\sqrt{2}}$ be an additive subgroup of $\C$ and let $\mu:(W,+)\rightarrow(\C^{\ast},\cdot)$ a group homomorphism given by
	$$\begin{array}{cc}
	\mu(1)=e^{-1}, & \mu(\sqrt{2})=e^{\sqrt{2}}.
	\end{array}
	$$
Consider the group $\Gamma=\Gamma_{W,\mu}$ defined as before. $\Gamma$ is a commutative group (see Proposition \ref{prop_conmutativo_c1_descripcion}). The rank of $W$ is 2 and if $\set{w_n}\subset W$ is a sequence of distinct elements such that $w_n=p_n+q_n\sqrt{2}\rightarrow 0$, then, without loss of generality assume that $p_n,-q_n\rightarrow \infty$, then
	$$\mu(w_n)=\mu(1)^{p_n}\mu(\sqrt{2})^{q_n}=e^{-p_n+q_n\sqrt{2}}\rightarrow 0,$$
and therefore $\Gamma$ is discrete (see Proposition \ref{prop_conmutativo_c1_rangoW}).\\

On the other hand, if $x\in\mu(W)$ then $x=\mu(p+q\sqrt{2})$ for some $p,q\in\Z$ and then
	$$\valorabs{x}=\valorabs{\mu(1)^p\mu(\sqrt{2})^q}=e^{-p+q\sqrt{2}}\neq 1.$$
Observe that, $\valorabs{x}=1$ if and only if $-p+q\sqrt{2}=0$ but this cannot happen because $set{1,\sqrt{2}}$ are rationally independent. This verifies that $\mu(W)$ doesn't have rotations and therefore belongs to the case C1.6.\\

Now we determine its Kulkarni limit set. As we see in Theorem \ref{thm_case1_kulkarni}, 
	$$\KulL(\Gamma)=\set{e_1}\cup\linproy{e_2,e_3},$$
unless condition \titemT{(F)} is satisfied. That is, unless there is a sequence $\set{w_n}\subset W$ such that $w_n\rightarrow \infty$, $\mu(w_n)\rightarrow 0$ and 
	$$w_n\mu(w_n)^3\rightarrow b\in\C^{\ast}.$$ 
Assume that this happens, since $w_n\rightarrow \infty$, there are the following possibilities for the sequences $\set{p_n}\text{, }\set{q_n}\subset\Z$:
\begin{enumerate}
\item If both sequences $\set{p_n},\set{q_n}$ are bounded, then there exists $R>0$ such that $\valorabs{p_n},\valorabs{q_n}<R$ and then
	$$\valorabs{p_n+q_n\sqrt{2}}<R(\sqrt{2}+1).$$  	
Therefore it's impossible that $w_n\rightarrow\infty$.
\item If $\mu(w_n)\rightarrow 0$ with $p_n\rightarrow\infty$ and $\set{q_n}$ bounded. Then
	\begin{align*}
	w_n\mu(w_n)^3 &= \parentesis{p_n+q_n\sqrt{2}}e^{-3p_n+3q_n\sqrt{2}}\\
	 &= p_ne^{-3p_n+3q_n\sqrt{2}} + q_n\sqrt{2} e^{-3p_n+3q_n\sqrt{2}}\\
	 &= \underbrace{p_n e^{-3p_n}}_{\rightarrow 0}\underbrace{e^{3q_n\sqrt{2}}}_{\text{bounded}} + \underbrace{q_n\sqrt{2}}_{\text{bounded}} \underbrace{e^{-3p_n+3q_n\sqrt{2}}}_{\rightarrow 0} \rightarrow 0
	\end{align*}

\item If $\mu(w_n)\rightarrow 0$ with $q_n\rightarrow-\infty$ and $\set{p_n}$ bounded. Then
	\begin{align*}
	w_n\mu(w_n)^3 &= \parentesis{p_n+q_n\sqrt{2}}e^{-3p_n+3q_n\sqrt{2}}\\
	 &= p_ne^{-3p_n+3q_n\sqrt{2}} + q_n\sqrt{2} e^{-3p_n+3q_n\sqrt{2}}\\
	 &= \underbrace{p_n}_{\text{bounded}}\underbrace{e^{-3p_n+3q_n\sqrt{2}}}_{\rightarrow 0} + \underbrace{q_n\sqrt{2}e^{3q_n\sqrt{2}}}_{\rightarrow 0} \underbrace{e^{-3p_n}}_{\text{bounded}} \rightarrow 0
	\end{align*}
\end{enumerate}
Then neither of these 3 cases occur and therefore condition \titemT{(F)} doesn't hold. Using Theorem \ref{thm_case1_kulkarni}, we conclude
	$$\KulL(\Gamma)=\set{e_1}\cup\linproy{e_2,e_3}.$$
\end{ejem}

\subsection{Case 2}

Finally, we study the case where $\Gamma$ is conjugate to a diagonal group. We start by describing the form of these groups.

\begin{prop}\label{prop_conmutativo_c3_descripcion}
Let $\Gamma\subset U_+$ be a commutative subgroup such that each element of $\Gamma$ has the form
		\begin{equation}\label{eq_prop_conmutativo_c3_descripcion_1}				
		\parentesis{\begin{array}{ccc}
		\alpha & 0 & 0 \\
		0 & \beta & 0 \\
		0 & 0 & \alpha^{-1}\beta^{-1}
		\end{array}}.
		\end{equation}
Then there exist two multiplicative subgroups $W_1,W_2\subset(\C^\ast,\cdot)$ such that
	\begin{equation}\label{eq_prop_conmutativo_c3_descripcion_2}
	\Gamma=\Gamma_{W_1,W_2}=\SET{\corchetes{\begin{array}{ccc}
	w_1 & 0 & 0 \\
	0 & w_2 & 0 \\
	0 & 0 & 1
	\end{array}}}{w_1\in W_1,\;w_2\in W_2}.	
	\end{equation}
\end{prop}

\begin{proof}
Let $\Gamma\subset U_+$ be a commutative group with the form given by (\ref{eq_prop_conmutativo_c3_descripcion_1}). Let $\gamma=\text{Diag}(\alpha,\beta,\alpha^{-1}\beta^{-1})\in\Gamma$, then $\gamma=\text{Diag}(\alpha^2\beta,\alpha\beta^2,1)$. Let $W_1,W_2\subset \C^\ast$ be the two multiplicative groups given by $W_1=\lambda_{13}(\Gamma)$ and $W_2=\lambda_{23}(\Gamma)$, then $\gamma=\text{Diag}(w_1,w_2,1)$, where
	\begin{align*}
	w_1 &:= \alpha^2\beta = \lambda_{13}(\gamma)\in W_1 \\
	w_2 &:= \alpha\beta^2 = \lambda_{23}(\gamma)\in W_2. 
 	\end{align*}
Then $\Gamma$ has the form given by (\ref{eq_prop_conmutativo_c3_descripcion_2}).
\end{proof}

\begin{prop}\label{prop_caso_diagonal_rango}
Let $\Gamma\subset U_+$ be a diagonal discrete group such that every element has the form
	 $$\gamma=\corchetes{\begin{array}{ccc}
	w_1 & 0 & 0 \\
	0 & w_2 & 0 \\
	0 & 0 & w_1^{-1}w_2^{-1}
	\end{array}}$$
then $\text{rank}(\Gamma)\leq 2$.	
\end{prop}

\begin{proof}
Recall the group morphisms $\lambda_{ij}$ defined in Section \ref{sec_commutative_triangular}. Let $\mu:\Gamma\rightarrow\R^2$ given by
	$$\mu(\gamma)=\parentesis{\log\valorabs{\lambda_{13}(\gamma)},\log\valorabs{\lambda_{23}(\gamma)}}.$$
Clearly $\mu$ is well defined and it's a group homomorphism between $\Gamma$ and the additive subgroup $\mu(\Gamma)\subset\R^2$. Furthermore,
	$$\kernel(\mu)=\SET{\gamma\in\Gamma}{\valorabs{\lambda_{13}(\gamma)}=\valorabs{\lambda_{23}(\gamma)}=1}=\set{\id}.$$
Then, $\mu:\Gamma\rightarrow\mu(\Gamma)$ is a group isomorphism. Since $\Gamma$ is discrete, then $\mu(\Gamma)$ is discrete and therefore $\text{rank}\parentesis{\mu(\Gamma)}\leq 2$, then $\text{rank}(\Gamma)\leq 2$.	
\end{proof}

\begin{prop}
Let $\Gamma=\Gamma_{W_1,W_2}$ be a subgroup of $U_+$ with the form given by Proposition \ref{prop_conmutativo_c3_descripcion}. Then 
	$$\Gamma\cong W_1\oplus W_2.$$ 
\end{prop}

\begin{proof}
Let $H_1,H_2\subset \Gamma$ be two subgroups given by
	$$H_1=\SET{\text{Diag}(w,1,1)}{w\in W_1},\;\;\;H_2=\SET{\text{Diag}(1,w,1)}{w\in W_2}.$$
Observe that $H_i\cong W_i$. Both $H_1$ and $H_2$ are normal in $\Gamma$ and $\Gamma=\prodint{H_1,H_2}$. Also, $H_1\cap H_2=\set{\id}$, therefore 
	$$\Gamma=H_1\oplus H_2\cong W_1\oplus W_2.$$ 
\end{proof}

The previous proposition, together with Proposition \ref{prop_caso_diagonal_rango}, imply that 
$$\rank(W_1)+\rank(W_2)\leq 2.$$
If $\rank(W_1)=1$, $\rank(W_2)=0$ or $\rank(W_1)=0$, $\rank(W_2)=1$, then $\Gamma$ is cyclic and therefore, its Kulkarni limit set is described in Section 4.2 of \cite{ckg_libro}.\\

The case $\rank(W_1)=2$, $\rank(W_2)=0$ is conjugated to the case $\rank(W_1)=0$, $\rank(W_2)=2$. However, both cases imply that $\Gamma$ is not discrete. Therefore we just have to describe the case:

	\begin{equation}\label{eq_forma_diagonal}
	\Gamma:=\Gamma_{\alpha,\beta}=\SET{\corchetes{\begin{array}{ccc}
	\alpha^{n} & 0 & 0 \\
	0 & \beta^{m} & 0 \\
	0 & 0 & 1
	\end{array}}}{n,m\in\Z}.
	\end{equation}

for some $\alpha,\beta\in \C^\ast$ such that $\valorabs{\alpha}\neq 1$ or $\valorabs{\beta}\neq 1$.\\

Consider a sequence of distinct elements $\set{\gamma_k}\subset \Gamma$ given by
	$$\corchetes{\begin{array}{ccc}
	\alpha^{n_k} & 0 & 0 \\
	0 & \beta^{m_k} & 0 \\
	0 & 0 & 1
	\end{array}}$$

In the following table we show all the possible quasi-projective limits of sequences in the diagonal case.

\begin{table}[H]
\begin{center}
  \begin{tabular}{ | l | c | c | c | c | }
    \hline
    Case & $\tau$ & Conditions & Ker($\tau$) & Im($\tau$) \\ \hline
    (i) & $\corchetes{\begin{array}{ccc}
	1 & 0 & 0\\
	0 & 0 & 0 \\
	0 & 0 & 0\\
	\end{array}}$ & \begin{tabular}{c}
	\small $\alpha^{n_k}\rightarrow\infty$, $\beta^{m_k}\rightarrow \infty$ and $\alpha^{n_k}\beta^{-m_k}\rightarrow \infty$ \\
	\small $\alpha^{n_k}\rightarrow\infty$ and $\beta^{m_k}\rightarrow b\in\C$	
	\end{tabular}
		
	 & $\linproy{e_2,e_3}$ & $\set{e_1}$  \\ \hline
    (ii) & $\corchetes{\begin{array}{ccc}
	0 & 0 & 0\\
	0 & 1 & 0 \\
	0 & 0 & 0\\
	\end{array}}$ & \begin{tabular}{c}
	\small $\alpha^{n_k}\rightarrow\infty$, $\beta^{m_k}\rightarrow \infty$ and $\alpha^{-n_k}\beta^{m_k}\rightarrow \infty$ \\
	\small $\alpha^{n_k}\rightarrow 0$ and $\beta^{m_k}\rightarrow \infty$	
	\end{tabular} & $\linproy{e_1,e_3}$ & $\set{e_2}$  \\ \hline
    (iii) & $\corchetes{\begin{array}{ccc}
	0 & 0 & 0\\
	0 & 0 & 0 \\
	0 & 0 & 1\\
	\end{array}}$ & \small $\alpha^{n_k}\rightarrow 0$ and $\beta^{m_k}\rightarrow 0$
   & $\linproy{e_1,e_2}$ & $\set{e_3}$  \\ \hline
    (iv) & $\corchetes{\begin{array}{ccc}
	b & 0 & 0\\
	0 & 0 & 0 \\
	0 & 0 & 1\\
	\end{array}}$ & $\alpha^{n_k}\rightarrow b\in\C^\ast$ and $\beta^{m_k}\rightarrow 0$ & $\set{e_2}$ & $\linproy{e_1,e_3}$ \\
    \hline
    (v) & $\corchetes{\begin{array}{ccc}
	b & 0 & 0\\
	0 & 1 & 0 \\
	0 & 0 & 0\\
	\end{array}}$ & $\alpha^{n_k}\rightarrow \infty$, $\beta^{m_k}\rightarrow \infty$ and $\alpha^{n_k}\beta^{-m_k}\rightarrow b\in\C^\ast$ & $\set{e_3}$ & $\linproy{e_1,e_2}$ \\
    \hline
    (vi) & $\corchetes{\begin{array}{ccc}
	0 & 0 & 0\\
	0 & b & 0 \\
	0 & 0 & 1\\
	\end{array}}$ & $\alpha^{n_k}\rightarrow 0$ and $\beta^{m_k}\rightarrow b\in\C^\ast$ & $\set{e_1}$ & $\linproy{e_2,e_3}$ \\
    \hline
  \end{tabular}
\caption{Quasi-projective limits of sequences of distinct elements in $\Gamma$.}
\label{fig_cdiag_casos_qp}
\end{center}
\end{table}


The following lemmas will be useful to determine the Kulkarni limit set of the groups $\Gamma_{\alpha,\beta}$.

\begin{lem}\label{lem_caso_diagonal_L0}
Let $\Gamma_{\alpha,\beta}\subset U_+$ be a discrete group containing loxodromic elements such that 
	\begin{equation}\label{eq_lem_caso_diagonal_L0_1}
	\alpha^n=\beta^m,
	\end{equation}		
	for some $n,m\in\Z$, then
	$$L_0(\Gamma)=\linproy{e_1,e_2}\cup\set{e_3}.$$
If (\ref{eq_lem_caso_diagonal_L0_1}) does not hold, then 
	$$L_0(\Gamma)=\set{e_1,e_2,e_3}.$$
\end{lem}

\begin{proof}
Let $\Gamma=\Gamma_{\alpha,\beta}$ as in the hypothesis of the lemma and let $z=\corchetes{z_1:z_2:z_3}\in L_0(\Gamma)$. Then $\corchetes{\alpha^p z_1: \beta^q z_2: z_3}= \corchetes{z_1:z_2:z_3}$ for an infinite number of $p,q\in\Z$.\\

If $z_3\neq 0$ then, either $\alpha$ and $\beta$ are rational rotations or $z_1=z_2=0$. If the former happens, then $\Gamma$ contains no loxodromic elements, contradicting the hypothesis. If the latter happens, then $z=e_3$.\\  

If $z_3=0$, we can assume without loss of generality that $z_1\neq 0$. If $z_2=0$, then $z=e_2$. If $z_2\neq 0$, then $\alpha^p=\beta^q$ for an infinite number of integers $p,q$ (this holds because (\ref{eq_lem_caso_diagonal_L0_1}) implies that $\alpha^{jn}=\beta^{jm}$ for any $j\in\Z$). Then $z\in \linproy{e_1,e_2}$. If (\ref{eq_lem_caso_diagonal_L0_1}) doesn't hold, then no point in $\linproy{e_1,e_2}\setminus\set{e_1,e_2}$ satisfies that $\corchetes{\alpha^p z_1: \beta^q z_2: z_3}= \corchetes{z_1:z_2:z_3}$ for an infinite number of $p,q\in\Z$.\\
 
All of these, together with Proposition \ref{prop_ei_en_L0} proves the lemma.
\end{proof}	

\begin{lem}\label{lem_caso_diagonal_L1}
Let $\Gamma_{\alpha,\beta}\subset U_+$ be a discrete group containing loxodromic elements, then
	\begin{itemize}
	\item $L_1(\Gamma)=\set{e_1,e_2}$, if $\valorabs{\alpha}>1>\valorabs{\beta}$ or $\valorabs{\alpha}<1<\valorabs{\beta}$.
	\item $L_1(\Gamma)=\set{e_1,e_3}$, if $\valorabs{\alpha}>\valorabs{\beta}>1$ or $\valorabs{\alpha}<\valorabs{\beta}<1$.
	\item $L_1(\Gamma)=\set{e_1}\cup\linproy{e_2,e_3}$, if $\beta$ is an irrational rotation.	
	\end{itemize}	 
\end{lem}

\begin{proof}
Let $\alpha,\beta\in\C^\ast$, as before, since $\Gamma=\Gamma_{\alpha,\beta}$ contains loxodromic elements, we can assume without loss of generality that $\valorabs{\alpha}\neq 1$. Let $\set{\gamma_k}\subset\Gamma$ be a sequence of distinct elements given by $\gamma_k=\text{Diag}(\alpha^k,\beta^k,1)$.\\

The set $L_0(\Gamma)$ is given by Lemma \ref{lem_caso_diagonal_L0} and let $z=\corchetes{z_1:z_2:z_3}\in \CP^2\setminus L_0(\Gamma)$. In any of the two possible outcomes for $L_0(\Gamma)$ described in Lemma \ref{lem_caso_diagonal_L0}, we have
	$$\gamma_k\corchetes{\begin{array}{c}
	z_1\\
	z_2\\
	z_3
\end{array}}=\corchetes{\begin{array}{c}
	\alpha^k z_1\\
	\beta^k z_2\\
	z_3
\end{array}}.$$
We have two essentially different sequences in $\Gamma$, $\set{\gamma^k}$ and $\set{\gamma^{-k}}$. Since $\valorabs{\alpha}\neq 1$, we can assume, without loss of generality that $\valorabs{\alpha}>1$, then $\alpha^k\rightarrow\infty$ as $k\rightarrow \infty$. We have three cases:
		\begin{enumerate}[(i)]
		\item $\valorabs{\beta}<1$, then $\alpha^k\rightarrow 0$ $k\rightarrow \infty$, then $\gamma_k z\rightarrow e_1$ as $k\rightarrow \infty$. Analogously, $\gamma_k z\rightarrow e_2$ as $k\rightarrow -\infty$.
		\item $\valorabs{\beta}=1$, $\beta$ cannot be a rational rotation because $\Gamma$ is torsion free, then $\beta$ is an irrational rotation. Then, there are subsequences (still denoted by $\set{\gamma_k}$) converging to any $b\in\C$, $\valorabs{b}=1$. Then $\gamma_k z\rightarrow e_1$ as $k\rightarrow \infty$ and $\gamma_k z\rightarrow \corchetes{0:b^{-1}z_2:z_3}$ as $k\rightarrow -\infty$. 
		\item $\valorabs{\beta}>1$, then $\alpha^k\rightarrow \infty$ as $k\rightarrow \infty$. We have three subcases:
			\begin{itemize}
			\item $\alpha^k\beta^{-k}\rightarrow\infty$ as $k\rightarrow \infty$, then $\gamma_k z\rightarrow e_1$ as $k\rightarrow \infty$. Also, $\gamma_k z\rightarrow e_3$ as $k\rightarrow -\infty$.
			\item $\alpha^k\beta^{-k}\rightarrow b\in\C^\ast$ as $k\rightarrow \infty$, then $\gamma_k z\rightarrow \corchetes{b:z_2:0}$ as $k\rightarrow \infty$, contradicting that $\Gamma$ is discrete.
			\item $\alpha^k\beta^{-k}\rightarrow 0$ as $k\rightarrow \infty$, $\gamma_k z\rightarrow e_2$ as $k\rightarrow \infty$ and $\gamma_k z\rightarrow e_2$ as $k\rightarrow \infty$, and $\gamma_k z\rightarrow e_3$ as $k\rightarrow -\infty$.
			\end{itemize}
		\end{enumerate}		 
This concludes the proof.
\end{proof}

If $\Gamma=\Gamma_{\alpha,\beta}\subset U_+$ is a discrete group containing loxodromic elements, then, putting together Lemmas \ref{lem_caso_diagonal_L0} and \ref{lem_caso_diagonal_L1} we have the following cases:
If $\alpha^n=\beta^m$ for some $n,m\in\Z$:
	\begin{itemize}
	\item[\titem{[D1]}] $L_0(\Gamma)\cup L_1(\Gamma)=\linproy{e_1,e_2}\cup\set{e_3}$, if $\valorabs{\alpha}>1>\valorabs{\beta}$ or $\valorabs{\alpha}<1<\valorabs{\beta}$.
	\item[\titem{[D2]}] $L_0(\Gamma)\cup L_1(\Gamma)=\linproy{e_1,e_2}\cup\set{e_3}$, if $\valorabs{\alpha}>\valorabs{\beta}>1$ or $\valorabs{\alpha}<\valorabs{\beta}<1$.
	\end{itemize}
If there are no integers $n,m$ such that $\alpha^n=\beta^m$:
	\begin{itemize}
	\item[\titem{[D3]}] $L_0(\Gamma)\cup L_1(\Gamma)=\set{e_1,e_2,e_3}$, if $\valorabs{\alpha}>1>\valorabs{\beta}$ or $\valorabs{\alpha}<1<\valorabs{\beta}$.
	\item[\titem{[D4]}] $L_0(\Gamma)\cup L_1(\Gamma)=\set{e_1,e_2,e_3}$, if $\valorabs{\alpha}>\valorabs{\beta}>1$ or $\valorabs{\alpha}<\valorabs{\beta}<1$.
	\item[\titem{[D5]}] $L_0(\Gamma)\cup L_1(\Gamma)=\linproy{e_1,e_2}\cup\linproy{e_2,e_3}$, if $\beta$ is an irrational rotation.
	\end{itemize}

The following proposition gives the full description of the Kulkarni limit set for groups in this case, we will use the notation described in the previous paragraph.

\begin{prop}\label{prop_kulkarni_diagonales}
Let $\Gamma_{\alpha,\beta}\subset U_+$ be a discrete group containing loxodromic elements, then
	\begin{enumerate}[(i)]
	\item $\KulL(\Gamma)=\linproy{e_1,e_2}\cup\set{e_3}$ in Cases \titem{[D1]} and \titem{[D2]}.
	\item $\KulL(\Gamma)=\set{e_1,e_2,e_3}$ in Cases \titem{[D3]} and \titem{[D4]}.
	\item $\KulL(\Gamma)=\linproy{e_1,e_2}\cup\linproy{e_2,e_3}$ in Case \titem{[D5]}.
	\end{enumerate}
\end{prop}

\begin{proof}
We consider a sequence $\set{\gamma^k}\subset\Gamma$. We will determine the quasi-projective limits of this sequence using Table \ref{fig_cdiag_casos_qp} and then, using the $\lambda$-lemma (Proposition \ref{prop_descripcion_Eq}) we determine the set $L_2(\Gamma)$ and therefore, $\KulL(\Gamma)$.
	\begin{itemize}
	\item In Case \titem{[D1]}, $\alpha^k\rightarrow\infty$ and $\beta^k\rightarrow 0$ and the quasi-projective limit is given by (i) of Table \ref{fig_cdiag_casos_qp}. Therefore, the orbits of compact subsets of $\CP\setminus L_0(\Gamma)\cup L_1(\Gamma)$ accumulate on $\set{e_1}$ and $\set{e_2}$ under the sequences $\set{\gamma^k}$ and $\set{\gamma^{-k}}$ respectively. Then $\KulL(\Gamma)=\linproy{e_1,e_2}\cup\set{e_3}$.
	\item In Case \titem{[D2]}, $\alpha^k\rightarrow\infty$, $\beta^k\rightarrow \infty$ and $\alpha^k\beta^{-k}\rightarrow\infty$, then the quasi-projective limit is also given by (i) of Table \ref{fig_cdiag_casos_qp}. We have again, $\KulL(\Gamma)=\linproy{e_1,e_2}\cup\set{e_3}$.
	\item In Cases \titem{[D3]} and \titem{[D4]}, the orbits of compact subsets of $\CP\setminus L_0(\Gamma)\cup L_1(\Gamma)$ accumulate on $\set{e_1}$ (and accumulate on $\set{e_1}$ under the sequence of $\set{\gamma^{-k}}$). Then, $\KulL(\Gamma)=\set{e_1,e_2,e_3}$.
	\item In Case \titem{[D5]}, we have $\valorabs{\beta}=1$ and we can assume without loss of generality that $\valorabs{\alpha}>1$. Then, the quasi-projective limit of the sequence $\set{\gamma_k}$ is $\tau=\text{Diag}(1,0,0)$ and therefore, the orbits of compact subsets of $\CP^2\setminus \kernel(\tau)=\CP^2\setminus\linproy{e_2,e_3}$ accumulate on $\text{Im}(\tau)=\set{e_1}$. If $K\subset\CP^2\setminus\parentesis{\linproy{e_1,e_2}\cup\linproy{e_2,e_3}}$, then $K\subset\CP^2\setminus\linproy{e_2,e_3}$. Therefore $L_2(\Gamma)=\set{e_1}$, we finally conclude that
		$$\KulL(\Gamma)=\linproy{e_1,e_2}\cup\linproy{e_2,e_3}.$$
	\end{itemize}
This concludes the proof.
\end{proof}

\section{Non-commutative triangular groups}\label{sec_non_commutative_triangular}

In this section we describe non-commutative upper triangular discrete subgroups of $\PSL$. The main theorem of this Section describes these groups (Theorem \ref{thm_descomposicion_caso_noconmutativo}). Before we state it and prove it, we will need to prove several propositions and lemmas.

\subsection{Restrictions on the elements of a non-commutative group}

The following propositions and corollaries hint that complex homotheties, irrational screws and irrational ellipto-parabolic elements impose strong restrictions on the groups they belong to, if we want these groups to be discrete. These propositions generalize the results presented in Chapter 5 of \cite{tesisadriana} (see for example, Propositions 5.12 and 5.13 of \cite{tesisadriana}, which are generalized in Propositions \ref{prop_HC_discreto_invariante} and \ref{prop_IS_discreto_invariante} respectively).\\

We start studying the case of complex homotheties.

\begin{prop}\label{prop_HC_discreto_invariante}
Let $\gamma\in \PSL$ be a complex homothety given by 
	$$\gamma=\text{Diag}(\lambda^{-2},\lambda,\lambda)$$ 
for some $\lambda\in\C^\ast$ with $\valorabs{\lambda}\neq 1$. Let $\alpha\in U_{+}\setminus\prodint{\gamma}$ such that $\alpha$ is neither a complex homothety or a screw, then the group $\prodint{\alpha,\gamma}$ is discrete if and only if $\alpha$ leaves $\KulL(\gamma)$ invariant. 
\end{prop}

\begin{proof}
We have that $\KulL(\gamma)=\set{e_1}\cup\linproy{e_2,e_3}$ (see Proposition 4.2.23 of \cite{ckg_libro}). Denote $\alpha=\corchetes{\alpha_{ij}}$, a straight forward calculation shows that $\alpha$ leaves $\linproy{e_2,e_3}$ invariant if and only if $\alpha_{12}=\alpha_{13}=0$.\\

First we prove that if $\prodint{\alpha,\gamma}$ is discrete then $\alpha$ leaves $\KulL(\gamma)$ invariant. Assume that $\alpha$ doesn't leave $\KulL(\gamma)$ invariant, this means that $\valorabs{\alpha_{12}}+\valorabs{\alpha_{13}}\neq 0$. Let $\set{g_n}\subset\prodint{\alpha,\gamma}$ be the sequence of distinct elements given by
	$$g_n:=\gamma^n \alpha \gamma^{-n}=\corchetes{\begin{array}{ccc}
	\alpha_{11} & \frac{\alpha_{12}}{\lambda^{3n}} & \frac{\alpha_{23}}{\lambda^{3n}}\\
	0 & \alpha_{22} & \alpha_{23} \\
	0 & 0 & \alpha_{33}\\
	\end{array}}.$$
Clearly, 
	$$g_n\rightarrow \corchetes{\begin{array}{ccc}
	\alpha_{11} & 0 & 0\\
	0 & \alpha_{22} & \alpha_{23} \\
	0 & 0 & \alpha_{33}\\
	\end{array}}\in\PSL.$$
Therefore $\prodint{\alpha,\gamma}$ is not discrete.\\

Now we prove that if $\alpha$ leaves $\KulL(\gamma)$ invariant then $\prodint{\alpha,\gamma}$ is discrete. We rewrite $\gamma=\text{Diag}(\lambda^{-3},1,1)$. Suppose that $\prodint{\alpha,\gamma}$ is not discrete, then there exists a sequence of distinct elements $\set{w_k}\subset\prodint{\alpha,\gamma}$ such that $w_k\rightarrow \id$. Since $\alpha$ leaves $\KulL(\gamma)$ invariant then
	$$\alpha=\corchetes{\begin{array}{ccc}
	\alpha_{11} & 0 & 0\\
	0 & \alpha_{22} & \alpha_{23} \\
	0 & 0 & \alpha_{33}\\
	\end{array}}$$
and therefore $\corchetes{\alpha,\gamma}=\id$. Then the sequence $w_k$ can be expressed as 
	\begin{equation}\label{eq_dem_prop_HC_discreto_invariante_1}
	w_k=\gamma^{i_k}\alpha^{j_k}.
	\end{equation}
Observe that $h\gamma h^{-1}=\gamma$ for any $h\in\PSL$ with the form
	$$h=\corchetes{\begin{array}{ccc}
	h_{11} & 0 & 0\\
	0 & h_{22} & h_{23} \\
	0 & 0 & h_{33}\\
	\end{array}}.$$
Furthermore, there exists an element $h\in\PSL$ with the previous form such that $h\alpha h^{-1}$ has one of the following forms:
	\begin{enumerate}
	\item 
		$$h\alpha h^{-1}=\corchetes{\begin{array}{ccc}
	\alpha_{22}^{-1}\alpha_{33}^{-1} & 0 & 0\\
	0 & \alpha_{22} & 0 \\
	0 & 0 & \alpha_{33}\\
	\end{array}},$$
	if $\Pi(\alpha)$ is loxodromic. In this case, from the previous equation and (\ref{eq_dem_prop_HC_discreto_invariante_1}) it follows
	$$w_k=\corchetes{\begin{array}{ccc}
	\lambda^{-3i_k}\alpha_{22}^{-j_k}\alpha_{33}^{-j_k} & 0 & 0\\
	0 & \alpha_{22}^{j_k} & 0 \\
	0 & 0 & \alpha_{33}^{j_k}\\
	\end{array}}\rightarrow \id.$$
	Therefore $\alpha_{22}^{j_k}$, $\alpha_{33}^{j_k}\rightarrow 1$, but then we cannot have 
		$$\lambda^{-3i_k}\alpha_{22}^{-j_k}\alpha_{33}^{-j_k}\rightarrow 1,$$
	since $\valorabs{\lambda}\neq 1$.
	\item $$h\alpha h^{-1}=\corchetes{\begin{array}{ccc}
	\alpha_{11} & 0 & 0\\
	0 & 1 & 1\\
	0 & 0 & 1\\
	\end{array}},$$
	if $\Pi(\alpha)$ is parabolic. Then we have
	$$w_k=\corchetes{\begin{array}{ccc}
	\lambda^{-3i_k}\alpha_{11}^{j_k} & 0 & 0\\
	0 & 1 & 1 \\
	0 & 0 & 1\\
	\end{array}}\rightarrow \id.$$
	But this is clearly impossible.	
	\item $$h\alpha h^{-1}=\corchetes{\begin{array}{ccc}
	\alpha_{11} & 0 & 0\\
	0 & \alpha_{22} & 0 \\
	0 & 0 & \alpha_{22}\\
	\end{array}},$$
	with $\valorabs{\alpha_{22}}=1$, if $\Pi(\alpha)$ is elliptic. But in this case $\alpha$ is a complex homothety or a screw.
	\end{enumerate}	 
Therefore neither of the three cases can occur and thus, $\prodint{\alpha,\gamma}$ is discrete.
\end{proof}

From the proof of the previous proposition, we have the following imme-diate consequence.

\begin{cor}\label{cor_forma_afin_HC}
If $\Gamma\subset U_+$ is a discrete subgroup and $\Gamma$ contains a complex homothety as in the previous proposition, then every element of $\Gamma$ has the form
	$$\alpha= \corchetes{\begin{array}{ccc}
	\alpha_{11} & 0 & 0\\
	0 & \alpha_{22} & \alpha_{23} \\
	0 & 0 & \alpha_{33}\\
	\end{array}}.$$ 
\end{cor}

\begin{prop}\label{prop_HC_conmutativo}
Let $\Gamma\subset U_+$ be a discrete group containing a complex homothety as in the previous proposition, if the control group $\Pi(\Gamma)$ is discrete then 
	\begin{enumerate}
	\item $\Pi(\Gamma)$ is purely parabolic or purely loxodromic.
	\item $\Gamma$ is commutative.
	\end{enumerate}
\end{prop}

\begin{proof}
Using Corollary \ref{cor_forma_afin_HC} we know that every element of $\Gamma$ has the form
	$$\alpha= \corchetes{\begin{array}{ccc}
	\alpha_{11} & 0 & 0\\
	0 & \alpha_{22} & \alpha_{23} \\
	0 & 0 & \alpha_{33}\\
	\end{array}}.$$ 
Consider the elements of the control group $\Sigma=\Pi(\Gamma)$:
\begin{enumerate}
\item If $\Sigma$ contains an elliptic element 
	$$\Pi(\alpha)=\corchetes{\begin{array}{cc}
	e^{2\pi i \theta_1} & 0\\
	0 & e^{2\pi i \theta_2} \\
	\end{array}},$$
then
	$$\alpha= \corchetes{\begin{array}{ccc}
	e^{-2\pi i (\theta_1+\theta_2)} & 0 & 0\\
	0 & e^{2\pi i \theta_1} & 0 \\
	0 & 0 & e^{2\pi i \theta_2}\\
	\end{array}}\in\Gamma$$
is elliptic. If $\alpha$ has infinite order, this contradicts that $\Gamma$ is discrete and if $\alpha$ has finite order, this contradicts that $\Gamma$ is torsion free. Therefore, $\Sigma$ has no elliptic elements.  	  
\item If $\Sigma$ has a parabolic element $\Pi(\alpha)$ then $\Sigma$ is purely parabolic. In order to verify this, assume that there exists $\beta\in\Gamma$ such that $\Pi(\beta)$ is loxodromic then, denoting $A=\Pi(\alpha)$ and $B=\Pi(\beta)$, a direct calculation shows that $[A,B]$ is parabolic and
	$$\text{Tr}(A)=\text{Tr}\parentesis{[A,B]}=2$$
then
	$$\valorabs{\text{Tr}(A)^2-4}+\valorabs{\text{Tr}\parentesis{[A,B]}-2}\geq 1.$$
Using the J\o rgensen inequality, it follows that $\prodint{A,B}\subset\Sigma$ is not discrete, contradicting that $\Sigma$ is discrete. This verifies that $\Sigma$ is purely parabolic. Then every element of $\Gamma$ has the form
	$$\alpha= \corchetes{\begin{array}{ccc}
	\alpha_{11} & 0 & 0\\
	0 & 1 & \alpha_{23} \\
	0 & 0 & 1\\
	\end{array}}.$$
A simple calculation shows that $\Gamma$ is commutative.  
\item If $\Sigma$ has a loxodromic element $\Pi(\alpha)$, the previous case shows that $\Sigma$ is purely loxodromic. Let $\Pi(\beta)\in\Sigma\setminus\set{\alpha}$ be another loxodromic element, observe that $\Pi(\alpha)$ and $\Pi(\beta)$ share at least one fixed point, since they both are upper triangular elements. We have two cases:
	\begin{enumerate}
	\item If $\Pi(\alpha)$ and $\Pi(\beta)$ share exactly one fixed point then they have the form
		$$\begin{array}{cc}
		\Pi(\alpha)=\corchetes{\begin{array}{cc}
	\alpha_{22} & \alpha_{23}\\
	0 & \alpha_{33} \\
	\end{array}}, & \Pi(\beta)=\corchetes{\begin{array}{cc}
	\beta_{22} & \beta_{23}\\
	0 & \beta_{33} \\
	\end{array}}
		\end{array}$$
	with $\alpha_{23},\beta_{23}\neq 0$ and 
	$$[\alpha_{23}:\alpha_{33}-\alpha_{22}]\neq[\beta_{23}:\beta_{33}-\beta_{22}].$$
Then the element $[\Pi(\alpha),\Pi(\beta)]\neq \id$ is parabolic, contradicting that $\Sigma$ is purely loxodromic. This shows that this case cannot happen.
	\item If $\fix(\Pi(\alpha))=\fix(\Pi(\beta))$, then every element of $\Gamma$ has the form
	$$\alpha= \corchetes{\begin{array}{ccc}
	\alpha_{11} & 0 & 0\\
	0 & 1 & \alpha_{23} \\
	0 & 0 & 1\\
	\end{array}},$$
with $\fix(\alpha_1)=\set{e_1,e_2}$ and $\fix(\alpha_2)=\set{e_2,p}$ for some point $p$. Then, using the calculations described in the proof of proposition \ref{prop_formas_conmutativas} we see that $\Gamma$ is commutative. 	
\end{enumerate} 
	
	\end{enumerate}
Then, under the hypothesis of the proposition, $\Sigma$ is purely parabolic or purely loxodromic and in both cases, $\Gamma$ is commutative.
\end{proof}

Last proposition describes the restrictions complex homotheties impose on discrete subgroups of $U_+$ whose control group is discrete. The following proposition analyzes the case of discrete subgroups of $U_+$ with non-discrete control group.

\begin{prop}\label{prop_HC_conmutativo_control_no_discreto}
Let $\Gamma\subset U_+$ be a discrete group such that $\Sigma=\Pi(\Gamma)$ is not discrete and $\Lambda_{Gr}(\Sigma)=\Ss^1$. Then $\Gamma$ cannot contain a complex homothety with the form $\text{Diag}(\lambda^{-2},\lambda,\lambda)$, $\valorabs{\lambda}\neq 1$.
\end{prop}

\begin{proof}
Suppose that $\Gamma$ contains such a complex homothety. Then, by Corollary \ref{cor_forma_afin_HC}, each element in $\Gamma$ has the form
	\begin{equation}\label{eq_dem_prop_HC_conmutativo_control_no_discreto_2}
	\alpha= \corchetes{\begin{array}{ccc}
	\alpha_{11} & 0 & 0\\
	0 & \alpha_{22} & \alpha_{23} \\
	0 & 0 & \alpha_{33}\\
	\end{array}}.
	\end{equation} 
Since $\Lambda_{Gr}(\Sigma)=\Ss^1$, then $\Gamma$ is a non-elemental, non-discrete group and therefore, $\Lambda_{Gr}(\Sigma)$ is the closure of fixed points of loxodromic elements of $\Sigma$ (Theorem \ref{thm_no_discreto_no_elem_cerradura}). This means that there are an infinite number of loxodromic elements in $\Sigma$ sharing exactly one fixed point. Let $f\in\Sigma$ such that $\Pi(f)$ is a loxodromic element and
	$$f= \corchetes{\begin{array}{ccc}
	\alpha^{-1} & 0 & 0\\
	0 & \alpha & 0 \\
	0 & 0 & 1\\
	\end{array}},$$
with $\valorabs{\alpha}<1$. Let $h_1,h_2\in\Gamma$ such that $\Pi(h_1),\Pi(h_2)$ are loxodromic elements, 
	\begin{equation}\label{eq_dem_prop_HC_conmutativo_control_no_discreto_1}
	\fix\parentesis{\Pi(h_1)}\neq \fix\parentesis{\Pi(h_2)}
	\end{equation}
and $h_i\neq f$, for $i=1,2$. As a consequence of (\ref{eq_dem_prop_HC_conmutativo_control_no_discreto_1}), $h:=\corchetes{h_1,h_2}\neq \id$ and then
	$$h=\corchetes{\begin{array}{ccc}
	1 & 0 & 0\\
	0 & 1 & h_0 \\
	0 & 0 & 1\\
	\end{array}}$$
is a parabolic element. Let $\set{f^k h f^{-k}}\subset\Gamma$ be the sequence given by
	$$f^k h f^{-k}=\corchetes{\begin{array}{ccc}
	1 & 0 & 0\\
	0 & 1 & \alpha^k h_0 \\
	0 & 0 & 1\\
	\end{array}}.$$
Then $f^k h f^{-k}\rightarrow\id$, contradicting that $\Gamma$ is discrete. If $\valorabs{\alpha}>1$, we consider the sequence $\set{f^{-k} h f^{k}}$ instead and get the same contradiction. This proves that $\Gamma$ cannot contain a complex homothety with the form $\text{Diag}(\lambda^{-2},\lambda,\lambda)$.
\end{proof}

\begin{prop}\label{prop_HC_conmutativo_control_no_discreto_grenlim_2ptos}
Let $\Gamma\subset U_+$ be a non-commutative discrete group such that $\Sigma=\Pi(\Gamma)$ is not discrete and $\valorabs{\grL(\Sigma)}=2$. Then $\Gamma$ cannot contain a complex homothety with the form $\text{Diag}(\lambda^{-2},\lambda,\lambda)$, $\valorabs{\lambda}\neq 1$.
\end{prop}

\begin{proof}
Suppose that $\Gamma$ contains such a complex homothety. Then, by Corollary \ref{cor_forma_afin_HC}, each element in $\Gamma$ has the form
	\begin{equation}\label{eq_dem_prop_HC_conmutativo_control_no_discreto_grenlim_2ptos_1}
	\alpha= \corchetes{\begin{array}{ccc}
	\alpha_{11} & 0 & 0\\
	0 & \alpha_{22} & \alpha_{23} \\
	0 & 0 & \alpha_{33}\\
	\end{array}}.
	\end{equation} 
If $\valorabs{\grL(\Sigma)}=2$, then, up to conjugation, every element of $\Sigma$ has the form $\text{Diag}(\beta,\delta)$ for some $\beta,\delta\in\C^\ast$. Then, using \ref{eq_dem_prop_HC_conmutativo_control_no_discreto_grenlim_2ptos_1}, it follows that each element in $\Gamma$ has the form
	$$\alpha= \corchetes{\begin{array}{ccc}
	\alpha_{11} & 0 & 0\\
	0 & \alpha_{22} & 0 \\
	0 & 0 & \alpha_{33}\\
	\end{array}}.$$
In other words, every element of $\Gamma$ is diagonal and hence, $\Gamma$ would be commutative. This contradiction proves the proposition.
\end{proof}

\begin{prop}\label{prop_HC_conmutativo_control_no_discreto_grenlim_1pto}
Let $\Gamma\subset U_+$ be a non-commutative discrete group such that $\Sigma=\Pi(\Gamma)$ is not discrete and $\valorabs{\grL(\Sigma)}=1$. Then $\Gamma$ cannot contain a complex homothety with the form $\text{Diag}(\lambda^{-2},\lambda,\lambda)$, $\valorabs{\lambda}\neq 1$.
\end{prop}

\begin{proof}
Suppose that there exists such a complex homothety in $\Gamma$, denote it by $\gamma$. Then, by Corollary \ref{cor_forma_afin_HC}, each element in $\Gamma$ has the form given by Equation (\ref{eq_dem_prop_HC_conmutativo_control_no_discreto_grenlim_2ptos_1}) of the proof of the previous proposition.\\

Since $\Sigma=\Pi(\Gamma)$ is not discrete and $\valorabs{\grL(\Sigma)}=1$, then Proposition \ref{prop_opciones_conj_lim_greenberg} implies that $\Sigma$ is a dense subgroup of $\epa$ containing parabolic elements. Then, every element of $\Sigma$ looks like 
	$$\Pi(\mu)= \corchetes{\begin{array}{cc}
	a & b \\
	0 & a^{-1} \\
	\end{array}},$$
for $\valorabs{a}=1$ and $b\in\C^\ast$. This, together with (\ref{eq_dem_prop_HC_conmutativo_control_no_discreto_grenlim_2ptos_1}), implies that every element of $\Gamma$ has the form
	$$\mu= \corchetes{\begin{array}{ccc}
	\alpha^{-2} & 0 & 0\\
	0 & \alpha a & \alpha b \\
	0 & 0 & \alpha a^{-1} \\
	\end{array}},$$
for some $\alpha\in\C^\ast$. Now, we have two possibilities:
\begin{enumerate}[(1)]
\item If every element of $\Sigma$ is parabolic, then every element of $\Gamma\setminus\set{\id}$ has the form
	$$\mu= \corchetes{\begin{array}{ccc}
	\alpha^{-2} & 0 & 0\\
	0 & \alpha & \alpha b \\
	0 & 0 & \alpha \\
	\end{array}}\text{, for some }\alpha\in\C^\ast.$$
A direct calculation shows that $\Gamma$ would be commutative in this case, this contradiction doesn't allow this case to happen.
\item If there is an elliptic element $\Pi(\gamma)\in\Sigma$ then
	$$\Pi(\gamma)= \corchetes{\begin{array}{cc}
	e^{2\pi i \theta} & 0 \\
	0 & e^{-2\pi i \theta} \\
	\end{array}}\text{, for some }\theta\nin\Z.$$
	Then 
	\begin{equation}\label{eq_dem_prop_HC_conmutativo_control_no_discreto_grenlim_1pto_1}
	\gamma= \corchetes{\begin{array}{ccc}
	\beta^{-2} & 0 & 0\\
	0 & \beta e^{2\pi i \theta} & 0 \\
	0 & 0 & \beta e^{-2\pi i \theta}\\
	\end{array}},
	\end{equation}		
	for some $\beta\in\C^{\ast}$. Observe that, if $\theta\in\Q$ then $\lambda_{23}(\gamma)=e^{4\pi i \theta}$ would be a torsion element in $\lambda_{23}(\Gamma)$. This contradictions implies that, if there exists an elliptic element in $\Sigma$, then $\theta\in\R\setminus\Q$.\\
	
	On the other hand, since $\Sigma$ is not discrete, there is a sequence of distinct elements $\set{\Pi(\mu_k)}\subset\Sigma$ such that $\Pi(\mu_k)\rightarrow\id$. We have the following two cases:
		\begin{enumerate}[(i)]
		\item If the sequence $\set{\Pi(\mu_k)}$ contains an infinite number of parabolic elements, then, considering an adequate subsequence, we can assume that $\set{\Pi(\mu_k)}$ is a sequence of distinct parabolic elements and we denote
			$$\Pi(\mu_k)= \corchetes{\begin{array}{cc}
	1 & b_k \\
	0 & 1 \\
	\end{array}},$$
with $\set{b_k}\subset\C^{\ast}$ a sequence of distinct elements such that $b_k\rightarrow 0$. Then,
	$$\mu_k = \corchetes{\begin{array}{ccc}
	\alpha_k^{-2} & 0 & 0\\
	0 & \alpha_k & \alpha_k b_k \\
	0 & 0 & \alpha_k\\
	\end{array}},$$
for some $\set{\alpha_k}\subset \C^{\ast}$. Let $\set{\xi_k}\subset\Gamma$ be the sequence of distinct elements given by
	$$\xi_k = \corchetes{\gamma,\mu_k}= \corchetes{\begin{array}{ccc}
	1 & 0 & 0\\
	0 & 1 & b_k\parentesis{1- e^{-4\pi i \theta}} \\
	0 & 0 & 1\\
	\end{array}},$$
where $\gamma$ is an elliptic element in $\Sigma$ given by (\ref{eq_dem_prop_HC_conmutativo_control_no_discreto_grenlim_1pto_1}). Then $\xi_k\rightarrow \id$, which is impossible, since $\Gamma$ is discrete.
		\item  If the sequence $\set{\Pi(\mu_k)}$ contains only a finite number of parabolic elements, then we can assume that the whole sequence $\set{\Pi(\mu_k)}$ is a sequence of irrational elliptic elements. We denote,
		$$\Pi(\mu_k)= \corchetes{\begin{array}{cc}
	e^{2\pi i \theta_k} & b_k \\
	0 & e^{-2\pi i \theta_k} \\
	\end{array}},$$
	with $b_k\rightarrow 0$. Since $\set{\theta_k}\subset\R\setminus\Q$, we can pick an adequate subsequence of $\set{\Pi(\mu_k)}$, still denoted in the same way, such that $\theta_k\rightarrow 0$ by distinct elements $\set{\theta_k}$. Let $\set{\Pi(\sigma_k)}\subset\Sigma$ the sequence of distinct elements given by
		\begin{align*}
		\Pi(\sigma_k) &= \corchetes{\Pi(\mu_k),\Pi(\mu_{k+1})}\\
		&=\corchetes{\begin{array}{cc}
	1 & e^{-4\pi i (\theta_k+\theta_{k+1})\parentesis{b_k e^{2\pi i \theta_k}\parentesis{1- e^{4\pi i \theta_{k+1}}}-b_{k+1} e^{2\pi i \theta_{k+1}}\parentesis{1- e^{4\pi i \theta_k}}}} \\
	0 & 1 \\
	\end{array}}.
		\end{align*}
		The sequence $\set{\Pi(\sigma_k)}$ is made up of distinct parabolic elements of $\Sigma$, and since $\theta_k\rightarrow 0$, we have that
			$$e^{-4\pi i (\theta_k+\theta_{k+1})\parentesis{b_k e^{2\pi i \theta_k}\parentesis{1- e^{4\pi i \theta_{k+1}}}-b_{k+1} e^{2\pi i \theta_{k+1}}\parentesis{1- e^{4\pi i \theta_k}}}}\rightarrow 0.$$	
		This implies that $\Pi(\sigma_k)\rightarrow\id$. Applying the same argument as in the last case (i), we get a contradiction. 
		\end{enumerate}
		The previous two contradictions imply that this case cannot happen either.
\end{enumerate}
	This concludes the proof.
\end{proof}

Together, Propositions \ref{prop_HC_conmutativo}, \ref{prop_HC_conmutativo_control_no_discreto}, \ref{prop_HC_conmutativo_control_no_discreto_grenlim_2ptos}, \ref{prop_HC_conmutativo_control_no_discreto_grenlim_1pto} and Lemma \ref{lem_control_nodiscreto_descartar_grL_vacio} imply the following conclusion.

\begin{cor}\label{cor_HC_no_hay_en_no_conmutativos}
Let $\Gamma\subset U_+$ be a non-commutative, torsion-free discrete subgroup, then $\Gamma$ cannot contain a complex homothety with the form $\text{Diag}(\lambda^{-2},\lambda,\lambda)$, $\valorabs{\lambda}\neq 1$.  
\end{cor}

The proof of Proposition \ref{prop_HC_discreto_invariante} can be repeated in a similar way to get the same conclusion but with a complex homothety with the form $\gamma=\text{Diag}(\lambda,\lambda,\lambda^{-2})$, $\valorabs{\lambda}\neq 1$. Repeating the proofs of Corollary \ref{cor_forma_afin_HC} and Propositions \ref{prop_HC_conmutativo}, \ref{prop_HC_conmutativo_control_no_discreto}, \ref{prop_HC_conmutativo_control_no_discreto_grenlim_2ptos} and \ref{prop_HC_conmutativo_control_no_discreto_grenlim_1pto}, we can prove the following corollary.

\begin{cor}\label{cor_HC3_no_hay_en_no_conmutativos}
Let $\Gamma\subset U_+$ be a non-commutative, torsion-free discrete subgroup, then $\Gamma$ cannot contain a complex homothety with the form $\text{Diag}(\lambda,\lambda,\lambda^{-2})$, $\valorabs{\lambda}\neq 1$.  
\end{cor}

\begin{obs}
Corollary \ref{cor_HC_no_hay_en_no_conmutativos} state that, the presence of certain diagonal complex homotheties in discrete subgroups of $U_+$ imply that $\Gamma$ is commutative. The same can be said about general complex homotheties with the same form in their eigenvalues. 
\end{obs}

\begin{lem}\label{lem_parte_parab_no_discreta}
Let $\Sigma$ be a non-discrete, upper triangular subgroup of $\psl$ such that $\Lambda_{Gr}(\Sigma)=\Ss^1$. Then the parabolic part of $\Sigma$ is a non-discrete group.
\end{lem}

\begin{proof}
Let $\Sigma_p$ be the parabolic part of $\Sigma$. It is straightforward to verify that $\Sigma_p$ is a subgroup of $\Sigma$. We will find a sequence $\set{f_k}\subset\Sigma_p$ such that $f_k\rightarrow\id$.\\

Let $g\in\Sigma$ be a loxodromic element such that 
	$$g=\corchetes{\begin{array}{cc}
	\alpha & 0 \\
	0 & \alpha^{-1} \\
	\end{array}},$$
for $\valorabs{\alpha}<1$. Let $h_1,h_2\in\Sigma$ loxodromic elements such that
	$$\fix(h_1)\neq\fix(h_2)\;\;\;\;\;\;\fix(h_i)\neq\fix(g).$$
Let $h=\corchetes{h_1,h_2}$, as in the proof of Proposition \ref{prop_HC_conmutativo_control_no_discreto}, $h\neq\id$ is a parabolic element in $\Sigma$,
	$$h=\corchetes{\begin{array}{cc}
	1 & h_0 \\
	0 & 1 \\
	\end{array}}.$$
Hence,
	$$f_k:=g^k h g^{-k}=\corchetes{\begin{array}{cc}
	1 & \alpha^{2k} h_0 \\
	0 & 1 \\
	\end{array}}$$
is the sequence of distinct elements in $\Sigma_p$ such that $f_k\rightarrow\id$. Again, if $\valorabs{\alpha}>1$, we take the sequence $f_k:=g^{-k} h g^{k}$ instead. The fact that one can take the elements $g,h_1,h_2$ is a consequence of 
	$$\Lambda_{Gr}(\Sigma)=\Ss^1.$$
This finishes the proof of the lemma.
\end{proof}

Now we look at the case of irrational screws.

\begin{prop}\label{prop_IS_discreto_invariante}
Let $\gamma\in U_+$ be an irrational screw given by
	$$\gamma=\corchetes{\begin{array}{ccc}
	\beta^{-2} & 0 & 0\\
	0 & \beta e^{4\pi i \theta} & 0 \\
	0 & 0 & \beta e^{2\pi i \theta}\\
	\end{array}},\;\;\; \valorabs{\beta}\neq 1,\;\;\;\theta\in\R\setminus\Q.$$
Let $\alpha\in U_+\setminus\prodint{\gamma}$, if $\prodint{\alpha,\gamma}$ is discrete then $\alpha$ is diagonal.
\end{prop}

\begin{proof}
Let $\gamma,\alpha\in U_+$ be as in the statement of the proposition, denote $\alpha=\corchetes{\alpha_{ij}}$. Suppose that $\alpha$ is not diagonal, then $\alpha_{12}\neq 0$, $\alpha_{13}\neq 0$ or $\alpha_{23}\neq 0$.\\

Since $\theta\in\R\setminus\Q$, there is a subsequence, still denoted by the index $k$, such that $e^{2 k \pi i \theta}\rightarrow 1$. If $\valorabs{\beta}>1$, consider the sequence $\set{\gamma^{-k}\alpha\gamma^{k}}\subset\prodint{\alpha,\gamma}$ given by
	$$\gamma^k\alpha\gamma^{-k}=\corchetes{\begin{array}{ccc}
	\alpha_{11} & \alpha_{12}\beta^{-3k}e^{-4k\pi i \theta} & \alpha_{13}\beta^{-3k}e^{-2k\pi i \theta}\\
	0 & \alpha_{22} & \alpha_{23}e^{2k\pi i\theta} \\
	0 & 0 & \alpha_{33}\\
	\end{array}}.$$
Since $\alpha_{12}\neq 0$, $\alpha_{13}\neq 0$ or $\alpha_{23}\neq 0$, this is a sequence of distinct elements such that 
	$$\gamma^{-k}\alpha\gamma^{k}\rightarrow\text{Diag}\parentesis{\alpha_{11},\alpha_{22},\alpha_{33}}.$$
Then $\prodint{\alpha,\gamma}$ is not discrete.\\

If $\valorabs{\beta}<1$, consider instead the sequence $\set{\gamma^{k}\alpha\gamma^{-k}}\subset\prodint{\alpha,\gamma}$ and we get the same conclusion. This proofs the proposition.
\end{proof}

Observe that the converse is not necessarily true. Let $\beta\in\C$ with $\valorabs{\beta}\neq 1$ and let $\gamma,\alpha\in U_+$ be two elements given by
	\begin{align*}
	\gamma &=\text{Diag}\parentesis{\beta^{-2\sqrt{2}}, \beta^{\sqrt{2}} e^{4\pi i \sqrt{2}}, \beta^{\sqrt{2}} e^{2\pi i \sqrt{2}}}\\
	\alpha &=\text{Diag}\parentesis{\beta^{-2\sqrt{3}},\beta^{\sqrt{3}} e^{4\pi i \sqrt{3}},\beta^{\sqrt{3}} e^{2\pi i \sqrt{3}}}.
\end{align*}

Then $\alpha$ is diagonal but $\prodint{\alpha,\gamma}$ is not discrete. To verify this observe that the additive group $\text{Span}_\Z\parentesis{\sqrt{2},\sqrt{3}}$ is not discrete and therefore there exist sequences $\set{p_k},\set{q_k}\subset\Z$ such that 
	$$p_k\sqrt{2}+q_k\sqrt{3}\rightarrow 0.$$
The sequence of distinct elements $\set{\gamma^{p_k}\alpha^{q_k}}\subset\prodint{\alpha,\gamma}$ satisfies
	$$\gamma^{p_k}\alpha^{q_k}\rightarrow\id.$$ 

We have the following consequences of proposition \ref{prop_IS_discreto_invariante}.

\begin{cor}\label{cor_IS_conmutativo}
Let $\Gamma\subset U_+$ be a discrete subgroup and $\gamma\in\Gamma$ an irrational screw as in proposition \ref{prop_IS_discreto_invariante}, then $\Gamma$ is commutative.   
\end{cor}

\begin{proof}
Let $\alpha\in U_+\setminus\prodint{\gamma}$ be an arbitrary element. Since $\Gamma$ is discrete, then $\prodint{\alpha,\gamma}$ is discrete and, by Proposition \ref{prop_IS_discreto_invariante}, $\alpha=\text{Diag}\parentesis{\alpha_{11},\alpha_{22},\alpha_{33}}$. Since both $\alpha$ and $\gamma$ are diagonal, $\corchetes{\gamma,\alpha}=\id$ and thus, $\Gamma$ is commutative.
\end{proof}

\begin{cor}\label{cor_IS_RotInf_prohibidas}
Let $\Gamma\subset U_+$ be a discrete subgroup such that $\Sigma=\Pi(\Gamma)$ is non-discrete and $\overline{\Sigma}=\rot$. Then $\Gamma$ is commutative.
\end{cor}

\begin{proof}
The elements of $\Sigma$ have the form
	$$\Pi(\gamma)=\corchetes{\begin{array}{cc}
	 e^{2\pi i \theta} & 0\\
	 0 & e^{-2\pi i \theta} \\
	\end{array}}.$$
Then every element of $\Gamma\setminus\set{\id}$ has the form
	$$\gamma=\corchetes{\begin{array}{ccc}
	\beta^{-2}e^{-2\pi i \theta} & \gamma_{12} & \gamma_{13}\\
	0 & \beta e^{2\pi i \theta} & 0 \\
	0 & 0 & \beta \\
	\end{array}},$$
with $\valorabs{\beta}\neq 1$. Then $\gamma$ is diagonalizable and therefore an irrational screw, then we can assume that $\gamma_{12}=\gamma_{13}=0$ and therefore, 
	$$\gamma=\corchetes{\begin{array}{ccc}
	\beta^{-2}e^{-2\pi i \theta} & 0 & 0\\
	0 & \beta e^{2\pi i \theta} & 0 \\
	0 & 0 & \beta \\
	\end{array}}=\corchetes{\begin{array}{ccc}
	\beta^{-2} & 0 & 0\\
	0 & \beta e^{4\pi i \theta} & 0 \\
	0 & 0 & \beta e^{2\pi i \theta} \\
	\end{array}}.$$
This means that $\Gamma$ contains only screws.\\

If there were a rational screw $\gamma\in\Gamma$ then $\lambda_{23}$ would be a torsion element in $\lambda_{23}(\Gamma)$, contradicting that $\lambda_{23}(\Gamma)$ is torsion free. Then $\Gamma$ contains only irrational screws, it follows from corollary \ref{cor_IS_conmutativo} that $\Gamma$ is commutative.
\end{proof}

\begin{prop}\label{prop_EPI_conmutativo}
Let $\Gamma\subset U_+$ be a discrete group containing an irrational ellipto-parabolic element $\gamma$ with one of the following two forms:
$$\begin{array}{cc}
\gamma = \corchetes{\begin{array}{ccc}
	e^{-4\pi i \theta} & \beta & \gamma \\
	0 & e^{2\pi i \theta} & \mu \\
	0 & 0 & e^{2\pi i \theta}\\
	\end{array}}, &
	\gamma = \corchetes{\begin{array}{ccc}
	e^{2\pi i \theta} & \beta & \gamma \\
	0 & e^{2\pi i \theta} & \mu \\
	0 & 0 & e^{-4\pi i \theta}\\
	\end{array}}\\
	\mu\neq 0 & \beta\neq 0
\end{array},$$
with $\theta\in\R\setminus\Q$. Then $\Gamma$ is commutative. 
\end{prop}

\begin{proof}
If $\gamma$ has the first form then $\lambda_{12}(\gamma)$ has infinite order in $\lambda_{12}(\Gamma)$. If $\gamma$ has the second form then $\lambda_{23}(\gamma)$ has infinite order in $\lambda_{23}(\Gamma)$. Since $\gamma$ is parabolic, it follows from Lemma 5.10 of \cite{ppar} that $\Gamma$ is abelian in both cases.
\end{proof}

\begin{lem}\label{lem_control_nodiscreto_descartar_grL_vacio}
Let $\Gamma\subset U_+$ be a torsion free, non-commutative, discrete group such that $\Sigma=\Pi(\Gamma)$ is not discrete. Then
	$$\grL(\Sigma)\neq \emptyset.$$
\end{lem}

\begin{proof}
Let us suppose that $\grL(\Sigma)=\emptyset$. Using Proposition \ref{prop_opciones_conj_lim_greenberg}, up to conjugation, we have the following two possibilities:
	\begin{enumerate}[(i)]
	\item $\overline{\Sigma}=\so$. Since $\Gamma\subset U_+$, $\Gamma$ is solvable and, by Proposition \ref{prop_cerradura_solvable}, $\overline{\Sigma}$ is solvable. On the other hand, $\so$ is not solvable. Therefore this case cannot occur.
	\item $\overline{\Sigma}=\dih$ or $\overline{\Sigma}=\rot$. By definition, 
		$$\dih=\prodint{\rot,\;z\mapsto -z}.$$
		Then $\Sigma\cap \rot\neq \emptyset$, otherwise $\Sigma$ would be discrete. Let $\Pi(\gamma)\in\Sigma\cap \rot$ and therefore
			$$\Pi(\gamma)=\corchetes{\begin{array}{cc}
	 e^{2\pi i \theta} & 0\\
	 0 & e^{-2\pi i \theta} \\
	\end{array}}$$
		and
			$$\gamma=\corchetes{\begin{array}{ccc}
			\beta^{-2}e^{-2\pi i \theta} & \gamma_{12} & \gamma_{13} \\
			0 & \beta e^{4\pi i \theta} & 0 \\
			0 & 0 & \beta\\
			\end{array}}.$$
		If $\theta\in\Q$ then $\lambda_{23}(\gamma)$ is a torsion element in $\lambda_{23}(\Gamma)$, contradicting that $\Gamma$ is torsion free. Then $\theta\in\R\setminus\Q$, hence $\gamma$ is an irrational screw. Conjugating $\gamma$ by a suitable element of $\PSL$ and applying Corollary \ref{cor_IS_RotInf_prohibidas}, we conclude that $\Gamma$ is commutative. This contradiction means that this case cannot occur either.
	\end{enumerate}
	These contradictions verify the lemma.
\end{proof}

Propositions \ref{prop_HC_conmutativo}, \ref{cor_IS_conmutativo} and \ref{prop_EPI_conmutativo} say that whenever a discrete subgroup of $U_+$ contains either a complex homothety, an irrational screw or an irrational ellipto-parabolic element, the group has to be commutative.\\

Before stating the main theorem of this section, we need to define and study an important purely parabolic subgroup of $\Gamma$, which determines the dynamics of $\Gamma$. Let us define 
	$$\core(\Gamma)=\kernel(\Gamma)\cap\kernel(\lambda_{12})\cap\kernel(\lambda_{23}).$$
We denote the elements of $\core(\Gamma)$ by $g_{x,y}$, with
	$$g_{x,y}=\corchetes{\begin{array}{ccc}
	1 & x & y \\
	0 & 1 & 0 \\
	0 & 0 & 1\\
	\end{array}}.$$

Lemma 7.11 of \cite{ppar} states that 
	$$\KulL\parentesis{\core(\Gamma)}=\bigcup_{g_{x,y}\in\core(\Gamma)} \linproy{e_1,[0:-y:x]}.$$
We denote this pencil of lines by $\mathcal{C}(\Gamma)=\KulL\parentesis{\core(\Gamma)}$.

\begin{prop}\label{prop_cono_invariante}
Let $\Gamma\subset U_+$ be a discrete group, then every element of $\Gamma$ leaves $\mathcal{C}(\Gamma)$ invariant.
\end{prop}

\begin{proof}
For $g_{x,y}\in\core(\Gamma)$, denote by 
		$$\ell_{x,y}=\linproy{e_1,[0:-y:x]}$$
	the line in $\mathcal{C}(\Gamma)$ determined by the element $g_{x,y}$. Let $\gamma=\corchetes{\gamma_{ij}}\in\Gamma$, observe that
		$$\gamma g_{x,y}\gamma^{-1}=\corchetes{\begin{array}{ccc}
	1 & \frac{\gamma_{11}}{\gamma_{22}}x & \frac{\gamma_{11}}{\gamma_{22}\gamma_{33}}(\gamma_{22} y - \gamma_{23} x) \\
	0 & 1 & 0 \\
	0 & 0 & 1\\
	\end{array}}\in \core(\Gamma).$$ 
	This element determines the line 
	\begin{align*}
	\linproy{e_1,\corchetes{0:-\frac{\gamma_{11}}{\gamma_{22}\gamma_{33}}(\gamma_{22} y - \gamma_{23} x):\frac{\gamma_{11}}{\gamma_{22}}x}}&=\linproy{e_1,\corchetes{0:\gamma_{23}x-\gamma_{22}y:\gamma_{33}x}}\\
		&=\ell_{\gamma_{33}x,\gamma_{22}y-\gamma_{23}x}.
	\end{align*}
	Therefore, this is a line in $\mathcal{C}(\Gamma)$.
	On the other hand, a direct calculation shows that
		$$\gamma(\ell_{x,y})=\ell_{\gamma_{33}x,\gamma_{22}y-\gamma_{23}x}.$$
	The last two equations show that $\gamma$ leaves $\mathcal{C}(\Gamma)$ invariant. 
\end{proof}

\begin{figure}[H]
\begin{center}	
\includegraphics[height=36mm]{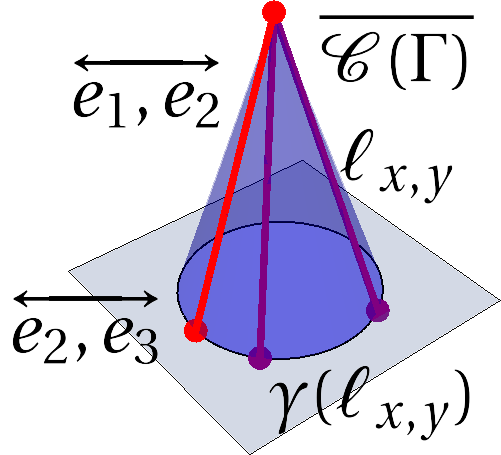}
\caption{The elements of $\Gamma$ moving the line $\ell_{x,y}$.}
\label{fig_dinamica_lineas_cono}
\end{center}	
\end{figure}

This last proposition says that the elements of $\core(\Gamma)$ determine a pencil of lines passing through $e_1$, whose closure is a cone in $\CP^2$), and this cone determines, in certain way, the dynamics of the group $\Gamma$. Every element $\gamma\in\Gamma$ moves the line $\ell_{x,y}$ to the line $\gamma(\ell_{x,y})$ according to the proof of the last proposition, the line $\linproy{e_1,e_2}$ is fixed by every element of $\Gamma$ (see figure \ref{fig_dinamica_lineas_cono}). 

In Subsetion \ref{subsection:representation} we will study in greater detail how a loxodromic element interacts with $\mathcal{C}(\Gamma)$ and leaves it invariant. As we will see, this invariance imposes strong restrictions on these loxodromic elements.\\

We say that the discrete group $\Gamma\subset U_+$ is \emph{conic} if $\overline{\mathcal{C}(\Gamma)}$ is a cone homeomorphic to the complement of $\C\times\parentesis{\Hh^+\cup\Hh^-}$. If $\mathcal{C}(\Gamma)$ is a line, we say that $\Gamma$ is \emph{non-conic}. In $\cite{ppar}$, conic groups are called \emph{irreducible} and non-conic groups, \emph{reducible}. This nomenclature is taken from \cite{benoistdiscrete}, however it was decided to change it in order to avoid confusion with reducible/irreducible groups, which are groups that leave a proper subspace invariant (or leave no proper subspace invariant respectively) (see Definition \ref{defn_irreducible}).\\

\begin{ejem}
If $\Gamma\subset U_+$ is such that $\core(\Gamma)=\prodint{g_{0,1},g_{1,0}}$ then the set $\mathcal{C}(\Gamma)\cap\linproy{e_2,e_3}$, viewed as a subset of $\hC\cong\linproy{e_2,e_3}$, is the set of rational numbers $\Q$.
\end{ejem}

The next proposition follows immediately from Propositions \ref{prop_HC_conmutativo} and \ref{prop_IS_discreto_invariante}. 

\begin{prop}\label{prop_kernel_finito_kernel_trivial}
Let $\Gamma\subset U_+$ be a discrete torsion free group. If $\kernel(\Gamma)$ is finite, then
	$$\kernel(\Gamma)=\set{\id}.$$
\end{prop}

\begin{prop}\label{prop_kernel_igual_core}
Let $\Gamma\subset U_+$ be a non-commutative discrete group such that one of the following hypothesis hold:
	\begin{itemize}
	\item Its control group $\Sigma=\Pi(\Gamma)$ is discrete.
	\item $\Sigma$ is not discrete and $\grL(\Sigma)=\Ss^1$.
	\item $\Sigma$ is not discrete and $\valorabs{\grL(\Sigma)}=2$.
	\end{itemize}
Then $\kernel(\Gamma)=\core(\Gamma)$.
\end{prop}

\begin{proof}
Clearly $\core(\Gamma)\subset\kernel(\Gamma)$, we only have to prove that $\kernel(\Gamma)\subset\core(\Gamma)$. Let $\gamma\in\kernel(\Gamma)$, then
	$$\gamma=\corchetes{\begin{array}{ccc}
	\alpha^{-2} & \gamma_{12} & \gamma_{13} \\
	0 & \alpha & 0 \\
	0 & 0 & \alpha\\
	\end{array}},$$
for some $\alpha\in\C^\ast$.
	\begin{itemize}
	\item If $\valorabs{\alpha}\neq 1$ then $\gamma$ is a complex homothety. 
		\begin{itemize}
		\item If $\Sigma$ is discrete, then by Proposition \ref{prop_HC_conmutativo}, $\Gamma$ would be commutative.
		\item If $\Sigma$ is not discrete and $\grL(\Sigma)=\Ss^1$, we get a contradiction using Proposition \ref{prop_HC_conmutativo_control_no_discreto}. 		
		\item $\Sigma$ is not discrete and $\valorabs{\grL(\Sigma)}=2$, $\Gamma$ would be commutative as a consequence of Proposition \ref{prop_HC_conmutativo_control_no_discreto_grenlim_2ptos} and we get a contradiction.
		\end{itemize}		 
	\item If $\valorabs{\alpha}=1$ but $\alpha\neq 1$ we have two possibilities:
		\begin{itemize}
		\item If $\alpha$ is a rational rotation then $\lambda_{12}(\gamma)$ is a torsion element in the torsion free group $\lambda_{12}(\Gamma)$.
		\item If $\alpha$ is an irrational rotation then $\gamma$ is a irrational ellipto-parabolic element, and then $\Gamma$ would be commutative, by using Proposition \ref{prop_EPI_conmutativo}.
		\end{itemize}
	\end{itemize}
Then $\alpha=1$ and then $\gamma\in\core(\Gamma)$.
\end{proof}

The proof of the following corollary is similar to the proof of the previous proposition.

\begin{cor}\label{cor_diagonal_es_1}
Under the hypothesis of the previous proposition, we can conclude that if 
	$$\gamma=\corchetes{\begin{array}{ccc}
	\alpha^{-2} & \gamma_{12} & \gamma_{13} \\
	0 & \alpha & \gamma_{23} \\
	0 & 0 & \alpha\\
	\end{array}}\in\Gamma,$$
then $\alpha=1$. 
\end{cor}

\begin{prop}\label{prop_control_nodiscreto_conjlim_2puntos_irreducible}
Let $\Gamma\subset U_+$ be a non-commutative, discrete group such that $\valorabs{\grL(\Pi(\Gamma))}=2$. Let $\ell$ be a line passing through $e_1$ such that $\ell\neq \linproy{e_1,e_2}$, $\ell\neq \linproy{e_1,e_3}$ and $\ell\subset\mathcal{C}(\Gamma)$, then $\Gamma$ is conic.  
\end{prop}

\begin{proof} 
If $\valorabs{\grL(\Pi(\Gamma))}=2$ then, up to conjugation, every element of $\Sigma$ has the form
	$$\gamma=\corchetes{\begin{array}{cc}
	\beta & 0 \\
	0 & \delta \\
	\end{array}},$$
with $\valorabs{\beta}\neq\valorabs{\delta}$. Then, every element of $\Gamma$ has the form
	\begin{equation}\label{eq_dem_prop_control_nodiscreto_conjlim_2puntos_irreducible_1}
	\gamma=\corchetes{\begin{array}{ccc}
	\gamma_{22}^{-1}\gamma_{33}^{-1} & \gamma_{12} & \gamma_{13} \\
	0 & \gamma_{22} & 0 \\
	0 & 0 & \gamma_{33}\\
	\end{array}}\text{, with }\valorabs{\gamma_{22}}\neq\valorabs{\gamma_{33}}.
	\end{equation}
Let us assume that $\Gamma$ is not conic, with $\mathcal{C}(\Gamma)=\ell$. Let $[0:-y:x]\in\linproy{e_2,e_3}$ such that $\ell=\ell_{x,y}$, then by hypothesis, $x\neq 0$ and $y\neq 0$. \\

Since $\valorabs{\kernel(\Gamma)}=\infty$, using Proposition \ref{prop_kernel_igual_core} it follows that $\kernel(\Gamma)=\core(\Gamma)$. Since $\kernel(\Gamma)$ is a normal subgroup of $\Gamma$, so is $\core(\Gamma)$. Therefore, if $\gamma=\corchetes{\gamma_{ij}}\in\Gamma$,
	$$\gamma g_{x,y} \gamma^{-1} =\corchetes{\begin{array}{ccc}
	1 & \frac{x}{\gamma_{22}^2\gamma_{33}} & \frac{y}{\gamma_{22}\gamma_{33}^2} \\
	0 & 1 & 0 \\
	0 & 0 & 1\\
	\end{array}}
	\in\core(\Gamma).$$ 
Furthermore, this element $\gamma g_{x,y} \gamma^{-1}$ determines the same line $\ell_{x,y}$ in $\mathcal{C}(\Gamma)$. Then 
	$$\corchetes{-\frac{y}{\gamma_{22}\gamma_{33}^2}:\frac{x}{\gamma_{22}^2\gamma_{33}}}=\corchetes{-y:x},$$
which means that $\gamma_{22}=\gamma_{33}$, contradicting (\ref{eq_dem_prop_control_nodiscreto_conjlim_2puntos_irreducible_1}). This contradiction proves the proposition.	
\end{proof}

We can re-state Lemma 3.4 of \cite{ppar} in the following way.

\begin{prop}\label{prop_proyectiv_adit}
Let $W\subset\C^2$ be an additive group, then $\corchetes{W\setminus\set{0}}$ is either a point, a circle or $\CP^1$. 
\end{prop}

\subsection{Decomposition of non-commutative discrete groups of $U_+$}

In this subsection we state and prove the main theorem of this section. For the sake of clarity, we will divide the theorem in two parts (Theorems \ref{thm_descomposicion_caso_noconmutativo} and \ref{thm_descomposicion_caso_noconmutativo2}). In the first part we give a decomposition of the group in four layers, the first two containing only parabolic elements and the last two containing only loxodromic elements. In the second part we bound the rank of the group. Finally, we give some consequences of these theorems.

\subsubsection*{First theorem} 

\begin{thm}\label{thm_descomposicion_caso_noconmutativo}
Let $\Gamma\subset U_{+}$ be a non-commutative, torsion free, complex Kleinian group, then $\Gamma$ can be written in the following way
	$$\Gamma = \core(\Gamma)\rtimes\prodint{\xi_1}\rtimes...\rtimes\prodint{\xi_r}\rtimes\prodint{\eta_1}\rtimes...\rtimes\prodint{\eta_m}\rtimes\prodint{\gamma_1}\rtimes...\rtimes\prodint{\gamma_n}$$
	where
	$$\begin{array}{rl}
	\lambda_{23}(\Gamma)=\prodint{\lambda_{23}(\gamma_1),...,\lambda_{23}(\gamma_n)}, & n=\text{rank}\parentesis{\lambda_{23}(\Gamma)}.\\
	\lambda_{12}\parentesis{\kernel(\lambda_{23})}=\prodint{\lambda_{12}(\eta_1),...,\lambda_{12}(\eta_m)}, & m=\text{rank}\parentesis{\lambda_{12}(\kernel(\lambda_{23})}.\\
	\Pi\parentesis{\kernel(\lambda_{12})\cap\kernel(\lambda_{23})}=\prodint{\Pi(\xi_1),...,\Pi(\xi_r)}, & r=\text{rank}\parentesis{\Pi\parentesis{\kernel(\lambda_{12})\cap\kernel(\lambda_{23})}}.
	\end{array}$$
	
\end{thm}
  
\begin{proof}
\underline{Part I. Decomposition of $\Gamma$ in terms of $\kernel(\lambda_{23})$.} Let $\Gamma\subset U_{+}$ be a torsion free complex Kleinian group. Since $\Gamma$ is triangular, it is solvable and therefore it is finitely generated (see Theorem \ref{thm_solubles_fg}), then $\lambda_{23}(\Gamma)$ is finitely generated. Let $n=\text{rank}\parentesis{\lambda_{23}(\Gamma)}$ and let $\set{\tilde{\gamma}_1,...,\tilde{\gamma}_n}\subset\C^\ast$ be a generating set for $\lambda_{23}(\Gamma)$. We choose elements $\gamma_1,...,\gamma_n\in\Gamma$ such that $\gamma_i\in\lambda_{23}^{-1}(\tilde{\gamma}_i)$. Observe that 
	\begin{equation}\label{eq_thm_descomposicion_caso_noconmutativo_1}
	\Gamma=\prodint{\text{Ker}(\lambda_{23}),\gamma_1,...,\gamma_n}.
	\end{equation}
Furthermore, we will prove that 
	\begin{equation}\label{eq_thm_descomposicion_caso_noconmutativo_2}
	\Gamma=\parentesis{\parentesis{\text{Ker}(\lambda_{23})\rtimes\prodint{\gamma_1}}\rtimes...}\rtimes\prodint{\gamma_n}.
	\end{equation}
Since $\lambda_{23}$ is a group homomorphism, $\text{Ker}(\lambda_{23})$ is a normal subgroup of $\Gamma$ and therefore, it is a normal subgroup of $\prodint{\text{Ker}(\lambda_{23}),\lambda_1}$.\\ 
Now, assume that $\text{Ker}(\lambda_{23})\cap\lambda_1$ is not trivial, then there exist $p\in\Z$ such that $\gamma_1^p\in\text{Ker}(\lambda_{23})$. If $\gamma_1=\corchetes{a_{ij}}$ then the previous assumption means that either $a_{22}=a_{33}$ or $a_{22}^p=a_{23}^p$. In the latter case, this means, without loss of generality that 
	\begin{equation}\label{eq_thm_descomposicion_caso_noconmutativo_3}
	a_{22}a^{-1}_{33}=\omega,
	\end{equation}	 
where $\omega$ is a $p$-th root of the unity, with $p>1$. On the other hand, since $\Gamma$ is torsion free, the group $\lambda_{23}(\Gamma)\subset\C^{\ast}$ is torsion free. However, it follows from (\ref{eq_thm_descomposicion_caso_noconmutativo_3}) that $a_{22}a^{-1}_{33}$ is a torsion element of $\lambda_{23}(\Gamma)$, contradicting that $\lambda_{23}(\Gamma)$ is torsion free; thus $a_{22}=a_{33}$. But if this is the case then $\gamma_1\in\text{Ker}(\lambda_{23})$ contradicting that $\lambda_{23}(\gamma_1)=\tilde{\gamma}_1$ belongs to a generating set for $\lambda_{23}(\Gamma)$.\\

Thus, we can form the semi-direct product $\text{Ker}(\lambda_{23})\rtimes\prodint{\gamma_1}$. Now we verify that we can form the semi-direct product $\parentesis{\parentesis{\text{Ker}(\lambda_{23})\rtimes\prodint{\gamma_1}}}\rtimes\prodint{\gamma_2}$.\\

First we verify that $\text{Ker}(\lambda_{23})\rtimes\prodint{\gamma_1}$ is a normal subgroup of $\prodint{\text{Ker}(\lambda_{23})\rtimes\prodint{\gamma_1},\gamma_2}$. Let $g\in\text{Ker}(\lambda_{23})$ and consider arbitrary elements $g\gamma_1^n\in \text{Ker}(\lambda_{23})\rtimes\prodint{\gamma_1}$ and $\gamma_2^m\in\prodint{\gamma_2}$, then
	\begin{align*}
	\lambda_2^m \parentesis{g \lambda_1^n} \lambda_2^{-m} &= \underbrace{\parentesis{\lambda_2^m g \lambda_2^{-m}}}_{\text{Ker}(\lambda_{23})\triangleleft \Gamma}\lambda_2^m \lambda_1^n\lambda_2^{-m}\\
		&= g_1\cdot \parentesis{\lambda_2^m \lambda_1^n\lambda_2^{-m}}
	\end{align*}
where $g_1=\lambda_2^m g \lambda_2^{-m}\in\text{Ker}(\lambda_{23})$. Since $\lambda_{23}\parentesis{g_1\lambda_2^m \lambda_1^n\lambda_2^{-m}}=\lambda_{23}\parentesis{\lambda_1^n}$, we know that there exists $g_2\in\text{Ker}(\lambda_{23})$ such that $g_1\lambda_2^m \lambda_1^n\lambda_2^{-m}=g_2 \lambda_1^n$. Thus 
	$$\lambda_2^m \parentesis{g \lambda_1^n} \lambda_2^{-m}=g_2 \lambda_1^n\in \text{Ker}(\lambda_{23})\rtimes\prodint{\gamma_1}.$$
Therefore $\text{Ker}(\lambda_{23})\rtimes\prodint{\gamma_1}$ is a normal subgroup of $\prodint{\text{Ker}(\lambda_{23})\rtimes\prodint{\gamma_1},\gamma_2}$.\\

Now we verify that $\parentesis{\text{Ker}(\lambda_{23})\rtimes\prodint{\gamma_1}}\cap \prodint{\gamma_2}=\set{\id}$. Assume that there exist $p,q\in\Z$ such that $\gamma_2^p = g\lambda_1^q$ for some $g\in\kernel(\lambda_{23})$ with $\gamma_2^p\neq\id$, then 
	\begin{equation}\label{eq_thm_descomposicion_caso_noconmutativo_4}
	\gamma_2^p\gamma_1^{-q}\in \text{Ker}(\lambda_{23})
	\end{equation}		 	
If we write $\lambda_1=\corchetes{\alpha_{ij}}$ and $\lambda_2=\corchetes{\beta_{ij}}$, then (\ref{eq_thm_descomposicion_caso_noconmutativo_4}) means that 
	$$\parentesis{\beta_{22}\alpha_{22}^{-\frac{q}{p}}}^p=\parentesis{\beta_{33}\alpha_{33}^{-\frac{q}{p}}}^p,$$
which means that 
	\begin{equation}\label{eq_thm_descomposicion_caso_noconmutativo_5}
	\beta_{22}\alpha_{22}^{-\frac{q}{p}}\parentesis{\beta_{33}\alpha_{33}^{-\frac{q}{p}}}^{-1}=\omega,
	\end{equation}
where $\omega$ is a $p$-th root of the unity. As before, $\lambda_{23}(\Gamma)\subset\C^\ast$ is a torsion free group, but (\ref{eq_thm_descomposicion_caso_noconmutativo_5}) gives a torsion element in $\lambda_{23}(\Gamma)$. This contradiction proves that  
	$$\parentesis{\text{Ker}(\lambda_{23})\rtimes\prodint{\gamma_1}}\cap \prodint{\gamma_2}=\set{\id}.$$
Therefore we can define the semi-direct product
	$$\parentesis{\text{Ker}(\lambda_{23})\rtimes\prodint{\gamma_1}}\rtimes\prodint{\gamma_2}.$$
In the same way we can form the semi-direct product
	$$\parentesis{\parentesis{\text{Ker}(\lambda_{23})\rtimes\prodint{\gamma_1}}\rtimes...}\rtimes\prodint{\gamma_n}.$$
For clarity we will just write $\text{Ker}(\lambda_{23})\rtimes\prodint{\gamma_1}\rtimes...\rtimes\prodint{\gamma_n}$ instead. Using (\ref{eq_thm_descomposicion_caso_noconmutativo_1}) we have proven (\ref{eq_thm_descomposicion_caso_noconmutativo_2}). \\

\underline{Part II. Decompose $\kernel(\lambda_{23})$ in terms of $\kernel(\lambda_{12})$.} Now consider the restriction
	\begin{equation}\label{eq_thm_descomposicion_caso_noconmutativo_6}
	\lambda_{12}:\kernel(\lambda_{23})\rightarrow \C^\ast.
	\end{equation}
Then, again, $\lambda_{12}\parentesis{\kernel(\lambda_{23})}$ is finitely generated and let $m$ be its rank. Let 
	$$\set{\tilde{\eta}_1,...,\tilde{\eta}_m}$$ 
be a generating set, we choose elements $\eta_i\in\lambda_{12}^{-1}\parentesis{\tilde{\eta}_i}$. Denote $$A=\kernel(\lambda_{12})\cap \kernel(\lambda_{23}),$$ observe that 
	\begin{equation}\label{eq_thm_descomposicion_caso_noconmutativo_7}
	\kernel(\lambda_{23})=\prodint{A,\eta_1,...,\eta_m}.	
	\end{equation}	 
Observe that every element of $A$ has the form
	$$\corchetes{\begin{array}{ccc}
	1 & x & y \\
	0 & 1 & z\\
	0 & 0 & 1\\ 
	\end{array}
	}.$$
Since $A$ is the kernel of the morphism defined in (\ref{eq_thm_descomposicion_caso_noconmutativo_6}), $A$ is a normal subgroup of $\prodint{A,\eta_1}$.\\ 

Now suppose that $A\cap\prodint{\eta_1}$ is not trivial, let $\eta_1^p\in A$ with $\eta_1^p\neq\id$. Denote $\eta_1=\corchetes{a_{ij}}$, since $\eta_1^p\in A$ it must hold $a_{11}^p=a_{22}^p=a_{33}^p$, which means that $a_{11}a_{22}^{-1}=\omega$, where either $\omega\neq 1$ is a $p$-th root of unity or $\omega=1$. In the former case, $a_{11}a_{22}^{-1}$ is a torsion element in $\lambda_{12}\parentesis{\kernel(\lambda_{23})}$, which is a torsion free group, so this case cannot happen. If $\omega=1$ then $a_{11}=a_{22}$ and since $\eta_1\in\kernel(\lambda_{23})$, $a_{22}=a_{33}$ it follows that $a_{11}=a_{22}=a_{33}$ and thus, $\eta_1\in A$ which contradicts that $\eta_1$ is part of the generating set. All this proves that $A\cap\prodint{\eta_1}= \emptyset$.\\

This guarantees that we can form the semi-direct product $A\rtimes\prodint{\eta_1}$.\\

Now we verify that we can make the semi-direct product with $\prodint{\eta_2}$. The same argument used in part I to prove normality when we added $\gamma_2$ to $\kernel(\lambda_{23})\rtimes\prodint{\gamma_1}$ can be applied in the same way now to prove that $A\rtimes\prodint{\eta_1}$ is a normal subgroup of $\prodint{A\rtimes\prodint{\eta_1},\eta_2}$ and that $\parentesis{A\rtimes\prodint{\eta_1}}\cap\prodint{\eta_2}$ is trivial. Using this argument for $\eta_3,...,\eta_m$ we get,
	\begin{equation}\label{eq_thm_descomposicion_caso_noconmutativo_8}
	\kernel(\lambda_{23})=A\rtimes\prodint{\eta_1}\rtimes...\rtimes\prodint{\eta_m}.
	\end{equation}
	
\underline{Part III. Decompose $A$ in terms of $\kernel(\Gamma)$.} Consider the restriction of the previously defined group homomorphism
	$$\Pi:A\rightarrow\psl$$
As before, $\core(\Gamma)$ it is a normal subgroup of $\Gamma$. Since $A$ is solvable, it is finitely generated and so $\Pi(A)\subset\psl$ is finitely generated, denote by $r$ its rank.  Let $\set{\tilde{\xi}_1,...,\tilde{\xi}_r}\subset\Pi(A)$ be a generating set for $\Pi(A)$, choose $\xi_i\in\tilde{\xi}_i^{-1}\subset A$, then
	\begin{equation}\label{eq_thm_descomposicion_caso_noconmutativo_10}
	A=\prodint{\core(\Gamma),\xi_1,...,\xi_r}
	\end{equation}		
$\core(\Gamma)$ is a normal subgroup of $\prodint{\core(\Gamma),\xi_1}$. Assume that $\core(\Gamma)\cap\xi_1$ is not trivial, in other words, assume that $\xi_1^p\in\core(\Gamma)$ for some $p\in\Z\setminus{0}$. Since $\xi_1\in A$, we write
	$$\xi_1=\corchetes{\begin{array}{ccc}
	1 & x & y \\
	0 & 1 & z\\
	0 & 0 & 1\\ 
	\end{array}
	}.$$
Then 
	$$\xi^p_1=\corchetes{\begin{array}{ccc}
	1 & x' & y' \\
	0 & 1 & pz\\
	0 & 0 & 1\\ 
	\end{array}
	}\text{, for some }x',y'\in\C.$$
Observe that $\xi_1^p\in\core(\Gamma)$ if and only if $z=0$, which contradicts that $\tilde{\xi}_1$ is a generator of $\Pi(A)$. Then we can form the semi-direct product
	$$\core(\Gamma)\rtimes\prodint{\xi_1}.$$
Observe that $\Pi(A)$ is a commutative subgroup of $\psl$, then we can apply the same argument we used when we formed the semi-direct product with $\lambda_2$ to conclude that $\core(\Gamma)\rtimes\prodint{\xi_1}$ is a normal subgroup of $\prodint{\core(\Gamma)\rtimes\prodint{\xi_1},\xi_2}$. This concludes the proof of the first part of the theorem.
\end{proof}

\subsubsection*{Second theorem}

Before proving Theorem \ref{thm_descomposicion_caso_noconmutativo2}, we state some results needed through the proof. The following result can be viewed in \cite{kapovich02} (see Theorem 1).

\begin{thm}\label{thm_obdim_1}
Let $\Gamma$ be a group acting properly and discontinuously on a contractible manifold of dimension $m$, then $\text{obdim}(\Gamma)\leq m$.
\end{thm}

In the statement of the previous Theorem, $\text{obdim}(\Gamma)$ is called the \emph{obstructor dimension}. Instead of giving its definition, we will determine its value using the hypothesis of our Theorem \ref{thm_descomposicion_caso_noconmutativo} and the following properties (see Corollary 27 of \cite{kapovich02} and \cite{bestvina02} respectively):
	\begin{itemize}
	\item If $\Gamma=H\rtimes Q$ with $H$ and $Q$ finitely generated and $H$ weakly convex, then
		$$\text{obdim}(\Gamma)\geq \text{obdim}(H)+\text{obdim}(Q).$$
	\item In the particular case when $\Gamma=\Z^n$, it holds
		$$\text{obdim}(\Z^n) = n.$$		
	\end{itemize}
	
Under the hypothesis of Theorem \ref{thm_descomposicion_caso_noconmutativo} and using the alternative decomposition described in Corollary \ref{cor_descomposicion_bloques_noconm}, the previous properties imply the following reformulation of Theorem \ref{thm_obdim_1}.

\begin{thm}\label{thm_obdim_2}
Let $\Gamma\subset U_{+}$ be a non-commutative, torsion free, complex Kleinian group acting properly and discontinuously on a simply connected domain $\Omega\subset\CP^2$, then 
	$$k+r+m+n\leq 4.$$
\end{thm}

The strategy to prove Theorem \ref{thm_descomposicion_caso_noconmutativo2} will be to find a simply connected domain $\Omega\subset\CP^2$ where $\Gamma$ acts properly and discontinuously and then apply Theorem \ref{thm_obdim_2}. In some cases, we will obtain the explicit decomposition of $\Gamma$ and verify that $\text{rank}(\Gamma)\leq 4$.\\

\begin{thm}\label{thm_descomposicion_caso_noconmutativo2}
Let $\Gamma\subset U_{+}$ be a non-commutative, torsion free, complex Kleinian group, then $\rank(\Gamma)\leq 4$. Using the notation of Theorem \ref{thm_descomposicion_caso_noconmutativo} it holds 
		$$k+r+m+n\leq 4.$$	
where $k=\text{rank}\parentesis{\core(\Gamma)}$.
\end{thm}

\begin{proof}
We denote the control group by $\Sigma=\Pi(\Gamma)$. If $k,r,m,n$ are defined as in the proof of Theorem \ref{thm_descomposicion_caso_noconmutativo}, then
	\begin{equation}\label{eq_thm_descomposicion_rango_caso_noconmutativo_1}
	\text{rank}(\Gamma)\leq k+r+m+n.
	\end{equation} 
We will divide the proof in the following cases:
	\begin{enumerate}[(i)]
	\item $\Sigma$ is discrete and $\kernel(\Gamma)$ is finite.
	\item $\Sigma$ is discrete and $\kernel(\Gamma)$ is infinite.
	\item $\Sigma$ is not discrete.	
	\end{enumerate}

Before dealing with each case, consider the diagram in Figure \ref{fig_Diagrama_Casos_Thm_2}. This diagram summarizes each subcase we will be considering. Each subcase is listed with the same name it appears on the proof.

\begin{figure}
	\begin{center}	
	\includegraphics[height=190mm]{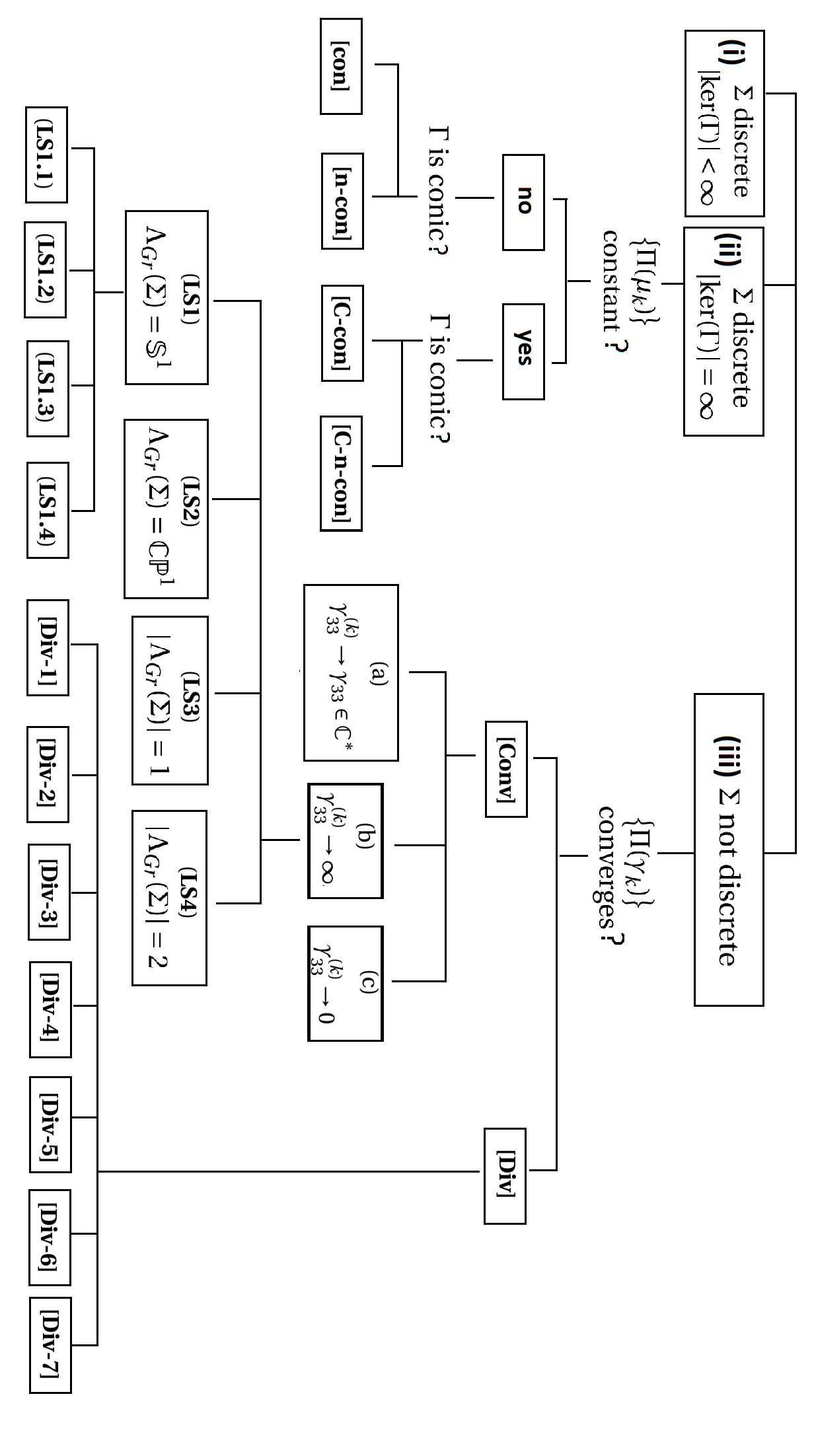}	
	\caption{Subcases for the proof of Theorem \ref{thm_descomposicion_caso_noconmutativo2}.}
	\label{fig_Diagrama_Casos_Thm_2}		
	\end{center}	
	\end{figure}  

Now we deal with each case.
	\begin{enumerate}[(i)]
	\item Assume that \textbf{$\Sigma$ is discrete and $\kernel(\Gamma)$ is finite}. Then by Theorem 5.8.2 of \cite{ckg_libro} we know that $\Gamma$ acts properly and discontinuously on 
		$$\Omega = \parentesis{\bigcup_{z\in\Omega(\Sigma)} \linproy{e_1,z}}\setminus \set{e_1},$$
	where $\Omega(\Sigma)=\CP^1\setminus\Lambda(\Sigma)$ denotes the discontinuity set of $\Sigma$. If the cardina-lity of the limit set $\valorabs{\Lambda(\Sigma)}$ is $0$, $1$ or $\infty$ then each connected component of $\Omega$ is simply connected, since they are respectively homeomorphic to $\CP^2$, $\C^2$ or the complement of a cone in $\CP^2$ (see the figure below). Using Theorem \ref{thm_obdim_2} it follows
		$$k+r+m+n\leq 4.$$
	
	\begin{figure}[H]
	\begin{center}	
	\includegraphics[height=35mm]{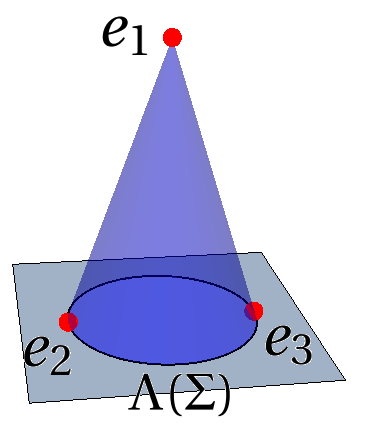}	
	\end{center}	
	\end{figure}	
	
	 If $\valorabs{\Lambda(\Sigma)}=2$, then $\Sigma$ is elemental and therefore it is a cyclic group generated by a loxodromic element (see Theorem 1.6 of \cite{series}). Hence, $\Sigma\cong\Z$. On the other hand, since $\Gamma$ is torsion free and $\kernel(\Gamma)$ is finite, $\kernel(\Gamma)=\set{\id}$ (Proposition \ref{prop_kernel_finito_kernel_trivial}) and therefore, $\Pi:\Gamma\rightarrow\Sigma$ is a group isomorphism and then, $\Gamma\cong \Sigma$, it follows that 
	 	$$\Gamma\cong\Z$$
	 and therefore, $\text{rank}(\Gamma)=1$. Then it holds trivially that $k+r+m+n\leq 4$.\\
	
	\item Now, let us assume that \textbf{$\kernel(\Gamma)$ is infinite and $\Sigma$ is discrete}. Since $\Sigma$ is discrete and torsion free, it follows from Proposition \ref{prop_kernel_igual_core} that $\core(\Gamma)=\kernel(\Gamma)$. Then $\core(\Gamma)$ is infinite, which means that there exist elements $g_{x,y}\in\core(\Gamma)$, with $g_{x,y}\neq\id$. Define $\mathcal{B}(\Gamma)=\pi\parentesis{\mathcal{C}(\Gamma)\setminus\set{e_1}}$, then $\overline{\mathcal{B}(\Gamma)}\cong\Ss^1$ or it is a single point (Proposition \ref{prop_proyectiv_adit}).\\
	
	On the other hand, consider a sequence of distinct elements $\set{\mu_k}\subset \Gamma$. Since $\Sigma$ is discrete, the sequence $\set{\Pi(\mu_k)}$ is either constant or converges by distinct elements to a quasi-projective map $\sigma\in\qp$. Let us assume first that $\set{\Pi(\mu_k)}$ diverges, then $\sigma$ has one of the following forms
	$$\sigma_1=\corchetes{\begin{array}{cc}
			1 & 0 \\
			0 & 0\\ 
			\end{array}},\;\;\;\sigma_2=\corchetes{\begin{array}{cc}
			0 & 1 \\
			0 & 0\\ 
			\end{array}},\;\;\;\sigma_3=\corchetes{\begin{array}{cc}
			0 & 0 \\
			0 & 1\\ 
			\end{array}},\;\;\;\sigma_4=\corchetes{\begin{array}{cc}
			1 & a\\
			0 & 0\\ 
			\end{array}},$$
	with $a\in\C^\ast$. Then
		\begin{align*}
		\kernel(\sigma_2)=\kernel(\sigma_3)&=\text{Im}(\sigma_1)=\text{Im}(\sigma_2)=\set{e_1}\\
		\kernel(\sigma_1)&=\text{Im}(\sigma_3)=\set{e_2}.
		\end{align*}		
	Since $\Sigma$ is triangular, it follows that $\Sigma$ is a solvable discrete subgroup of $\psl$ and therefore, it is elemental. Using Theorem 1.6 of \cite{series}, we have three possibilities:
		\begin{itemize}
		\item[\titem{[e1]}] $\Sigma=\prodint{h}$, with $h\in\psl$ a loxodromic element. Then, the quasi-projective limit of sequences in $\Sigma$ can only be $\sigma_1$ or $\sigma_3$.
		\item[\titem{[e2]}] $\Sigma=\prodint{h}$, with $h\in\psl$ a parabolic element. 
		\item[\titem{[e3]}] $\Sigma=\prodint{g,h}$, with $g,h\in\psl$ parabolic elements with different translation directions.	\\
		\end{itemize}		 
				
		Now, considering whether $\Gamma$ is conic or not, we have the two possible cases:
		
		\begin{itemize}
		\item[\titem{[con]}] $\Gamma$ is conic, then $\mathcal{B}(\Gamma)\cong\Ss^1$. Define
			$$\Omega(\Gamma)=\CP^2\setminus \overline{\mathcal{C}(\Gamma)}.$$
		We will verify that $\Gamma$ acts properly and discontinuously on $\Omega(\Gamma)$. Let $K\subset \Omega(\Gamma)$ be a compact set, denote 
			$$\mathfrak{K}=\SET{\gamma\in\Gamma}{\gamma(K)\cap K\neq \emptyset},$$		
		 and assume that $\valorabs{\mathfrak{K}}=\infty$. Then, we can write $\mathfrak{K}=\set{\gamma_1,\gamma_2,...}$. Consider the sequence of distinct elements $\set{\gamma_k}\subset\Gamma$, and the sequence $\set{\Pi(\gamma_k)}\subset\Sigma$. Let $\sigma\in\qp$ be the quasi-projective limit of $\set{\Pi(\gamma_k)}$, then
		 	$$\kernel(\sigma),\text{Im}(\sigma)\subset\mathcal{B}(\Gamma).$$
		By definition, $\pi(K)\cap\mathcal{B}(\Gamma)=\emptyset$ and then, by Proposition \ref{prop_convergencia_qp}, $\Pi(\gamma_k)(\pi(K))$ accumulates on $e_1$ (observe that we are considering this $e_1$ as a point in $\CP^1\cong \linproy{e_2,e_3}$, but we are actually referring to $\pi(e_2)$).
		On the other hand, since $\valorabs{\mathfrak{K}}=\infty$, 
			$$\valorabs{\SET{\alpha\in\Sigma}{\alpha(\pi(K))\cap \pi(K)\neq \emptyset}}=\infty.$$
			
		\begin{figure}[H]
		\begin{center}	
		\includegraphics[height=40mm]{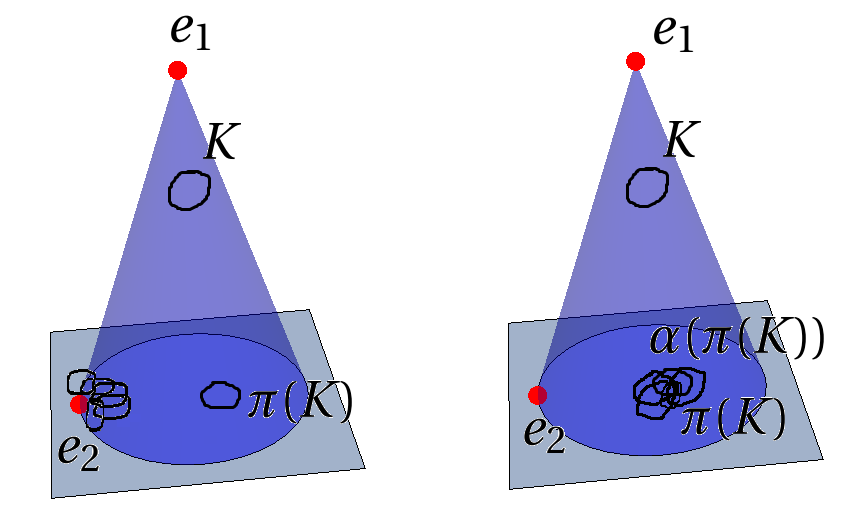}
		\caption{Left: Orbits of $\pi(K)$ accumulate on $e_2$ according to Proposition \ref{prop_convergencia_qp}. Right: Orbits of $\pi(K)$ accumulate away from $e_2$ when we assume $\valorabs{\mathfrak{K}}=\infty$.}
		\label{fig_contradiccion_acum_cono}
		\end{center}	
		\end{figure}	
		This contradicts the fact that $\Pi(\gamma_k)(\pi(K))$ accumulates on $e_1$ (see Figure \ref{fig_contradiccion_acum_cono}). Therefore, $\valorabs{\mathfrak{K}}<\infty$ and then, $\Gamma$ acts properly and discontinuously on each connected component of $\Omega(\Gamma)$.\\
		
		Finally, observe that 
			$$\Omega(\Gamma)\cong \C\times\parentesis{\Hh\cup\Hh^-}.$$
		Then, each connected component of $\Omega(\Gamma)$ is simply connected. U-sing Theorem \ref{thm_obdim_2}, yields
		$$k+r+m+n\leq 4.$$
		
		\item[\titem{[n-con]}] $\Gamma$ is not conic, then $\mathcal{B}(\Gamma)=\set{p}$ for some $p\in\CP^1$. It follows from proposition \ref{prop_control_nodiscreto_conjlim_2puntos_irreducible}, that if $p\neq e_1$ and $p\neq e_2$, then $\Gamma$ is conic. Therefore, we can assume that $p=e_1$ or $p=e_2$. We have the following cases:
			\begin{itemize}
			\item[\titem{[n-con-1]}] If $\Sigma$ has the form \titemT{[e2]} or \titemT{[e3]} then $\Lambda(\Sigma)=\set{e_1}$. Let 
				$$\Omega(\Gamma)=\CP^2\setminus\linproy{e_1,e_2}.$$
			The same argument used in the previous case \titemT{[con]} shows that $\Gamma$ acts properly and discontinuously on $\Omega(\Gamma)\cong \C^2$. Since $\Omega(\Gamma)$ is simply connected, it follows from Theorem \ref{thm_obdim_2} that
		$$k+r+m+n\leq 4.$$ 	
			\item[\titem{[n-con-2]}] If $\Sigma$ has the form \titemT{[e1]} then $\Sigma$ is cyclic and therefore,
			\begin{equation}\label{eq_dem_eq_thm_descomposicion_caso_noconmutativo_indices_INFDISC_1}
			\Sigma\cong \Z.
			\end{equation}
			Since $\pi$ is a group morphism, then
				\begin{equation}\label{eq_dem_eq_thm_descomposicion_caso_noconmutativo_indices_INFDISC_0}
				\Gamma\cong \kernel(\Gamma)\rtimes \text{Im}(\Gamma)=\core(\Gamma)\rtimes \Sigma.
				\end{equation}			
			Now, since $\core(\Gamma)$ is not conic,
				\begin{equation}\label{eq_dem_eq_thm_descomposicion_caso_noconmutativo_indices_INFDISC_2}
				\core(\Gamma)\cong \Z.
				\end{equation}
			Using (\ref{eq_dem_eq_thm_descomposicion_caso_noconmutativo_indices_INFDISC_1}) and (\ref{eq_dem_eq_thm_descomposicion_caso_noconmutativo_indices_INFDISC_2}) in (\ref{eq_dem_eq_thm_descomposicion_caso_noconmutativo_indices_INFDISC_0}) yields $\Gamma\cong\Z\rtimes\Z$. Then $\text{rank}(\Gamma)=2\leq 4$.\\
			\end{itemize}						
		\end{itemize}		 
		Now, let us assume that $\set{\Pi(\mu_k)}$ is a constant sequence. Then, there exists a sequence $\set{g_k}\subset\ker(\Gamma)$ such that
		$$\gamma_k=\gamma_0 g_k.$$		
		It follows that $\set{g_k}\subset\core(\Gamma)$ (see Proposition \ref{prop_kernel_igual_core}). Let $\tau\in\QP$ be the quasi-projective limit of $\set{g_k}$, since 
			$$\text{Im}(\tau)\nsubseteq \kernel(\tau),$$
		then $\gamma_0 g_k=\gamma_k\rightarrow \gamma_0\tau$. It is straightforward to verify that $\text{Im}(\tau)=\text{Im}(\gamma_0\tau)$ and $\kernel(\tau)=\kernel(\gamma_0 \tau)$. We now have two possibilities:
		
		\begin{itemize}
		\item[\titem{[C-con]}] $\Gamma$ is conic, then $\core(\Gamma)=\prodint{g_{x_1,y_1},g_{x_2,y_2}}$ for some $x_1,x_2,y_1,y_2\in\C$. Denote
			$$g_k=g_{n_k x_1 + m_k x_2 , n_k y_1 + m_k y_2},$$
		for some sequences $\set{n_k},\set{m_k}\subset\Z$. Then
			$$\tau=\corchetes{\begin{array}{ccc}
		0 & x & y \\
		0 & 0 & 0\\
		0 & 0 & 0\\ 
		\end{array}}\text{, for some }x,y\in\C,$$
		and then $\kernel(\tau)=\ell_{x,y}\in\overline{\mathcal{C}(\Gamma)}$. Let
			$$\Omega=\CP^2\setminus\overline{\mathcal{C}(\Gamma)}\cong \C\times\parentesis{\Hh^+\cup\Hh^-}.$$
		Let $K\subset\Omega$ be a compact set, then $K\subset\CP^2\setminus\ker(\tau)$. By Lemma \ref{prop_convergencia_qp}, $\Gamma K$ accumulates on $\text{Im}(\tau)=\set{e_1}\nin\Omega$. This means that $\Gamma$ acts properly and discontinuously on each component of $\Omega$. Since each component of $\Omega$ is simply connected, Theorem \ref{thm_obdim_2} implies $k+r+m+n\leq 4$.\\
		
		\item[\titem{[C-n-con]}] $\Gamma$ is not conic, then $\core(\Gamma)=\prodint{g_{x_0,y_0}}$ for some $x_0,y_0\in\C$. Now the sequence $\set{g_k}$ has the form
			$$g_k=\corchetes{\begin{array}{ccc}
		1 & n_k x_0 & n_k y_0 \\
		0 & 1 & 0\\
		0 & 0 & 1\\ 
		\end{array}}\text{, for some }\set{n_k}\subset\Z.$$
		Then, the quasi-projective limit has the form
			$$\tau=\corchetes{\begin{array}{ccc}
		0 & x_0 & y_0 \\
		0 & 0 & 0\\
		0 & 0 & 0\\ 
		\end{array}},$$
		then $\kernel(\tau)=\ell_{x_0,y_0}\in\mathcal{C}(\Gamma)$. Now we apply the same argument as in the previous case \titem{[con]} to conclude $k+r+m+n\leq 4$.			 		
		\end{itemize}				
		
		This proves the second statement of Theorem \ref{thm_descomposicion_caso_noconmutativo} in the case of discrete subgroups of $U_+$ with infinite kernel and discrete control group.
	
	\item Assume that \textbf{$\Sigma$ is not discrete}. Let $\set{\gamma_k}\subset\Gamma$ be a sequence of distinct elements, denote $\gamma_k=\corchetes{\gamma_{ij}^{(k)}}$. A direct calculations shows that 
		\begin{equation}\label{eq_thm_descomposicion_caso_noconmutativo_indices_0}
		\gamma_k^{-1}=\corchetes{\begin{array}{ccc}
		{\gamma_{11}^{(k)}}^{-1} & -\frac{\gamma_{12}^{(k)}}{\gamma_{11}^{(k)}\gamma_{22}^{(k)}} & \frac{\gamma_{12}^{(k)}\gamma_{23}^{(k)}-\gamma_{13}^{(k)}\gamma_{22}^{(k)}}{\gamma_{11}^{(k)}\gamma_{22}^{(k)}\gamma_{33}^{(k)}} \\
		0 & {\gamma_{22}^{(k)}}^{-1} & -\frac{\gamma_{23}^{(k)}}{\gamma_{22}^{(k)}\gamma_{33}^{(k)}}\\
		0 & 0 & {\gamma_{33}^{(k)}}^{-1}\\ 
		\end{array}}.
		\end{equation}				
	Consider the sequence $\set{\Pi(\gamma_k)}\subset\Sigma$. We have two cases:
		\begin{itemize}
		\item[\titem{[Conv]}] The sequence $\set{\Pi(\gamma_k)}\subset\Sigma$ converges to some $\alpha\in\psl$ where
			$$\alpha=\corchetes{\begin{array}{cc}
			\gamma_{22} & \gamma_{23} \\
			0 & 1\\ 
			\end{array}},$$
		with $\gamma_{22}\in\C^\ast$ and 
			$$\frac{\gamma_{22}^{(k)}}{\gamma_{33}^{(k)}}\rightarrow \gamma_{22},\;\;\;\frac{\gamma_{23}^{(k)}}{\gamma_{33}^{(k)}}\rightarrow \gamma_{23}.$$
		 Consider the sequence $\set{\gamma_{33}^{(k)}}\subset\C^{\ast}$, we have three cases:
			\begin{enumerate}
			\item $\gamma_{33}^{(k)}\longrightarrow \gamma_{33}\in\C^{\ast}$. Since $\Pi(\gamma_k)$ converges then
				$$\Pi(\gamma_k)\rightarrow\corchetes{\begin{array}{cc}
				\gamma_{22} & \gamma_{23} \\
				0 & 1\\ 
				\end{array}},$$
				for some $\gamma_{22}\in\C^\ast$ and $\gamma_{23}\in\C$ such that
				\begin{equation}\label{eq_thm_descomposicion_noconmutativo_iii_conv_a_1}
				\frac{\gamma_{22}^{(k)}}{\gamma_{33}^{(k)}}\longrightarrow\gamma_{22},\;\;\;\frac{\gamma_{23}^{(k)}}{\gamma_{33}^{(k)}}\longrightarrow\gamma_{23}.
				\end{equation}		
				Since $\gamma_{11}^{(k)}\gamma_{22}^{(k)}\gamma_{33}^{(k)}=1$, then
				\begin{equation}\label{eq_thm_descomposicion_noconmutativo_iii_conv_a_2}
				\frac{\gamma_{11}^{(k)}}{\gamma_{33}^{(k)}}\longrightarrow \frac{1}{\gamma_{22}\gamma_{33}^3}.
				\end{equation}
				If $\frac{\gamma_{12}^{(k)}}{\gamma_{33}^{(k)}}\rightarrow \gamma_{12}\in\C$ and $\frac{\gamma_{13}^{(k)}}{\gamma_{33}^{(k)}}\rightarrow \gamma_{13}\in\C$ then (\ref{eq_thm_descomposicion_noconmutativo_iii_conv_a_1}) and (\ref{eq_thm_descomposicion_noconmutativo_iii_conv_a_2}) yields 
				$$\gamma_k\longrightarrow\corchetes{\begin{array}{ccc}
				\frac{1}{\gamma_{22}\gamma_{33}^3} & \gamma_{12} & \gamma_{13} \\
				0 & \gamma_{22} & \gamma_{23}\\
				0 & 0 & 1\\ 
				\end{array}}\in\PSL,$$
				contradicting that $\Gamma$ is discrete. Then, $\frac{\gamma_{12}^{(k)}}{\gamma_{33}^{(k)}}\rightarrow \infty$ or $\frac{\gamma_{13}^{(k)}}{\gamma_{33}^{(k)}}\rightarrow \infty$. Thus,	
				\begin{equation}\label{eq_thm_descomposicion_noconmutativo_iii_conv_a_3}
				\gamma_k\longrightarrow\tau_{a,b}:=\corchetes{\begin{array}{ccc}
				0 & a & b \\
				0 & 0 & 0\\
				0 & 0 & 0\\ 
				\end{array}},
				\end{equation}	
				for some $a,b\in\C$, $\valorabs{a}+\valorabs{b}\neq 0$.	Observe that 
				\begin{align}
				\kernel(\tau_{a,b}) &= \linproy{e_1,\corchetes{0:-b:a}},\label{eq_thm_descomposicion_noconmutativo_iii_conv_a_4}\\ 
				\text{Im}(\tau_{a,b}) &=\set{e_1}. \nonumber
				\end{align}
				For each quasi-projective limit $\tau_{a,b}\in\QP$ of sequences $\set{\gamma_k}\subset\Gamma$ like the ones studied in this case, we consider the point $\corchetes{0:-b:a}\in\linproy{e_2,e_3}$ and the horocycle determined by this point and $e_2$, then we consider the pencil of lines passing through $e_1$ and each point of the horocycle (see figure \ref{fig_region_horociclos} below). Denote by $\Omega$ to the complement of this pencils of lines in $\CP^2$. Each connected component of $\Omega$ is simply connected.  
				
				\begin{figure}[H]
				\begin{center}
				\includegraphics[height=40mm]{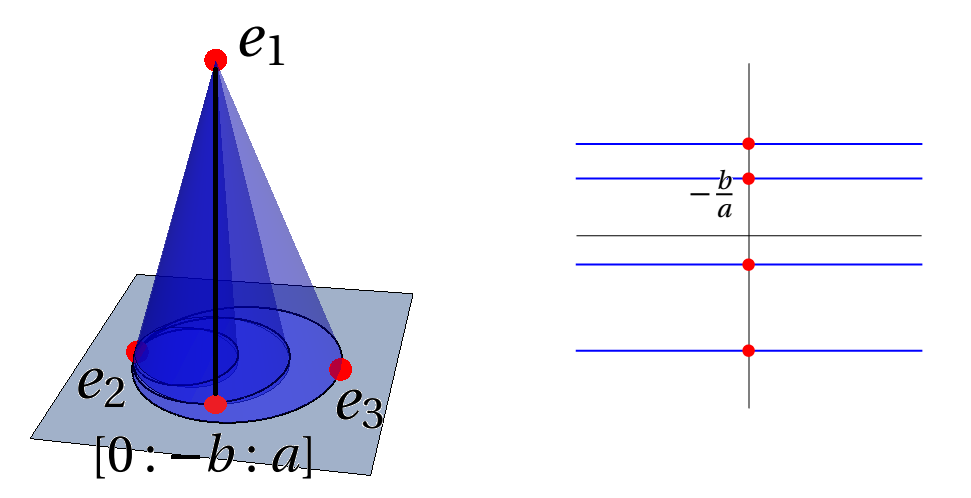}
				\caption{Left: Diferents cones determined by different horocycles. Right: The horocycles when we regard $\protect\linproy{e_2,e_3}$ as $\hC$ with $e_2$ corresponding to $\infty$.}
				\label{fig_region_horociclos}
				\end{center}	
				\end{figure}
				
				As a consequence of Proposition $\ref{prop_convergencia_qp}$ and (\ref{eq_thm_descomposicion_noconmutativo_iii_conv_a_4}), the accumulation points of orbits of compact subsets of $\Omega$ is $\set{e_1}$ and therefore, the action of $\Gamma$ in $\Omega$ is proper and discontinuous.
	 
			\item $\gamma_{33}^{(k)}\rightarrow \infty$. Since $\Pi(\gamma_k)$ converges then
				$$\Pi(\gamma_k)\longrightarrow\corchetes{\begin{array}{cc}
				\gamma_{22} & \gamma_{23} \\
				0 & 1\\ 
				\end{array}},$$
				for some $\gamma_{22},\gamma_{33}\in\C^\ast$ and $\gamma_{23}\in\C$ such that
				\begin{equation}\label{eq_thm_descomposicion_noconmutativo_iii_conv_b_1}
				\frac{\gamma_{22}^{(k)}}{\gamma_{33}^{(k)}}\longrightarrow\gamma_{22},\;\;\;\frac{\gamma_{23}^{(k)}}{\gamma_{33}^{(k)}}\longrightarrow\gamma_{23}.
				\end{equation}								
				Since $\gamma_{22}\in\C^\ast$ and $\gamma_{33}^{(k)}\rightarrow \infty$, then from (\ref{eq_thm_descomposicion_noconmutativo_iii_conv_b_1}) it follows $\gamma_{22}^{(k)}\rightarrow\infty$. Given that 
				$$\gamma_{11}^{(k)}\gamma_{22}^{(k)}\gamma_{33}^{(k)}=1$$
				for all $k\in\N$, then $\gamma_{11}^{(k)}\rightarrow 0$. Then
				$$\frac{\gamma_{11}^{(k)}}{\gamma_{33}^{(k)}}\longrightarrow 0$$
				and then
				\begin{equation}\label{eq_thm_descomposicion_noconmutativo_iii_conv_b_2}
				\gamma_k\longrightarrow\corchetes{\begin{array}{ccc}
				0 & \gamma_{12} & \gamma_{13} \\
				0 & \gamma_{22} & \gamma_{23}\\
				0 & 0 & 1\\ 
				\end{array}},
				\end{equation}
				where 
				$$\frac{\gamma_{12}^{(k)}}{\gamma_{33}^{(k)}}\rightarrow\gamma_{12}\in\C,\;\;\;\frac{\gamma_{13}^{(k)}}{\gamma_{33}^{(k)}}\rightarrow\gamma_{13}\in\C.$$
				If 
				$$\frac{\gamma_{12}^{(k)}}{\gamma_{33}^{(k)}}\rightarrow\infty\;\;\;\text{or}\;\;\;\frac{\gamma_{13}^{(k)}}{\gamma_{33}^{(k)}}\rightarrow\infty,$$
				then the quasi-projective limit $\tau$ has the form (\ref{eq_thm_descomposicion_caso_noconmutativo_indices_3}) with $a=0$ and $b,c\in\C$. We will deal with this case in case (c). Therefore we can assume that $\tau$ has the form (\ref{eq_thm_descomposicion_noconmutativo_iii_conv_b_2}).\\
				
				We can conjugate $\tau$ by an adequate element of $\PSL$ in such way that
				\begin{equation}\label{eq_thm_descomposicion_caso_noconmutativo_indices_001}							
				\tau=\corchetes{\begin{array}{ccc}
				0 & 0 & 0 \\
				0 & \gamma_{22} & \gamma_{23}\\
				0 & 0 & 1\\ 
				\end{array}}.
				\end{equation}
				Consider the Greenberg limit set of $\Sigma$, using proposition \ref{prop_opciones_conj_lim_greenberg} and Lemma \ref{lem_control_nodiscreto_descartar_grL_vacio}, we have the following possible cases:
				\begin{itemize}
				\item[{\scriptsize(\textbf{LS1})}] $\Lambda_{Gr}(\Sigma)=\Ss^1$.
				\item[{\scriptsize(\textbf{LS2})}] $\Lambda_{Gr}(\Sigma)=\CP^1$. 
				\item[{\scriptsize(\textbf{LS3})}] $\valorabs{\Lambda_{Gr}(\Sigma)}=1$.
				\item[{\scriptsize(\textbf{LS4})}] $\valorabs{\Lambda_{Gr}(\Sigma)}=2$.\\				 				
				\end{itemize}								
				
				Let us consider the case {\footnotesize(\textbf{LS1})}. There are four possibilities for the convergence of $\set{\Pi(\gamma_k)}$:
				\begin{itemize}
				\item[{\scriptsize(\textbf{LS1.1})}] The sequence $\set{\Pi(\gamma_k)}$ converges to the identity in $\Sigma$.
				\item[{\scriptsize(\textbf{LS1.2})}] The sequence $\set{\Pi(\gamma_k)}$ converges to a elliptic element in $\Sigma$.
				\item[{\scriptsize(\textbf{LS1.3})}] The sequence $\set{\Pi(\gamma_k)}$ converges to a parabolic element in $\Sigma$.
				\item[{\scriptsize(\textbf{LS1.4})}] The sequence $\set{\Pi(\gamma_k)}$ converges to a loxodromic element in $\Sigma$.\\
				\end{itemize}					
				
				Consider the case {\footnotesize(\textbf{LS1.1})}, then $\gamma_{22}=1$ and $\gamma_{23}=0$, denote the quasi-projective limit of $\set{\gamma_k}$  
				$$\tau_{id}=\corchetes{\begin{array}{ccc}
				0 & 0 & 0 \\
				0 & 1 & 0\\
				0 & 0 & 1\\ 
				\end{array}}.$$
				Since $\Lambda_{Gr}(\Sigma)=\Ss^1$, we can guarantee the existence of parabolic elements in $\Sigma$ (see the proof of proposition \ref{prop_HC_conmutativo_control_no_discreto}). Let $g=\corchetes{g_{ij}}\in\Gamma$ such that $\Pi(g)$ is a parabolic element in $\Sigma$, then $g_{22}=g_{33}=\lambda\in\C^\ast$ and $g_{23}\neq 0$ and therefore $g_{11}=\lambda^{-2}$. Then, denoting $\gamma_k g\gamma_k^{-1}=\corchetes{f_{ij}^{(k)}}$, we have
				\begin{equation}\label{eq_thm_descomposicion_caso_noconmutativo_indices_4}
				\gamma_k g\gamma_k^{-1}=\corchetes{\begin{array}{ccc}
				\lambda^{-2} & \frac{\gamma_{12}^{(k)}(\lambda-\lambda^{-2})+\gamma_{11}^{(k)}g_{12}}{\gamma_{22}^{(k)}} & f_{13}^{(k)} \\
				0 & \lambda & \frac{\gamma_{22}^{(k)}}{\gamma_{33}^{(k)}}g_{23}\\
				0 & 0 & \lambda\\ 
				\end{array}}
				\end{equation}
				where
				\begin{align*}
				f_{13}^{(k)} &= \frac{\parentesis{\gamma_{12}^{(k)}\gamma_{23}^{(k)}-\gamma_{13}^{(k)}\gamma_{22}^{(k)}}\lambda^{-2}-\gamma_{23}^{(k)}\parentesis{\gamma_{12}^{(k)}\lambda+\gamma_{11}^{(k)}g_{12}}}{\gamma_{22}^{(k)}\gamma_{33}^{(k)}}\\
				&+\frac{\gamma_{12}^{(k)}g_{23}+\gamma_{13}^{(k)}\lambda+\gamma_{11}^{(k)}g_{12}}{\gamma_{33}^{(k)}},
				\end{align*}
				and since $\gamma_k\rightarrow\tau_{id}$, we have
				\begin{equation}\label{eq_thm_descomposicion_caso_noconmutativo_indices_5}
				\frac{\gamma_{11}^{(k)}}{\gamma_{33}^{(k)}},\frac{\gamma_{12}^{(k)}}{\gamma_{33}^{(k)}},\frac{\gamma_{13}^{(k)}}{\gamma_{33}^{(k)}},\frac{\gamma_{23}^{(k)}}{\gamma_{33}^{(k)}}\rightarrow 0,\;\;\;\frac{\gamma_{22}^{(k)}}{\gamma_{33}^{(k)}}\rightarrow 1.
				\end{equation}
				From (\ref{eq_thm_descomposicion_caso_noconmutativo_indices_4}) and (\ref{eq_thm_descomposicion_caso_noconmutativo_indices_5}) it follows
				$$\gamma_k g\gamma_k^{-1}\rightarrow \corchetes{\begin{array}{ccc}
				\lambda^{-2} & 0 & 0 \\
				0 & \lambda & g_{23}\\
				0 & 0 & \lambda\\ 
				\end{array}}\in\PSL$$
				which contradicts that $\Gamma$ is discrete, unless the sequence $\set{\gamma_k g\gamma_k^{-1}}$ is eventually constant. Without loss of generality, we can assume that this sequence is constant and therefore,
				$$f_{23}^{(k)}=g_{23},\text{ for all }k.$$
				This means that 
				$$\gamma_{22}^{(k)}=\gamma_{33}^{(k)}=\xi_k\in\C^\ast$$
				and then 
				$$\gamma_k = \corchetes{\begin{array}{ccc}
				\xi_k^{-2} & \gamma_{12}^{(k)} & \gamma_{13}^{(k)} \\
				0 & \xi_k & \gamma_{23}^{(k)}\\
				0 & 0 & \xi_k\\ 
				\end{array}}.$$
				Using corollary \ref{cor_diagonal_es_1}, we get $\xi_k=1$ for all $k$ and then
				\begin{equation}\label{eq_thm_descomposicion_caso_noconmutativo_indices_9}
				\gamma_k = \corchetes{\begin{array}{ccc}
				1 & \gamma_{12}^{(k)} & \gamma_{13}^{(k)} \\
				0 & 1 & \gamma_{23}^{(k)}\\
				0 & 0 & 1\\ 
				\end{array}}.
				\end{equation}
				But then it is not possible that $\gamma_k\rightarrow\tau_{id}$. This means that there cannot be a sequence $\set{\gamma_k}\subset\Gamma$ converging to $\tau_{id}$.\\
				
				Now we deal with the case {\footnotesize(\textbf{LS1.2})}. Let us assume that $\Pi(\gamma_k)$ converges to an elliptic element of $\Sigma$ with the form
				$$\corchetes{\begin{array}{cc}
				e^{2\pi i \theta} & 0\\
				0 & 1\\
				\end{array}}.$$				
				As before, using the adequate conjugation we can assume that
				\begin{equation}\label{eq_thm_descomposicion_caso_noconmutativo_indices_6}
				\gamma_k\rightarrow\tau_\theta=\corchetes{\begin{array}{ccc}
				0 & 0 & 0 \\
				0 & e^{2\pi i \theta} & 0\\
				0 & 0 & 1\\ 
				\end{array}}.				
				\end{equation}
				Now we verify that $\gamma_k^n\rightarrow\tau_\theta^n$ for any $n\in\N$.\\
				
				Decomposing any matrix corresponding to $\gamma_k$ in its canonical Jordan form it follows that
				\begin{equation}\label{eq_thm_descomposicion_caso_noconmutativo_indices_7}
				\gamma_k^n=\corchetes{\begin{array}{ccc}
				{\gamma_{11}^{(k)}}^n & \gamma_{12}^{(k)}\frac{{\gamma_{11}^{(k)}}^n-{\gamma_{22}^{(k)}}^n}{\gamma_{11}^{(k)}-\gamma_{22}^{(k)}} & \alpha_{k,n} \\
				0 & e^{2\pi i n\theta} & \gamma_{23}^{(k)}\frac{{\gamma_{22}^{(k)}}^n-{\gamma_{33}^{(k)}}^n}{\gamma_{22}^{(k)}-\gamma_{33}^{(k)}}\\
				0 & 0 & 1\\ 
				\end{array}},				
				\end{equation}
				where 
				\begin{align*}
				\alpha_{k,n} &= \frac{{\gamma_{11}^{(k)}}^n\parentesis{\gamma_{11}^{(k)}\gamma_{13}^{(k)}-\gamma_{13}^{(k)}+\gamma_{12}^{(k)}\gamma_{23}^{(k)}}}{\parentesis{\gamma_{11}^{(k)}-\gamma_{22}^{(k)}}\parentesis{\gamma_{11}^{(k)}-\gamma_{33}^{(k)}}}-\frac{\gamma_{12}^{(k)}\gamma_{22}^{(k)}\gamma_{23}^{(k)}}{\parentesis{\gamma_{11}^{(k)}-\gamma_{22}^{(k)}}\parentesis{\gamma_{22}^{(k)}-\gamma_{33}^{(k)}}}\\
				&+\frac{{\gamma_{33}^{(k)}}^n\parentesis{\gamma_{12}^{(k)}\gamma_{23}^{(k)}-\gamma_{13}^{(k)}\gamma_{22}^{(k)}+\gamma_{13}^{(k)}\gamma_{33}^{(k)}}}{\parentesis{\gamma_{11}^{(k)}-\gamma_{33}^{(k)}}\parentesis{\gamma_{22}^{(k)}-\gamma_{33}^{(k)}}}.				
				\end{align*}						
				From (\ref{eq_thm_descomposicion_caso_noconmutativo_indices_6}) it follows that
				\begin{equation}\label{eq_thm_descomposicion_caso_noconmutativo_indices_8}
				\frac{\gamma_{11}^{(k)}}{\gamma_{33}^{(k)}},\frac{\gamma_{12}^{(k)}}{\gamma_{33}^{(k)}},\frac{\gamma_{13}^{(k)}}{\gamma_{33}^{(k)}},\frac{\gamma_{23}^{(k)}}{\gamma_{33}^{(k)}}\rightarrow 0,\;\;\;\frac{\gamma_{22}^{(k)}}{\gamma_{33}^{(k)}}\rightarrow e^{2\pi i \theta}.
				\end{equation}
				Dividing each entry of $\gamma_k^n$ by ${\gamma_{33}^{(k)}}^n$ and using (\ref{eq_thm_descomposicion_caso_noconmutativo_indices_7}) and (\ref{eq_thm_descomposicion_caso_noconmutativo_indices_8}), we get
				\begin{align*}
				\frac{{\gamma_{11}^{(k)}}^n}{{\gamma_{33}^{(k)}}^n} &\rightarrow 0\\
				\frac{\gamma_{12}^{(k)}}{{\gamma_{33}^{(k)}}^n}\frac{{\gamma_{11}^{(k)}}^n-{\gamma_{22}^{(k)}}^n}{\gamma_{11}^{(k)}-\gamma_{22}^{(k)}} &= \frac{\gamma_{12}^{(k)}}{{\gamma_{33}^{(k)}}}\parentesis{\frac{{\gamma_{11}^{(k)}}^{n-1}+{\gamma_{11}^{(k)}}^{n-2}{\gamma_{22}^{(k)}}+\cdots+{\gamma_{11}^{(k)}}{\gamma_{11}^{(k)}}^{n-2}+{\gamma_{22}^{(k)}}^{n-1}}{{\gamma_{11}^{(k)}}^{n-1}}}\\
				&= \frac{\gamma_{12}^{(k)}}{{\gamma_{33}^{(k)}}}\parentesis{\frac{{\gamma_{11}^{(k)}}^{n-1}}{{\gamma_{11}^{(k)}}^{n-1}}+\frac{{\gamma_{11}^{(k)}}^{n-2}}{{{\gamma_{11}^{(k)}}^{n-2}}}\frac{\gamma_{22}^{(k)}}{\gamma_{11}^{(k)}}+\cdots+\frac{\gamma_{11}^{(k)}}{\gamma_{11}^{(k)}}\frac{{\gamma_{11}^{(k)}}^{n-2}}{{{\gamma_{11}^{(k)}}^{n-2}}}+\frac{{\gamma_{22}^{(k)}}^{n-1}}{{{\gamma_{11}^{(k)}}^{n-1}}}}\\
				&\rightarrow 0\\
				\frac{\gamma_{23}^{(k)}}{{\gamma_{33}^{(k)}}^n}\frac{{\gamma_{22}^{(k)}}^n-{\gamma_{33}^{(k)}}^n}{\gamma_{22}^{(k)}-\gamma_{33}^{(k)}} &\rightarrow 0\\
				\frac{{\gamma_{22}^{(k)}}^n}{{\gamma_{33}^{(k)}}^n} &\rightarrow e^{2\pi i n\theta}\\
				\frac{\alpha_{k,n}}{{\gamma_{33}^{(k)}}^n} &\rightarrow 0.
				\end{align*}							 
				From this, it follows that $\gamma_k^n\rightarrow\tau_\theta^n$ for any $n\in\N$.\\
				
				If $\theta\in\Q$, then there exist $p\in\Z$ such that $p\theta\in\Z$ and then $\tau_{\theta}^p=\tau_{\id}$. Then we have a sequence $\set{\gamma_k^p}\subset \Gamma$ such that $\gamma_k^p\rightarrow\tau_{\id}$, which cannot happen, as we have already proven.\\
				
				If $\theta\in\R\setminus\Q$, consider the sequence $\set{e^{2\pi i n \theta}}\subset\Ss^1$. Since $\Ss^1$ is compact then there is a subsequence $\set{e^{2\pi i n_j \theta}}$ such that
				$$e^{2\pi i n_j \theta}\rightarrow 1.$$
				Then $\tau_\theta^{n_j}\rightarrow \tau_\id$ as $j\rightarrow\infty$. Consider the diagonal sequence $\set{\gamma_k^{n_k}}\subset\Gamma$, then
				$$\gamma_k^{n_k}\rightarrow \tau_\id$$
				which cannot happen. All of the above proves that $\Pi(\gamma_k)$ cannot converge to an elliptic element of $\Sigma$.\\
				
				Next, consider the case {\footnotesize(\textbf{LS1.3})}. Assume that 
				\begin{equation}\label{eq_thm_descomposicion_caso_noconmutativo_indices_10}
				\gamma_k\rightarrow\tau_b = \corchetes{\begin{array}{ccc}
				0 & 0 & 0\\
				0 & 1 & b\\
				0 & 0 & 1				
				\end{array}}.
				\end{equation}
				Therefore, the sequence $\set{\Pi(\gamma_k)}$ converges to the parabolic e-lement 
				$$\corchetes{\begin{array}{cc}
				1 & b\\
				0 & 1\\
				\end{array}}.$$	
				Since $\Sigma$ is not discrete, it follows from lemma \ref{lem_parte_parab_no_discreta} that there exists a sequence $\set{h_k}\subset\Gamma$ such that $\Pi(h_k)$ is parabolic for all $k$, and 
				$$\Pi(h_k)=\corchetes{\begin{array}{cc}
				1 & \varepsilon_k\\
				0 & 1\\
				\end{array}}\rightarrow \id,$$
				for some sequence $\varepsilon_k\rightarrow 0$. Then, using the same reasoning as in the case {\footnotesize(\textbf{LS1.1})} (see (\ref{eq_thm_descomposicion_caso_noconmutativo_indices_9})), we know that 
				$$h_k=\corchetes{\begin{array}{ccc}
				1 & h_{12} & h_{13}\\
				0 & 1 & \varepsilon_k\\
				0 & 0 & 1				
				\end{array}}.$$
				Consider the sequence
				$$\gamma_k h_k \gamma_k^{-1} h_k^{-1} = \corchetes{\begin{array}{ccc}
				1 & h_{12} \parentesis{\frac{\gamma_{11}^{(k)}}{\gamma_{22}^{(k)}}-1} & f_{13}^{(k)}\\
				0 & 1 & \parentesis{\frac{\gamma_{22}^{(k)}}{\gamma_{33}^{(k)}}-1}\varepsilon_k\\
				0 & 0 & 1				
				\end{array}},$$
				where
				$$f_{13}^{(k)}=h_{13}\parentesis{\frac{\gamma_{11}^{(k)}}{\gamma_{33}^{(k)}}-1}+\frac{\gamma_{12}^{(k)}}{\gamma_{33}^{(k)}}\varepsilon_k + h_{12}\varepsilon_k-h_{12}\parentesis{\frac{\gamma_{11}^{(k)}}{\gamma_{23}^{(k)}}\frac{\gamma_{23}^{(k)}}{\gamma_{33}^{(k)}}+\frac{\gamma_{11}^{(k)}}{\gamma_{22}^{(k)}}\varepsilon_k}.$$
				From (\ref{eq_thm_descomposicion_caso_noconmutativo_indices_10}) it follows
				\begin{equation}\label{eq_thm_descomposicion_caso_noconmutativo_indices_11}
				\frac{\gamma_{11}^{(k)}}{\gamma_{33}^{(k)}},\frac{\gamma_{12}^{(k)}}{\gamma_{33}^{(k)}},\frac{\gamma_{13}^{(k)}}{\gamma_{33}^{(k)}}\rightarrow 0,\;\;\;\frac{\gamma_{23}^{(k)}}{\gamma_{33}^{(k)}}\rightarrow b,\;\;\;\frac{\gamma_{22}^{(k)}}{\gamma_{33}^{(k)}}\rightarrow 1,
				\end{equation}
				and therefore 
				$$\gamma_k h_k \gamma_k^{-1} h_k^{-1}\rightarrow \corchetes{\begin{array}{ccc}
				1 & h_{12} & 0\\
				0 & 1 & 0\\
				0 & 0 & 1				
				\end{array}}\in\PSL,$$
				contradicting that $\Gamma$ is discrete. We conclude that this case cannot occur.\\
				
				Finally, we deal with the case {\footnotesize(\textbf{LS1.4})}. Let us suppose that 
				\begin{equation}\label{eq_thm_descomposicion_caso_noconmutativo_indices_10}
				\gamma_k\rightarrow\tau_\alpha = \corchetes{\begin{array}{ccc}
				0 & 0 & 0\\
				0 & \alpha & 0\\
				0 & 0 & 1				
				\end{array}}.
				\end{equation}
				Since $\Lambda_{Gr}(\Sigma)=\Ss^1$, there is an element $h=\corchetes{h_{ij}}\in\Gamma$ such that $h_{23}\neq 0$ and 
				$$h = \corchetes{\begin{array}{ccc}
				1 & h_{12} & h_{13}\\
				0 & \lambda & h_{23}\\
				0 & 0 & \lambda^{-1}				
				\end{array}},$$
				with $\valorabs{\lambda}\neq 1$. Consider the sequence
				$$f_k:= h\gamma_k h^{-1} \gamma_k^{-1} = \corchetes{\begin{array}{ccc}
				1 & f_{12}^{(k)} & f_{13}^{(k)}\\
				0 & 1 & f_{23}^{(k)} \\
				0 & 0 & 1				
				\end{array}},$$
				where
				\begin{align*}
				f_{12}^{(k)} &= \gamma_{12}^{(k)}\frac{1-\lambda}{\lambda\gamma_{22}^{(k)}} + h_{12}\frac{\gamma_{22}^{(k)}-\gamma_{11}^{(k)}}{\lambda\gamma_{22}^{(k)}},\\
				f_{13}^{(k)} &= -\gamma_{33}^{(k)}\frac{\gamma_{12}^{(k)}+h_{12}\parentesis{\gamma_{12}^{(k)}-\gamma_{11}^{(k)}}}{\lambda\gamma_{22}^{(k)}\gamma_{33}^{(k)}} - \frac{h_{23}\gamma_{12}^{(k)}+\gamma_{13}^{(k)}+h_{12}h_{23}\parentesis{\gamma_{22}^{(k)}-\gamma_{11}^{(k)}}}{\gamma_{33}^{(k)}}  \\
				&+ \frac{\gamma_{12}^{(k)}\gamma_{23}^{(k)}}{\gamma_{22}^{(k)}\gamma_{33}^{(k)}} + \lambda\frac{\gamma_{13}^{(k)}+h_{12}\gamma_{23}^{(k)}+h_{13}\parentesis{\gamma_{33}^{(k)}-\gamma_{11}^{(k)}}}{\gamma_{33}^{(k)}}\\
				f_{23}^{(k)} &= \gamma_{23}^{(k)}\frac{\lambda^2-1}{\gamma_{33}^{(k)}} + \lambda h_{23}\frac{\gamma_{33}^{(k)}-\gamma_{22}^{(k)}}{\gamma_{33}^{(k)}}.
				\end{align*}
				from (\ref{eq_thm_descomposicion_caso_noconmutativo_indices_10}), it follows
				$$f_k\rightarrow  \corchetes{\begin{array}{ccc}
				1 & \frac{h_{12}}{\lambda} & h_{13}\parentesis{\lambda-h_{12}\alpha}\\
				0 & 1 & \lambda h_{23}(1-\alpha) \\
				0 & 0 & 1				
				\end{array}},$$
				contradicting that $\Gamma$ is discrete.\\
				
				We have proved that sequences in $\Gamma$ cannot have a quasi projective limit with the form
				$$\tau=\corchetes{\begin{array}{ccc}
				0 & \gamma_{12} & \gamma_{13} \\
				0 & \gamma_{22} & \gamma_{23}\\
				0 & 0 & 1\\ 
				\end{array}},$$
				when $\Lambda_{Gr}(\Sigma)=\Ss^1$.\\
				
				Observe that the case {\footnotesize(\textbf{LS2})} cannot happen, otherwise the action of $\Gamma$ would be nowhere proper and discontinuous.\\
				
				Consider the case {\footnotesize(\textbf{LS3})}, that is, $\valorabs{\Lambda_{Gr}(\Sigma)}=1$. Using proposition \ref{prop_opciones_conj_lim_greenberg} we know that $\Sigma$ is conjugated to a subgroup of $\epa$. Let $h\in\Gamma$ such that $\Pi(h)$ is parabolic, that is,
				$$\Pi(h)=\corchetes{\begin{array}{cc}
				1 & h_{23}\\
				0 & 1\\
				\end{array}}$$
				for some $h_{23}\neq 0$. Then
				$$h=\corchetes{\begin{array}{ccc}
				\lambda^{-2} & h_{12} & h_{13} \\
				0 & \lambda & \lambda h_{23} \\
				0 & 0 & \lambda
				\end{array}},$$
				for some $\lambda\in\C^\ast$ and $h_{12},h_{13}\in\C$. Consider the sequence of distinct elements $\set{\eta_k}\subset \Gamma$ given by
				$$\eta_k:=h\gamma_k h^{-1} \gamma_k^{-1}=\corchetes{\begin{array}{ccc}
				1 & \frac{\gamma_{12}^{(k)}}{\gamma_{22}^{(k)}}\parentesis{\gamma^{-3}-1}+\frac{h_{12}}{\lambda}\parentesis{1-\frac{\gamma_{11}^{(k)}}{\gamma_{22}^{(k)}}} & \eta_{13}^{(k)} \\
				0 & 1 & h_{23}\parentesis{1-\frac{\gamma_{22}^{(k)}}{\gamma_{33}^{(k)}}} \\
				0 & 0 & 1
				\end{array}},$$
				where $\eta_{13}^{(k)}$ is a sequence depending on all the terms $h_{ij}$, $\gamma_{ij}^{(k)}$ and $\lambda$ (we don't write the explicit expression for the sake of clarity, however we will write its limit). Using the limits calculated at the beginning of this case (b) yields
				\begin{align*}
				\frac{\gamma_{12}^{(k)}}{\gamma_{22}^{(k)}}\parentesis{\gamma^{-3}-1}+\frac{h_{12}}{\lambda}\parentesis{1-\frac{\gamma_{11}^{(k)}}{\gamma_{22}^{(k)}}} & \longrightarrow \frac{\gamma_{12}}{\gamma_{22}}\parentesis{\gamma^{-3}-1}+ \frac{h_{12}}{\lambda} \\
				\eta_{13}^{(k)} &\longrightarrow \lambda^{-3}\parentesis{\gamma_{13}-\frac{\gamma_{12}\gamma_{23}}{\gamma_{22}\gamma_{33}}-\gamma_{12}h_{23}}+\\
				&\frac{h_{23}}{\lambda}\parentesis{1-h_{12}\gamma_{22}}-\gamma_{13}+\frac{\gamma_{12}\gamma_{23}}{\gamma_{22}\gamma_{33}} \\
				h_{23}\parentesis{1-\frac{\gamma_{22}^{(k)}}{\gamma_{33}^{(k)}}}  &\longrightarrow h_{23}(1-\gamma_{22})
				\end{align*}
				and therefore, $\eta_k$ converges to an element of $\PSL$, contradicting that $\Gamma$ is discrete. This proves that the case {\footnotesize(\textbf{LS3})} cannot happen.\\				
				
				Now, finally, consider the case {\footnotesize(\textbf{LS4})}, that is, $\valorabs{\Lambda_{Gr}(\Sigma)}=2$. Using proposition \ref{prop_opciones_conj_lim_greenberg} we know that $\Sigma$ is conjugated to a subgroup of $\autc$. This means that, conjugating by a suitable element of $\PSL$, every element of $\Sigma$ has the form
				$$\Pi(g)=\corchetes{\begin{array}{cc}
				\alpha & 0\\
				0 & \alpha^{-1}\\
				\end{array}}$$
				for some $\alpha\in\C^\ast$. Then,
				\begin{equation}\textsc{\label{eq_thm_descomposicion_caso_noconmutativo_indices_401}}							
				g=\corchetes{\begin{array}{ccc}
				\lambda^{-2} & g_{12} & g_{13} \\
				0 & \lambda\alpha & 0\\
				0 & 0 & \lambda\alpha^{-1}\\ 
				\end{array}}.
				\end{equation}
				Let $g_1,g_2\in\Gamma$ be two elements such that $\corchetes{g_1,g_2}\neq\id$, then
				$$h=\corchetes{g_1,g_2}=\corchetes{\begin{array}{ccc}
				1 & h_{12} & h_{13} \\
				0 & 1 & 0\\
				0 & 0 & 1\\ 
				\end{array}},$$
				with $\valorabs{h_{12}}+\valorabs{h_{13}}\neq 0$. The sequence $\set{\gamma_k}$ has the same form described in (\ref{eq_thm_descomposicion_caso_noconmutativo_indices_401}),
				$$\gamma_k=\corchetes{\begin{array}{ccc}
				\gamma_{11}^{(k)} & \gamma_{12}^{(k)} & \gamma_{13}^{(k)} \\
				0 & \gamma_{22}^{(k)} & 0\\
				0 & 0 & \gamma_{33}^{(k)}\\ 
				\end{array}}.$$
				Then the quasi-projective limit $\tau$, given by (\ref{eq_thm_descomposicion_caso_noconmutativo_indices_001}) has the form
				$$\tau=\corchetes{\begin{array}{ccc}
				0 & 0 & 0 \\
				0 & \gamma_{22} & 0\\
				0 & 0 & 1\\ 
				\end{array}},$$
				and then
				\begin{equation}\label{eq_thm_descomposicion_caso_noconmutativo_indices_402}
				\frac{\gamma_{11}^{(k)}}{\gamma_{33}^{(k)}},\frac{\gamma_{12}^{(k)}}{\gamma_{33}^{(k)}},\frac{\gamma_{13}^{(k)}}{\gamma_{33}^{(k)}}\rightarrow 0,\;\;\;\frac{\gamma_{22}^{(k)}}{\gamma_{33}^{(k)}}\rightarrow \gamma_{22}.
				\end{equation}
				Consider the sequence $\set{f_k}\in\Gamma$ given by
				$$f_k=h \gamma_k h^{-1}\gamma_k^{-1}=\corchetes{\begin{array}{ccc}
				1 & h_{12}\parentesis{1-\frac{\gamma_{11}^{(k)}}{\gamma_{22}^{(k)}}} & h_{13}\parentesis{1-\frac{\gamma_{11}^{(k)}}{\gamma_{33}^{(k)}}} \\
				0 & 1 & 0\\
				0 & 0 & 1\\ 
				\end{array}}.$$
				From (\ref{eq_thm_descomposicion_caso_noconmutativo_indices_402}) it follows that
				$$f_k\rightarrow \corchetes{\begin{array}{ccc}
				1 & h_{12} & h_{13} \\
				0 & 1 & 0\\
				0 & 0 & 1\\ 
				\end{array}}\in\PSL,$$
				contradicting that $\Gamma$ is discrete. Therefore the case {\footnotesize(\textbf{LS4})} cannot happen when $\gamma_{33}^{(k)}\rightarrow\infty$.				 
			\item $\gamma_{33}^{(k)}\rightarrow 0$. From this, it follows that $\gamma_{22}^{(k)}\rightarrow 0$, otherwise $\alpha\nin\psl$. Since $\gamma_{11}^{(k)}\gamma_{22}^{(k)}\gamma_{33}^{(k)}=1$, we have that $\gamma_{11}^{(k)}\rightarrow \infty$. Also, $\gamma_{22}^{(k)}\rightarrow a\in\C$. All of this yields
			\begin{equation}\label{eq_thm_descomposicion_caso_noconmutativo_indices_3}
			\gamma_k\rightarrow\tau=\corchetes{\begin{array}{ccc}
		a & b & c \\
		0 & 0 & 0\\
		0 & 0 & 0\\ 
		\end{array}},
			\end{equation}
			with $a\neq 0$, unless
				\begin{equation}\label{eq_thm_descomposicion_caso_noconmutativo_indices_1}
				\frac{\gamma_{11}^{(k)}}{\gamma_{12}^{(k)}}\rightarrow 0,\;\;\;\text{or}\;\;\;\frac{\gamma_{11}^{(k)}}{\gamma_{13}^{(k)}}\rightarrow 0.
				\end{equation}
			But then we are in case (b). Since (\ref{eq_thm_descomposicion_caso_noconmutativo_indices_1}) does not happen, then 
				\begin{equation}\label{eq_thm_descomposicion_caso_noconmutativo_indices_2}
				\gamma_{12}^{(k)}\parentesis{\gamma_{11}^{(k)}}^{-1}\rightarrow b,\;\;\;\text{and}\;\;\;\gamma_{13}^{(k)}\parentesis{\gamma_{11}^{(k)}}^{-1}\rightarrow c,
				\end{equation}
			for some $b,c\in\C$.
				
			Dividing each entry of (\ref{eq_thm_descomposicion_caso_noconmutativo_indices_0}) by the entry, ${\gamma_{33}^{(k)}}^{-1}$, we have
			\begin{align*}
			\frac{{\gamma_{11}^{(k)}}^{-1}}{{\gamma_{33}^{(k)}}^{-1}} = \frac{\gamma_{33}^{(k)}}{\gamma_{11}^{(k)}} & \rightarrow 0. \\
			-\frac{\gamma_{12}^{(k)}}{\gamma_{11}^{(k)}\gamma_{22}^{(k)}{\gamma_{33}^{(k)}}^{-1}} = -\frac{\gamma_{12}^{(k)}}{\gamma_{11}^{(k)}}\frac{\gamma_{33}^{(k)}}{\gamma_{22}^{(k)}}&\rightarrow - \frac{b}{\gamma_{22}}.\\
			\frac{\gamma_{12}^{(k)}\gamma_{23}^{(k)}-\gamma_{13}^{(k)}\gamma_{22}^{(k)}}{\gamma_{11}^{(k)}\gamma_{22}^{(k)}\gamma_{33}^{(k)}{\gamma_{33}^{(k)}}^{-1}} = \frac{\gamma_{12}^{(k)}\gamma_{23}^{(k)}}{\gamma_{11}^{(k)}\gamma_{22}^{(k)}}-\frac{\gamma_{13}^{(k)}}{\gamma_{11}^{(k)}} &\rightarrow b\frac{\gamma_{23}}{\gamma_{22}}-c.\\
			\frac{{\gamma_{22}^{(k)}}^{-1}}{{\gamma_{33}^{(k)}}^{-1}} = \frac{\gamma_{33}^{(k)}}{\gamma_{22}^{(k)}} & \rightarrow \frac{1}{\gamma_{22}}. \\
			-\frac{\gamma_{23}^{(k)}}{\gamma_{22}^{(k)}\gamma_{33}^{(k)}{\gamma_{33}^{(k)}}^{-1}} = - \frac{\gamma_{23}^{(k)}}{\gamma_{22}^{(k)}} &\rightarrow -\frac{\gamma_{23}}{\gamma_{22}}.		
			\end{align*}
			And therefore, the sequence of the inverses satisfies
			$$
			\gamma_{k}^{-1}\rightarrow \mu = \corchetes{\begin{array}{ccc}
		0 &  - \frac{b}{\gamma_{22}} & b\frac{\gamma_{23}}{\gamma_{22}}-c \\
		0 & \frac{1}{\gamma_{22}} & -\frac{\gamma_{23}}{\gamma_{22}}\\
		0 & 0 & 1\\ 
		\end{array}}.			
			$$
			This means that, if we consider the sequence of the inverses instead, this case is the same as case (b).			\end{enumerate}
			
		\item[\titem{[Div]}] The sequence $\set{\Pi(\gamma_k)}\subset\Sigma$ diverges. Proposition \ref{prop_opciones_conj_lim_greenberg} implies that one of the following cases occur:
			\begin{itemize}
			\item[\titem{[Div-1]}] $\overline{\Sigma}=\so$. This case cannot happen since $\Sigma$ is solvable but $\so$ is not (see example \ref{ej_solvable}).
			\item[\titem{[Div-2]}] $\overline{\Sigma}=\rot$. This case cannot happen, otherwise $\grL(\Sigma)=\emptyset$, contradicting Lemma \ref{lem_control_nodiscreto_descartar_grL_vacio}.
			\item[\titem{[Div-3]}] $\overline{\Sigma}=\dih$. This case cannot happen by the same argument of case \titem{[Div-2]}.
			\item[\titem{[Div-4]}] $\overline{\Sigma}=\psl$. This case cannot happen since $\psl$ is not solvable (as it contains Schottky groups, see Theorem \ref{teo_tits2}).
			\item[\titem{[Div-5]}] The group $\overline{\Sigma}$ is a subgroup of the affine group $\epa$. First, observe that there cannot be elliptic elements in $\Sigma$. If there is an elliptic element 
				$$\Pi(\gamma)=\corchetes{\begin{array}{cc}
				e^{i\theta} & a\\
				0 & e^{-i\theta}\\
				\end{array}}\in\Sigma,$$
			where $\theta\neq 0$, then, either $\theta\in\Q$ or $\theta\in\R\setminus\Q$. If the former occurs then $\lambda_{23}(\gamma)$ would be a torsion element in the torsion-free group $\lambda_{23}(\Gamma)$. If the latter happens, then $\gamma$ would be an irrational screw, contradicting that $\Gamma$ is not commutative (see Corollary \ref{cor_IS_conmutativo}). Then every element of $\Sigma$ has the form
				$$\Pi(\gamma)=\corchetes{\begin{array}{cc}
				1 & a\\
				0 & 1\\
				\end{array}},$$
			where $a\in A$, for some additive group $A\subset(\C,+)$. Furthermore, $A\cong \Sigma$, and therefore, $A$ is not discrete. As a consequence of Proposition \ref{prop_forma_subgrupos_aditivos}, either $\overline{A}\cong \R$ or $\overline{A}\cong \R\oplus\Z$ (see Figure \ref{fig_OpcionesBandas_Div5}).
				\begin{figure}[H]
				\begin{center}
				\includegraphics[height=35mm]{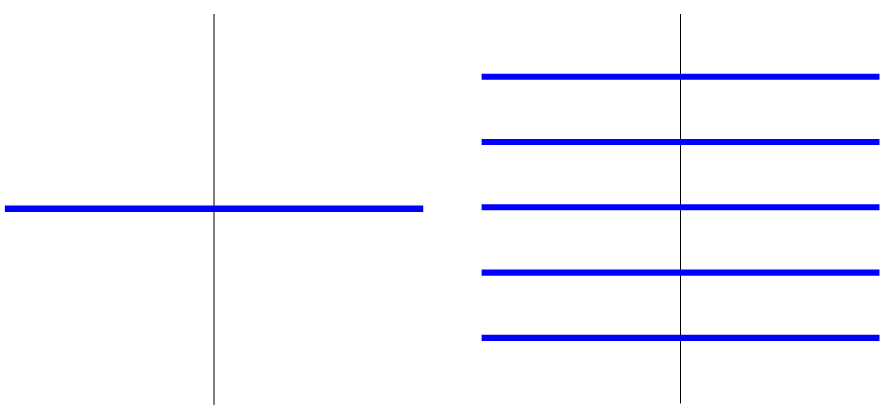}
				\caption{The two possibilities for the set $\overline{A}\subset\C$.}
				\label{fig_OpcionesBandas_Div5}
				\end{center}	
				\end{figure}
				Let us define the following union of pencils of lines passing through $e_1$,				
					$$\Lambda=\set{e_1}\cup\parentesis{\bigcup_{p\in\overline{A}}\pi(p)^{-1}},$$
				and $\Omega=\CP^2\setminus \Lambda$ (see Figure \ref{fig_OpcionesConos_Div5}). 
				\begin{figure}[H]
				\begin{center}
				\includegraphics[height=35mm]{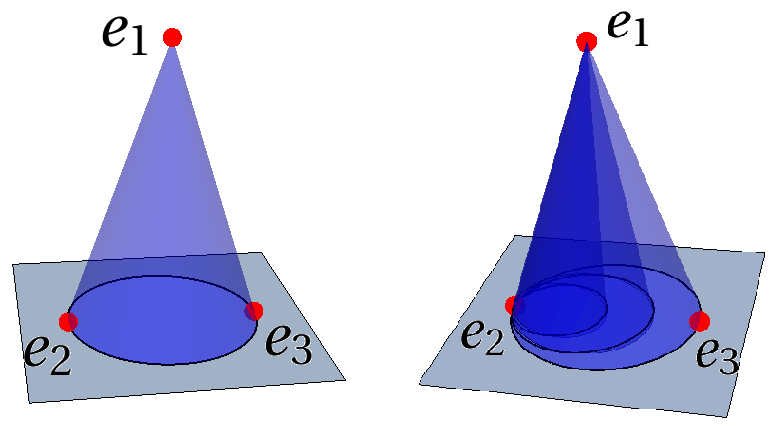}
				\caption{The two possibilities for the set $\Omega\subset\CP^2$.}
				\label{fig_OpcionesConos_Div5}
				\end{center}	
				\end{figure}				
				The $\Sigma$-orbits of compact subsets of $\hC\setminus A$ accumulate on $\infty$, which correspond to the point $\pi(e_2)$. Then, the orbits of compact subsets of $\Omega$ accumulate on $\set{e_1}$, and therefore, the action of $\Gamma$ on each connected component of $\Omega$ is proper and discontinuous. Since each connected component of $\Omega$ is simply connected, this concludes the case \titem{[Div-5]}.	
				
			\item[\titem{[Div-6]}] The group $\overline{\Sigma}$ is a subgroup of the group $\autc$. Then $\Sigma$ is a purely loxodromic group and then, up to conjugation, each element has the form
				\begin{equation}\label{eq_thm_descomposicion_caso_noconmutativo_div6_1}
				\Pi(\gamma)=\corchetes{\begin{array}{cc}
				\alpha & 0\\
				0 & 1\\
				\end{array}},
				\end{equation}
				for some $\valorabs{\alpha}\neq 1$. Let $G=\lambda_{23}(\Gamma)$, from (\ref{eq_thm_descomposicion_caso_noconmutativo_div6_1}) it follows that $G\cong\Sigma$, and since $\Sigma$ is not discrete, then $G$ is not discrete. Let us write each $\alpha\in G$ as
				$$\begin{array}{ll}
				\alpha = re^{i\theta}, & r\in\R^+,\; r\in A\subset (\C^\ast,\cdot)\\
				& \theta\in\R,\; \theta\in B\subset (\C,+),
				\end{array}$$
			for some multiplicative group $A$ and some additive group $B$. Then $G\cong A\times B$ and, since $G$ is not discrete, then $A$ is not discrete or $B$ is not discrete (it is not possible that both $A$ and $B$ are not discrete, otherwise the action of $\Gamma$ would be nowhere proper and discontinuous). Therefore we have two possibilities:
			\begin{itemize}
			\item $A$ is discrete and $B$, not discrete. Since $A$ is a discrete multiplicative subgroup of $\C$, then $\rank(A)=1$ and therefore, $A=\SET{r^n}{n\in\Z}$ for some $r\in\C^\ast$. Hence, 
				$$\overline{G}=\overline{\lambda_{23}(\Gamma)}=\SET{r^ne^{i\theta}}{\theta\in\left[0,2\pi\right),\;n\in\Z}.$$
				In particular, 
				$$\SET{r^ne^{i\theta}}{\theta\in\left[0,2\pi\right)}\subset\overline{\lambda_{23}(\Gamma)}.$$
				This means that there exists an element $\gamma\in\Gamma$ such that $\lambda_{23}=e^{i\theta}$ for some $\theta\in\left[0,2\pi\right)$. If $\theta\in\Q$, then $\lambda_{23}(\gamma)$ is a torsion e-lement, contradicting that $\lambda_{23}(\Gamma)$ is torsion free. If $\theta\in\R\setminus\Q$, then $\gamma$ is an irrational screw, contradicting that $\Gamma$ is non-commutative (see Corollary \ref{cor_IS_conmutativo}). This dismisses this possibility.
			\item $A$ is not discrete and $B$, discrete. Since $B$ is discrete, then $B$ is finite and $B\cong\Z_k$, for some $k\in\N$. Therefore $\overline{G}\subset\C$ looks like the the set depicted in Figure \ref{fig_G_Anodisc_Bdisc}.
				\begin{figure}[H]
				\begin{center}
				\includegraphics[height=40mm]{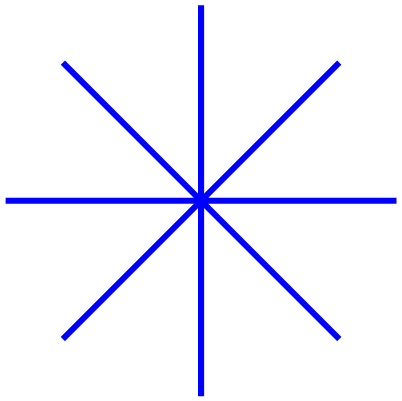}
				\caption{The set $\overline{G}\subset\C$.}
				\label{fig_G_Anodisc_Bdisc}
				\end{center}	
				\end{figure}
				We know there is an homomorphism $\varphi:\linproy{e_2,e_3}\rightarrow\hC$. For the sake of simplifying notation, we will denote indistinctly a point $x\in\overline{G}\subset\C$ and the point $\varphi^{-1}(x)\in \linproy{e_2,e_3}$. Let $\Lambda\subset\CP^2$ be given by
					$$\Lambda=\set{e_1}\cup\parentesis{\bigcup_{p\in\overline{G}}\pi(p)^{-1}},$$
				and $\Omega=\CP^2\setminus\Lambda$ (see Figure \ref{fig_cono_lenticular}).
				\begin{figure}[H]
				\begin{center}
				\includegraphics[height=40mm]{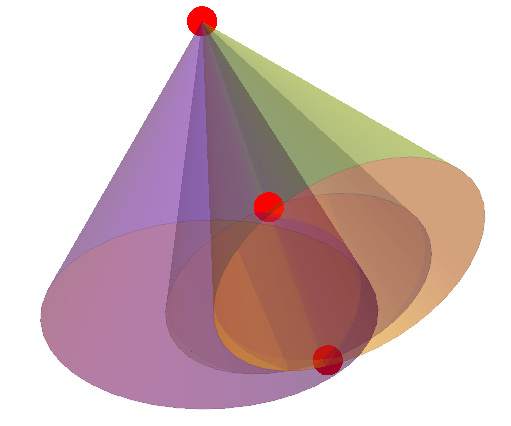}
				\caption{The set $\Lambda\subset\CP^2$.}
				\label{fig_cono_lenticular}
				\end{center}	
				\end{figure}
				Since the orbits of compact subsets of $\C\setminus\overline{G}$ accumulate on $\pi(e_2)$ and $\pi(e_3)$, then the orbits of compact subsets of $\Omega$ accumulate on $\set{e_1}$. This means that the action of $\Gamma$ on $\Omega$ is proper and discontinuous. Besides, each connected component of $\Omega$ is simply connected (see Figure \ref{fig_cono_lenticular}). This concludes the case \titem{[Div-6]}, since we have constructed a simply connected region where $\Gamma$ acts properly and discontinuously.	
			\end{itemize}			 
			\item[\titem{[Div-7]}] The group $\overline{\Sigma}$ is a subgroup of the group $\pslr$. Then, $\grL(\Sigma)\cong\hat{\R}$ (Proposition \ref{prop_opciones_conj_lim_greenberg}). Then, up to a suitable conjugation, the orbits of compact subsets of $\CP^1\setminus \grL(\Sigma)$ accumulate on $\hat{\R}$, we regard the points $\pi(e_2)$ and $\pi(e_3)$ as the points $0$ and $\infty$ in this euclidian circle. We define the pencil of lines passing through $e_1$,
				$$\Lambda=\set{e_1}\cup\parentesis{\bigcup_{p\in\hat{\R}}\pi(p)^{-1}}.$$ 
				\begin{figure}[H]
				\begin{center}
				\includegraphics[height=35mm]{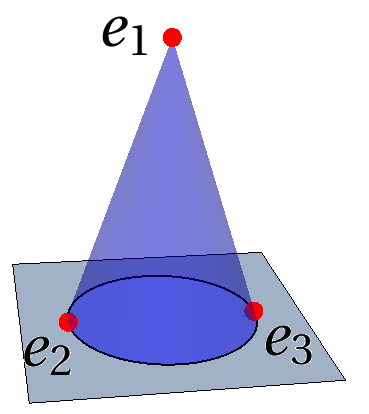}
				\caption{The set $\Lambda\subset\CP^2$.}
				\label{fig_cono_sencillo_Div7}
				\end{center}	
				\end{figure} 
			Then its complement $\Omega=\CP^2\setminus\Lambda$ is homeomorphic to a cone $\C\times\parentesis{\Hh^+\cup\Hh^-}$ and therefore, each of its connected components are simply connected (see Figure \ref{fig_cono_sencillo_Div7}). The action of $\Gamma$ on each of these connected components is proper and discontinuous, this completes this case \titem{[Div-7]}.
				
			\end{itemize}
		\end{itemize}
	\end{enumerate}
\end{proof} 

In each subcase of the previous proof, for every divergent sequence $a:=\set{\gamma_k}\subset\Gamma$, we have constructed an open subset $U_a\subset\CP^2$, such that the the orbits of every compact set $K\subset U_a$ accumulate on $\CP^2\setminus U_a$. Thus, we can define a limit set for the action of $\Gamma$ in the following form
	$$\Lambda_S(\Gamma):=\overline{\bigcup_{a}\parentesis{\CP^2\setminus U_a}}.$$  
This limit set describes the dynamics of non-commutative upper triangular subgroups of $\PSL$.

\subsection{Consequences of the theorem of decomposition of non-commutative groups}

Theorem \ref{thm_descomposicion_caso_noconmutativo} gives a decomposition of the group $\Gamma$ in four layers, the first two layers are made of parabolic elements and the last two layers are made of loxodromic elements (see Corollary \ref{cor_lox_nivel34} and Table \ref{fig_noconmutativo_capas}). The description of these four layers are given in the proof of previous theorem and are summarized in the next table.

\begin{table}[H]
\begin{center}
  \begin{tabular}{  c  c  c  c  }
    
    \multicolumn{2}{c}{Parabolic} & \multicolumn{2}{c}{Loxodromic}\\
    \multicolumn{2}{c}{$\overbrace{\hspace{25ex}}$} & \multicolumn{2}{c}{$\overbrace{\hspace{45ex}}$}\\
    $\corchetes{\begin{array}{ccc}
	1 & x & y\\
	0 & 1 & 0 \\
	0 & 0 & 1\\
	\end{array}}$ & $\corchetes{\begin{array}{ccc}
	1 & x & y\\
	0 & 1 & z\\
	0 & 0 & 1\\
	\end{array}}$ & $\corchetes{\begin{array}{ccc}
	\alpha & x & y\\
	0 & \beta & z\\
	0 & 0 & \beta\\
	\end{array}}$ & $\corchetes{\begin{array}{ccc}
	\alpha & x & y\\
	0 & \beta & z\\
	0 & 0 & \gamma\\
	\end{array}}$\\ 
     & $z\neq 0$ & $\alpha\neq\beta$, $z\neq 0$ & $\beta\neq\gamma$     \\
     & & \begin{tabular}{c}
     \small{Loxo-parabolic} \\
     \end{tabular} & \begin{tabular}{c}
     \small{Loxo-parabolic} \\
     \small{Complex homothety} \\
     \small{Strongly loxodromic} \\
     \end{tabular} \\  
    $\core(\Gamma)$ & $A\setminus\kernel(\Gamma)$ & $\text{Ker}(\lambda_{23})\setminus A$ & \\
	\multicolumn{2}{c}{$\underbrace{\hspace{25ex}}$} & & \\    
    \multicolumn{2}{c}{$A=\text{Ker}(\lambda_{13})\cap\text{Ker}(\lambda_{23})$} & & \\
    \multicolumn{3}{c}{$\underbrace{\hspace{50ex}}$} & \\
    \multicolumn{3}{c}{$\text{Ker}(\lambda_{23})$} & $\Gamma\setminus \text{Ker}(\lambda_{23})$\\
	\multicolumn{4}{c}{$\underbrace{\hspace{75ex}}$} \\    
    \multicolumn{4}{c}{$\Gamma$} \\
    
  \end{tabular}
\end{center}
\caption{\small The decomposition of a non-commutative subgroup of $U_+$ in four layers.}
\label{fig_noconmutativo_capas}
\end{table}

Let $F_1,F_2,F_3\subset U_+$ be the pairwise disjoint subsets defined as
	\begin{align*}
	F_1&=\SET{\alpha=\corchetes{\begin{array}{ccc}
	\alpha_{11} & 0 & \alpha_{13}\\
	0 & \alpha_{22} & 0\\
	0 & 0 & \alpha_{33}
	\end{array}}}{\alpha\in U_+} \\
	F_2&=\SET{\alpha=\corchetes{\begin{array}{ccc}
	\alpha_{11} & \alpha_{12} & \alpha_{13}\\
	0 & \alpha_{22} & 0\\
	0 & 0 & \alpha_{33}
	\end{array}}}{\alpha_{12}\neq 0\text{, }\alpha\in U_+} \\
	F_3&=\SET{\alpha=\corchetes{\begin{array}{ccc}
	\alpha_{11} & 0 & \alpha_{13}\\
	0 & \alpha_{22} & \alpha_{23}\\
	0 & 0 & \alpha_{33}
	\end{array}}}{\alpha_{23}\neq 0\text{, }\alpha\in U_+}	\\
	F_4&=\SET{\alpha=\corchetes{\begin{array}{ccc}
	\alpha_{11} & \alpha_{12} & \alpha_{13}\\
	0 & \alpha_{22} & \alpha_{23}\\
	0 & 0 & \alpha_{33}
	\end{array}}}{\alpha_{12},\alpha_{23}\neq 0,\alpha\in U_+}	
	\end{align*}	

These four subsets classify the elements of $U_+$ depending on whether they have zeroes in positions 12 and 23. We need this classification because, as we will see in the next proposition, two elements of $U_+$ commute only if they have the same \emph{form} given by the set $F_i$ they belong to. This argument will be key to prove Corollary \ref{cor_descomposicion_bloques_noconm}.

\begin{prop}\label{prop_formas_conmutativas}
Let $\Gamma\subset U_+$ be a subgroup and let $\alpha=\corchetes{\alpha_{ij}},\beta=\corchetes{\beta_{ij}}\in\Gamma\setminus\set{\id}$. If $\corchetes{\alpha,\beta}=\id$ then $\alpha,\beta\in F_i$ for some $i=1,2,3,4$.
\end{prop}

\begin{proof}
Let $\alpha=\corchetes{\alpha_{ij}}$ and $\beta=\corchetes{\beta_{ij}}$ two elements in $\Gamma$, suppose that they commute. Denote by $\alpha_1$ and $\alpha_2$ (resp. $\beta_1$ and $\beta_2$) the upper left and bottom right $2\times 2$ blocks of $\alpha$ (resp. $\beta$)
	$$\alpha_1=\corchetes{\begin{array}{cc}
	\alpha_{11} & \alpha_{12} \\
	0 & \alpha_{22} 
	\end{array}},
	\;\;\;\;
	\alpha_2=\corchetes{\begin{array}{cc}
	\alpha_{22} & \alpha_{23} \\
	0 & \alpha_{33} 
	\end{array}}.
	$$
A direct calculation shows that, since $\alpha$ and $\beta$ commute, $\alpha_i$ and $\beta_i$ commute, for $i=1,2$. Considering $\alpha_i$ and $\beta_i$ as elements of $\psl$ we observe that 
	\begin{align*}
	\fix(\alpha_1)&=\set{e_1,\corchetes{\alpha_{12}:\alpha_{22}-\alpha_{11}}}\\
	\fix(\alpha_2)&=\set{e_1,\corchetes{\alpha_{23}:\alpha_{33}-\alpha_{22}}}.
	\end{align*}
Similar expressions hold for $\fix(\beta_i)$. It follows that
	\begin{align*}
	\corchetes{\alpha_1,\beta_1}=\id &\Leftrightarrow \alpha_{11}\beta_{12}+\alpha_{12}\beta_{22}=\alpha_{12}\beta_{11}+\alpha_{22}\beta_{12}\\
	&\Leftrightarrow \corchetes{\alpha_{12}:\alpha_{22}-\alpha_{11}}=\corchetes{\beta_{12}:\beta_{22}-\beta_{11}}\\
	&\Leftrightarrow \fix(\alpha_1)=\fix(\beta_1).
	\end{align*}
If $\alpha_{12}=0$ and $\alpha_{11}\neq\alpha_{22}$, then $\fix(\alpha_1)=\set{e_1,e_2}$. In this case, the previous calculation shows that $\beta_{12}=0$, then $\fix(\beta_1)=\set{e_1,e_2}=\fix(\alpha_1)$.\\

If $\alpha_{11}=\alpha_{22}$ and $\alpha_{12}\neq 0$, then $\fix(\alpha_1)=\set{e_1}$. Again, the previous calculation shows that $\beta_{11}=\beta_{22}$, then $\fix(\beta_1)=\set{e_1}=\fix(\alpha_1)$. All this shows that $\alpha_1$ and $\beta_1$ commute if and only if $\fix(\alpha_1)=\fix(\beta_1)$.\\
	
A similar conclusion holds for $\alpha_2$ and $\beta_2$. From this we can conclude that, if $\alpha$ and $\beta$ commute then $\alpha_{12}$ and $\beta_{12}$ are zero or non-zero simultaneously, the same holds for $\alpha_{23}$ and $\beta_{23}$. This proves the proposition.
\end{proof}

This previous proposition says that a necessary condition for two elements of $U_+$ to commute is that they both have the same \emph{form} given by the sets $F_1,F_2,F_3$. And this is equivalent to say that the two $2\times 2$ sub-blocks (see the proof of the previous proposition) of one element share the same fixed points with the corresponding sub-block of the other element (see the next figure).

\begin{figure}[H]
\begin{center}	
\includegraphics[height=45mm]{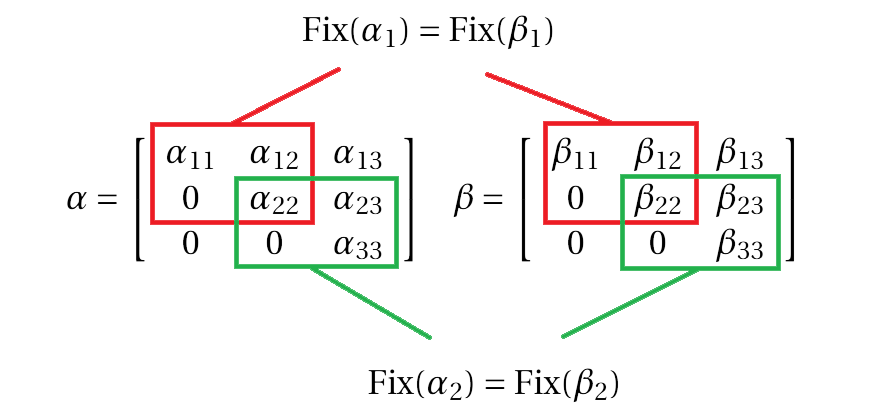}	
\end{center}	
\end{figure}

\begin{prop}\label{prop_transitividad_conmutatividad}
Let $\alpha,\beta\in F_i$, $i=1,2,3$. Then, the relation defined by $\alpha\sim\beta$ if and only if $\corchetes{\alpha,\beta}=\id$, is an equivalence relation on each subset $F_i$, $i=1,2,3$. 
\end{prop}

\begin{proof}
It is direct to verify that the relation is reflexive and symmetric. We now verify that it is also transitive, let $\alpha,\beta,\gamma\in F_i$, for some $i=1,2,3$ such that $\corchetes{\alpha,\beta}=\id$ and $\corchetes{\beta,\gamma}=\id$. Denote $\alpha=\corchetes{\alpha_{ij}}$, $\beta=\corchetes{\beta_{ij}}$ and $\gamma=\corchetes{\gamma_{ij}}$. Then $\corchetes{\alpha,\beta}=\id$ if and only if
	\begin{equation}\label{eq_prop_transitividad_conmutatividad_1}
	\beta_{13}(\alpha_{11}-\alpha_{33})+\alpha_{13}(\beta_{11}-\beta_{33})=\alpha_{23}\beta_{12}-\alpha_{12}\beta_{23}.
	\end{equation}	 
If both $\alpha,\beta\in F_i$, $i=1,2,3$, then $\alpha_{12}=\beta_{12}=0$ or $\alpha_{23}=\beta_{23}=0$, which implies that the right side of (\ref{eq_prop_transitividad_conmutatividad_1}) is zero and then 
	\begin{equation}\label{eq_prop_transitividad_conmutatividad_1}
	\beta_{13}(\alpha_{11}-\alpha_{33})+\alpha_{13}(\beta_{11}-\beta_{33})=0.
	\end{equation}	
From this, analogously to the proof of Proposition \ref{prop_formas_conmutativas}, it follows that $\fix(\alpha)=\fix(\beta)$. Analogously, $\corchetes{\beta,\gamma}=\id$ implies that $\fix(\beta)=\fix(\gamma)$, then $\fix(\alpha)=\fix(\gamma)$ and this implies that  $\corchetes{\alpha,\gamma}=\id$. This proves that the relation is transitive and thus, it is an equivalence relation on $F_i$, $i=1,2,3$.
\end{proof}

Propositions \ref{prop_formas_conmutativas} and \ref{prop_transitividad_conmutatividad} state that, unlike the purely parabolic case (see \cite{ppar}), commutativity does not define an equivalence relation on $U_+$, this equivalence relation only happens on each $F_i$, $i=1,2,3$ separately. On $F_4$, commutativity does not define an equivalence relation.\\

We can now prove the following corollary, which simplifies the decomposition described in Theorem \ref{thm_descomposicion_caso_noconmutativo}.

\begin{cor}\label{cor_descomposicion_bloques_noconm}
Under the same hypothesis and notation of Theorem \ref{thm_descomposicion_caso_noconmutativo}, the group $\Gamma$ can also be written as
	$$\Gamma\cong\Z^{r_0}\rtimes ...\rtimes \Z^{r_m},$$
for integers $r_0,...,r_m\geq 1$ satisfying
	$$r_0+\cdots+r_m\leq 4.$$
\end{cor}

\begin{proof}
Using the notation of Theorem \ref{thm_descomposicion_caso_noconmutativo} and Corollary \ref{cor_lox_nivel34}, we know that the group 
	$$A=\core(\Gamma)\rtimes\prodint{\xi_1}\rtimes...\rtimes \prodint{\xi_r}$$
is purely parabolic and therefore, by Lemma 7.9 of \cite{ppar} we can write
	\begin{equation}\label{eq_dem_cor_descomposicion_bloques_noconm_1}
	A\cong \Z^{k_0}\rtimes ...\rtimes \Z^{k_{n_1}}.
	\end{equation}
For some integers $k_0,...,k_{n_1}$ such that $k_0+...+k_{n_1}\leq 4$. We denote 
	\begin{equation}\label{eq_dem_cor_descomposicion_bloques_noconm_5}
	r_i=k_i,\;\;\;\text{for }i=0,...,n_1.
	\end{equation}

Let us re-order the elements $\set{\eta_1,...,\eta_m}$ in the third layer, in such way that if $i<j$ then $\eta_i\in F_{s_1}$, $\eta_j\in F_{s_2}$ with $s_1\leq s_2$. We re-order the elements $\set{\gamma_1,...,\gamma_n}$ in the same way.

Define the relation in 
	$$\Gamma_1 := \prodint{\eta_1,...,\eta_m}\cap\parentesis{F_1\cup F_2\cup F_3}$$
given by $\alpha\sim \beta$ if and only if $\corchetes{\alpha,\beta}=\id$, this is an equivalence relation (Proposition \ref{prop_transitividad_conmutatividad}). Denote by $A_1,...,A_{n_2}$ the equivalence classes in $\Gamma_1$. Let $B_i=\prodint{A_i}$, clearly $B_i$ is a commutative and torsion free group, denoting $p_i=\text{rank}(B_i)$ and using Propostion \ref{prop_grupos_conm_sintorsion_Zk}, we have
	$$B_i\cong \Z^{p_i}.$$ 
Then 
	\begin{equation}\label{eq_dem_cor_descomposicion_bloques_noconm_2}
	\Gamma_1\cong \Z^{p_1}\rtimes ...\rtimes \Z^{p_{n_2}}.
	\end{equation}
Denote by $\eta_{\tilde{p}_i}$ the remaining elements of the third layer. That is, 
	$$\prodint{\eta_1,...,\eta_m}\cap F_4=\set{\eta_{\tilde{p}_1},...,\eta_{\tilde{p}_{n_3}}}.$$
Then, it follows from (\ref{eq_dem_cor_descomposicion_bloques_noconm_2}),
	\begin{equation}\label{eq_dem_cor_descomposicion_bloques_noconm_3}
	\prodint{\eta_1,...,\eta_m}\cong \Z^{p_1}\rtimes ...\rtimes \Z^{p_{n_2}}\rtimes\prodint{\eta_{\tilde{p}_1}}\rtimes...\rtimes \prodint{\eta_{\tilde{p}_{n_3}}}.
	\end{equation}
Let us denote
	\begin{align}
	r_{n_1+i}=p_i &, \;\;\;\text{for }i=1,...,n_2. \label{eq_dem_cor_descomposicion_bloques_noconm_6} \\
	r_{n_1+n_2+i}= 1 &,\;\;\;\text{for }i=1,...,n_3. \nonumber
	\end{align}
	
Applying the same argument to the elements of the fourth layer $\set{\gamma_1,...,\gamma_n}$ we have  
	\begin{equation}\label{eq_dem_cor_descomposicion_bloques_noconm_4}
	\prodint{\gamma_1,...,\gamma_n}\cong \Z^{q_1}\rtimes ...\rtimes \Z^{q_{n_4}}\rtimes\prodint{\gamma_{\tilde{q}_1}}\rtimes...\rtimes \prodint{\gamma_{\tilde{q}_{n_5}}}.
	\end{equation}
Again, we denote
\begin{align}
	r_{n_1+n_2+n_3+i}=q_i &, \;\;\;\text{for }i=1,...,n_4. \label{eq_dem_cor_descomposicion_bloques_noconm_7} \\
	r_{n_1+n_2+n_3+n_4+i}= 1 &,\;\;\;\text{for }i=1,...,n_5. \nonumber
	\end{align}
Putting together (\ref{eq_dem_cor_descomposicion_bloques_noconm_1}), (\ref{eq_dem_cor_descomposicion_bloques_noconm_3}) and (\ref{eq_dem_cor_descomposicion_bloques_noconm_4}) we prove the corollary. The indices $r_0,...,r_m$ are given by (\ref{eq_dem_cor_descomposicion_bloques_noconm_5}), (\ref{eq_dem_cor_descomposicion_bloques_noconm_6}) and (\ref{eq_dem_cor_descomposicion_bloques_noconm_7}) and $m=n_1+...+n_5$.
\end{proof}

The following corollary describes the type of elements found in each layer of the decomposition of a non-commutative upper triangular discrete group.

\begin{cor}\label{cor_lox_nivel34}
Let $\Gamma\subset U_+$ be a non-commutative discrete subgroup. Consider the decomposition in four layers described in the proof of Theorem \ref{thm_descomposicion_caso_noconmutativo}, and summarized in Table \ref{fig_noconmutativo_capas}.\\
The first two layers $\core(\Gamma)$ and $A\setminus\core(\Gamma)$ are purely parabolic and the last two layers $\kernel(\lambda_{23})\setminus A$ and $\Gamma\setminus\kernel(\lambda_{23})$ are made up entirely of loxodromic elements. Furthermore,
	\begin{enumerate}[(i)]
	\item The third layer $\kernel(\lambda_{23})\setminus A$ contains only loxo-parabolic elements.
	\item The fourth layer $\Gamma\setminus\kernel(\lambda_{23})$ contains only loxo-parabolic and strongly loxodromic elements or complex homotheties with the form $\text{Diag}(\lambda,\lambda^{-2},\lambda)$. 
	\end{enumerate}
\end{cor}

\begin{proof}
From the definition of $\core(\Gamma)$ and $A=\kernel(\lambda_{12})\cap\kernel(\lambda_{23})$ is clear that this two subgroups are purely parabolic.\\

Now we deal with the third layer $\kernel(\lambda_{23})\setminus A$. If there were a elliptic element $\gamma$ in this layer, we have two cases:
	\begin{itemize}
	\item If $\gamma$ has infinite order then $\Gamma$ cannot be discrete. 
	\item If $\gamma$ has finite order $p>0$, but we are assuming that $\Gamma$ is torsion free.
	\end{itemize}

If there were a parabolic element $\gamma$ in this layer, then it must have exactly two repeated eigenvalues (if it had 3, then $\gamma\in A$). Furthermore, all of its eigenvalues must be unitary (they can't be 1, because then $\gamma\in A$). Then,
	$$\gamma=\corchetes{\begin{array}{ccc}
	e^{-4\pi i \theta} & x & y\\
	0 & e^{2\pi i \theta} & z\\
	0 & 0 & e^{2\pi i \theta}\\
	\end{array}}$$
with $z\neq 0$ and $\theta\in\R\setminus\Q$ (otherwise $\lambda_{12}(\gamma)$ would be a torsion element in the torsion free group $\lambda_{12}(\Gamma)$). This means that $\gamma$ is an irrational ellipto-parabolic element, and by Proposition \ref{prop_EPI_conmutativo}, $\Gamma$ would be commutative. All of these arguments prove that the third layer $\kernel(\lambda_{23})\setminus A$ is purely loxodromic. Finally, since $\Gamma$ is not commutative, it cannot contain a complex homothety as a consequence of Corollary \ref{cor_HC_no_hay_en_no_conmutativos}.\\

Now we deal with the fourth layer $\Gamma\setminus\kernel(\lambda_{23})$. Using the same argument as in the third layer, there cannot be elliptic elements. Now assume that there is a parabolic element $\gamma\in\Gamma\setminus\kernel(\lambda_{23})$. In the same way as before, $\gamma$ must have exactly two distinct eigenvalues and neither of them are equal to 1. Since $\gamma\nin\kernel(\lambda_{23})$ then
	$$\gamma=\corchetes{\begin{array}{ccc}
	e^{2\pi i \theta} & x & y\\
	0 & e^{2\pi i \theta} & z\\
	0 & 0 & e^{-4\pi i \theta}\\
	\end{array}}$$
with $x\neq 0$ and $\theta\in\R\setminus\Q$. Then $\gamma$ is an irrational ellipto-parabolic element, and by Proposition \ref{prop_EPI_conmutativo}, $\Gamma$ is commutative. Then the fourth layer $\Gamma\setminus\kernel(\lambda_{23})$ is purely loxodromic. Inspecting the form of these elements, they can be strongly loxodromic, complex homotheties of the form $\text{Diag}(\lambda,\lambda^{-2},\lambda)$ or a complex homothety of the form $\text{Diag}(\lambda,\lambda,\lambda^{-2})$. This latter possibility is dismissed by Proposition \ref{cor_HC3_no_hay_en_no_conmutativos}.
\end{proof}

We can summarize the results and ideas given in proofs of Theorems \ref{thm_case1_kulkarni}, \ref{thm_kulkarni_diagonales}, \ref{thm_descomposicion_caso_noconmutativo}, \ref{thm_descomposicion_caso_noconmutativo2} and Corollaries \ref{cor_descomposicion_bloques_noconm}, \ref{cor_lox_nivel34} in the following theorem.

\begin{thm}\label{thm_main_solvable}
Let $\Gamma\subset\PSL$ be a solvable complex Kleinian group such that its Kulkarni limit set does not contain exactly four lines in general position. Then, there exists a non-empty open region $\Omega_\Gamma\subset \CP^2$ such that
	\begin{enumerate}[(i)]
	\item $\Omega_\Gamma$ is the maximal open set where the action is proper and discontinuous.
	\item $\Omega_\Gamma$ is homeomorphic to one of the following regions: $\C^2$, $\C^2\setminus\set{0}$, $\C\times\parentesis{\Hh^+\cup\Hh^-}$ or $\C\times\C^\ast$.
	\item $\Gamma$ is finitely generated and $\rank(\Gamma)\leq 4$.
	\item The group $\Gamma$ can be written as
		$$\Gamma=\Gamma_p\rtimes \Gamma_{LP}\rtimes\Gamma_{L},$$
	where $\Gamma_p$ is the subgroup of $\Gamma$ consisting of all the parabolic elements of $\Gamma$, $\Gamma_{LP}$ is the subgroup of $\Gamma$ consisting of all the loxo-parabolic elements of $\Gamma$ and $\Gamma_{L}$ is the subgroup of $\Gamma$ consisting of all the strongly loxodromic and  complex homotheties of $\Gamma$.
	\item The group $\Gamma$, up to a finite index subgroup, leaves a full flag invariant. 
	\end{enumerate}
\end{thm}

\subsection{Representations of non-commutative triangular groups}
\label{subsection:representation}

In this section we give an explicit description of the non-commutative triangular groups using Theorem \ref{thm_descomposicion_caso_noconmutativo}. In order to do this we need several results from \cite{ppar}, these results will describe the parabolic part of any complex Kleinian group $\Gamma\subset U_+$. We can rewrite Theorem 0.1 of \cite{ppar} in the following way.

\begin{thm}\label{thm_descripcion_parte_parabolica}
Let $\tilde{\Gamma}\subset\PSL$ be a complex Kleinian group without loxodromic elements, then there exists a finite index subgroup $\Gamma\subset\tilde{\Gamma}$ conjugate to one of the following groups.
	\begin{itemize}
	\item If $\Gamma$ is commutative:	
	\begin{enumerate}[(1)]
	\item 
		$$\Gamma_{W,\mu}=\SET{\corchetes{\begin{array}{ccc}
	\mu(w) & w\mu(w) & 0\\
	0 & \mu(w) & 0\\
	0 & 0 & \mu(w)^{-2}\\
	\end{array}}}{w\in W}$$
	where $W\subset\C$ is a discrete subgroup and $\mu:(\C,+)\rightarrow(\C^\ast,\cdot)$ is a group morphism. 
	\item 
		$$\Gamma_{W}^\ast=\SET{\corchetes{\begin{array}{ccc}
	1 & 0 & x\\
	0 & 1 & y\\
	0 & 0 & 1\\
	\end{array}}}{(x,y)\in W}$$
	where $W\subset\C^2$ is an additive subgroup.
	\item 
		$$\Gamma_{R,L,W}=\SET{\corchetes{\begin{array}{ccc}
	1 & x & L(x)+\frac{x^2}{2}+w\\
	0 & 1 & x\\
	0 & 0 & 1\\
	\end{array}}}{x\in R, w\in W}$$
	where $W\subset\C$ is an additive discrete subgroup, $R\subset\C$ is an additive subgroup and $L:R\rightarrow\C$ is an additive function such that
		\begin{itemize}
		\item If $R$ is discrete, then $\rank(W)+\rank(R)\leq 4$.
		\item If $R$ is not discrete, then $\rank(W)\leq 1$, $\rank(W)+\rank(R)\leq 4$ and
			$$\lim_{n\rightarrow\infty} L(x_n) +w_n=\infty$$
			for any sequence $\set{w_n}\subset W$, and any sequence $\set{x_n}\subset R$ such that $x_n\rightarrow 0$. 		
		\end{itemize}		 
	\item 
		$$\Gamma_{W}^\ast=\SET{\corchetes{\begin{array}{ccc}
	1 & x & y\\
	0 & 1 & 0\\
	0 & 0 & 1\\
	\end{array}}}{(x,y)\in W}$$
	where $W\subset\C^2$ is a discrete additive subgroup with $\rank(W)\leq 2$.
	\end{enumerate}	
	\item If $\Gamma$ is not commutative:
	\begin{enumerate}[(1)]
	\setcounter{enumi}{4}
	\item 
		$$\Gamma_{\mathbf{w}}=\SET{\corchetes{\begin{array}{ccc}
	1 & k+nc+mx & nd+m(k+nc)+{m\choose 2}x+my\\
	0 & 1 & m\\
	0 & 0 & 1\\
	\end{array}}}{k,m,n\in\Z},$$
	where $\mathbf{w}=(x,y,p,q,r)$; $x,y\in\C$, $p,q,r\in\Z$ such that $p,q$ are co-primers, $q^2$ divides $r$ and $c=pq^{-1}$, $d=r^{-1}$.
	\item 
		$$\Gamma_{W,a,b,c}=\SET{\corchetes{\begin{array}{ccc}
	1 & 0 & w\\
	0 & 1 & 0\\
	0 & 0 & 1\\
	\end{array}}\corchetes{\begin{array}{ccc}
	1 & 1 & 0\\
	0 & 1 & 1\\
	0 & 0 & 1\\
	\end{array}}^n\corchetes{\begin{array}{ccc}
	1 & a+c & b\\
	0 & 1 & c\\
	0 & 0 & 1\\
	\end{array}}^m}{\begin{array}{c}
	m,n\in\Z\\
	w\in W\\
	\end{array}}$$
	with $W\subset\C$ an additive discrete group, and $a,b,c\in\C$ are such that $a\in W\setminus\set{0}$ and, either:
		\begin{itemize}
		\item $\set{1,c}$ is a $\R$-linearly independent set, or
		\item $\set{1,c}$ is a $\R$-linearly dependent set, but is a $\Z$-linearly independent set.
		\end{itemize}
	\end{enumerate}
	\end{itemize}		
\end{thm}

Given a torsion free, non-commutative complex Kleinian group $\Gamma\subset U_+$, we denote by $\Gamma_p\subset\Gamma$ the parabolic part of $\Gamma$. This subgroup is conjugated to one of the six options given by Theorem \ref{thm_descripcion_parte_parabolica}. Observe that, according to Theorem \ref{thm_descomposicion_caso_noconmutativo}, $\Gamma_p$ correspond to the semidirect product of elements of the first two layers. Therefore, in order to give the full description of $\Gamma$, we have to check whether we can add one, two or three loxodromic elements of the third or fourth layer, taking into account that the rank of $\Gamma$ has to be at most, four.\\

In Table \ref{table:casos_rangos_abc} we summarize all these possible cases.

\begin{center}
\begin{table}
  \begin{tabular}{ | c | c | c | c || c | c | c | c |}
    \hline
    Case & $\rank(\Gamma_p)$ & $\rank(N_3)$ & $\rank(N_4)$ & Case & $\rank(\Gamma_p)$ & $\rank(N_3)$ & $\rank(N_4)$\\ \hline
    400 & 4 & 0 & 0 & 200 & 2 & 0 & 0 \\ \hline
    310 & 3 & 1 & 0 & 121 & 1 & 2 & 1 \\ \hline
    301 & 3 & 0 & 1 & 120 & 1 & 2 & 0 \\ \hline
	300 & 3 & 0 & 0 & 112 & 1 & 1 & 2 \\ \hline    
	220 & 2 & 2 & 0 & 111 & 1 & 1 & 1 \\ \hline    
	211 & 2 & 1 & 1 & 110 & 1 & 1 & 0 \\ \hline    
	210 & 2 & 1 & 0 & 102 & 1 & 0 & 2 \\ \hline
	202 & 2 & 0 & 2 & 101 & 1 & 0 & 1 \\ \hline    
	201 & 2 & 0 & 1 & 100 & 1 & 0 & 0 \\ \hline
	
  \end{tabular}
\caption{All the possibilities for the rank of each layer of $\Gamma$.}
\label{table:casos_rangos_abc}
\end{table}
\end{center}

$N_3$ and $N_4$ are the groups generated by the elements of the third and fourth layer respectively (see Theorem \ref{thm_descomposicion_caso_noconmutativo} and Table \ref{fig_noconmutativo_capas} for a description of these elements).\\

Observe that the description of cases 400, 300, 200 and 100 is already given in \cite{ppar}, since they are purely parabolic groups.\\

In each one of the cases described in Table \ref{table:casos_rangos_abc}, the parabolic part of $\Gamma$ can be conjugated to one of six possibilities given by Theorem \ref{thm_descripcion_parte_parabolica}. Given $m,n$ not simultaneously $0$, consider the case $kmn$ (that is, $\rank(\Gamma_p)=k$, $\rank(N_3)=m$ and $\rank(N_4)=n$), we will denote by $kmn(j)$ the subcase in which $\Gamma_p$ is conjugated to the group given by the case (j) in Theorem \ref{thm_descripcion_parte_parabolica}. So, the subcases $kmn(1),...,kmn(4)$ cover the case when $\Gamma_p$ is commutative and the rest, when $\Gamma_p$ is not commutative. 

\subsubsection*{$\Gamma_p$ is commutative}

As a consequence of Theorem \ref{thm_descripcion_parte_parabolica} we know that, in each case of Table \ref{table:casos_rangos_abc}, $\Gamma_p$ can be conjugated to a group of the form $\Gamma_{W,\mu}$, $\Gamma_{W}$, $\Gamma_{R,L,W}$ or $\Gamma_{W}^{\ast}$.\\

The following lemmas and propositions will be describing or discarding each subcase separately.\\

The following lemma dismisses all the subcases kmn(1), since they are commutative and therefore, have been already described in Section \ref{sec_commutative_triangular}.

\begin{lem}
Let $\Gamma\subset U_+$ be a complex Kleinian group such that $\Gamma_p$ is conjugated to the group $\Gamma_{W,\mu}$ for some discrete group $W$ and group morphism $\mu$ as in the case (1) of Theorem \ref{thm_descripcion_parte_parabolica}. Then $\Gamma$ is commutative.
\end{lem}

\begin{proof}
We can assume that $\Gamma_p=\Gamma_{W,\mu}$, and since $\Gamma_p$ is a purely parabolic group, then $\valorabs{\mu^{-3}(w)}=1$ for any $w\in W$. Then $\lambda_{23}(\gamma)$ has infinite order for any $\gamma\in\Gamma_p$, it follows from Lemma 5.10 of \cite{ppar} that $\Gamma$ is commutative. 
\end{proof}

The following lemma discards the subcases 110(2), 111(2), 112(2), 120(2), 121(2), 130(2), 210(2), 211(2), 220(2) and 310(2). 

\begin{lem}\label{lem_casos_kmn2_noexisten}
Let $\Gamma\subset U_+$ be a complex Kleinian group such that $\Gamma_p$ is conjugated to the group $\Gamma_W$ for some discrete group $W$ as in the case (2) of Theorem \ref{thm_descripcion_parte_parabolica}. Then $\Gamma$ cannot contain elements in its third layer.
\end{lem}

\begin{proof}
Let us assume that $\Gamma_p=\Gamma_W$ with $W\subset\C^2$ an additive subgroup given by
	$$W=\text{Span}_\Z\parentesis{(x_1,y_1),...,(x_n,y_n)},$$
where $\set{(x_1,y_1),...,(x_n,y_n)}$ is a $\R$-linearly independent set and $n\leq 3$. If $(x,y)\in W$, denote 
	$$h_{x,y}=\corchetes{\begin{array}{ccc}
	1 & 0 & x\\
	0 & 1 & y\\
	0 & 0 & 1\\
	\end{array}}.$$
Let us suppose first that there exists a loxodromic element $\gamma\in N_3$, 
	$$\gamma=\corchetes{\begin{array}{ccc}
	\alpha^{-2} & \gamma_{12} & \gamma_{13}\\
	0 & \alpha & \gamma_{23}\\
	0 & 0 & \alpha\\
	\end{array}}$$
for some $\alpha\in\C^\ast$ such that $\valorabs{\alpha}\neq 1$ and $\gamma_{12},\gamma_{13}\in\C$ and $\gamma_{23}\in\C^\ast$. Conjugating $\Gamma$ by
	$$A=\corchetes{\begin{array}{ccc}
	1 & \frac{\gamma_{12}}{\alpha^{-2}-\alpha} & \gamma_{13}\\
	0 & 1 & 0\\
	0 & 0 & \gamma_{23}\\
	\end{array}},$$
we can assume that $\gamma_{12}=0$ and $\gamma_{23}=1$. Now, assuming without loss of generality that $\valorabs{\alpha}>1$, consider the sequence of distinct elements
	$$\gamma^k h_{x,y}\gamma^{-k}=\corchetes{\begin{array}{ccc}
	1 & 0 & \alpha^{-3k}x\\
	0 & 1 & y\\
	0 & 0 & 1\\
	\end{array}}.$$
It is clear that $\gamma^k h_{x,y}\gamma^{-k}$ converges to an element of $\PSL$, contradicting that $\Gamma$ is discrete.\\
\end{proof}

The following lemma discards the subcases 101(2), 102(2) and 103(2).

\begin{lem}
Let $\Gamma\subset U_+$ be a complex Kleinian group such that $\Gamma_p$ is conjugated to the group $\Gamma_W$ for some discrete group $W$ as in the case (2) of Theorem \ref{thm_descripcion_parte_parabolica}. If $\rank(\Gamma_p)=1$ then $\Gamma$ contains no loxodromic elements in its fourth layer, unless
	$$\Gamma_W=\prodint{\corchetes{\begin{array}{ccc}
	1 & 0 & x\\
	0 & 1 & 0\\
	0 & 0 & 1\\
	\end{array}}},$$
for some $x\in\C^\ast$. In this latter case, the combinations of this group are described in Lemma \ref{lem_kmn4_combinaciones_Gamma5}. 
\end{lem}

\begin{proof}
Assume that $\Gamma_p=\prodint{h}$,
	$$h=\corchetes{\begin{array}{ccc}
	1 & 0 & x\\
	0 & 1 & y\\
	0 & 0 & 1\\
	\end{array}}$$
for some $x,y\in\C$, $\valorabs{x}+\valorabs{y}\neq 0$. Let us suppose that there exists an element $\gamma=\corchetes{\gamma_{ij}}\in N_4$, conjugating by
	$$A=\corchetes{\begin{array}{ccc}
	1 & \frac{\gamma_{12}}{\alpha-\beta} & 0\\
	0 & 1 & 0\\
	0 & 0 & 1\\
	\end{array}}$$
if necessary, we can assume that the loxodromic element $\gamma$ has the form
	$$\gamma=\corchetes{\begin{array}{ccc}
	\alpha & 0 & \gamma_{13}\\
	0 & \beta & \gamma_{23}\\
	0 & 0 & \alpha^{-1}\beta^{-1}\\
	\end{array}},$$
with $\alpha,\beta\in\C^\ast$, $\alpha\beta^2\neq 1$.
Since $\prodint{h}$ is a normal subgroup of $\prodint{h,\gamma}$, we have that $\gamma h\gamma^{-1}=g^n$ for some $n\in\Z$. Comparing the entries 13 and 23 in the previous equation yields
	$$\alpha^2\beta x = nx\text{ and }\alpha\beta^2 y = ny.$$
If $x,y\neq 0$, the equation above implies that $\alpha=\beta$, which in turn implies that $\gamma$ is a complex homothety, contradicting Corollary \ref{cor_HC3_no_hay_en_no_conmutativos}.\\

If $x=0$ and $y\neq 0$, the sequence
	$$\tau_k:=\gamma^k h \gamma^{-k} = \corchetes{\begin{array}{ccc}
	1 & 0 & 0\\
	0 & 1 & (\alpha\beta^2)^k y\\
	0 & 0 & 1\\
	\end{array}}$$
contains a subsequence of distinct elements $\set{\tau_{k_j}}$ such that, either $\tau_{k_j}\rightarrow\id$ or $\tau_{k_j}\rightarrow h$ (depending on whether $\valorabs{\alpha\beta^2}\neq 1$ or $\valorabs{\alpha\beta^2}=1$).\\

If $x\neq 0$ and $y=0$, then $\Gamma_p$ is the group described in Lemma \ref{lem_kmn4_combinaciones_Gamma5} and its behaviour is studied in the aforementioned lemma.	

\end{proof}

The following lemma describes the case 301(2).

\begin{lem}
Let $\Gamma\subset U_+$ be a complex Kleinian group such that $\Gamma_p$ is conjugated to the group $\Gamma_W$ for some discrete group $W$ as in the case (2) of Theorem \ref{thm_descripcion_parte_parabolica}. If $\rank(\Gamma_p)=3$ and $\rank(N_4)=1$ then one generator of $N_4$ is 
	$$\gamma=\corchetes{\begin{array}{ccc}
	\alpha^2\beta & 0 & 0\\
	0 & \alpha\beta^2 & 0\\
	0 & 0 & 1\\
	\end{array}}$$
with $\alpha\neq 1$ and $\alpha\beta^2\nin\R$. Moreover $\Gamma$ is the fundamental group of an Inoue surface.
\end{lem}

\begin{proof}
One can assume that $\Gamma_p=\Gamma_W$. Denote by $h_1,h_2,h_3$ three generators of $\Gamma_W$ given by
	$$h_i=\corchetes{\begin{array}{ccc}
	1 & 0 & x_i\\
	0 & 1 & y_i\\
	0 & 0 & 1\\
	\end{array}}.$$
Let $\gamma\in N_4$ be a loxodromic generator of $N_4$, up to a suitable conjugation, one can assume that
	$$\gamma=\corchetes{\begin{array}{ccc}
	\alpha^2\beta & 0 & 0\\
	0 & \alpha\beta^2 & 0\\
	0 & 0 & 1\\
	\end{array}}$$
for some $\alpha,\beta\in\C^\ast$. Since $\Gamma_p$ is normal in $\Gamma$ we have that $\gamma h_i\gamma^{-1}=h_1^{k_{i,1}}h_2^{k_{i,2}}h_3^{k_{i,3}}$ for some $k_{i,j}\in\Z$, $i=1,2,3$. Let us define $A=\corchetes{k_{i,j}}\subset \text{PSL}(3,\Z)$, then
	$$A\corchetes{\begin{array}{c}
	x_1\\
	x_2\\
	x_3\\
	\end{array}}=\alpha^2\beta\corchetes{\begin{array}{c}
	x_1\\
	x_2\\
	x_3\\
	\end{array}},\;\;\;\;A\corchetes{\begin{array}{c}
	y_1\\
	y_2\\
	y_3\\
	\end{array}}=\alpha\beta^2\corchetes{\begin{array}{c}
	y_1\\
	y_2\\
	y_3\\
	\end{array}}.$$
This means that $\alpha^2\beta$ and $\alpha\beta^2$ are eigenvalues of $A$, with respective eigenvectors $\corchetes{x_1:x_2:x_3}$ and $\corchetes{y_1:y_2:y_3}$. A direct calculation shows that $\alpha\neq 1$ and $\alpha\beta^2\nin\R$. Using (1) of Theorem 4.4 of \cite{cs2014} finishes the proof. 
\end{proof}

The following lemma discards all subcases kmn(3).

\begin{lem}\label{lem_Gamma_RLW_no_tiene_lox}
Let $\Gamma\subset U_+$ be a complex Kleinian group such that $\Gamma_P$ is conjugated to the group $\Gamma_{R,L,W}$ for some $R,L,W$ as in the case (3) of Theorem \ref{thm_descripcion_parte_parabolica}. Then $\Gamma$ cannot contain loxodromic elements.
\end{lem}

\begin{proof}
Assume that $\Gamma_p=\Gamma_{R,L,W}$. Let us suppose that $\Gamma$ contain a loxodromic element $\gamma\in\Gamma$. If $\gamma\in N_3$ then
	$$\gamma=\corchetes{\begin{array}{ccc}
	\alpha^{-2} & \gamma_{12} & \gamma_{13}\\
	0 & \alpha & \gamma_{23}\\
	0 & 0 & \alpha\\
	\end{array}}$$
for some $\alpha\in\C^\ast$ such that $\valorabs{\alpha}\neq 1$ and $\gamma_{12},\gamma_{13}\in\C$ and $\gamma_{23}\in\C^\ast$. Since $\Gamma_p$ is a normal subgroup of $\Gamma$, we have $\gamma g \gamma^{-1}\in\Gamma_p$ for any $g\in\Gamma_p$. If 
	\begin{equation}\label{eq_dem_lem_Gamma_RLW_no_tiene_lox_1}		
	g=\corchetes{\begin{array}{ccc}
	1 & x & L(x)+\frac{x^{2}}{2}+w\\
	0 & 1 & x\\
	0 & 0 & 1\\
	\end{array}},
	\end{equation}
then,
	$$\gamma g\gamma^{-1}=\corchetes{\begin{array}{ccc}
	1 & \alpha^{-3}x & \alpha\parentesis{L(x)+\frac{x^{2}}{2}+w}-\frac{\gamma_{23}x}{1-\alpha^3}\\
	0 & 1 & x\\
	0 & 0 & 1\\
	\end{array}},$$
and, since $\gamma g \gamma^{-1}\in\Gamma_p$, then $\alpha^{-3}x=x$, which is impossible, given that $\valorabs{\alpha}\neq 1$.\\

Now, if $\gamma\in N_4$, then
	$$\gamma=\corchetes{\begin{array}{ccc}
	\alpha & \gamma_{12} & \gamma_{13}\\
	0 & \beta & \gamma_{23}\\
	0 & 0 & \alpha^{-1}\beta^{-1}\\
	\end{array}}$$
with $\alpha,\beta\in\C^\ast$, $\alpha^{-1}\beta^{-1}\neq \beta$ (or equivalently, $\alpha\beta^{2}\neq 1$). Again, since $\Gamma_p$ is a normal subgroup of $\Gamma$, we have $\gamma g \gamma^{-1}\in\Gamma_p$ for any $g\in\Gamma_p$. If $g$ is given by (\ref{eq_dem_lem_Gamma_RLW_no_tiene_lox_1}), then
	$$\gamma g\gamma^{-1}=\corchetes{\begin{array}{ccc}
	1 & \alpha\beta^{-1}x & \alpha^{2}\beta\parentesis{L(x)+\frac{x^{2}}{2}+w}+\alpha\gamma_{23}x(\beta -\alpha)\\
	0 & 1 & \alpha\beta^{2}x\\
	0 & 0 & 1\\
	\end{array}},$$
and, since $\gamma g \gamma^{-1}\in\Gamma_p$, then $\alpha\beta^{-1}=1$ and $\alpha\beta^{2}=1$, but this last condition is a contradiction.\\

This two contradictions imply that $\Gamma$ contains no elements in the third and fourth layer, proving the lemma. 
\end{proof}

Now we deal with \textbf{subcases kmn(4)}, i.e. the cases when $\Gamma_p$ is conjugated to the group $\Gamma^{\ast}_W$. The following proposition delimits the rank of $W$.

\begin{prop}\label{prop_casos_kmn4_rangoW_menor_2}
Let $\Gamma\subset U_+$ be a complex Kleinian group such that $\Gamma_p$ is conjugated to a group $\Gamma^{\ast}_W$ for some discrete additive subgroup $W\subset\C^2$, then $\rank(W)\leq 2$.
\end{prop}

Now, the following proposition describes the simplest forms we can assume for $\Gamma_p$, up to conjugation.

\begin{prop}\label{prop_casos_kmn4_conjugaciones}
Let $\Gamma\subset U_+$ be a complex Kleinian group such that $\Gamma_p$ is conjugated to a group $\Gamma^{\ast}_W$ for some discrete additive subgroup $W\subset\C^2$, given by
	$$\Gamma^{\ast}_W=\SET{\corchetes{\begin{array}{ccc}
	1 & n x_1 + m x_2 & n y_1 + m y_2\\
	0 & 1 & 0\\
	0 & 0 & 1\\
	\end{array}}}{n,m\in\Z}$$
and $W=\text{Span}_\Z\parentesis{(x_1,y_1),(x_2,y_2)}$. Then $\Gamma_p$ is conjugated to one of the following groups:
\begin{enumerate}[(i)]
\item
	$$\Gamma_1=\prodint{\corchetes{\begin{array}{ccc}
	1 & 1 & 0\\
	0 & 1 & 0\\
	0 & 0 & 1\\
	\end{array}},\corchetes{\begin{array}{ccc}
	1 & 0 & 1\\
	0 & 1 & 0\\
	0 & 0 & 1\\
	\end{array}}}.$$
\item
	$$\Gamma_2=\prodint{\corchetes{\begin{array}{ccc}
	1 & 0 & 1\\
	0 & 1 & 0\\
	0 & 0 & 1\\
	\end{array}},\corchetes{\begin{array}{ccc}
	1 & 0 & y\\
	0 & 1 & 0\\
	0 & 0 & 1\\
	\end{array}}},$$
	with $y\nin\R$.
\item
	$$\Gamma_3=\prodint{\corchetes{\begin{array}{ccc}
	1 & 1 & 0\\
	0 & 1 & 0\\
	0 & 0 & 1\\
	\end{array}},\corchetes{\begin{array}{ccc}
	1 & x & 0\\
	0 & 1 & 0\\
	0 & 0 & 1\\
	\end{array}}},$$
	with $x\nin\R$.

\item
	$$\Gamma_4=\prodint{\corchetes{\begin{array}{ccc}
	1 & 1 & 0\\
	0 & 1 & 0\\
	0 & 0 & 1\\
	\end{array}}}.$$
\item 
	$$\Gamma_5=\prodint{\corchetes{\begin{array}{ccc}
	1 & 0 & 1\\
	0 & 1 & 0\\
	0 & 0 & 1\\
	\end{array}}}.$$
\end{enumerate}
\end{prop}

\begin{proof}
Using an appropriate conjugation, we can assume that $\Gamma_p=\Gamma^{\ast}_W$, with $W$ as in the hypothesis of the proposition. If $\rank(W)=2$, denote by
	$$g_1=\corchetes{\begin{array}{ccc}
	1 & x_1 & y_1\\
	0 & 1 & 0\\
	0 & 0 & 1\\
	\end{array}},\;\;\;g_2=\corchetes{\begin{array}{ccc}
	1 & x_2 & y_2\\
	0 & 1 & 0\\
	0 & 0 & 1\\
	\end{array}}$$
the generators of $\Gamma_p$. If $x_1\neq 0$ and $x_1 y_2 - x_2y_1\neq 0$, let us define
	$$A_1=\corchetes{\begin{array}{ccc}
	\frac{1}{\sqrt{x_1}} & 0 & 0\\
	0 & \sqrt{x_1} & \frac{y_1}{\sqrt{x_1}}\\
	0 & \frac{x_2}{\sqrt{x_1}} & \frac{y_2}{\sqrt{x_1}}\\
	\end{array}}.$$
Then, conjugating each generator of $\Gamma_p$ by $A_1$ we get
	$$\Gamma_1=A_1\Gamma_p A_1^{-1}=\prodint{\corchetes{\begin{array}{ccc}
	1 & 1 & 0\\
	0 & 1 & 0\\
	0 & 0 & 1\\
	\end{array}},\corchetes{\begin{array}{ccc}
	1 & 0 & 1\\
	0 & 1 & 0\\
	0 & 0 & 1\\
	\end{array}}}.$$
If $x_1=0$ and $x_2\neq 0$, conjugating $\Gamma_p$ by 
	$$A=\corchetes{\begin{array}{ccc}
	\frac{1}{\sqrt{x_2}} & 0 & 0\\
	0 & \sqrt{x_2} & \frac{y_2}{\sqrt{x_2}}\\
	0 & 0 & \frac{y_1}{\sqrt{x_2}}\\
	\end{array}},$$
we conclude that $\Gamma_p$ is conjugated to the same group $\Gamma_1$.\\
	
If $x_1=x_2=0$, then $y_1,y_2\neq 0$ and conjugating $\Gamma_p$ by $A_2=\text{Diag}(1,1,y_1)$ we conclude that $\Gamma_p$ is conjugated to
	$$\Gamma_2=A_2\Gamma_p A_2^{-1}=\prodint{\corchetes{\begin{array}{ccc}
	1 & 0 & 1\\
	0 & 1 & 0\\
	0 & 0 & 1\\
	\end{array}},\corchetes{\begin{array}{ccc}
	1 & 0 & \frac{y_2}{y_1}\\
	0 & 1 & 0\\
	0 & 0 & 1\\
	\end{array}}}.$$
Furthermore, $y:=y_2 y_1^{-1}\nin\R$; otherwise, $W$ would not be discrete.\\

Finally, if $x_1\neq 0$ and $x_1y_2 - x_2y_1=0$, conjugating $\Gamma_p$ by 
	$$A_3=\corchetes{\begin{array}{ccc}
	\frac{1}{\sqrt{x_1}} & 0 & 0\\
	0 & \sqrt{x_1} & \frac{y_1}{\sqrt{x_1}}\\
	0 & 0 & 1\\
	\end{array}},$$
we conclude that $\Gamma_p$ is conjugated to
	$$\Gamma_3=A_3\Gamma_p A_3^{-1}=\prodint{\corchetes{\begin{array}{ccc}
	1 & 1 & 0\\
	0 & 1 & 0\\
	0 & 0 & 1\\
	\end{array}},\corchetes{\begin{array}{ccc}
	1 & \frac{x_2}{x_1} & 0\\
	0 & 1 & 0\\
	0 & 0 & 1\\
	\end{array}}}.$$
Again, $x:=x_2 x_1^{-1}\nin\R$; otherwise, $W$ would not be discrete.\\

Now, if $\rank(W)=1$, we have 
	$$\Gamma_p=\prodint{\corchetes{\begin{array}{ccc}
	1 & x & y\\
	0 & 1 & 0\\
	0 & 0 & 1\\
	\end{array}}},$$
for some $x,y\in\C$ with $\valorabs{x}+\valorabs{y}\neq 0$. If $x\neq 0$, conjugating $\Gamma_p$ by 
	$$A_4=\corchetes{\begin{array}{ccc}
	\frac{1}{\sqrt{x}} & 0 & 0\\
	0 & \sqrt{x} & \frac{y}{\sqrt{x}}\\
	0 & 0 & 1\\
	\end{array}},$$
we conclude that $\Gamma_p$ is conjugated to
	$$\Gamma_4=A_4\Gamma_p A_4^{-1}=\prodint{\corchetes{\begin{array}{ccc}
	1 & 1 & 0\\
	0 & 1 & 0\\
	0 & 0 & 1\\
	\end{array}}}.$$
If $x=0$, then $y\neq 0$ and conjugating $\Gamma_p$ by $A_5=\text{Diag}(1,1,y)$ we conclude that $\Gamma_p$ is conjugated to
	$$\Gamma_5=A_5\Gamma_p A_5^{-1}=\prodint{\corchetes{\begin{array}{ccc}
	1 & 0 & 1\\
	0 & 1 & 0\\
	0 & 0 & 1\\
	\end{array}}}.$$
This concludes the proof.
\end{proof}

The next lemmas study the different combinations of $\Gamma_p$ and elements of the third and fourth layer of $\Gamma$, for the subcases kmn(4). For the sake of clarity, we summarize this combinations and conclusions in observation \ref{obs_resumen_casos_kmn4}. In all of these lemmas, $\Gamma$ is a complex Kleinian subgroup of $U_+$.\\

\begin{lem}\label{lem_kmn4_combinaciones_Gamma1}
If $\Gamma$ is such that $\Gamma_p$ is conjugated to 
	\begin{equation}\label{eq_lem_kmn4_combinaciones_Gamma1_1}
	\Gamma_1=\prodint{\corchetes{\begin{array}{ccc}
	1 & 1 & 0\\
	0 & 1 & 0\\
	0 & 0 & 1\\
	\end{array}},\corchetes{\begin{array}{ccc}
	1 & 0 & 1\\
	0 & 1 & 0\\
	0 & 0 & 1\\
	\end{array}}},
	\end{equation}
Then:
	\begin{enumerate}[(i)]
	\item If $\rank(N_3)=1$ then one generator of $N_3$ is given by
		\begin{equation}\label{eq_lem_kmn4_combinaciones_Gamma1_2}
		\gamma=\corchetes{\begin{array}{ccc}
	p & \gamma_{12} & \gamma_{13}\\
	0 & 1 & \frac{q}{p}\\
	0 & 0 & 1\\
	\end{array}},
		\end{equation}
	for $p\in\Z\setminus\set{-1,0,1}$, $q\in\Z\setminus\set{0}$, $\gamma_{12},\gamma_{13}\in\C$.
	\item If $\rank(N_3)=2$ then one generator of $N_3$ is given by (\ref{eq_lem_kmn4_combinaciones_Gamma1_2}) and the other is given by
		\begin{equation}\label{eq_lem_kmn4_combinaciones_Gamma1_4}\mu=\corchetes{\begin{array}{ccc}
	p_2 & \mu_{12} & \mu_{13}\\
	0 & 1 & \frac{q_2}{p_2}\\
	0 & 0 & 1\\
	\end{array}},
	\end{equation}
	with $p_2\in\Z\setminus\set{-1,0,1}$, $q_2\in\Z\setminus\set{0}$, $\mu_{12},\mu_{13}\in\Q\setminus\set{0}$ and
		\begin{align*}
		\mu_{12} &= -\frac{p_1 p_2(m+jp_1p_2)}{(1-p_1)(p_1q_2+p_2q_1)}\\
		\mu_{13} &= \frac{p_2\parentesis{mq_1+jp_1^2(p_2q_1-(1-p_1)q_2)}}{(1-p_1)^2(p_1q_2+p_2q_1)}
		\end{align*}
		where $j,m\in\Z$ are integers such that $p_1 p_2+p_1q_2 \mid m+jp_1p_2$ and $p_1 p_2+p_1q_2 \mid mp_1q_2 -jp_1p_2^2q_1$.
		
	\item If $\rank(N_3)=0$ and $\rank(N_4)=1$, then one generator of $N_4$ is given by
		\begin{equation}\label{eq_lem_kmn4_combinaciones_Gamma1_3}
		\gamma=\corchetes{\begin{array}{ccc}
	pq & \gamma_{12} & \gamma_{13}\\
	0 & q & r\\
	0 & 0 & p\\
	\end{array}},
		\end{equation}
	for $p,q\in\Z\setminus\set{0}$ such that $\gamma$ is loxodromic, $r\in\Z$, $\gamma_{12},\gamma_{13}\in\C$.			
	\item It's not possible to have $\rank(N_4)=2$.
	\item If $\rank(N_3)=1$, then there are no elements in the fourth layer $N_4$.
	\end{enumerate}
\end{lem}

\begin{proof}
We can assume that $\Gamma_p=\Gamma_1$. Let $g\in\Gamma_p$ be an element with the form
	$$g=\corchetes{\begin{array}{ccc}
	1 & n & m\\
	0 & 1 & 0\\
	0 & 0 & 1\\
	\end{array}}$$
We prove each conclusion separately:
	\begin{enumerate}[(i)]
	\item Let $\gamma\in N_3$ be a generator of the third layer of $\Gamma$. Then
	$$\gamma=\corchetes{\begin{array}{ccc}
	\alpha^{-2} & \gamma_{12} & \gamma_{13}\\
	0 & \alpha & \gamma_{23}\\
	0 & 0 & \alpha\\
	\end{array}}\in N_3,$$
	for some $\gamma_{23},\alpha\in\C^\ast$ with $\valorabs{\alpha}\neq 1$. Since $\Gamma_p$ is a normal subgroup of $\Gamma$, then $\gamma g \gamma^{-1}\in\Gamma_p$. 
		$$\gamma g \gamma^{-1}=\corchetes{\begin{array}{ccc}
	1 & \alpha^{-3}n & \alpha^{-4}(m\alpha-n\gamma_{23})\\
	0 & 1 & 0\\
	0 & 0 & 1\\
	\end{array}}\in\Gamma_p,$$
	which means that $\alpha^{-3}n,\alpha^{-4}(m\alpha-n\gamma_{23})\in\Z$ for any $n,m\in\Z$. In particular, for $n=0,m=1$ we get $\alpha^{-3}=p$ for some $p\in\Z\setminus\set{0}$; and for $n=1,m=0$ we get $\gamma_{23}=q p^{-\frac{4}{3}}$ for some $q\in\Z\setminus\set{0}$. Then
		$$\gamma=\corchetes{\begin{array}{ccc}
	p^{\frac{2}{3}} & \gamma_{12} & \gamma_{13}\\
	0 & p^{-\frac{1}{3}} & qp^{-\frac{4}{3}}\\
	0 & 0 & p^{-\frac{1}{3}}\\
	\end{array}}=\corchetes{\begin{array}{ccc}
	p & p^{\frac{1}{3}}\gamma_{12} & p^{\frac{1}{3}}\gamma_{13}\\
	0 & 1 & \frac{q}{p}\\
	0 & 0 & 1\\
	\end{array}},$$
and since $\gamma_{12}$ and $\gamma_{13}$ are arbitrary, we get (\ref{eq_lem_kmn4_combinaciones_Gamma1_2}).
	\item Conjugating $\Gamma$ by 
		$$A=\corchetes{\begin{array}{ccc}
	1 & -\frac{\gamma_{12}}{1-p_1} & \frac{q_1\gamma_{12}-p_1\gamma_{13}(1-p_1)}{p_1(1-p_1)^2}\\
	0 & 1 & \frac{q}{p}\\
	0 & 0 & 1\\
	\end{array}}$$
	if necessary, we can assume that 
		$$\gamma=\corchetes{\begin{array}{ccc}
	p_1 & 0 & 0\\
	0 & 1 & \frac{q_1}{p_1}\\
	0 & 0 & 1\\
	\end{array}}.$$
	Using the proof of (i), the second generator of $N_3$ has the form (\ref{eq_lem_kmn4_combinaciones_Gamma1_4}). Observe that if $\corchetes{\mu,\gamma}=\id$ then $\mu_{12}=\mu_{13}=0$ and then $\Gamma$ would be commutative. Then, we can assume that $\corchetes{\mu,\gamma}\in\Gamma_p\setminus\set{\id}$, comparing the entry 13 of the last expression yields
		\begin{equation}\label{eq_dem_lem_kmn4_combinaciones_Gamma1_4}
		\frac{q\mu_{12}+(1-p_1)p_1\mu_{13}}{p_1^2p_2}=j
		\end{equation}
	for some $j\in\Z\setminus\set{0}$. Since $\Gamma_p\rtimes\prodint{\gamma}$ is normal in $\Gamma$ we have
		$$\mu\gamma\mu^{-1}\in \Gamma_p\rtimes\prodint{\gamma}=\SET{\corchetes{\begin{array}{ccc}
	p_1^k & n & m+\frac{knq_1}{p_1}\\
	0 & 1 & k\frac{q_1}{p_1}\\
	0 & 0 & 1\\
	\end{array}}}{m,n\in\Z}.$$		
		Comparing the entries 12 and 13 in the last expression yields		
		\begin{align}
		\mu_{12}(1-p_1)&=n\label{eq_dem_lem_kmn4_combinaciones_Gamma1_2}\\
		\mu_{12}\parentesis{\frac{q_1}{p_1}-\frac{q_2}{p_2}}+\mu_{12}\frac{p_1q_2}{p_2}+\mu_{13}(1-p_1)&=m+n\frac{q_1}{p_1}\label{eq_dem_lem_kmn4_combinaciones_Gamma1_3}
		\end{align}
	for some $m,n\in\Z$. Combining both (\ref{eq_dem_lem_kmn4_combinaciones_Gamma1_2}) and (\ref{eq_dem_lem_kmn4_combinaciones_Gamma1_3}) we get
		\begin{equation}\label{eq_dem_lem_kmn4_combinaciones_Gamma1_5}
		mp_2+\parentesis{(1-p_1)q_2-p_2q_1}\mu_{12}+p_2\mu_{13}(1-p_1)=0.
		\end{equation}		 
		Solving (\ref{eq_dem_lem_kmn4_combinaciones_Gamma1_4}) and (\ref{eq_dem_lem_kmn4_combinaciones_Gamma1_5}) we get the desired expressions for $\mu_{12}$ and $\mu_{13}$. 
	\item If $\gamma\in N_4$ then 
		\begin{equation}\label{eq_dem_lem_kmn4_combinaciones_Gamma1_1}
		\gamma=\corchetes{\begin{array}{ccc}
	\alpha & \gamma_{12} & \gamma_{13}\\
	0 & \beta & \gamma_{23}\\
	0 & 0 & \alpha^{-1}\beta^{-1}\\
	\end{array}}
	\end{equation}				
	for some $\alpha,\beta\in\C^{\ast}$ such that $\gamma$ is loxodromic and $\alpha\neq\beta^2$, $\gamma_{12},\gamma_{13},\gamma_{23}\in\C$. Since $\Gamma_p$ is normal in $\Gamma$, it holds
		$$\gamma g \gamma^{-1}=\corchetes{\begin{array}{ccc}
	1 & \frac{\alpha}{\beta}n & \alpha^{2}(m\beta-n\gamma_{23})\\
	0 & 1 & 0\\
	0 & 0 & 1\\
	\end{array}}\in\Gamma_p.$$
	This means that $\alpha\beta^{-1},\alpha^{2}(m\beta-n\gamma_{23})\in\Z$ for any $m,n\in\Z$. If $n=0,m=1$ then $\alpha^2\beta=q$ for some $q\in\Z$. If $n=1,m=0$ then $\alpha\beta^{-1}=p\in\Z$ and $\alpha^2\gamma_{23}=r\in\Z$. All this together yields (\ref{eq_lem_kmn4_combinaciones_Gamma1_3}).
	\item Let us assume that $\rank(N_4)=2$ and let $\gamma$, $\mu$ be two generators of $N_4$. Using an adequate conjugation we can assume that 
		$$\gamma=\corchetes{\begin{array}{ccc}
	q_1p_1 & 0 & \gamma_{13}\\
	0 & q_1 & 0\\
	0 & 0 & p_1\\
	\end{array}},\;\;\;\;\mu=\corchetes{\begin{array}{ccc}
	q_2p_2 & \mu_{12} & \mu_{13}\\
	0 & q_2 & \mu_{23}\\
	0 & 0 & p_2\\
	\end{array}}.$$
	Then, using normality we have $\mu\gamma\mu^{-1}\in\Gamma_p\rtimes\prodint{\gamma}$. Comparing the entries 12, we have
		$$\frac{q_1\mu_{12}(1-p_1)}{q_2}=-q_1\mu_{12}(1-p_1)$$
	then $q_2=-1$. Then $\lambda_{13}(\mu)$ is a torsion element in $\lambda_{13}$, this contradiction proves that it is impossible to have $\rank(N_4)=2$.
	\item Let $\gamma$ be a generator of $N_3$, using an adequate conjugation we can assume that 
		$$\gamma =\corchetes{\begin{array}{ccc}
	p & 0 & \gamma_{13}\\
	0 & 1 & \frac{q}{p}\\
	0 & 0 & 1\\		
	\end{array}},$$
with $p,q\in\Z\setminus\set{0}$, and $\valorabs{p}\neq 1$. Let us assume that there is an element $\mu\in N_4$, then
		$$\mu =\corchetes{\begin{array}{ccc}
	\alpha & \mu_{12} & \mu_{13}\\
	0 & \beta & \mu_{23}\\
	0 & 0 & \alpha^{-1}\beta^{-1}\\		
	\end{array}},$$
with $\alpha\beta^{2}\neq 1$. Since $\corchetes{\gamma,\mu}\in\Gamma_p$, then comparing the entries 23 yields $q(1-\alpha\beta^2)p^{-1}\alpha^{-1}\beta^{-2}=0$. This means that, either $q=0$ or $\alpha\beta^2=1$. Both conclusions contradict the hypotheses. Then, there are no elements in $N_4$.

	\end{enumerate}
\end{proof}

\begin{lem}\label{lem_kmn4_combinaciones_Gamma2}
If $\Gamma$ is such that $\Gamma_p$ is conjugated to 
	\begin{equation}\label{eq_lem_kmn4_combinaciones_Gamma2_1}
	\Gamma_2=\prodint{\corchetes{\begin{array}{ccc}
	1 & 0 & 1\\
	0 & 1 & 0\\
	0 & 0 & 1\\
	\end{array}},\corchetes{\begin{array}{ccc}
	1 & 0 & y\\
	0 & 1 & 0\\
	0 & 0 & 1\\
	\end{array}}},
	\end{equation}
with $y\nin\R$. Then:
	\begin{enumerate}[(i)]
	\item If $\rank(N_3)=1$ then one generator of $N_3$ is given by
		\begin{equation}\label{eq_lem_kmn4_combinaciones_Gamma2_2}
		\gamma=\corchetes{\begin{array}{ccc}
	\parentesis{p+q\text{Re}(y)}+iq\text{Im}(y) & \gamma_{12} & \gamma_{13}\\
	0 & 1 & \gamma_{23}\\
	0 & 0 & 1\\
	\end{array}},
		\end{equation}
	for $p,q\in\Z$ such that $\valorabs{p}+\valorabs{q}\neq 0$ and $\valorabs{\parentesis{p+q\text{Re}(y)}+iq\text{Im}(y)}\neq 1$, $\gamma_{12},\gamma_{23}\in\C$ and $\gamma_{23}\in\C^\ast$.
	\item If $\rank(N_3)=0$ and $\rank(N_4)=1$ then one generator of $N_4$ is given by
		\begin{equation}\label{eq_lem_kmn4_combinaciones_Gamma2_3}
		\gamma=\corchetes{\begin{array}{ccc}
	\alpha & \gamma_{12} & \gamma_{13}\\
	0 & \alpha^{-2}z_{p,q} & \gamma_{23}\\
	0 & 0 & \alpha z_{p,q}^{-1}\\
	\end{array}},
		\end{equation}
	for $\beta\in\C^{\ast}$, $z_{p,q}=\parentesis{p+q\text{Re}(y)}+iq\text{Im}(y)$, $\valorabs{p}+\valorabs{q}\neq 0$, $\gamma_{12},\gamma_{23},\gamma_{23}\in\C$. Finally, $\gamma$ must be loxodromic. 
	\end{enumerate}
\end{lem}

\begin{proof}
We can assume that $\Gamma_p=\Gamma_2$. Let $g\in\Gamma_p$ be an element with the form
	$$g=\corchetes{\begin{array}{ccc}
	1 & 0 & n+my\\
	0 & 1 & 0\\
	0 & 0 & 1\\
	\end{array}}$$
We prove each conclusion separately:
	\begin{enumerate}[(i)]
	\item Let $\gamma\in N_3$ be a generator of the third layer of $\Gamma$. Then
	$$\gamma=\corchetes{\begin{array}{ccc}
	\alpha^{-2} & \gamma_{12} & \gamma_{13}\\
	0 & \alpha & \gamma_{23}\\
	0 & 0 & \alpha\\
	\end{array}}\in N_3,$$
	for some $\gamma_{23},\alpha\in\C^\ast$ with $\valorabs{\alpha}\neq 1$. Since $\Gamma_p$ is a normal subgroup of $\Gamma$, then $\gamma g \gamma^{-1}\in\Gamma_p$. 
		$$\gamma g \gamma^{-1}=\corchetes{\begin{array}{ccc}
	1 & 0 & \alpha^{-3}(n+my)\\
	0 & 1 & 0\\
	0 & 0 & 1\\
	\end{array}}\in\Gamma_p$$	
	This means that, for any $n,m\in\Z$,
		\begin{equation}\label{eq_dem_lem_kmn4_combinaciones_Gamma2_1}
		\alpha^{-3}(n+my)=p+qy
		\end{equation}		 
		for some $p,q\in\Z$. Denote $y=a+bi$ and $\alpha^{-3}=r(\cos\theta+i\sin\theta)$. In particular, (\ref{eq_dem_lem_kmn4_combinaciones_Gamma2_1}) holds for $n=1$, $m=0$. Then, (\ref{eq_dem_lem_kmn4_combinaciones_Gamma2_1}) becomes 
		$$r\cos\theta=p+qa,\;\;\;\; r\sin\theta= q b.$$
		Which means that $\text{Re}\parentesis{\alpha^{-3}}=p+qa$ and $\text{Im}\parentesis{\alpha^{-3}}=qb$ and therefore, 
		$$\alpha^{-3}=\parentesis{p+q\text{Re}(y)}+iq\text{Im}(y)$$
		for $p,q\in\Z$ such that $\valorabs{p}+\valorabs{q}\neq 0$ (otherwise, $\alpha=0$) and $\valorabs{\parentesis{p+q\text{Re}(y)}+iq\text{Im}(y)}\neq 1$ (otherwise, $\gamma$ would not be loxodromic). Since we can re-write $\gamma$ as
			$$\gamma=\corchetes{\begin{array}{ccc}
	\alpha^{-3} & \alpha^{-1}\gamma_{12} & \alpha^{-1}\gamma_{13}\\
	0 & 1 & \alpha^{-1}\gamma_{23}\\
	0 & 0 & 1\\
	\end{array}},$$
	we get (\ref{eq_lem_kmn4_combinaciones_Gamma2_2}).
	
	\item If $\gamma\in N_4$, then $\gamma$ has the form given by (\ref{eq_dem_lem_kmn4_combinaciones_Gamma1_1}). Since $\Gamma_p$ is normal in $\Gamma$, $\gamma g\gamma^{-1}\in \Gamma$, and this implies that
		$$\alpha^2\beta(n+my)=\tilde{n}+\tilde{m}y,$$
		for some $\tilde{n},\tilde{m}\in\Z$. The same argument and calculations used in the proof of (i) above show that $\alpha^2\beta=\parentesis{p+q\text{Re}(y)}+iq\text{Im}(y)$. Substituting this in the form of $\gamma$ yields (\ref{eq_lem_kmn4_combinaciones_Gamma2_3}).		  
	\end{enumerate}
\end{proof}

\begin{lem}\label{lem_kmn4_combinaciones_Gamma3}
If $\Gamma$ is such that $\Gamma_p$ is conjugated to 
	\begin{equation}\label{eq_lem_kmn34_combinaciones_Gamma3_1}
	\Gamma_3=\prodint{\corchetes{\begin{array}{ccc}
	1 & 1 & 0\\
	0 & 1 & 0\\
	0 & 0 & 1\\
	\end{array}},\corchetes{\begin{array}{ccc}
	1 & x & 0\\
	0 & 1 & 0\\
	0 & 0 & 1\\
	\end{array}}},
	\end{equation}
with $x\nin\R$. Then:
	\begin{enumerate}[(i)]
	\item $\Gamma$ cannot contain elements in its third layer.
	\item If $\rank(N_4)=1$ then one generator of $N_4$ is given by
		\begin{equation}\label{eq_lem_kmn4_combinaciones_Gamma3_2}
		\gamma=\corchetes{\begin{array}{ccc}
	p+qx & \gamma_{12} & \gamma_{13}\\
	0 & 1 & 0\\
	0 & 0 & (p+qx)^{-1}\beta^{-3}\\
	\end{array}},
		\end{equation}
	where $p,q\in\Z$ such that $\valorabs{p}+\valorabs{q}\neq 0$, $\beta\in\C^\ast$ such that $\gamma$ is loxodromic and $\gamma_{12}, \gamma_{13}\in\C$. Also
		\begin{equation}\label{eq_lem_kmn4_combinaciones_Gamma3_3}
		p-q\valorabs{x}^2 \in\Z,\;\;\;\;p+q+2q\text{Re}(x)\in\Z.
		\end{equation}
	\item If $\rank(N_4)=2$ then one generator is given by (\ref{eq_lem_kmn4_combinaciones_Gamma3_2}) and the other is given by
		\begin{equation}\label{eq_lem_kmn4_combinaciones_Gamma3_4}
		\mu=\corchetes{\begin{array}{ccc}
	p_2+q_2x & \mu_{12} & 0\\
	0 & 1 & 0\\
	0 & 0 & (p_2+q_2x)^{-1}\beta_2^{-3}\\
	\end{array}},
		\end{equation}
	where $p_2,q_2\in\Z$ such that $\valorabs{p_2}+\valorabs{q_2}\neq 0$ and both integers satisfy (\ref{eq_lem_kmn4_combinaciones_Gamma3_3}), $\beta_2\in\C^\ast$ is such that $\gamma$ is loxodromic and $p_2-q\valorabs{x}^2\in\Z$. Furthermore, $\mu_{12}\in\C$ such that
		\begin{equation}\label{eq_lem_kmn4_combinaciones_Gamma3_5}
		\mu_{12}=\frac{n+mx-(x+1)(p_2+q_2 x)}{1-(p_1+q_1 x)}
		\end{equation} 
		for some $n,m\in\Z$ such that
		\begin{equation}\label{eq_lem_kmn4_combinaciones_Gamma3_6}
		\mu_{12}\frac{p_1+q_1 x-1}{(p_1+q_1 x)(p_2+q_2 x)}=k_1+k_2 x,
		\end{equation}
		for some $k_1,k_2\in\Z$.
	\end{enumerate}
\end{lem}

\begin{proof}
We can assume that $\Gamma_p=\Gamma_2$. Let $g\in\Gamma_p$ be an element with the form
	$$g=g_{n,m}=\corchetes{\begin{array}{ccc}
	1 & n+mx & 0\\
	0 & 1 & 0\\
	0 & 0 & 1\\
	\end{array}}$$
We prove each conclusion separately:
	\begin{enumerate}[(i)]
	\item The argument is the same used in the first conclusion of \ref{lem_kmn4_combinaciones_Gamma4}.
	\item If $\gamma\in N_4$, then $\gamma$ has the form given by (\ref{eq_dem_lem_kmn4_combinaciones_Gamma1_1}). Since $\Gamma_p$ is normal in $\Gamma$, $\gamma g\gamma^{-1}\in \Gamma_p$, and this implies that
		$$\frac{\alpha}{\beta}(n+mx)=\tilde{n}+\tilde{m}x,\;\;\;\text{and}\;\;\; -\alpha^2\gamma_{23}(n+mx)=0,$$
	for some $\tilde{n},\tilde{m}\in\Z$. Then $\gamma_{23}=0$, and using the same argument and calculations as in the first part of the proof of lemma \ref{lem_kmn4_combinaciones_Gamma2} we have that
		$$\frac{\alpha}{\beta}=\parentesis{p+q\text{Re}(x)}+iq\text{Im}(x).$$
	Substituting this in the form of $\gamma$ yields (\ref{eq_lem_kmn4_combinaciones_Gamma3_2}). However, these conditions alone do not guarantee that $\gamma g\gamma^{-1}\in \Gamma_p$ for $\gamma$ as in (\ref{eq_lem_kmn4_combinaciones_Gamma3_2}), since 
		$$\gamma g\gamma^{-1}=\corchetes{\begin{array}{ccc}
	1 & p + (p+q)x + qx^2 & 0\\
	0 & 1 & 0\\
	0 & 0 & 1\\
	\end{array}},$$
	and it does not necessarily hold 
		\begin{equation}\label{eq_dem_lem_kmn4_combinaciones_Gamma3_1}
		p + (p+q)x + qx^2=p'+q'x,
		\end{equation}				
		 for some $p',q'\in\Z$. In order for (\ref{eq_dem_lem_kmn4_combinaciones_Gamma3_1}) to hold we solve (\ref{eq_dem_lem_kmn4_combinaciones_Gamma3_1}), by separating $x=\text{Re}(x)+i\text{Im}(x)$ and assuming $\text{Im}(x)\neq 0$. From the equation corresponding to the imaginary part we get $(p+q)+2q\text{Re}(x)\in\Z$, substituting this into the equation corresponding to the real part we get $p-q\valorabs{x}^2\in\Z$. This verifies (\ref{eq_lem_kmn4_combinaciones_Gamma3_3}).
		 
		 \item Let $\gamma$ be a generator of $N_4$ given by (ii) and let $\mu$ be the second generator of $\mu$, then $\mu$ has the same form given by (ii), that is
		 	$$\mu=\corchetes{\begin{array}{ccc}
			p_2+q_2x & \mu_{12} & 0\\
			0 & 1 & 0\\
			0 & 0 & (p_2+q_2x)^{-1}\beta_2^{-3}\\
			\end{array}}.$$
		 We can conjugate $\Gamma$ by a suitable $A\in\PSL$ such that
		 	$$A\gamma A^{-1}=\corchetes{\begin{array}{ccc}
				p_1+q_1 x & 0 & 0\\
				0 & 1 & 0\\
				0 & 0 & \alpha^{-3}(p+q x)^{-1}\\
			\end{array}}.$$
			We have that, either $\corchetes{\mu,\gamma}=0$ or $\corchetes{\mu,\gamma}=g_{k_1,k_2}$ for some $k_1,k_2\in\Z$. Both possibilities imply that $\mu_{13}=0$.\\
			
			If $\corchetes{\mu,\gamma}=0$ then $\mu_{12}=0$ and then $\mu$ is diagonal. Since $\gamma$ is diagonal and the parabolic part is commutative, then $\Gamma$ would be commutative. Then, $\corchetes{\mu,\gamma}=g_{k_1,k_2}$ for some $k_1,k_2\in\Z$, this yields (\ref{eq_lem_kmn4_combinaciones_Gamma3_6}).\\
			
			On the other hand, since $\Gamma_p\rtimes\prodint{\gamma}$ is normal in $\Gamma$, we have that $\mu g_{1,1} \gamma \mu^{-1}=g_{n,m} \gamma^k$ for some $n,m,k\in\Z$. A direct calculation shows that $k=1$, the remaining expression yields (\ref{eq_lem_kmn4_combinaciones_Gamma3_5}).
			 	 
	\end{enumerate}
\end{proof}

\begin{lem}\label{lem_kmn4_combinaciones_Gamma4}
If $\Gamma$ is such that $\Gamma_p$ is conjugated to 
	\begin{equation}\label{eq_lem_kmn4_combinaciones_Gamma4_1}
	\Gamma_4=\prodint{\corchetes{\begin{array}{ccc}
	1 & 1 & 0\\
	0 & 1 & 0\\
	0 & 0 & 1\\
	\end{array}}},
	\end{equation}
then:
	\begin{enumerate}[(i)]
	\item $\Gamma$ cannot contain elements of the third layer.
	\item If $\rank(N_4)=1$ then one generator of $N_4$ is given by
		\begin{equation}\label{eq_lem_kmn4_combinaciones_Gamma4_2}
		\gamma=\corchetes{\begin{array}{ccc}
	p\alpha & \gamma_{12} & \gamma_{13}\\
	0 & \alpha & 0\\
	0 & 0 & p^{-1}\alpha^{-2}\\
	\end{array}},
		\end{equation}
	where $p\in\Z\setminus\set{0,1}$, $\alpha\in\C^{\ast}$ such that $\gamma$ is loxodromic, $\valorabs{\beta}\neq 1$ and $\gamma_{12},\gamma_{13}\in\C$.
	\item If $\rank(N_4)=2$, one generator of $N_4$ is given by (\ref{eq_lem_kmn4_combinaciones_Gamma4_2}) and the other satisfies one of the following conditions:
		\begin{enumerate}
		\item If $p^2\alpha^3=1$ then 
			$$\mu=\corchetes{\begin{array}{ccc}
	q\beta & \frac{\beta jpq}{1-p} & \mu_{13}\\
	0 & \beta & 0\\
	0 & 0 &q^{-1}\beta^{-2}\\
	\end{array}},$$
		where $\mu_{13}\in\C$.
		\item If $p^2\alpha^3\neq 1$ then 
			$$\mu=\corchetes{\begin{array}{ccc}
	q\beta & \frac{\beta jpq}{1-p} & 0\\
	0 & \beta & 0\\
	0 & 0 &q^{-1}\beta^{-2}\\
	\end{array}}.$$
		\end{enumerate}
		In both cases, $\delta\in\C^\ast$, $j,q\in\Z$ and $q\neq 0$ such that $\mu$ is loxodromic.
	\item If $\rank(N_4)=3$, two generators of $N_4$ are given by (ii) and (iii) respectively and the third generator $\eta$ is
		\begin{equation}\label{eq_lem_kmn4_combinaciones_Gamma4_3}
		\eta=\corchetes{\begin{array}{ccc}
	r\delta & \frac{\delta k pr}{1-p} & \eta_{13}\\
	0 & \delta & 0\\
	0 & 0 &q^{-1}\delta^{-2}\\
	\end{array}},
		\end{equation}
	where $\delta\in\C^\ast$, $r\in\Z\setminus\set{-1,0,1}$, $r\delta^3\neq 1$ and
		\begin{align*}
		\eta_{13} &= \frac{p\alpha^2\parentesis{\gamma_{13}(1-r^2\delta^3)-pq\alpha\beta^2\mu_{13}(1+r^2\delta^3)}}{r\delta^2(1-p^2q^2\alpha^3\beta^3)} \\
		k &= \frac{(1-p)\parentesis{\alpha(n+r)-\gamma_{12}(1-r)}-j_2p^2q\alpha(1-r)}{pr\alpha(1-pq)} \\
		\end{align*}
	for $n\in\Z$. Also, either
		\begin{align*}
		k r(1-q) &= jq(1-r) \text{, or}\\
		qr(1-p) & \left.\right|\; p\parentesis{kr(1-q)-jq(1-r)},
		\end{align*}			 
	depending on whether $\mu$ and $\eta$ commute.
	\end{enumerate}		 
\end{lem}

\begin{proof}
We can assume that $\Gamma_p=\Gamma_4$. We denote by $g$ the generator of $\Gamma_4$ given in (\ref{eq_lem_kmn4_combinaciones_Gamma4_1}). We prove each conclusion separately:
	\begin{enumerate}[(i)]
	\item Assume that $\Gamma$ contain an element of the third layer, 
		\begin{equation}\label{eq_dem_lem_kmn4_combinaciones_Gamma4_2}
		\gamma=\corchetes{\begin{array}{ccc}
	\alpha^{-2} & \gamma_{12} & \gamma_{13}\\
	0 & \alpha & \gamma_{23}\\
	0 & 0 & \alpha\\
	\end{array}}\in N_3,
	\end{equation}				
	for some $\gamma_{23},\alpha\in\C^\ast$ with $\valorabs{\alpha}\neq 1$. Since $\Gamma_p$ is a normal subgroup of $\Gamma$, then $\gamma g \gamma^{-1}\in\Gamma_p$. 
		$$\gamma g \gamma^{-1}=\corchetes{\begin{array}{ccc}
	1 & \alpha^{-3} & -\alpha^{-4}\gamma_{23}\\
	0 & 1 & 0\\
	0 & 0 & 1\\
	\end{array}}\in\Gamma_p$$
	implies that $\gamma_{23}=0$, contradicting that $\gamma\in N_3$.
	\item If $\rank(N_4)=1$, let $\gamma\in N_4$ be an element given by
		\begin{equation}\label{eq_dem_lem_kmn4_combinaciones_Gamma4_1}
		\gamma=\corchetes{\begin{array}{ccc}
	\beta & \gamma_{12} & \gamma_{13}\\
	0 & \alpha & \gamma_{23}\\
	0 & 0 & \alpha^{-1}\beta^{-1}\\
	\end{array}},
		\end{equation}				
	with $\alpha^2\beta\neq 1$, $\alpha\neq\beta$ such that $\gamma$ is loxodromic. Since $\Gamma_p$ is a normal subgroup of $\Gamma$, then $\gamma g \gamma^{-1}\in\Gamma_p$. This means that
		$$\gamma g \gamma^{-1}=\corchetes{\begin{array}{ccc}
	1 & \beta\alpha^{-1}n & -\beta^{-2}\gamma_{23}n\\
	0 & 1 & 0\\
	0 & 0 & 1\\
	\end{array}}\in\Gamma_p,$$
	and this implies that $\gamma_{23}=0$ and that, for $n=1$ in particular, $\beta\alpha^{-1}\in\Z$, in other words, $\beta=p\alpha$ for any $p\in\Z$. All of these conditions together with (\ref{eq_dem_lem_kmn4_combinaciones_Gamma4_1}) imply (\ref{eq_lem_kmn4_combinaciones_Gamma4_2}).
	\item Let $\gamma\in N_4$ be one generator of $N_4$ given by (\ref{eq_lem_kmn4_combinaciones_Gamma4_2}). Up to a suitable conjugation we can assume that 
		$$\gamma=\corchetes{\begin{array}{ccc}
	p\alpha & 0 & 0\\
	0 & \alpha & 0\\
	0 & 0 & p^{-1}\alpha^{-2}\\
	\end{array}},$$
	with $p\in\Z\setminus\set{0}$ and $\alpha\in\C^\ast$. Using part (ii), we know that 
		$$\mu=\corchetes{\begin{array}{ccc}
	q\beta & \mu_{12} & \mu_{13}\\
	0 & \beta & 0\\
	0 & 0 & q^{-1}\beta^{-2}\\
	\end{array}}.$$
	Since $\corchetes{\mu,\gamma}\in\prodint{g}$, comparing the respective entries we have $\mu_{13}(1-p^2\alpha^3)=0$ and $\mu_{12}(1-p)=\beta jpq$ for some $j\in\Z$. Therefore
	\begin{align}
	\mu_{13} &= 0\text{ or }p^2\alpha^3=1 \label{eq_dem_lem_kmn4_combinaciones_Gamma4_03} \\
	\mu_{12} &= \frac{\beta jpq}{1-p}. \label{eq_dem_lem_kmn4_combinaciones_Gamma4_04}
	\end{align}		
	A direct calculation shows that these two conditions are equivalent to prove that $\Gamma_p\rtimes\prodint{\gamma}$ is normal in $\Gamma_p\rtimes\prodint{\gamma}\rtimes\prodint{\mu}$. Therefore the only restrictions for $\mu$ are (\ref{eq_dem_lem_kmn4_combinaciones_Gamma4_03}) and (\ref{eq_dem_lem_kmn4_combinaciones_Gamma4_04}). This proves this part of the lemma.
	\item If $\rank(N_4)=3$, we denote by $\gamma$ and $\mu$ the first two generators of $N_4$ and, by the proof of (iii), we know the form of each element (see (ii) and (iii)). Furthermore, $\eta$ has the form (\ref{eq_lem_kmn4_combinaciones_Gamma4_3}). By the normality of $\Gamma_p\rtimes\prodint{\gamma}\rtimes\prodint{\mu}$ in $\Gamma$ we have $\eta g \gamma\mu\eta^{-1}=g^{n}\gamma^{m_1}\mu^{m_2}$ for some integers $n,m_1,m_2\in\Z$. Comparing entries 12 and 13 in both sides in the previous equation and solving for $eta_{13}$ and $k$ yields the expressions for these variables.\\	
	On the other hand, $\corchetes{\mu,\eta}\in\Gamma_p$, let $\xi$ be the entry 13 of $\corchetes{\mu,\eta}$, then 
		$$\xi=\frac{q\beta\mu_{13}\parentesis{1-r^2\delta^3}}{r\delta^2\parentesis{1-q^2\beta^3}}.$$
	If $\corchetes{\mu,\eta}=\id$ then $\xi=0$, if $\corchetes{\mu,\eta}\neq\id$ the $\xi\in\Z$. This completes the proof. 
	\end{enumerate}
\end{proof}

\begin{lem}\label{lem_kmn4_combinaciones_Gamma5}
If $\Gamma$ is such that $\Gamma_p$ is conjugated to 
	\begin{equation}\label{eq_lem_kmn4_combinaciones_Gamma5_1}
	\Gamma_5=\prodint{\corchetes{\begin{array}{ccc}
	1 & 0 & 1\\
	0 & 1 & 0\\
	0 & 0 & 1\\
	\end{array}}},
	\end{equation}
then:
	\begin{enumerate}
	\item If $\rank(N_3)=1$, then one generator of $N_3$ is the element
		\begin{equation}\label{eq_lem_kmn4_combinaciones_Gamma5_2}
		\gamma=\corchetes{\begin{array}{ccc}
	p & \gamma_{12} & \gamma_{13}\\
	0 & 1 & \gamma_{23}\\
	0 & 0 & 1\\
	\end{array}},
	\end{equation}				

	for $p\in\Z\setminus\set{-1,0,1}$, $\gamma_{12},\gamma_{13}\in\C$ and $\gamma_{23}\in\C^\ast$.
	\item It's not possible to have $\rank(N_3)=2$.
	\item If $\rank(N_3)=0$ and $\rank(N_4)=1$, then one generator of $N_4$ is the element
		\begin{equation}\label{eq_lem_kmn4_combinaciones_Gamma5_0}
		\gamma=\corchetes{\begin{array}{ccc}
		\alpha & \gamma_{12} & \gamma_{13}\\
		0 & p\alpha^{-2} & \gamma_{23}\\
		0 & 0 & p^{-1}\alpha\\
		\end{array}},
		\end{equation}
	with $\alpha\in\C^{\ast}$, $p\in\Z\setminus\set{0}$ such that $p^2\neq \alpha^3$ and $\gamma$ is loxodromic.
	\item If $\rank(N_3)=1$, then $\Gamma$ cannot contain elements in its fourth layer.	
	\item If $\rank(N_3)=0$ and $\rank(N_4)=2$, then one generator of $N_4$ is given by (\ref{eq_lem_kmn4_combinaciones_Gamma5_0}) and the other by
		\begin{equation}\label{eq_lem_kmn4_combinaciones_Gamma5_00}
		\mu=\corchetes{\begin{array}{ccc}
		\beta & 0 & \frac{jp}{1-p}\beta\\
		0 & q\beta^{-2} & 0\\
		0 & 0 & q^{-1}\beta\\
		\end{array}},
		\end{equation}
		for $\beta\in\C^\ast$, $q,j\in\Z\setminus\set{0}$ such that $q^2\neq\beta^3$ and $\mu$ is loxodromic.
	\item It is not possible to have $\rank(N_4)=3$.
	\end{enumerate}
\end{lem}

\begin{proof}
Assume that $\Gamma_p=\Gamma_5$ and let $g\in\Gamma_p$ given by
	$$g=\corchetes{\begin{array}{ccc}
	1 & 0 & n\\
	0 & 1 & 0\\
	0 & 0 & 1\\
	\end{array}}.$$
We now verify each conclusion:
\begin{enumerate}[(i)]
	\item Let $\gamma\in N_3$, then $\gamma$ has the form given by (\ref{eq_dem_lem_kmn4_combinaciones_Gamma4_2}). As before, the normality of $\Gamma_p$ in $\Gamma$ implies that $\gamma g \gamma^{-1}\in\Gamma_p$, and this implies that $\alpha^{-3}n\in \Z$ for all $n\in\Z$. Taking $n=1$, yields $\alpha=p^{-\frac{1}{3}}$. Substituting this in the form of $\gamma$ verifies (\ref{eq_lem_kmn4_combinaciones_Gamma5_2}).
	
	\item Let us suppose that $\rank(N_2)=2$ with $N_2=\prodint{\gamma,\mu}$ where $\gamma$ has the form (\ref{eq_lem_kmn4_combinaciones_Gamma5_2}). Then $\mu$ must have the same form and then we denote
	  $$\mu=\corchetes{\begin{array}{ccc}
	q & \mu_{12} & \mu_{13}\\
	0 & 1 & \mu_{23}\\
	0 & 0 & 1\\
	\end{array}},$$
	for some $q\in\Z$. Observe that
	\begin{equation}\label{eq_dem_lem_kmn4_combinaciones_Gamma5_1}
	\corchetes{\gamma,\mu}=\corchetes{\begin{array}{ccc}
	1 & \frac{\gamma_{12}(1-q)-\mu_{12}(1-p)}{pq} & \frac{\gamma_{13}(1-q)-\mu_{13}(1-p)+\gamma_{12}\mu_{23}-\gamma_{23}\mu_{12}}{pq}\\
	0 & 1 & 0\\
	0 & 0 & 1\\
	\end{array}}.
	\end{equation}
This element $\corchetes{\gamma,\mu}$ is either parabolic (and therefore, $\corchetes{\gamma,\mu}\in\Gamma_p$) or $\corchetes{\gamma,\mu}=\id$. In either case, the entry 12 of $\corchetes{\gamma,\mu}$ satisfies
	\begin{equation}\label{eq_dem_lem_kmn4_combinaciones_Gamma5_2}
	\frac{\gamma_{12}(1-q)-\mu_{12}(1-p)}{pq} = 0.	
	\end{equation}	 
This means that 
	$$\fix\parentesis{\corchetes{\begin{array}{ccc}
	q & \mu_{12} \\
	0 & 1 \\
	\end{array}}}=\fix\parentesis{\corchetes{\begin{array}{ccc}
	p & \gamma_{12} \\
	0 & 1 \\
	\end{array}}},$$
and therefore, conjugating $\Gamma$ by 
	$$A=\corchetes{\begin{array}{ccc}
	1 & -\frac{\gamma_{12}}{1-p} & 0\\
	0 & 1 & 0\\
	0 & 0 & 1\\
	\end{array}}$$
and using (\ref{eq_dem_lem_kmn4_combinaciones_Gamma5_2}), we can assume that
	$$\gamma=\corchetes{\begin{array}{ccc}
	p & 0 & \gamma_{13}\\
	0 & 1 & \gamma_{23}\\
	0 & 0 & 1\\
	\end{array}}\text{, and }\mu=\corchetes{\begin{array}{ccc}
	q & 0 & \mu_{13}\\
	0 & 1 & \mu_{23}\\
	0 & 0 & 1\\
	\end{array}},$$
for some $p,q\in\Z\setminus\set{0}$. With this simplified form, the commutator (\ref{eq_dem_lem_kmn4_combinaciones_Gamma5_1}) gives the conditions:
	\begin{align}
	\text{If }\corchetes{\gamma,\mu}=\id &\Rightarrow \frac{\gamma_{13}}{1-p}=\frac{\mu_{13}}{1-q}.\label{eq_dem_lem_kmn4_combinaciones_Gamma5_3}\\
	\text{If }\corchetes{\gamma,\mu}\neq\id &\Rightarrow \frac{\gamma_{13}{1-q}-\mu_{13}(1-p)}{pq}\in\Z.\label{eq_dem_lem_kmn4_combinaciones_Gamma5_4}
	\end{align}
Since $\Gamma_p\rtimes\prodint{\gamma}$ is normal in $\Gamma$, then $\mu g \gamma \mu^{-1}=g^n\gamma^k$ for some $n,k\in\Z$. Comparing the entries 23 of both expressions yields $k=1$, and then using the entries 13 we have
	\begin{equation}\label{eq_dem_lem_kmn4_combinaciones_Gamma5_5}
	q(1+\gamma_{13})+\mu_{13}(1-p)=1+\gamma_{13}.
	\end{equation}
If $\corchetes{\gamma,\mu}=\id$, then (\ref{eq_dem_lem_kmn4_combinaciones_Gamma5_3}) yield $q=1$, which cannot happen. If $\corchetes{\gamma,\mu}=\id$, then (\ref{eq_dem_lem_kmn4_combinaciones_Gamma5_3}) implies 
	$$\frac{\gamma_{13}{1-q}-\mu_{13}(1-p)}{pq}=j,$$
for some $j\in\Z$. Combining this equation with (\ref{eq_dem_lem_kmn4_combinaciones_Gamma5_5}), we have
	$$q=\frac{1}{1+jp}\in\Z$$	
	which cannot happen. Therefore there cannot be an element in $N_3\setminus\prodint{\gamma}$.
	\item The proof is similar to the proof of (i) above.
	\item We know, by (i), that if $\rank(N_3)=1$, then one generator of $N_3$ is the element $\gamma$ given by (\ref{eq_lem_kmn4_combinaciones_Gamma5_2}). Let us assume that there is an element $\mu\in N_4$, given by
		$$\mu=\corchetes{\begin{array}{ccc}
	\alpha & \mu_{12} & \mu_{13}\\
	0 & \beta & \mu_{23}\\
	0 & 0 & \alpha^{-1}\beta^{-1}\\
	\end{array}},$$
	where $\alpha,\beta\in\C^{\ast}$, such that $\alpha\beta^{-2}\neq 1$ and $\mu$ is loxodromic. Observe that 
		$$\corchetes{\gamma,\mu}=\corchetes{\begin{array}{ccc}
	\alpha & \frac{\gamma_{12}(\beta-\alpha)+\mu_{12}(1-p)}{\beta} & \xi_{13}\\
	0 & \beta & \frac{\gamma_{23}}{\alpha\beta^2}(1-\alpha\beta^2)\\
	0 & 0 & \alpha^{-1}\beta^{-1}\\
	\end{array}}\in\Gamma_p,$$
	then $\frac{\gamma_{23}}{\alpha\beta^2}(1-\alpha\beta^2)=0$. This means that, either $\gamma_{23}=0$ or $\alpha\beta^2=1$. Both conditions are forbidden by hypothesis, then $\Gamma$ cannot contain elements in $N_4$.

	\item If $\rank(N_4)=2$ and $\gamma$ and $\mu$ are two generators then $\gamma$ has the form given by (\ref{eq_lem_kmn4_combinaciones_Gamma5_0}), observe that we can conjugate $\Gamma$ by 
		$$A=\corchetes{\begin{array}{ccc}
	1 & -\frac{\gamma_{12}}{(p\alpha^{-2}-\alpha)} & -\frac{p\parentesis{\gamma_{13}(p-\alpha^3)-\alpha^2\gamma_{12}\gamma_{23}}}{\alpha(1-p)(p-\alpha^3)}\\
	0 & 1 & \frac{p\alpha^2\gamma_{23}}{p^2-\alpha^3}\\
	0 & 0 & 1\\
	\end{array}}$$
	and then we can assume that $\gamma_{12}=\gamma_{13}=\gamma_{23}=0$ in (\ref{eq_lem_kmn4_combinaciones_Gamma5_0}). Since $\corchetes{\mu,\gamma}\in\Gamma_p=\Gamma_5$ then $\mu$ has the form
		$$\mu=\corchetes{\begin{array}{ccc}
	\beta & 0 & \mu_{13}\\
	0 & q\beta^{-2} & 0\\
	0 & 0 & q^{-1}\beta\\
	\end{array}}$$
for $\beta\in\C^\ast$ and $q\in\Z\setminus\set{0}$ such that $q^2\neq\beta^3$ and $\mu$ is loxodromic. If $\corchetes{\mu,\gamma}=\id$ then $\Gamma$ would be commutative, therefore we can assume that $\corchetes{\mu,\gamma}\neq\id$ and then, comparing the entries 13 of $\corchetes{\mu,\gamma}$ and $g$ yields 
	\begin{equation}\label{eq_dem_lem_kmn4_combinaciones_Gamma5_6}
	\frac{1-p}{p\beta}\mu_{13}=j,
	\end{equation}		
for some $j\in\Z\setminus\set{0}$. Then, $\mu_{13}=jp\beta(1-p)^{-1}$, this verifies (\ref{eq_lem_kmn4_combinaciones_Gamma5_00}).
\end{enumerate}
\end{proof}

\begin{obs}\label{obs_resumen_casos_kmn4}
We now summarize all subcases kmn(4), using Lemmas \ref{lem_kmn4_combinaciones_Gamma1}, \ref{lem_kmn4_combinaciones_Gamma2}, \ref{lem_kmn4_combinaciones_Gamma3}, \ref{lem_kmn4_combinaciones_Gamma4} and \ref{lem_kmn4_combinaciones_Gamma5}.
\begin{itemize}
\item Subcase 101(4):
	\begin{enumerate}
	\item 
	$$\Gamma=\SET{\corchetes{\begin{array}{ccc}
	1 & m & 0\\
	0 & 1 & 0\\
	0 & 0 & 1\\
	\end{array}}\corchetes{\begin{array}{ccc}
	p\beta & \gamma_{12} & \gamma_{13}\\
	0 & \beta & 0\\
	0 & 0 & p^{-1}\beta^{-2}\\
	\end{array}}^k}{m,k\in\Z},$$
	where $p\in\Z\setminus\set{0,1}$, $\beta\in\C^{\ast}$ such that the second generator is loxodromic, $\valorabs{\beta}\neq 1$ and $\gamma_{12},\gamma_{13}\in\C$.
	\item 
	$$\Gamma=\SET{\corchetes{\begin{array}{ccc}
	1 & 0 & m\\
	0 & 1 & 0\\
	0 & 0 & 1\\
	\end{array}}\corchetes{\begin{array}{ccc}
	\alpha & \gamma_{12} & \gamma_{13}\\
	0 & p\alpha^{-2} & \gamma_{23}\\
	0 & 0 & p^{-1}\alpha\\
	\end{array}}^k}{m,k\in\Z},$$
	with $\alpha\in\C^{\ast}$, $p\in\Z\setminus\set{0}$ such that $p^2\neq \alpha^3$ and the second generator is loxodromic.
	\end{enumerate}
\item The two possible subcases 102(4):
	\begin{enumerate}
	\item $$\Gamma=\SET{\corchetes{\begin{array}{ccc}
	1 & 0 & k\\
	0 & 1 & 0\\
	0 & 0 & 1\\
	\end{array}}\corchetes{\begin{array}{ccc}
	p\alpha & \gamma_{12} & \gamma_{13}\\
	0 & \alpha & 0\\
	0 & 0 & p^{-1}\alpha^{-2}\\
	\end{array}}^m \corchetes{\begin{array}{ccc}
	q\beta & \frac{\beta jpq}{1-p} & \mu_{13}\\
	0 & \beta & 0\\
	0 & 0 & q^{-1}\beta^{-2}\\
	\end{array}}^n}{k,m,n\in\Z},$$
	for $p,q\in\Z\setminus\set{-1,0,1}$, $j\in\Z$, $\gamma_{12},\gamma_{13},\mu_{13}\in\C$ and $\alpha,\beta\in\C^\ast$ such that the second and third generators are loxodromic. If $p^2\alpha^3\neq 1$, then $\mu_{13}=0$.

	\item 
	$$\Gamma=\SET{\corchetes{\begin{array}{ccc}
	1 & k & 0\\
	0 & 1 & 0\\
	0 & 0 & 1\\
	\end{array}}\corchetes{\begin{array}{ccc}
	\alpha & \gamma_{12} & \gamma_{13}\\
	0 & p\alpha^{-2} & \gamma_{23}\\
	0 & 0 & p^{-1}\alpha\\
	\end{array}}^m \corchetes{\begin{array}{ccc}
	\beta & 0 & \frac{jp\beta}{1-p}\\
	0 & q\beta^{-2} & 0\\
	0 & 0 & q^{-1}\beta\\
	\end{array}}^n}{k,m,n\in\Z},$$
	for $p,q\in\Z\setminus\set{-1,0,1}$, $\gamma_{ij}\in\C$, $\alpha,\beta\in\C^\ast$ such that $p^2\neq\alpha^3$, $q^2\neq\beta^3$ and the second and third generators are loxodromic.
	\end{enumerate}

\item Subcase 103(4):	
	$$\Gamma=\SET{\corchetes{\begin{array}{ccc}
	1 & 0 & n_1\\
	0 & 1 & 0\\
	0 & 0 & 1\\
	\end{array}}\gamma^{n_2}\mu^{n_3}\eta^{n_4}}{n_1,n_2,n_3,n_4\in\Z},$$
	where $\gamma,\mu,\eta$ have the form defined in (ii), (iii) and (iv) of Lemma \ref{lem_kmn4_combinaciones_Gamma4}.

\item Subcase 110(4):
	$$\Gamma=\SET{\corchetes{\begin{array}{ccc}
	1 & 0 & m\\
	0 & 1 & 0\\
	0 & 0 & 1\\
	\end{array}}\corchetes{\begin{array}{ccc}
	p & \gamma_{12} & \gamma_{13}\\
	0 & 1 & \gamma_{23}\\
	0 & 0 & 1\\
	\end{array}}^n}{m,n\in\Z},$$
	for $p\in\Z\setminus\set{-1,0,1}$, $\gamma_{12},\gamma_{13}\in\C$ and $\gamma_{23}\in\C^\ast$.

\item The 3 possibilities for subcase 201(4):
	\begin{enumerate}
	\item
	$$\Gamma=\SET{\corchetes{\begin{array}{ccc}
	1 & n & m\\
	0 & 1 & 0\\
	0 & 0 & 1\\
	\end{array}}\corchetes{\begin{array}{ccc}
	\alpha q & \gamma_{12} & \gamma_{13}\\
	0 & \alpha^{-2}q^{2} & \alpha^5 r\\
	0 & 0 & \alpha\\
	\end{array}}^k}{m,n,k\in\Z}$$
	for $\alpha\in\C^\ast$ such that $\gamma$ is loxodromic, $q\in\Z\setminus\set{0}$, $r\in\Q$, $\gamma_{12},\gamma_{13}\in\C$.
	\item 
	$$\Gamma=\SET{\corchetes{\begin{array}{ccc}
	1 & 0 & n+my\\
	0 & 1 & 0\\
	0 & 0 & 1\\
	\end{array}}\corchetes{\begin{array}{ccc}
	\alpha & \gamma_{12} & \gamma_{13}\\
	0 & z_{p,q}\alpha^{-2} & \gamma_{23}\\
	0 & 0 & z_{p,q}^{-1}\alpha\\
	\end{array}}^k}{m,n,k\in\Z}$$
	where $y\nin\R$, $z_{p,q}=\parentesis{p+q\text{Re}(x)}+iq\text{Im}(x)$, $p,q\in\Z$ such that $\valorabs{p}+\valorabs{q}\neq 0$, $\alpha\in\C^\ast$ and $\gamma_{12},\gamma_{13},\gamma_{23}\in\C$. The second generator must be loxodromic.
	\item 
	$$\Gamma=\SET{\corchetes{\begin{array}{ccc}
	1 & n+mx & 0\\
	0 & 1 & 0\\
	0 & 0 & 1\\
	\end{array}}\corchetes{\begin{array}{ccc}
	z_{p,q} & \gamma_{12} & \gamma_{13}\\
	0 & 1 & 0\\
	0 & 0 & z_{p,q}^{-1}\beta^{-3}\\
	\end{array}}^k}{m,n,k\in\Z}$$
	where $x\nin\R$, $z_{p,q}=\parentesis{p+q\text{Re}(x)}+iq\text{Im}(x)$, $p,q\in\Z$ such that $\valorabs{p}+\valorabs{q}\neq 0$, $\beta\in\C^\ast$ and $\gamma_{12}, \gamma_{13}\in\C$. The second generator must be loxodromic.
	\end{enumerate}	

\item The two possibilities for subcase 210(4):
	\begin{enumerate}
	\item
	$$\Gamma=\SET{\corchetes{\begin{array}{ccc}
	1 & 0 & n+my\\
	0 & 1 & 0\\
	0 & 0 & 1\\
	\end{array}}\corchetes{\begin{array}{ccc}
	\parentesis{p+q\text{Re}(y)}+iq\text{Im}(y) & \gamma_{12} & \gamma_{13}\\
	0 & 1 & \gamma_{23}\\
	0 & 0 & 1\\
	\end{array}}^k}{m,n,k\in\Z}$$
	for $y\nin\R$, $p,q\in\Z$ such that $\valorabs{p}+\valorabs{q}\neq 0$ and $\valorabs{\parentesis{p+q\text{Re}(y)}+iq\text{Im}(y)}\neq 1$, $\gamma_{12},\gamma_{23}\in\C$ and $\gamma_{23}\in\C^\ast$.
	\item 	
	$$\Gamma=\SET{\corchetes{\begin{array}{ccc}
	1 & n & m\\
	0 & 1 & 0\\
	0 & 0 & 1\\
	\end{array}}\corchetes{\begin{array}{ccc}
	p & \gamma_{12} & \gamma_{13}\\
	0 & 1 & r\\
	0 & 0 & 1\\
	\end{array}}^k}{m,n,k\in\Z}$$
	for $p\in\Z\setminus\set{0}$, $p\in\Q\setminus\set{0}$, $\gamma_{12},\gamma_{13}\in\C$ and $\gamma_{23}\in\C^\ast$.
	\end{enumerate}

\item Subcase 202(4):
	\begin{enumerate}
	\item
		\begin{multline}
		\Gamma=\left\lbrace\corchetes{\begin{array}{ccc}
	1 & n+mx & 0\\
	0 & 1 & 0\\
	0 & 0 & 1\\
	\end{array}}\corchetes{\begin{array}{ccc}
	p_1+q_1x & 0 & 0\\
	0 & 1 & 0\\
	0 & 0 & (p_1+q_1x)^{-1}\alpha^{-3}\\
	\end{array}}\right.\\
	\left.\corchetes{\begin{array}{ccc}
	p_2+q_2x & \mu_{k_1,k_2} & 0\\
	0 & 1 & 0\\
	0 & 0 & (p_2+q_2x)^{-1}\beta^{-3}\\
	\end{array}}\middle| n,m,k_1,k_2\in\Z\right\rbrace,
	\end{multline}
	The conditions are described in Lemma \ref{lem_kmn4_combinaciones_Gamma3}.
	\end{enumerate}

\item Subcases 111(4), 120(4), 121(4), 130(4), 310(4) and 301(4) are impossible.	

\end{itemize}
\end{obs}

\subsubsection*{$\Gamma_p$ is not commutative}

In this part, we now describe the non-commutative upper triangular groups such that its parabolic part is not commutative. This parabolic part is given by part (5) and (6) of Theorem \ref{thm_descripcion_parte_parabolica}.\\

Observe that $\rank\parentesis{\Gamma_{\mathbf{w}}}=3$. This dismisses the subcases 1mn(5) and 2mn(5), where $m,n\in\set{0,1,2,3}$. In the following lemma, we dismiss the remaining subcases 310(5) and 301(5).

\begin{lem}\label{lem_kmn5_combinaciones}
If $\Gamma$ is such that $\Gamma_p$ is conjugated to $\Gamma_{\mathbf{w}}$, for some $\mathbf{w}$ as in the case (5) of Theorem \ref{thm_descripcion_parte_parabolica}. If $\rank\parentesis{\Gamma_{\mathbf{w}}}=3$, then $\Gamma$ cannot contain loxodromic elements.
\end{lem}

\begin{proof}
We can assume that $\Gamma_p=\Gamma_{\mathbf{w}}$. Observe that
	$$\Gamma_{\mathbf{w}}=\prodint{\corchetes{\begin{array}{ccc}
	1 & 1 & 0\\
	0 & 1 & 0\\
	0 & 0 & 1\\
	\end{array}},\corchetes{\begin{array}{ccc}
	1 & c & d\\
	0 & 1 & 0\\
	0 & 0 & 1\\
	\end{array}},\corchetes{\begin{array}{ccc}
	1 & x & y\\
	0 & 1 & 1\\
	0 & 0 & 1\\
	\end{array}}},$$
denote respectively $g_1,g_2,g_3$ these three generators. Considering $g_1$ and $g_2$ and since $d\neq 0$, it follows from Proposition \ref{prop_casos_kmn4_conjugaciones} that $\Gamma_{\mathbf{w}}$ is conjugated to 
	$$\Gamma_{\mathbf{w'}}=\prodint{\corchetes{\begin{array}{ccc}
	1 & 1 & 0\\
	0 & 1 & 0\\
	0 & 0 & 1\\
	\end{array}},\corchetes{\begin{array}{ccc}
	1 & 0 & 1\\
	0 & 1 & 0\\
	0 & 0 & 1\\
	\end{array}},\corchetes{\begin{array}{ccc}
	1 & x-\frac{c}{d}y & \frac{1}{d}y\\
	0 & 1-\frac{c}{d} & \frac{1}{d}\\
	0 & -\frac{1}{d^2}c & 1+\frac{c}{d}\\
	\end{array}}},$$
and therefore we can consider these generators for $\Gamma_p$ instead.\\ 

Let us suppose first that there is an element $\gamma\in N_3$. Observe that $\prodint{g_1,g_2}$ is normal in $\Gamma$ and therefore $\gamma g_1 g_2 \gamma^{-1}\in\prodint{g_1,g_2}$, then, by lemma \ref{lem_kmn4_combinaciones_Gamma1}
	$$\gamma=\corchetes{\begin{array}{ccc}
	p & \gamma_{12} & \gamma_{13}\\
	0 & 1 & \frac{q}{p}\\
	0 & 0 & 1\\
	\end{array}},$$
for some $p\in\Z\setminus\set{-1,0,1}$, $q\in\Z\setminus\set{0}$. Since $\Gamma_p$ is normal in $\Gamma_p\rtimes\prodint{\gamma}$ then $\gamma g_3 \gamma^{-1}=g_1^k g_2^n g_3^m$ for some $k,n,m\in\Z$. Comparing the entries 32 yields $m=1$, substituting in the equation resulting of comparing the entries 22 we have 
	$$\frac{c^2q}{dp}=0,$$
then $q=0$, which is a contradiction.\\
	
Now assume that there is a loxodromic element $\gamma\in N_4$. Using the same argument as in the case of an element in the third layer we know that 
	$$\gamma=\corchetes{\begin{array}{ccc}
	qp & \gamma_{12} & \gamma_{13}\\
	0 & q & r\\
	0 & 0 & p\\
	\end{array}},$$
for some $p,q,r\in\Z$, $p,q\neq 0$. Consider the commutator
	$$\corchetes{\gamma,g_1}=\corchetes{\begin{array}{ccc}
	1 & 1-\frac{1}{p} & -\frac{r}{pq}\\
	0 & 1 & 0\\
	0 & 0 & 1\\
	\end{array}}\in\prodint{g_1,g_2}$$
and therefore $1-\frac{1}{p}\in\Z$, then $p=\pm 1$. If $p=1$ then $\gamma$ is a complex homothety, which cannot happen by Corollary \ref{cor_HC3_no_hay_en_no_conmutativos}. If $p=-1$ then $\lambda_{12}(\gamma)$ is a torsion element in the torsion free group $\lambda_{12}(\Gamma)$. This concludes the proof.

\end{proof}

Observe that $3\leq \rank\parentesis{\Gamma_{W,a,b,c}}\leq 4$ by the definition of $\Gamma_{W,a,b,c}$. This means that the subcases 1mn(6) and 2mn(6) cannot happen, if $m,n\in\set{0,1,2,3}$. Therefore, we just need to describe the subcases 310(6) and 301(6), but these cases cannot happen, as we prove in the following lemma.

\begin{lem}\label{lem_kmn6_combinaciones}
If $\Gamma$ is such that $\Gamma_p$ is conjugated to $\Gamma_{W,a,b,c}$, for some $W,a,b,c$ as in the case (6) of Theorem \ref{thm_descripcion_parte_parabolica}. If $\rank\parentesis{\Gamma_{W,a,b,c}}=3$, then $\Gamma$ cannot contain loxodromic elements.
\end{lem} 

\begin{proof}
Assume that $\Gamma_p=\Gamma_{W,a,b,c}$, with $W,a,b,c$ as in the hypothesis of the lemma. Since $\rank\parentesis{\Gamma_{W,a,b,c}}=3$ then $\rank(W)=1$ and then we write $W=\SET{k w_0}{k\in \Z}$. Denote by $g_1,g_2,g_3$ to the generators of $\Gamma_p$ in this way:
	$$g_1=\corchetes{\begin{array}{ccc}
	1 & 0 & w_0\\
	0 & 1 & 0\\
	0 & 0 & 1\\
	\end{array}},\; g_2=\corchetes{\begin{array}{ccc}
	1 & 1 & 0\\
	0 & 1 & 1\\
	0 & 0 & 1\\
	\end{array}},\; g_3=\corchetes{\begin{array}{ccc}
	1 & a+c & b\\
	0 & 1 & c\\
	0 & 0 & 1\\
	\end{array}}.$$

Let us assume first that there is a loxodromic element $\gamma\in N_3$ in the third layer of $\Gamma$, then
	\begin{equation}\label{eq_dem_lem_kmn6_combinaciones1}
	\gamma=\corchetes{\begin{array}{ccc}
	\alpha^{-2} & \alpha_{12} & \gamma_{13}\\
	0 & \alpha & \gamma_{23}\\
	0 & 0 & \alpha\\
	\end{array}}
	\end{equation}		
for some $\alpha\in\C^\ast$ such that $\valorabs{\alpha}\neq 1$ and $\gamma_{23}\neq 0$. A direct calculation shows that
	\begin{equation}\label{eq_dem_lem_kmn6_combinaciones2}
	\Gamma_p = \SET{\corchetes{\begin{array}{ccc}
	1 & mn(a+c) & \psi_{k,m,n}\\
	0 & 1 & cmn\\
	0 & 0 & 1\\
	\end{array}}}{k,m,n\in\Z}, 
	\end{equation}
we are not interested in the value of $\psi_{k,m,n,w}$ for now. Since $\Gamma_p$ is normal in $\Gamma$, then 
	$$\gamma g_1 g_2 g_3 \gamma^{-1}=\corchetes{\begin{array}{ccc}
	1 & \alpha^{-3}(a+c+1) & \xi\\
	0 & 1 & c+1\\
	0 & 0 & 1\\
	\end{array}}\in \Gamma_p.$$  
Comparing the entry 23 with the entry 23 of (\ref{eq_dem_lem_kmn6_combinaciones2}), we get that $c+1=cmn$ for some $m,n\in\Z$. If $c\neq 0$, this means that
	\begin{equation}\label{eq_dem_lem_kmn6_combinaciones3}
	1+\frac{1}{c}\in\Z
	\end{equation}
then $c\in\Q$ and therefore $\set{1,c}$ is $\R$-linearly dependent, which means that $\set{1,c}$ is $\Z$-linearly independent (by hypothesis). From (\ref{eq_dem_lem_kmn6_combinaciones3}) we get $\frac{1}{c}\in\Z$, and then $-1+\frac{1}{c}c=0$ implies that $\set{1,c}$ is $\Z$-linearly dependent. This contradiction shows that $\Gamma$ cannot contain elements in its third layer.\\

Now assume that $\gamma\in N_4$ is a loxodromic element in the fourth layer, we denote
	$$\gamma=\corchetes{\begin{array}{ccc}
	\alpha & \alpha_{12} & \gamma_{13}\\
	0 & \beta & \gamma_{23}\\
	0 & 0 & \alpha^{-1}\beta^{-1}\\
	\end{array}}$$	 
with $\alpha,\beta\in\C^\ast$. Since $\Gamma_p$ is normal in $\Gamma$, in particular 
	$$\gamma g_1 g_3 \gamma^{-1}=\corchetes{\begin{array}{ccc}
	1 & \frac{\alpha}{\beta}(a+c) & \xi\\
	0 & 1 & c\alpha\beta\\
	0 & 0 & 1\\
	\end{array}}\in \Gamma_p.$$ 
Comparing the entries 12 and 23 with the same entries in (\ref{eq_dem_lem_kmn6_combinaciones2}), we get
	\begin{align}
	\frac{\alpha}{\beta} &= mn\text{, if }a+c\neq 0 \label{eq_dem_lem_kmn6_combinaciones4}\\
	\alpha\beta^2 &= mn\text{, if }c\neq 0 \label{eq_dem_lem_kmn6_combinaciones5}
	\end{align}
for some $m,n\in\Z$. From (\ref{eq_dem_lem_kmn6_combinaciones4}) and (\ref{eq_dem_lem_kmn6_combinaciones5}) we get that $\beta^3=1$. This implies that $\alpha=p\beta$ for some $p\in\Z$. Then, 
	\begin{equation}\label{eq_dem_lem_kmn6_combinaciones6}
	\gamma=\corchetes{\begin{array}{ccc}
	p\beta & \alpha_{12} & \gamma_{13}\\
	0 & \beta & \gamma_{23}\\
	0 & 0 & p^{-1}\beta^{-2}\\
	\end{array}},
	\end{equation}
for $\beta$ a cubic root of the unity. In the proof of proposition 7.15 of \cite{ppar} we see that $\prodint{g_1}$ is a normal subgroup of $\Gamma$ and therefore $\gamma g_1 \gamma^{-1}\in \prodint{g_1}$. As we proved in lemma \ref{lem_kmn4_combinaciones_Gamma5}, this implies that 
	\begin{equation}\label{eq_dem_lem_kmn6_combinaciones7}
	\gamma=\corchetes{\begin{array}{ccc}
	\alpha & \alpha_{12} & \gamma_{13}\\
	0 & q\alpha^{-2} & \gamma_{23}\\
	0 & 0 & q^{-1}\alpha^{-1}\\
	\end{array}},
	\end{equation}	
for some $q\in\Z$ and $\alpha\in\C^\ast$ such that $\lambda$ is loxodromic. Comparing the diagonal entries of both (\ref{eq_dem_lem_kmn6_combinaciones6}) and (\ref{eq_dem_lem_kmn6_combinaciones7}), and since $\beta^3=1$ we get 
	$$p=\pm 1,\;\;\; q=1,\;\;\; \alpha=p,\;\;\; \beta=1.$$
Then $\gamma$ is either parabolic or induces a torsion element in $\lambda_{12}(\Gamma)$, both are impossible. This finishes the proof.
\end{proof}
\chapter{The Frances limit set}
\label{chapter_frances}

In this chapter we propose a new limit set for the action of discrete subgroups of $\PSLN$. In \cite{frances}, Charles Frances defines a limit set for the action of discrete subgroups of the Lorentzian M\"obius group $O(2,n)$. We apply some of these ideas to define the Frances limit set, this limit set seems to be the \emph{right} concept of limit set for complex dimension $n\geq 3$. In general, this set is smaller that the Kulkarni limit set.

\section{The Cartan decomposition}

For any element $\mathbf{g}\in\SLN$ we have the singular value decomposition

	$$\mathbf{g}=U \mu(\mathbf{g}) V^\ast$$
	
where $U,V\in U(n+1)$ are complex unitary matrices and $\mu(\mathbf{g})\in\mathcal{M}_n(\R^+)$ is a diagonal matrix, $V^\ast$ denotes the conjugate transpose of the matrix $V$. The columns of $U$ are the eigenvectors of $\mathbf{g}\mathbf{g}^\ast$, the columns of $V$ are the eigenvectors of $\mathbf{g}^\ast\mathbf{g}$.\\

We call $\mu(\mathbf{g})$ the \emph{Cartan projection} of $\mathbf{g}$.\\

\subsection{The diagonal case}

From now on we restrict ourselves to transformations $g\in\PSLN$ given by lifts of the form

\begin{equation}\label{eq_forma_diagonal}
\mathbf{g}=\parentesis{\begin{array}{ccc}
\lambda_{1} & & \\
 & \ddots & \\
 & & \lambda_{n+1}
\end{array}}
\end{equation}

with $\lambda_1\geq ...\geq \lambda_{n+1}$. Let $\Gamma\subset\PSLN$ be a discrete subgroup and $\set{g_n}\in\Gamma$ a sequence of distinct elements of the form

\begin{equation}\label{eq_forma_diagonal_suc}
\mathbf{g}_k=\parentesis{\begin{array}{ccc}
\lambda_{1,k} & & \\
 & \ddots & \\
 & & \lambda_{n+1,k}
\end{array}}
\end{equation}

If there is no ambiguity we will denote indistinctly by $g$ both the transformation in $\PSLN$ and any of its lifts in $\SLN$. Let $\mathcal{B}=\set{e_1,...,e_{n+1}}$ be the canonical base of $\C^{n+1}$.\\

We say that $e_i\sim e_{i+1}$ are related in the matrix $g_k$ if their corresponding eigenvalues $\lambda_{i,k}$ and $\lambda_{i+1,k}$ satisfy:

	$$\lim_{k\rightarrow\infty} \frac{\lambda_{i,k}}{\lambda_{i+1,k}}\in\R^+.$$

It is clear that $\sim$ defines an equivalence relation in $\set{\lambda_{1,k},...,\lambda_{n+1,k}}$. Factoring in each equivalence class the first element we have

\begin{equation}\label{eq_forma_bloques}
\mathbf{g}_k=\parentesis{\begin{array}{ccc}
\alpha_{1,k} D_{1,k}& & \\
 & \ddots & \\
 & & \alpha_{m,k} D_{m,k}
\end{array}}
\end{equation}	

where $\alpha_{1,k}\geq...\geq \alpha_{m,k}$ and $D_{i,k}$ are diagonal matrices.\\

We denote by $\mathcal{A}_i\subset\CP^n$ the projectivization of the vector subspace $\tilde{\mathcal{A}_i}\subset\C^{n+1}$ generated by the columns of the matrix $D_{i,k}$. Clearly the spaces $\mathcal{A}_j$ are invariant under the sequence $\set{g_k}$.\\ 

We call the spaces $\mathcal{A}_1$ and $\mathcal{A}_m$ the \emph{attracting space} and \emph{repelling space} respectively. We denote 
	$$c_i=\text{dim}_{\C}\parentesis{\tilde{\mathcal{A}}_i}.$$

Each matrix $D_{i,k}$ has the form
	$$
	D_{i,k}=\parentesis{\begin{array}{ccc}
	\beta^{(i)}_{1,k} & & \\
	 & \ddots & \\
	 & & \beta^{(i)}_{c_i,k}
	\end{array}}
	$$

\begin{obs}
\mbox{}
	\begin{enumerate}[(i)]
	\item For $Z\in\CP^n$, we regroup the homogeneous coordinates $Z=\corchetes{\zeta_1:...:\zeta_{n+1}}$ as
		$Z=\corchetes{z_1:...:z_m}$ where
			$$z_i=\corchetes{\zeta_{c_1+...+c_{i-1}+1}:...:\zeta_{c_1+...+ c_i}}.$$
	\item Let $Z=\corchetes{z_1:...:z_m}\in\CP^n$. If there's no ambiguity we denote as well by $z_i$ the point $\corchetes{\overline{0}:...:\overline{0}:z_i:\overline{0}:...:\overline{0}}\in\mathcal{A}_i$, where $\overline{0}$ is the projectivization of the vector $0\in\tilde{\mathcal{A}}_i$ for $i=1,...,m$.
	\end{enumerate}

\end{obs}

\begin{prop}\label{prop_bloques_tienden_a_matrices_invertibles}
The sequence of matrices $\set{D_{i,k}}$ converges to an invertible matrix $D_i$, for $i=1,...,m$.
\end{prop}

\begin{proof}
Suppose that $D_i$ is non-invertible for some $i=1,...,m$. Since $D_i$ is dia-gonal, some element of the diagonal satisfies $a^{(i)}_{j,k}\overset{k\rightarrow\infty}{\rightarrow} 0$. However
	$$a^{(i)}_{j,k}=\frac{\lambda_{q,k}}{\alpha_{i,k}}$$
for some $q=1,...,n+1$. By construction of the blocks we have that
	$$\lim_{k\rightarrow\infty}\frac{\alpha_{i,k}}{\lambda_{q,k}}\in\R^+$$
which contradicts that $a^{(i)}_{j,k}\overset{k\rightarrow\infty}{\rightarrow} 0$. 
\end{proof}

The proof of the following proposition is analogous to the proof of the previous proposition.

\begin{prop}\label{prop_cocientes_coeficientes_a_infinito}
For some $i=1,...,m-1$, it holds
	$$\lim_{k\rightarrow\infty}\frac{\alpha_{i,k}}{\alpha_{i+1,k}}=\infty$$
\end{prop}

\begin{defn}
Let $\set{g_n}\subset \Gamma$ be a divergent sequence, for each $z\in\CP^n$ we define the set
	$$D_{\set{g_n}}(z)=\bigcup_{z_n\rightarrow z}\set{\text{acummulation points of }\set{g_n(z_n)}}$$
where the union is taken over all the sequences $\set{z_n}$ converging to $z$. For $E\subset\CP^n$ we define
	$$D_{\set{g_n}}(E)= \bigcup_{z\in E}D_{\set{g_n}}(z).$$
If we denote by $\text{Div}(\Gamma)$ to the set of divergent sequences of distinct elements in $\Gamma$ then
	$$D_\Gamma(E)=\bigcup_{\set{g_n}\in \text{Div}(\Gamma)}D_{\set{g_n}}(E).$$ 
\end{defn}

\begin{prop}
Let $K\subset\CP^n$ be a compact set, the set $D_{\Gamma}(K)$ coincides with the set of accumulation points of $\Gamma$-orbits of $K$.
\end{prop}

\begin{proof}
Let $K\subset\CP^n$ be a compact set. Assume that $w\in D_{\Gamma}(K)$, then for some $z\in K$, some sequence $\set{z_k}$ converging to $z$ and some sequence $\set{\gamma_k}\subset\Gamma$ it holds that $w$ is an accumulation point of $\set{\gamma_k(z_k)}$. For any $k$, $\gamma_k(z_k)\in\gamma_k(K)$, then $w$ is an accumulation point of the sequence $\set{\gamma_k(K)}$ and therefore, of $\set{\gamma(K)}_{\gamma\in\Gamma}$.\\

Now assume that $w$ is an accumulation point of $\set{\gamma(K)}_{\gamma\in\Gamma}$, then there is a sequence of points $\set{z_k}\subset K$ and a sequence $\set{\gamma_k}\subset\Gamma$ such that $\gamma_k(z_k)\rightarrow
 w$. Since $K$ is compact, there is a subsequence of $\set{z_k}$, $\set{\zeta_j=z_{k_j}}$ converging to a point $\hat{\zeta}\in K$. Thus, $w$ is an accumulation point of $\set{g_j(\zeta_j)}$ where $g_j=\gamma_{k_j}$ and then $w\in D_{\Gamma}(K)$.
\end{proof}

Let $W\subset\CP^n$ be a projective subspace and $z\in\CP^n\setminus W$, let $\mathbf{z}\in\C^{n+1}$ y $\mathbf{W}\subset\C^{n+1}$ be lifts of $z$ and $W$ respectively. Let $\mathbf{\pi_W}:\C^{n+1}\rightarrow \mathbf{W}$ be the orthogonal projection over $\mathbf{W}$. Then we define the \emph{projection} $\pi_{W}$ of $z$ on $W$ as the composition
	$$z\mapsto \pi_{W}(z)=\corchetes{\mathbf{\pi_W}(\mathbf{z})}.$$

Let $A,B\subset\CP^n$ be two projective subspaces and let $\mathbf{A},\mathbf{B}\subset\C^{n+1}$ be lifts of $A$ and $B$. The \emph{projective subspace generated} by $A$ and $B$, denoted by $\prodint{A,B}$, is defined as 
	$$\prodint{A,B}=\corchetes{\prodint{\mathbf{A},\mathbf{B}}}.$$

\begin{defn}\label{defn_simplemente_infinito}
Let $\Gamma\subset \PSLN$ be a discrete subgroup and $\set{g_j}\subset \Gamma$ a sequence of distinct elements. We denote the Cartan decomposition of each element of the sequence as
	$$g_j=\kappa^{(1)}_j \mu(g_j) \kappa^{(2)}_j.$$
  We say that $\set{g_j}$ \emph{tends simply to infinity} if both of the following two conditions hold:
  \begin{enumerate}[(i)]
  \item $\kappa^{(1)}_j\longrightarrow \kappa^{(1)}\in U(n+1)$ y $\kappa^{(2)}_j\rightarrow \kappa^{(2)}\in U(n+1)$ when $j\rightarrow\infty$. 
  \item $\mu(g_j)$ tends to infinity.
  \end{enumerate}
\end{defn}

It is clear that every sequence of distinct elements $\set{g_j}\subset\Gamma$ has a subsequence tending simply to infinity.

\begin{prop}
Let $\Gamma\subset \PSLN$ be a discrete subgroup and $\set{g_j}\subset \Gamma$ a sequence tending simply to infinity. Then there are quasi-projective maps $\tau$ and $\upsilon$ such that 
	$$g_j\rightarrow \tau\text{ and }g^{-1}_j\rightarrow \upsilon,$$
when $j\rightarrow\infty$.
\end{prop}

\begin{proof}
We can write every element of the sequence $\set{g_j}$ as 
	$$g_j=\corchetes{\kappa^{(1)}_j \parentesis{\begin{array}{ccc}
\alpha_{1,k} D_{1,k}& & \\
 & \ddots & \\
 & & \alpha_{m,k} D_{m,k}
\end{array}} \kappa^{(2)}_j}.$$

From Definition \ref{defn_simplemente_infinito} we know that $\kappa^{(1)}_j\longrightarrow \kappa^{(1)}\in U(n+1)$ and $\kappa^{(2)}_j\rightarrow \kappa^{(2)}\in U(n+1)$ when $j\rightarrow\infty$. From Proposition \ref{prop_bloques_tienden_a_matrices_invertibles} it follows that $D_{i,j}$ tends to an invertible matrix $D_i$ when $j\rightarrow\infty$. Then
	$$g_j\rightarrow \tau= \corchetes{\kappa^{(1)} \parentesis{\begin{array}{cccc}
 D_1 & & & \\
 & 0 & & \\ 
 & & \ddots & \\
 & & & 0
\end{array}} \kappa^{(2)}},$$

since

$$\corchetes{\kappa^{(1)}_j \parentesis{\begin{array}{ccc}
\alpha_{1,k} D_{1,k}& & \\
 & \ddots & \\
 & & \alpha_{m,k} D_{m,k}
\end{array}} \kappa^{(2)}_j}=\corchetes{\kappa^{(1)}_j \parentesis{\begin{array}{ccc}
\frac{\alpha_{1,k}}{\alpha_{1,k}} D_{1,k}& & \\
 & \ddots & \\
 & & \frac{\alpha_{m,k}}{\alpha_{1,k}} D_{m,k}
\end{array}} \kappa^{(2)}_j},$$

and by Proposition \ref{prop_cocientes_coeficientes_a_infinito}, 
	$$\frac{\alpha_{i,k}}{\alpha_{1,k}}\rightarrow 0$$
when $j\rightarrow \infty$, for $i=2,...,m$.\\

In the same way 

$$g^{-1}_j\rightarrow \upsilon= \corchetes{\kappa^{(1)} \parentesis{\begin{array}{cccc}
 0 & & & \\
 & \ddots & & \\ 
 & & 0 & \\
 & & & D_m
\end{array}} \kappa^{(2)}}.$$

Clearly $\tau,\upsilon\in\QPN$. 
\end{proof}

\section{The polygon associated to a sequence}

Before we describe the dynamics of sequences in $\Gamma$, we will present a graphic way to visualize their dynamics.\\ 

Consider a sequence $\set{g_k}\subset\Gamma$ separated in $m$ blocks, as in (\ref{eq_forma_bloques}). Consider a polygon with $m$ vertices and weights on each vertex. Each vertex $v_i$ represents a projective subspace $\mathcal{A}_i$. Each edge represents the projective subspace generated by the two subspaces associated to the vertices the edge is joining. In general, given a subset of vertices $\set{v_{i_1},...,v_{i_j}}\subset\set{v_1,...,v_m}$, its convex hull represents the subspace $\prodint{\mathcal{A}_{i_1},...,\mathcal{A}_{i_j}}\subset\CP^n$. In each vertex we put the dimension of the corresponding subspace as its weight. For example, if

$$
g_k=\corchetes{\begin{array}{ccccc}
\boxed{\begin{array}{ccc}
7^{k+1} & & \\
 & 7^k & \\
 & & 7^{k-1}\end{array}} & & & &\\
& \boxed{5^k} & & & \\
& & \boxed{\begin{array}{cc}
 3^{k+1} & \\
  & 3^k\end{array}} & & \\
 & & & \boxed{\begin{array}{cc}
 2^{k+1} & \\
  & 2^k\end{array}} & \\
& & & & \boxed{5^{-k}3^{-2k-1}2^{-2k-1}}
\end{array}}.$$

We associate to $\set{g_k}$ the polygon

\begin{figure}[H]
\begin{center}	
\includegraphics[height=35mm]{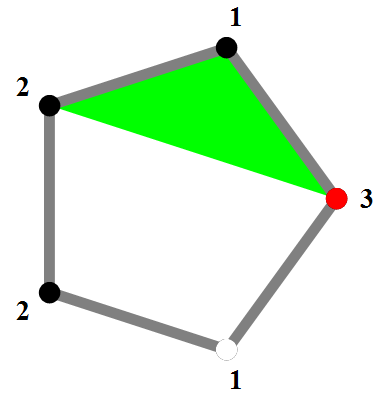}	
\end{center}	
\end{figure}

The red vertex represents the subspace $\mathcal{A}_1=\prodint{e_1,e_2,e_3}$ associated to the first block, the attractive block. The white vertex represents the subspace $\mathcal{A}_5=\set{e_9}$, associated to the last block, the repelling block.\\

The green polygon represents the subspace generated by the subspaces $\mathcal{A}_1$, $\mathcal{A}_2$ and $\mathcal{A}_3$.

\section{Description of the dynamics of sequences}

Let $\set{g_k}\subset\Gamma$ be a sequence. Define the flags $\set{V_i}$ and $\set{W_i}$ in the following way:

\begin{itemize}
\item Let $V_i=\prodint{\mathcal{A}_1,...,\mathcal{A}_i}$ for $i=1,...,m$.
\item Let $W_i=\prodint{\mathcal{A}_{i+1},...,\mathcal{A}_m}\setminus\mathcal{A}_m$ for $i=0,...,m-2$.
\end{itemize}

Let $z=[z_1:...:z_m]\in W_{i-1}\setminus W_i$, and let 
	$$y_z=\corchetes{\lim_{k\rightarrow\infty} D_{i,k}z_i}\in\mathcal{A}_i.$$

We define
	$$V^{(i)}_z=\SET{\corchetes{\zeta_1:...:\zeta_{i-1}:y_z:\overline{0}:...:\overline{0}}}{\zeta_1,...,\zeta_{i-1}\in\C}\subset V_i.$$

\begin{obs}
The set $V^{(i)}_z$ is a union of complex projective lines,
	$$V^{(i)}_z=\bigcup_{\zeta\in V_{i-1}}\overleftrightarrow{\zeta,y_z}.$$
\end{obs}

\begin{obs}
It is clear that 
	$$V^{(i)}_z\subset V_i.$$
\end{obs}

With the following proposition we describe the dynamics of sequences $\set{g_k}\subset\Gamma$ by describing the sets $D_{\set{g_k}}(K)$ for different compact sets $K\subset\CP^n$.

\begin{prop}\label{prop_descripcion_dinamica_frances_compactos}
Let $\set{g_k}\subset\Gamma$ be a sequence.

\begin{enumerate}[(i)]
\item If $z\in W_{i-1}\setminus W_i$, then $D_{\set{g_k}}(z)=V^{(i)}_z$, for $i=1,...,m-2$.
\item If $z\in \mathcal{A}_m$, then $D_{\set{g_k}}(z)=V^{(m)}_z$. 
\end{enumerate}

\end{prop}

\begin{proof}
Let $i\in\set{1,...,m-1}$. Let $\set{g_k}\subset\Gamma$ be a sequence given in the form (\ref{eq_forma_bloques}) and $Z\in W_{i-1}\setminus W_i$. Observe that $Z$ has the form
	$$\corchetes{\overline{0}:...:\overline{0}:Z_i:...:Z_m}$$
with $Z_i\neq \overline{0}$. 

\begin{enumerate}[(i)]
\item We will prove the double set inclusion. Let $w\in D_{\set{g_k}}(Z)$, then there exists a sequence $\set{X_k}\subset \CP^n$, such that $X_k\rightarrow Z$ and $g_k X_k\rightarrow w$.\\

Observe that $g_k X_k$ has the form

\begin{equation}
\begin{split}
\left[\underset{\text{Case 1}}{\underbrace{\frac{\alpha_{1,k}}{\alpha_{i,k}}D_{1,k}X_{1,k}:...:\frac{\alpha_{i-1,k}}{\alpha_{i,k}}D_{i-1,k}X_{i-1,k}}}:D_{i,k}X_{i,k}: \right. \\
\left. :\underset{\text{Case 3}}{\underbrace{\frac{\alpha_{i+1,k}}{\alpha_{i,k}}D_{i+1,k}X_{i+1,k}:...:\frac{\alpha_{m,k}}{\alpha_{i,k}}D_{m,k}X_{m,k}}}\right]
\end{split}
\end{equation}

We have three cases:

\begin{enumerate}[1.]
\item For $j=1,...,i-1$, since $X_{j,k}\rightarrow \overline{0}$ and $\frac{\alpha_{j,k}}{\alpha_{i,k}}\rightarrow \infty$, then $\frac{\alpha_{j,k}}{\alpha_{j,k}}D_{j,k}X_{j,k}$ can converge to any arbitrary $\zeta_j\in\C$.
\item $X_{i,k}\rightarrow Z_i$, then $D_{i,k}X_{i,k}\rightarrow y_z$ by definition.
\item Finally, for $j=i+1,...,m$ it holds $\frac{\alpha_{j,k}}{\alpha_{i,k}}\rightarrow 0$ and $X_{j,k}\rightarrow Z_j$, then 
	$$\frac{\alpha_{j,k}}{\alpha_{j,k}}D_{j,k}X_{j,k}\rightarrow \overline{0}.$$ 
\end{enumerate}

All this means that $w\in V^{(i)}_z$.\\

Now let $w\in V^{(i)}_z$ and suppose it has the form
	$$w=\corchetes{\zeta_1:...:\zeta_{i-1}:y_z:\overline{0}:...:\overline{0}}.$$
Let $i\in\set{1,...,m}$. We want to find a sequence $\set{X_k}\subset\CP^n$ such that $X_k\rightarrow Z$ and $g_k X_k\rightarrow w$. We define

$$X_{j,k}=\begin{cases}
\frac{\alpha_{i,k}}{\alpha_{j,k}}D_{j,k}^{-1}\zeta_j & \text{, } j=1,...,i-1\\
Z_j &\text{, } j=i,...,m.
\end{cases}$$

We can verify directly that this sequence has the desired properties. This means that $V^{(i)}_z\subset D_{\set{g_k}}(z)$ and therefore $D_{\set{g_k}}(z)=V^{(i)}_z$.

\item The proof of this second statement is identical to the proof of the first one.
\end{enumerate}

\end{proof}

\begin{center}
\begin{figure}[h]
\label{fig_ejemplo}
\includegraphics[width=100mm]{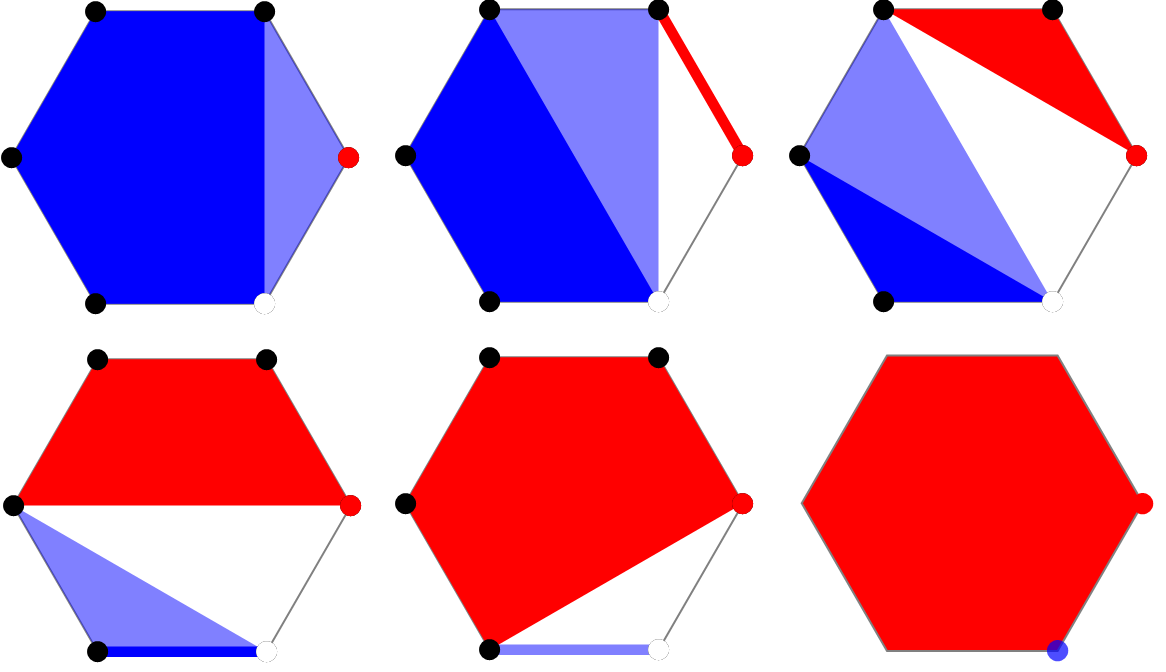}
\caption{Dynamics of $\set{g_k}$ in each element of the flag.}
\end{figure}
\end{center}

In Figure \ref{fig_ejemplo} we can visualize the dynamics described by the previous proposition using the polygons defined in the previous section, in this case $m=6$. In each one of the figures, the subspace $V_i$ is shown in red, the subspace $W_i$ in dark blue and the subspace $W_{i-1}$ in light blue. Proposition \ref{prop_descripcion_dinamica_frances_compactos} states that if we take a compact set $K$ in the light blue region, $D_{\set{g_k}}(K)$ lies inside the red region.\\

\section{The Frances limit set}

Using Proposition \ref{prop_descripcion_dinamica_frances_compactos} we can now define the Frances limit set.

\begin{defn}\label{def_conj_limite_frances}
For $a=\set{g_k}\in\text{Div}(\Gamma)$, let $s_0\in\set{1,...,m}$ such that $e_{\lfloor\frac{n+1}{2}\rfloor}\in \mathcal{A}_{s_0}$. We say that $\mathcal{A}_{s_0}$ is the \emph{middle subspace}. Let
	$$L_a=V_{s_0}$$
where $\set{V_s}$ is the flag defined in Proposition \ref{prop_descripcion_dinamica_frances_compactos}. We define the \emph{Frances limit set} of $\Gamma$, denoted by $\FraL(\Gamma)$ as
	$$\FraL(\Gamma)=\overline{\bigcup_{a\in\text{Div}(\Gamma)}L_a}.$$
The complement is denoted by
	$$\FraD(\Gamma)=\CP^n\setminus\FraL(\Gamma).$$
\end{defn}

The following theorem shows that the Frances limit set it is, in fact, a \emph{good} notion of limit set for the action of a discrete subgroup of $\PSLN$.

\begin{thm}\label{prop_frances_prop_disc}
The action of $\Gamma$ on $\FraD(\Gamma)$ is proper and discontinuous. 
\end{thm}

\begin{proof}
Since $\FraL(\Gamma)$ is closed, $\FraD(\Gamma)$ is open. In order to show that $\Gamma$ acts properly and discontinuously on $\FraD(\Gamma)$ we will show that there cannot be two points $z,w\in\FraD(\Gamma)$ such that $z\in D_\Gamma(w)$.\\

Assume that there exists such pair of points $z,w\in \FraD(\Gamma)$ with $z\in D_\Gamma(w)$. Since $D_\Gamma(w)=\cup_{a\in \text{Div}(\Gamma)}D_a(w)$ then $z\in D_b(w)$ for some sequence $b\in \text{Div}(\Gamma)$.\\ 

Besides, for any sequence $a=\set{g_k}\in \text{Div}(\Gamma)$, it holds that $a^{-1}=\set{g^{-1}_k}\in \text{Div}(\Gamma)$. On the other hand, by definition we have
	$$\FraD(\Gamma)=\bigcap_{a\in \text{Div}(\Gamma)}\parentesis{\CP^n\setminus L_a}.$$
Since $z,w\in \FraD(\Gamma)$,  
	\begin{equation}\label{eq_dem_propydisc_1}
	z,w\in\bigcap_{a\in \text{Div}(\Gamma)}\parentesis{\CP^n\setminus L_a}\subset \parentesis{\CP^n\setminus L_b}\cap\parentesis{\CP^n\setminus L_{b^{-1}}} \subset \CP^n\setminus\parentesis{L_b\cup L_{b^{-1}}}.
	\end{equation}
	
Observe that, for the sequences $b$ and $b^{-1}$, the decomposition in the spaces $\mathcal{A}_1,...,\mathcal{A}_m$ is the same, just in the opposite order. Therefore, in both cases the middle space is the same.\\

Denote $E=\CP^n\setminus\parentesis{L_b\cup L_{b^{-1}}}$. Then 
	\begin{equation}\label{eq_dem_propydisc_2}
	E=\parentesis{E\cap W_0\setminus W_1}\cup...\cup\parentesis{E\cap W_{s_0-2}\setminus W_{s_0-1}}.
	\end{equation}
To verify this, we prove the double set inclusion. It is clear that 
	$$\parentesis{E\cap W_0\setminus W_1}\cup...\cup\parentesis{E\cap W_{s_0-2}\setminus W_{s_0-1}}\subset E.$$ 
To verify the opposite inclusion take $Z=[z_1:...:z_m]\in E$. There is an index $j\in\set{1,...,s_0-1}$ such that $z_j\neq\overline{0}$, otherwise 
	$$Z=\corchetes{\overline{0}:...:\overline{0}:Z_{s_0}:...:Z_m}$$
and then $Z\in L_{b^{-1}}$, contradicting that $Z\in E$. Then $Z\in W_{j-1}\setminus W_j$, which proves (\ref{eq_dem_propydisc_2}).\\

From (\ref{eq_dem_propydisc_2}) we know there are two indices $j_1,j_2$ such that 
	\begin{equation}\label{eq_dem_propydisc_4}
	z\in W_{j_1-1}\setminus W_{j_1}\text{, } w\in W_{j_2-1}\setminus W_{j_2}.
	\end{equation}		
Now we verify that $j_1,j_2<s_0$. In order to do this, assume that 
	\begin{equation}\label{eq_dem_propydisc_3}
	j_1\geq s_0.
	\end{equation}
By hypothesis, $z\in \CP^n\setminus\parentesis{L_b\cup L_{b^{-1}}}$ and then $z\nin L_{b^{-1}}$, but 
	$$L_{b^{-1}}=\prodint{\mathcal{A}_{s_0},...,\mathcal{A}_m}=W_{s_0-1}.$$
Then $z\nin W_{s_0-1}$.\\
On the other hand, from (\ref{eq_dem_propydisc_4}) it follows that $z\in W_{j_1-1}$. From (\ref{eq_dem_propydisc_3}) we have $j_1-1\geq s_0-1$, and then
	$$W_{j_1-1}\subset W_{s_0-1}$$
which contradicts that $z\in W_{j_1-1}$ y $z\nin W_{s_0-1}$. This verifies that $j_1<s_0$. Using the same argument for $j_2$ it follows that $j_2<s_0$.\\

Now, since $z\in E=\CP^n\setminus\parentesis{L_b\cup L_{b^{-1}}}$, it follows that 

\begin{equation}\label{eq_dem_propydisc_5}
z\nin L_b=\prodint{\mathcal{A}_1,...,\mathcal{A}_{s_0}}=V_{s_0}.
\end{equation}

Applying Proposition \ref{prop_descripcion_dinamica_frances_compactos} in (\ref{eq_dem_propydisc_4}) it follows
	\begin{equation}\label{eq_dem_propydisc_6}
	D_\Gamma(w)\subset V_{j_2}.
	\end{equation}
By hypothesis, $z\in D_\Gamma(w)$ and then, from (\ref{eq_dem_propydisc_6}) we have that 
	\begin{equation}\label{eq_dem_propydisc_7}
	z\in V_{j_2}.
	\end{equation} 

Since $j_2<s_0$, it holds
	$$V_{j_2}\subset V_{s_0}$$
which contradicts (\ref{eq_dem_propydisc_5}) and (\ref{eq_dem_propydisc_7}). This contradiction proves that the action of $\Gamma$ on $\FraD(\Gamma)$ is proper and discontinous.
\end{proof}

As a consequence of the previous proposition, we call $\FraD(\Gamma)$ the \emph{Frances region of discontinuity}.\\

The following examples describe some relations between the Frances limit set and other limit sets.

\begin{ejem}
The Frances limit set is smaller than the Kulkarni limit set in general. This immediately follows from the definition of the Frances limit set. Each point in the Frances limit set is in the closure of accumulation points of orbits of compact sets and therefore, lies in the Kulkarni limit set. Therefore,
	$$\FraL(\Gamma)\subset\KulL(\Gamma),$$
this, together with Proposition \ref{prop_eq_in_kuld}, yields
	$$\Eq(\Gamma)\subset\KulD(\Gamma)\subset\FraD(\Gamma).$$
This is one of the advantages of the Frances limit set over the Kulkarni limit set, especially in complex dimension greater than 2.
\end{ejem}

\begin{ejem}
If $\Gamma\subset\PSL$ is strongly irreducible then
	$$\FraL(\Gamma)=\KulL(\Gamma).$$
\end{ejem}

\begin{ejem}
If $\Gamma\subset\PSL$ is not strongly irreducible then
	$$\FraD(\Gamma)=\Eq(\Gamma).$$
\end{ejem}

The following theorem directly follows from the definition of $\FraL(\Gamma)$ and Proposition \ref{prop_frances_prop_disc}.

\begin{thm}
Let $\Gamma\subset \PSLN$ be a strongly irreducible complex Kleinian group, then $\FraL(\Gamma)$ is the smallest limit set such that the group acts properly and discontinuously on the complement.
\end{thm}

Proposition \ref{prop_frances_CG} gives a relation between the Frances limit set and the Chen-Greenberg limit set for discrete subgroups of $\text{PU}(1,n)$. Before we can state and prove this proposition, we will need some definitions and results.\\

In the following we will consider $\C^{n+1}$ with the hermitian metric given by 
$$\prodint{v,w}=v_1\overline{w_{n+1}}+v_{n+1}\overline{w_{1}}+ \sum_{j=2}^{n}v_j\overline{w_j}$$
clearly this metric has signature (1,n), so in the following we will consider
$$U(1,n)=\SET{g\in \GLN}{\prodint{gu,gv}=\prodint{u,v} \textrm{ for any } u,v\in \C^{n+1}}.$$

Finally, we denote by $\text{PU}(1,n)$ to the projectivization of $\text{U}(1,n)$, this projectivization preserves the unitary complex ball 
	$$\Hh_\C^n=\SET{\corchetes{w}\in\CP^n}{\prodint{w,w}<0}.$$	  

We have the following version of the Cartan decomposition

\begin{thm}[Cartan decomposition]\label{thm_cartan_alternativo}
There is a compact  group $K$ of $\text{PU}(1,n)$ such that for every $\gamma\in U(1,n)$ there are elements $k_1,k_2\in K$ and  a unique  element $\mu(\gamma)\in U(1,n)$, called the Cartan projection,  such that $\gamma=k_1\mu(\gamma)k_2$ and
$$
\mu(\gamma)=
\left (
\begin{array}{lllll}
e^{\alpha}\\
                 & 1\\
                 &          & \ddots\\
                 &          &            &1\\
                 &          &            &     & e^{-\alpha}
\end{array} 
\right ).$$
\end{thm}

In order to provide a description of the Frances limit set for discrete subgroups of $\text{PU}(1,n)$ the following easy lemma will be useful.

\begin{lem} [See \cite{cano2010equicontinuity}] \label{l:eqcs}
Let $\Gamma\subset PU(1,n)$ be a discrete group, then if $\set{\gamma_m}\subset \Gamma$ is a sequence of distinct elements such that $\set{\gamma_m}$ converges to a pseudo-projective transformation $\tau$ and $\set{\gamma_m^{-1}}$ converges to a pseudo-projective transformation $\theta$, then $\ker(\tau)$ and $\ker(\theta)$ are hyperplanes tangent to $\partial \Hh^n_\C$ at points in the Chen-Greenberg limit set of $\Gamma$.
\end{lem}

\begin{prop}\label{prop_frances_CG}
Let $\Gamma\subset PU(1,n)$ be a discrete group, then the Frances limit set for $\Gamma$ are the hyperplanes tangent to $\partial \Hh^n_\C$ at points in the Chen-Greenberg limit set of $\Gamma$, {\it i. e.}

$$\FraL(\Gamma)=\bigcup_{p\in \Lambda_{CG}(\Gamma)}p^\bot.$$
 
\end{prop}

\begin{proof}
Let $\set{\gamma_m}\subset \Gamma$ be a sequence of distinct elements, since the space of pseudo-projective elements is compact, we can assume that there is a pseudo-projective transformations $\gamma$ and $\vartheta$ such that $\set{\gamma_m}$ converges to $\gamma$ and $\set{\gamma_m^{-1}}$ converges to $\vartheta$. On the other hand, by the Cartan decomposition theorem (Theorem \ref{thm_cartan_alternativo}) we can assume that there are convergent sequences $\set{\alpha_m}\in \R^+$, $\set{\kappa_m},\set{\tau_m}\in K$, such that 
\begin{equation} \label{e:dec}
\gamma=
\corchetes{
\kappa_m
\left (
\begin{array}{lllll}
e^{\alpha_m}\\
                 & 1\\
                 &          & \ddots\\
                 &          &            &1\\
                 &          &            &     & e^{-\alpha_m}
\end{array}
\right )
\tau_m}.
\end{equation}
If $\set{\kappa_m}$ converges to $\kappa$ and $\set{\tau_m}$ converges to $\tau$, then Equation \ref{e:dec} shows that   the projective subspaces in the Frances limit set introduced by the sequence $\set{\gamma_m}$ are $\mathcal{H}_1=\kappa \corchetes{\prodint{e_1,\ldots ,e_n }}$ and  $\mathcal{H}_2=\tau^{-1} \corchetes{\prodint{e_2,\ldots,e_{n+1}}}$. A straight-forward calculation shows that

$$\gamma=
\corchetes{\kappa
\parentesis{\begin{array}{lllll}
1   \\
                 & 0\\
                 &          & \ddots\\
                 &          &            &0\\
\end{array}}\tau};\;\;\;\; 
\vartheta=
\corchetes{
\tau^{-1}
\parentesis{\begin{array}{lllll}
 0\\
                           & \ddots\\
                           &       & 0\\
                           &       &            &1\\
\end{array}}\kappa^{-1}},$$

therefore $\ker(\gamma)=\tau^{-1}\corchetes{\prodint{e_2,\ldots,e_{n+1}}}$ and  $\ker(\vartheta)=\kappa\corchetes{\prodint{e_1,\ldots,e_{n}}}$. Using Lemma \ref{l:eqcs}, the proof is complete.
\end{proof}

\section{Purely dimensional sets}

\begin{defn}\label{defn_puramente_dimensional}
Let $\Gamma\subset \PSLN$ be a complex Kleinian group and $F\subset\CP^n$ a closed subset. Let $k>0$, we say that $F$ is \emph{purely $k$-dimensional} if $\exists$ $W\in Gr(k,n)$ such that 
	\begin{enumerate}
	\item $\overline{\Gamma(W)}=F$, and
	\item For any $k'>k$, there is no $W'\in Gr(k',n)$ such that $W'\subset F$. 
	\end{enumerate}	   
We call $k$ the \emph{pure dimension} of $\FraL(\Gamma)$.
\end{defn}

The previous definition describes the idea of subsets of $\CP^n$ made up only of projective subspaces of the same dimension. The Frances limit set will be a purely dimensional set under certain hypothesis (see Theorem \ref{thm_frances_puramente_dimensional}).

\begin{ejem}
If $\Gamma\subset\PSL$ is a Kleinian complex group with 4 lines in general position in $\KulL(\Gamma)$, then $\KulL(\Gamma)$ is purely 1-dimensional (see Theorem 6.3.3 of \cite{ckg_libro}).
\end{ejem}

\begin{thm}\label{thm_frances_puramente_dimensional}
If $\Gamma\subset\PSLN$ is a strongly irreducible discrete subgroup then $\FraL(\Gamma)$ is purely $k$-dimensional, where $k$ satisfies
	$$k\geq \lfloor\frac{n}{2}\rfloor.$$
\end{thm}

\begin{proof}
First we prove that, if $\FraL(\Gamma)$ is purely $k$-dimensional with 
	$$k\geq \lfloor\frac{n}{2}\rfloor.$$
	
From Definitions \ref{def_conj_limite_frances} and \ref{defn_puramente_dimensional} it follows that, if $\FraL(\Gamma)$ is purely $k$-dimensional then 
	\begin{equation}\label{eq_dem_puramente_dimensional_1}
	k\geq \text{dim}(L_a)
	\end{equation} 
for any sequence $a\in \text{Div}(\Gamma)$. If this doesn't hold, there is a sequence $\tilde{a}\in\text{Div}(\Gamma)$ such that $k<\text{dim}(L_{\tilde{a}})$. $L_{\tilde{a}}$ is a projective subspace of $\CP^n$ and thus $L_{\tilde{a}}\in Gr(\text{dim}(L_{\tilde{a}}),n)$; besides $L_{\tilde{a}}\subset \FraL(\Gamma)$. This contradicts that $\FraL(\Gamma)$ is purely $k$-dimensional.\\

On the other hand,
	\begin{equation}\label{eq_dem_puramente_dimensional_2}
	\text{dim}(L_{a})\geq \floor{\frac{n}{2}},
	\end{equation}

for any sequence $a\in\text{Div}(\Gamma)$. This follows from the fact that $e_{\floor{\frac{n+1}{2}}}\in L_a$.\\

Combining (\ref{eq_dem_puramente_dimensional_1}) and (\ref{eq_dem_puramente_dimensional_2}) it follows
	\begin{equation}\label{eq_dem_puramente_dimensional_3}
	k\geq \lfloor\frac{n}{2}\rfloor.
	\end{equation}		
	
We now prove that $\FraL(\Gamma)$ is purely $k$-dimensional, where $k$ satisfies (\ref{eq_dem_puramente_dimensional_3}). Let $a_0\in\text{Div}(\Gamma)$ such that $\text{dim}(L_{a_0})\geq \text{dim}(L_a)$ for any $a\in\text{Div}(\Gamma)$, this is possible since $\text{dim}(L_a)\leq n$. As we showed before, $k:=\text{dim}(L_{a_0})\geq \lfloor\frac{n}{2}\rfloor$. It is straight-forward to verify that (ii) of Definition \ref{defn_puramente_dimensional} holds. On the other hand, since $\Gamma$ is strongly irreducible $L_{a_0}$ doesn't have a finite $\Gamma$-orbit. This, together with the fact that $\gamma L_a = L_{\gamma a}$ for $a\in\text{Div}(\Gamma)$ and $\gamma\in\Gamma$, proves (i) of Definition \ref{defn_puramente_dimensional}.  
\end{proof}



With the following example we state that the Frances limit set is not stable under deformation.

\begin{ejem} For $0\leq\epsilon<1$, consider the element

\begin{equation}\label{eq_ej_no_estable_1}
A_\epsilon=
\corchetes{\begin{array}{cccc}
2& 0 &0 & 0\\
0 & 1+\epsilon^2 & 0 &0 \\
0 & 0 & (1+\epsilon^2)^{-1}&0\\ 
0 & 0& 0& 2^{-1}
\end{array}}.
\end{equation}

Then the family of groups $\Gamma_\epsilon:=\prodint{ A_\epsilon}$ show that the Frances limit set cannot be stable under deformations, since 
	$$\FraL\parentesis{\Gamma_\epsilon}=\begin{cases}
	\prodint{e_1,e_2}, & 0<\epsilon<1\\
	\prodint{e_1,e_2,e_3}, & \epsilon=0  	
	\end{cases}.$$

We can also give a similar example for strongly irreducible groups. Recalling (\ref{eq_ej_no_estable_1}), observe that $A_0\in \text{PU}(1,n)$. Let $b\in \text{PU}(1,n)$ such that
	$$b(e_1^\perp \cap e_4^\perp)\cap e_1^\perp \cap  e_4^\perp=\emptyset.$$
Denote $\Gamma_{\epsilon,k}=\prodint{ A_\epsilon^k, bA_\epsilon^k b^{-1}}$, then for a large enough $k$, $\Gamma_{0,k}$ is a strongly irreducible Schottky subgroup of $\text{PU}(1,n)$ and therefore the Frances limit set is a collection of lines tangent to $\Hh_\C^3$. For a small enough $\epsilon>0$, $\Gamma_{\epsilon,k}$ is a strongly irreducible Schottky-like group. Since the Cartan decomposition varies smoothly for smooth deformations, we conclude that the Frances limit set of $\Gamma_{\epsilon,k}$ for $\epsilon>0$ is a Cantor set of lines.
\end{ejem} 
\chapter{Measures on limit sets of complex Kleinian groups}
\label{chapter_measures}

Patterson-Sullivan measures are a family of probability measures associated to discrete groups $\Gamma\subset \psl$ (or in general, subgroups of isometries of the hyperbolic real space $\Hh^n_\R$), they are supported on the limit set of $\Gamma$ and are quasi-invariant under the action of $\Gamma$. These measures give the proportion of elements of any orbit that accumulates in a given region of $\Ss^2$ (or in general, $\Ss^{n-1}$).\\

We want to define similar measures for complex Kleinian groups acting on $\CP^2$. In this chapter we lay the foundations upon which we hope these measures can be constructed. Therefore, we present some partials results and give some ideas about the program we think should be followed. Because of the nature of this chapter, the exposition will be less rigorous.\\ 

In Section \ref{sec_parameter_space} we will describe the space of arrays of 5 complex projective lines in general position in $\CP^2$. Each of these configurations of lines determines a domain in $\CP^2$ (the complement in $\CP^2$ of the array), the existence of the measure ultimately depends on whether the entropy volume of the Kobaya-shi metric is finite on this domain.\\

In Section \ref{sec_volumen_entropia} we propose a way to estimate the entropy volume for the Kobayashi metric. This part of the work is not complete, however several advances are shown.\\

Finally, in Section \ref{sec_construction} we describe the needed steps to construct the Patterson-Sullivan measures, once we assume that the entropy volume of the Kobayashi (or Eisenman) metric is finite. We only remark some steps that are important to the construction in the case of complex dimension two, the remaining steps are the same as in the case of complex dimension 1 (see \cite{patterson} and \cite{nicholls} for these details).\\

\section{The space of configurations of 4 lines in general position in $\C^2$}
\label{sec_parameter_space}

In this section we want to determine the configurations of 4 complex lines in general position in $\C^2$ for which the complement has finite entropy volume with respect to the Kobayashi metric. The first step to achieve this is to parame-trize the space configurations of 4 complex lines in general position in $\C^2$. Once this space is described, we want to determine the arrays of lines such that any group acting properly and discontinuously on the complement of the array admits a measure supported on the limit set.\\

Before we work on the parametrization of this space we give some definitions and properties of complex lines in $\C^2$. Following the notation of \cite{shabat}, we define a complex line in the following way.

\begin{defn}\label{defn_linea_compleja}
A complex line $\ell$ in $\C^2$ is a subset
	$$\ell = \SET{Z\in\C^2}{Z=Z_0+\xi W_0\text{, donde }\xi\in\C}$$
where $Z_0\in\C^2$ is a point in $\ell$ and $W_0\in\C^2\setminus{(0,0)}$ is the \emph{complex direction} of $\ell$.
\end{defn}

In the previous definition is clear that, if $W_0,W'_0\in\C^2$ are points such that $W_1=\lambda W_2$ for some $\lambda\in\C\setminus\set{0}$ then the complex lines $\ell$ and $\ell'$ defined by the parametric equations $Z_0+\xi W_0$ and $Z_0+\xi W'_0$ respectively, are the same set. In other words, given a complex line, its complex direction is uniquely determined up to multiplication by a non-zero complex number.\\

For a complex line $\ell$, we denote its parametric equation by the same symbol
	$$\ell(\xi)=Z+\xi W,$$
if there is no ambiguity.

\begin{obs}\label{obs_relacion_lineas_complejas}
In Definition \ref{defn_linea} we define a \emph{complex projective line}, meanwhile in the previous Definition \ref{defn_linea_compleja} we use the term \emph{complex line}. However, there is a relationship between both concepts.\\

A complex projective line is a complex line in $\C^2$ compactified with a point lying in the intersection of the complex line and the complex projective line at infinity. Equivalently, a complex line is an affine part of a complex projective line.
\end{obs}

\begin{defn}
Let $\ell_1$ and $\ell_2$ be two complex lines with complex directions $W_1$ and $W_2$. If $W_1=\lambda W_2$ for some $\lambda\in\C^\ast$ we say that $\ell_1$ and $\ell_2$ have the same complex direction.
\end{defn}

The following proposition states that, unlike the case of complex projectives lines, not all complex lines have non-empty intersection.

\begin{prop}\label{prop_interseccion_lineas_complejas}
Let $\ell_1$ and $\ell_2$ two complex lines with different complex directions, then $\ell_1\cap\ell_2$ is exactly one point in $\C^2$.
\end{prop}

\begin{proof}
Suppose that the parametric equation of the complex line $\ell_i$ is given by 
	$$\ell_i(\zeta)=Z_i+\zeta W_i$$
for $i=1,2$, where $Z_i=(z^i_1,z^i_2)$ and $W_i=(w^i_1,w^i_2)$.\\

The intersection of both complex lines are the points that can be parametrized by $\zeta,\xi\in\C$ such that $\ell_1(\zeta)=\ell_2(\xi)$. In order for a point $\ell_1(\zeta)=\ell_2(\xi)$ to be in the intersection, it has to satisfy the equation
	\begin{equation}\label{eq_dem_interseccion_lineas_complejas_1}		
	\parentesis{\begin{array}{cc}
	w^1_1 & -w^2_1\\
	w^1_2 & -w^2_2
	\end{array}}\parentesis{\begin{array}{c}
	\zeta\\
	\xi
	\end{array}}=\parentesis{\begin{array}{c}
	z^2_1-z^1_1\\
	z^2_2-z^1_2
	\end{array}}.
	\end{equation}
Let 
	$$A=\parentesis{\begin{array}{cc}
	w^1_1 & -w^2_1\\
	w^1_2 & -w^2_2
	\end{array}}.$$
The equation (\ref{eq_dem_interseccion_lineas_complejas_1}) has a unique solution if $\det(A)\neq 0$. Suppose that the intersection $\ell_1\cap\ell_2$ is not a single point, then $\det(A)=0$, but this only happens if 
	\begin{equation}\label{eq_dem_interseccion_lineas_complejas_2}
	w^2_1 w^1_2 = w^1_1 w^2_2.
	\end{equation}
From (\ref{eq_dem_interseccion_lineas_complejas_2}) we see that if $w^2_1=0$ then either $w^1_1=0$ or $w^2_2=0$. If $w^2_2=0$ then $W_2=(0,0)$, contradicting that $W_2$ is the complex direction of a complex line. If $w^1_1=0$ then $W_1=(0,w^1_2)$ and $W_2=(0,w^2_2)$ are complex multiples, contradicting the hypothesis. Therefore $w^2_1\neq 0$ and thus,
	\begin{equation}\label{eq_dem_interseccion_lineas_complejas_3}
	w^1_2=\frac{w^1_1}{w^2_1}w^2_2.
	\end{equation}    
If $w^1_1=0$, from (\ref{eq_dem_interseccion_lineas_complejas_2}) it follows that $w^1_2=0$ since $w^2_1\neq 0$ and then $W_1=(0,0)$, thus $w^1_1\neq 0$. Then $\lambda:=\frac{w^1_1}{w^2_1}\in\C\setminus\set{0}$, furthermore
	\begin{align*}
	w^1_2 &= \lambda w^2_2\\
	w^1_1 &= \lambda w^2_1.
	\end{align*}
This means that $W_1=\lambda W_2$ contradicting that $\ell_1$ and $\ell_2$ have different complex direction. This proves that $\ell_1\cap\ell_2$ is exactly one point in $\C^2$.
\end{proof}

\begin{prop}\label{prop_unica_linea_compleja_por_dos_puntos}
Let $z,w\in\C^2$, there exists a unique complex line containing both $z$ and $w$.
\end{prop}

\begin{proof}
We know there is a unique projective complex line $\hat{\ell}$ containing the points $z$ and $w$. There is some affine part of $\hat{\ell}$ containing both $z$ and $w$.\\

If there were two different complex lines $\ell_1$ and $\ell_2$ containing $z$ and $w$, then their compactifications are two projective complex lines $\hat{\ell}_1$ and $\hat{\ell}_2$ containing the points $z$ and $w$, which contradicts the fact that there is only one projective complex line containing those points.
\end{proof}

As a consequence of the previous proposition, given two points $z,w\in\C^2$, we can denote by
	$$\ell=\overleftrightarrow{z,w}$$
to the complex line containing both points. This notation is similar to the notation used in the case of projective complex lines, we will use it for both cases as long as there is no ambiguity.

\begin{obs}\label{obs_posicion_gral_no_paralelas}
Let $\mathcal{L}$ be an array of complex projective lines in general position in $\CP^2$, we assume that one of these lines is the line at infinity $\ell_\infty$. Consider the subarray $\mathcal{L}'=\mathcal{L}\setminus\set{\ell_\infty}$, seen as a collection of complex lines in $\C^2$. There cannot be a pair of lines with the same complex direction $\mathcal{L}'$, otherwise such pair of lines would intersect at a point belonging to $\ell_\infty$, contradicting that the array is in general position in $\CP^2$.
\end{obs}

The previous observation states that two lines with the same complex direction cannot be in general position in $\CP^2$, together with the line at infinity.

\subsection{The parameter space}

In this subsection we determine the parameter space of the configurations of 4 complex lines in general position in $\C^2$.\\

In order to so, we will prove that $\PSL$ acts 4-transitively on the space of complex projective lines in general position in $\CP^2$. With this in mind, given an array of 5 complex projective lines in general position, we can arbitrarily pick 4 of them and describe the sought parameter space as the parameter space of configurations of 4 complex lines in general position.\\

Equivalently, choosing one of the 4 complex projective lines as the line at infinity, we can parameterize the configurations of 4 complex lines in general position as the space of all complex lines in general position with the 3 previously chosen complex lines.\\  

The following is a well known result.

\begin{thm}\label{thm_4_transitividad_puntos}
For any pair of $n+2$ points $p_1,...,p_{n+2}$ and $q_1,...,q_{n+2}$ in $\CP^n$ there is a unique transformation $\gamma\in\PSLN$ such that $\gamma(p_i)=q_i$ for $i=1,...,n+2$, assuming no $n+1$ points $p_i$ lay in the same hyperplane in $\CP^n$.
\end{thm}

The following proposition can be directly verified.

\begin{prop}
Any element of $\PSLN$ takes projective subspaces of dimension $k$ of $\CP^n$ in projective subspaces of dimension $k$ of $\CP^n$.
\end{prop}

The following theorem is proven in \cite{green} for the general case of projective hyperplanes in $\CP^n$.

\begin{thm}\label{thm_complemento_hiperplanos}
Let $\ell_1,...,\ell_5$ be 5 projective complex lines in general position, then $\CP^2\setminus\parentesis{\ell_1\cup...\cup\ell_5}$ is complete Kobayashi hyperbolic. 
\end{thm}

Using the argument in Observation \ref{obs_relacion_lineas_complejas} we can restate the previous theorem in the following way: If $\ell_1,...,\ell_4$ are 4 complex lines in general position then $\C^2\setminus\parentesis{\ell_1\cup...\cup\ell_4}$ is complete Kobayashi hyperbolic.

\begin{prop}
The group $\PSL$ acts sharply 4-transitively on the space of complex projective lines in general position in $\CP^2$.
\end{prop}

\begin{proof}
Let $\ell^1_1,...,\ell^1_4$ and $\ell^2_1,...,\ell^2_4$ be two arrays of 4 complex lines in general position in $\C^2$. From Observation \ref{obs_posicion_gral_no_paralelas} we know that neither array contain a pair of lines with the same complex direction and therefore, as a consequence of Proposition \ref{prop_interseccion_lineas_complejas}, each pair of lines intersect. Let us denote then 
\begin{equation}\label{eq_dem_4_transitivo_lineas_1}
\begin{array}{cccc}
p^i_1=\ell^i_1\cap \ell^i_3 & p^i_2=\ell^i_2\cap \ell^i_3 & p^i_3=\ell^i_2\cap \ell^i_4 & p^i_4=\ell^i_1\cap \ell^i_4\\
\end{array}
\end{equation}
for $i=1,2$. Theorem \ref{thm_4_transitividad_puntos} implies that there is a unique element $\gamma\in\PSL$ such that 
	$$p^1_i\overset{\gamma}{\mapsto}p^2_i.$$
for $i=1,...,4$. Observe that $\ell^i_1=\overleftrightarrow{p^i_1,p^i_4}$, this is because (\ref{eq_dem_4_transitivo_lineas_1}) implies that $p^i_1,p^i_4\in\ell^i_1$ and there is a unique line containing both points. In the same way, it follows that 
$$
\begin{array}{ccc}
\ell^i_2=\overleftrightarrow{p^i_2,p^i_3}, & \ell^i_3=\overleftrightarrow{p^i_1,p^i_2}, & \ell^i_4=\overleftrightarrow{p^i_3,p^i_4}
\end{array}
$$
for $i=1,2$. Then 
	\begin{equation}\label{eq_dem_4_transitivo_lineas_2}
	\gamma(\ell^1_j)=\ell^2_j
	\end{equation} 
for $j=1,...,4$.\\

Finally, we verify the uniqueness of $\gamma$. If there was another $\mu\in\PSL$ satisfying (\ref{eq_dem_4_transitivo_lineas_2}) then $\mu$ should take the intersections $p^1_i$ to the intersections $p^2_i$. Then Theorem \ref{thm_4_transitividad_puntos} implies $\mu=\gamma$.
\end{proof}

The previous result states that we can choose 3 arbritrary complex lines out of the 4 lines in general position we are considering in $\C^2$. For the rest of the chapter we choose 3 complex lines $\ell_0,\ell_1,\ell_2$ and three points $v_0,v_1,v_2$ in $\C^2$ in the following way:

$$\begin{array}{ccc}
v_0=(0,0) & v_1=(1,0) & v_2=(0,1)\\
 & \text{and} & \\
\ell_0=\overleftrightarrow{v_0,v_1} & \ell_1=\overleftrightarrow{v_1,v_2} & \ell_2=\overleftrightarrow{v_2,v_0}.\\ 
\end{array}$$

We will work with the following parametrizations for the 3 complex lines we have chosen:

	\begin{equation}\label{eq_defn_lineas_escogidas}
	\begin{array}{cl}
	\ell_0: & \ell_0(\xi)=\xi(1,0)\\
	\ell_1: & \ell_1(\xi)=(1,0)+\xi(-1,1)\\
	\ell_2: & \ell_2(\xi)=\xi(0,1).
	\end{array}	
	\end{equation}		

\begin{figure}[H]
\begin{center}	
\includegraphics[height=45mm]{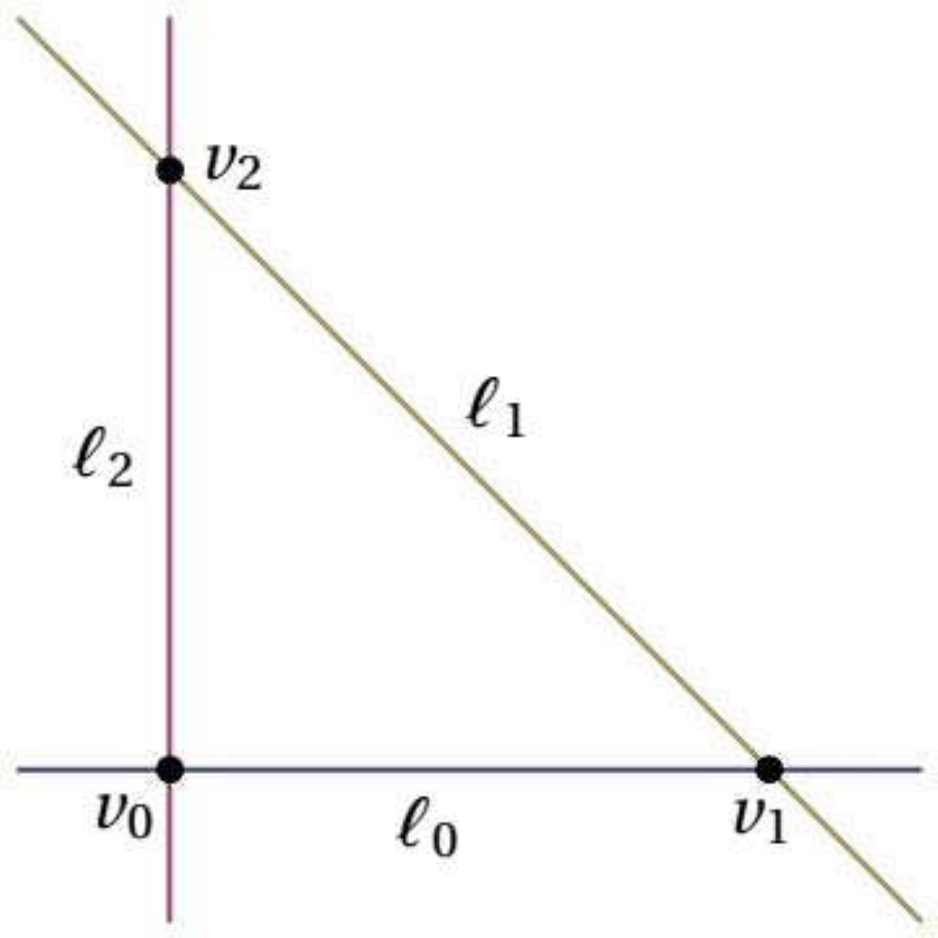}	
\end{center}	
\end{figure}

If we see the lines $\ell_0$, $\ell_1$ y $\ell_2$ as complex projective lines, then we will consider the following parametrizations:
	\begin{align*}
	\ell_0&=\SET{\corchetes{\alpha:0:\beta}}{(\alpha,\beta)\in\C^2\setminus\set{(0,0)}}\\
	\ell_1&=\SET{\corchetes{\beta-\alpha:\alpha:\beta}}{(\alpha,\beta)\in\C^2\setminus\set{(0,0)}}\\
	\ell_2&=\SET{\corchetes{0:\alpha:\beta}}{(\alpha,\beta)\in\C^2\setminus\set{(0,0)}}.
	\end{align*}	
We have to keep in mind that we have chosen a fourth line $\ell_4$, the line at infinity. We will consider the following parametrization for $\ell_4$:
	$$\ell_4=\SET{\corchetes{\alpha:\beta:0}}{(\alpha,\beta)\in\C^2\setminus\set{(0,0)}}$$
In order to describe each complex line $\ell$ in general position with the 3 chosen lines $\ell_0$, $\ell_1$, $\ell_2$, we will use two points in $\ell$, the first one will be $\ell\cap\ell_0$ and the second, $\ell\cap\ell_2$. Each of these points is determined by a complex parameter. Therefore, the space of all complex lines in general position with the 3 previously chosen lines is a subset of $\C^2$. This space will be described in Proposition \ref{prop_descripcion_parametros}.\\

In Proposition \ref{prop_descripcion_parametros}, $\zeta_i\in\C$ denotes the complex parameter determining the point $z_i:=\ell\cap\ell_{i-1}$, for $i=1,2$.\\

For convention, the complex line $\ell_3$ determined by the parameter $\Upsilon=(\zeta_1,\zeta_2)\in\C^2$ will have the parametrization
	
	$$\ell^{\Upsilon}_3(\xi)=\ell_0(\zeta_1)+\xi\parentesis{\ell_2(\zeta_2)-\ell_0(\zeta_1)}.$$  

Using the parametrizations for $\ell_0$ y $\ell_2$ given in (\ref{eq_defn_lineas_escogidas}), it follows

	$$\ell^{\Upsilon}_3(\xi)=\parentesis{(1-\xi)\zeta_1,\xi\zeta_2}.$$
	
Considering $\ell^{\Upsilon}_3$ as a complex projective line, its parametrization is given by
	
	$$\ell^{\Upsilon}_3=\SET{\corchetes{(\beta-\alpha)\zeta_1:\alpha\zeta_2:\beta}}{(\alpha,\beta)\in\C^2\setminus\set{(0,0)}}.$$

As long as there is no confusion, we will omit the dependence on $\Upsilon$ and just write $\ell_3$ instead of $\ell_3^\Upsilon$.

\begin{prop}\label{prop_descripcion_parametros}
The space of arrays of 4 complex lines in general position in $\C^2$ is given by 
$$
\mathcal{P}=\SET{(\zeta_1,\zeta_2)\in\C^2}{\zeta_1\in\C\setminus\set{0,1}\text{, }\zeta_2\in \C\setminus\set{0,1,\zeta_1}}.$$
\end{prop}

\begin{proof}
We will prove the double set inclusion.\\
\begin{itemize}
\item Let $\ell$ be the complex line determined by the parameter $(\xi_1,\xi_2)\in\mathcal{P}$. Since $\xi_1\neq 0,1$, $z_1\neq v_0,v_1$, and since $\xi_2\neq 1$, then $z_2\neq v_2$ and hence, $\ell$ doesn't contain any intersection of the lines $\ell_0$, $\ell_1$ and $\ell_2$.\\

Since $\xi_2\neq 0$, then $\ell\neq \ell_0$. Since $\xi_2\neq \xi_1$, then the line $\ell$ doesn't have the same complex direction of $\ell_1$. Since $\ell$ intersects $\ell_2$ and $\ell_1$, then $\ell$ doesn't have the same complex direction of neither $\ell_1$ nor $\ell_2$.\\ 

All of these considerations implies that $\ell$ is in general position with $\set{\ell_0,\ell_1,\ell_2}$.

\item Let $\ell$ be a line in general position with the array of lines $\set{\ell_0,\ell_1,\ell_2}$. From Observation \ref{obs_posicion_gral_no_paralelas} we know that $\ell$ cannot have the same complex direction of $\ell_0$ and therefore they should intersect $\ell_0$, but this intersection cannot be $v_0$ or $v_1$, since the four lines are in general position, then $z_1\neq v_0,v_1$ and then $\xi_1\in\C\setminus\set{0,1}$.\\

In the same way, $\ell$ intersects $\ell_2\setminus\set{v_0,v_2}$ and then $\xi_2\neq 0,1$. Finally, $\ell$ doesn't have the same complex direction of $\ell_1$ and then $\xi_2\neq\xi_1$. Hence, $\xi_2\in\C\setminus\set{0,1,\xi_1}$. We finally conclude that $(\xi_1,\xi_2)\in\mathcal{P}$.
\end{itemize}
\end{proof}

\begin{obs}
The set $\mathcal{P}\subset\C^2$ can be regarded as $\parentesis{\C\setminus\set{0,1}}\times\parentesis{\C\setminus\set{0,1}}$ minus a complex line going through the origin.
\end{obs}

Now we describe the space of the arrays of 4 complex lines in general position mutually tangent to a ball in $\C^2$. This parameter space is important because we know that Patterson-Sullivan measures exist for subgroups of $\pu$, this is an immediate consequence of the work of D. Sullivan (see \cite{sullivan79}).

\begin{prop}
The subspace $\mathcal{P}_1\subset\mathcal{P}$ parametrizing arrays of complex lines in general position, tangent to a ball, is given by

\begin{equation}\label{eq_descripcion_mut_tangentes}
	\mathcal{P}_1=\SET{\parentesis{\xi_1,\xi_2}\in\C^2}{\xi_1\in\C\setminus\set{0,1,\alpha_0,3-2\sqrt{2},2-\sqrt{2}}\text{, }\xi_2=\varphi(\xi_1)}
	\end{equation}

where $\varphi:\C\rightarrow\C$ is the function satisfying that, if $(\xi,0)\in\ell_0$ then $\parentesis{0,\varphi(\xi)}\in\ell_2$ is the intersection of the line tangent to the ball with boundary $S$, passing through $(\xi,0)$, with the line $\ell_2$.
\end{prop}

\begin{proof}
The lines $\ell_0$, $\ell_1$ and $\ell_2$ are mutually tangent to the ball with boundary
	$$S=\SET{(z_1,z_2)\in\C^2}{\valorabs{z_1-\alpha_0}^2+\valorabs{z_2-\alpha_0}^2=\alpha_0^2}$$
where $\alpha_0=\frac{1}{\sqrt{2}\parentesis{1+\sqrt{2}}}$.\\

Any point $(z_1,z_2)\in S$ determines a complex line tangent to $S$ passing through $(z_1,z_2)$. This line intersects $\ell_0$ in the point $(\xi,0)$, where
	$$\xi=\frac{\alpha_0\overline{z}_1+\alpha_0\overline{z}_2-\parentesis{\valorabs{z_1}^2+\valorabs{z_2}^2}}{\alpha_0-\overline{z}_1}.$$

Each one of these lines in $S$ is in general position with the 3 lines $\ell_0$, $\ell_1$ and $\ell_2$, unless:
\begin{enumerate}[(i)]
\item The point $(z_1,z_2)$ is the tangency point with the lines $\ell_0$, $\ell_1$, $\ell_2$. This happens if 
	$$\begin{array}{ccc}
	(z_1,z_2)=(\alpha_0,0), & (z_1,z_2)=(\frac{1}{2},\frac{1}{2}), & (z_1,z_2)=(0,\alpha_0)
	\end{array}$$ 
	respectively.\\

In the first case, the line would be $\ell_0$. In the other two cases, these tangent lines determine the intersections in $\ell_0$ given by the parameters $\xi=1$ and $\xi=0$ respectively.
 
\item The tangent line determined has the same complex direction that any of the lines $\ell_0$, $\ell_1$, $\ell_2$. This happens if 
	$$\begin{array}{ccc}
	(z_1,z_2)=(\alpha_0,2\alpha_0), & (z_1,z_2)=(\frac{3-2\sqrt{2}}{2},\frac{3-2\sqrt{2}}{2}), & (z_1,z_2)=(2\alpha_0,\alpha_0)
	\end{array}$$
respectively.\\

In the first case the line would have the same complex direction of $\ell_0$. In the other two cases, these tangent lines determine the intersections in $\ell_0$ given by the parameters $\xi=3-2\sqrt{2}$ and $\xi=2-\sqrt{2}$ respectively.
\end{enumerate}

All of this implies (\ref{eq_descripcion_mut_tangentes}).
\end{proof}

It can be seen in (\ref{eq_descripcion_mut_tangentes}) that $\mathcal{P}_1$ is a complex curve in $\C^2$ omitting $5$ points.\\

\section{Estimating the entropy volume}
\label{sec_volumen_entropia}

For the rest of the chapter we will use the Kobayashi and Eisenman metric (and volume) indistinctly since they share the same needed properties (Proposition \ref{prop_dist_contractante} and \ref{prop_kobayashi_cubriente}), for a quick review see section \ref{sec_kobayashi} and for more details see \cite{kobayashi}, \cite{kobayashi05} and \cite{graham1985some}.\\

In this section we propose a way to bound the entropy volume of the Koba-yashi metric on a domain, which is the complement of an array of 5 complex projective lines in general position in $\CP^2$. In the first subsection we study what happens when we consider one line $L_\eta$ intersecting an array of lines in general position.\\

In the second subsection, the line $L_\eta$ will be regarded as a \emph{diameter} of a Kobayashi ball which is homeomorphic to a domain in $\CP^1$. We will consider the Kobayashi metric on this domain and give an estimate of the volume of balls with respect to that metric.\\

In the third subsection we give some ideas on how to estimate the entropy volume we're interested in, using the estimate found in each \emph{diameter} of a given Kobayashi ball.\\

\begin{defn}
Let $M$ be a complex manifold. For $z\in M$ we define the \emph{entropy volume} as
$$e(M,z)=\lim_{r\rightarrow\infty}\frac{\log\parentesis{\text{Vol}\parentesis{\bola{r}{\Omega}{z}}}}{r}.$$
\end{defn}

From now on we assume that $\Gamma\subset\PSL$ is a strongly irreducible complex Kleinian group. By Theorem \ref{thm_complemento_kob_hyp}, $\KulD(\Gamma)$ is complete Kobayashi hyperbolic. We define the orbital counting function in the same way as in \cite{nicholls}.

\begin{defn}
Let $r>0$ and $z,w\in\KulD(\Gamma)$, we define the orbital counting function as
	$$N(r,z,w)=\valorabs{\SET{\gamma\in\Gamma}{d\parentesis{z,\gamma(w)<r}}}.$$
where $d$ is the Kobayashi metric on $\KulD(\Gamma)$.
\end{defn} 

$N$ counts the elements of the $\Gamma$-orbit of $w$ lying inside the open ball $\bola{R}{\KulD(\Gamma)}{z}$.

\begin{prop}
$N(r,z,w)=N(r,w,z)$.
\end{prop}

\begin{proof}
Let $A:=\SET{\gamma\in\Gamma}{\rho\parentesis{z,\gamma(w)}<r}$ and $B:=\SET{\gamma\in\Gamma}{\rho\parentesis{w,\gamma(z)}<r}$. To each $\psi\in A$ we can assign $\psi^{-1}\in B$, to do this observe that, if $\psi\in A$ then $\rho(z,\psi w)<r$ and it follows
	$$\rho(w,\psi^{-1} z)=\rho(\psi w, z)< r,$$
this is a consequence of Proposition \ref{prop_kob_isometrias}. This mapping $\psi\mapsto\psi^{-1}$ is bijective and thus $\valorabs{A}=\valorabs{B}$. This means, $N(r,z,w)=N(r,w,z)$.
\end{proof}

In the following lemma we prove that, if a region has finite entropy volume then the region resulting of removing one additional line has finite entropy volume as well.

\begin{lem}
Let $\Omega\subset\C^2$ be a domain given by 
	$$\Omega=\C^2\setminus\bigcup_{i=1}^k \ell_i$$
where $\ell_i$ are complex lines in $\C^2$ in general position. Let $\ell$ be a complex line in general position with the lines $\ell_1,...,\ell_k$. If $\Omega$ has finite entropy volume, then $\Omega\setminus\ell$ has finite entropy volume.
\end{lem}

\begin{proof}
We know that $\Omega$ is complete Kobayashi hyperbolic if $k\geq 4$ (Theorem \ref{thm_complemento_hiperplanos}), then $\Omega':=\Omega\setminus\ell$ is complete Kobayashi hyperbolic (see Proposition 4.2 of \cite{kobayashi05}). Let us denote by $d_{\Omega}$ and $d_{\Omega'}$ the Kobayashi metrics on $\Omega$ and $\Omega'$ respectively. Using the inclusion $\mathsf{I}:\Omega'\rightarrow\Omega$ in Proposition \ref{prop_dist_contractante} it follows that
	\begin{equation}\label{eq_dem_lema_regiones_mas_chicas_1}
	d_{\Omega}(z,w)\leq d_{\Omega'}(z,w)
	\end{equation}		
for any $z,w\in\Omega'$. In the same way, if $F_{\Omega}$ and $F_{\Omega'}$ are the local forms of the Kobayashi metric on $\Omega$ and $\Omega'$ respectively, then
	$$F_{\Omega'}(z,\zeta)\leq F_{\Omega}(z,\zeta)$$
for any $z\in\C^2$ and any tangent vector $\zeta\in\C^2$ on $z$.\\

Let us denote, for $z\in\Omega$ and $R>0$
	$$B_R^\Omega(z)=\SET{w\in\Omega}{d_\Omega(z,w)\leq R}$$
and analogously, $B_R^{\Omega'}(z)$. $\text{Vol}(U)$ denotes the Kobayashi volume of the region $U$.\\

Let $w\in B_R^{\Omega'}(z)$. Then $d_{\Omega'}(w,z)\leq R$, from this and (\ref{eq_dem_lema_regiones_mas_chicas_1}) it follows that $d_{\Omega}(w,z)\leq R$. Hence, 
	\begin{equation}\label{eq_dem_lema_regiones_mas_chicas_2}
	B_R^{\Omega'}(z)\subset B_R^{\Omega}(z)
	\end{equation}
for $z\in \Omega'$ and $R>0$.\\

As a consequence of this, 
	\begin{equation}\label{eq_dem_lema_regiones_mas_chicas_3}	
	\text{Vol}(B_R^{\Omega'}(z))\leq \text{Vol}(B_R^{\Omega}(z)).
	\end{equation}
By hyphothesis, the entropy volume of $\Omega$ is finite, that is, $e(\Omega,z)<\infty$. Then, from (\ref{eq_dem_lema_regiones_mas_chicas_3}) it follows that the entropy volume of $\Omega'$ es finite.
\end{proof}

\subsection*{Interaction of one line with an array of lines}

Let $\Upsilon=(\zeta_1,\zeta_2)\in \mathcal{P}$ be a parameter determining an array of 5 complex projective lines in general position, $\mathcal{L}_\Upsilon$. Denote $\Omega_{\Upsilon}=\C^2\setminus \mathcal{L}_\Upsilon$. For now, we will omit the dependence on the parameter $\Upsilon$ and just write $\Omega$, for the sake of legibility.\\ 

Let $z=(z_1,z_2)\in\Omega$ and $R>0$. Consider a complex line $L$ passing through $z$, if $L$ doesn't have the same complex direction of $\ell_0$ then $L\cap\ell_0\neq \emptyset$ (see Observation \ref{prop_interseccion_lineas_complejas}) and then, we can parametrize all the complex lines passing through $z$ using the complex number $\eta$ such that $\ell_0(\eta)=L\cap\ell_0$. \\

For $\eta\in\C$, we will denote by $L_\eta$ to the complex line passing through $z$ and $\ell_0(\eta)=(\eta,0)$. Finally we denote 
 $$\hat{L}_{R,\eta}=L_\eta\cap \bola{r}{\Omega}{z}.$$

If there are is ambiguity, we will omit the dependence on $R$ and just write $\hat{L}_\eta$.\\

Consider the parametrization of $L_\eta$ given by the function $h_\eta:\C\rightarrow \C^2$ defined as
	$$h_\eta(\xi)= \parentesis{z_1+\xi (\eta-z_1),(1-\xi)z_2}.$$

Now we need to define several lines determined by the intersections of the lines in the array. We denote $q_{i,j}=\ell_i\cap\ell_j$ for integers $i<j$, $0\leq i,j\leq 4$. Let $\mathcal{L}=\ell_0\cup...\cup\ell_3$; since the lines in $\mathcal{L}$ are in general position, then the points $\set{q_{i,j}}$ are pairwise distinct.\\

Using the parametrizations we have chosen for the lines $\ell_0,...,\ell_4$, it follows 

$$\begin{array}{llll}
q_{0,1} = (1,0) & q_{1,2} = (0,1) & q_{2,3} = (0,\zeta_2) & q_{3,4} = \corchetes{-\zeta_1:\zeta_2:0} \\ 
q_{0,2} = (0,0) & q_{1,3} = \parentesis{\frac{\zeta_1(1-\zeta_2)}{\zeta_1-\zeta_2},\frac{\zeta_2(\zeta_1-1)}{\zeta_1-\zeta_2}} & q_{2,4} = \corchetes{0:1:0} &  \\
q_{0,3} = (\zeta_1,0) & q_{1,4} = \corchetes{-1:1:0} &  &  \\
q_{0,4} = \corchetes{1:0:0} &  &  &  
\end{array}$$

\begin{figure}[H]
\begin{center}	
\includegraphics[height=55mm]{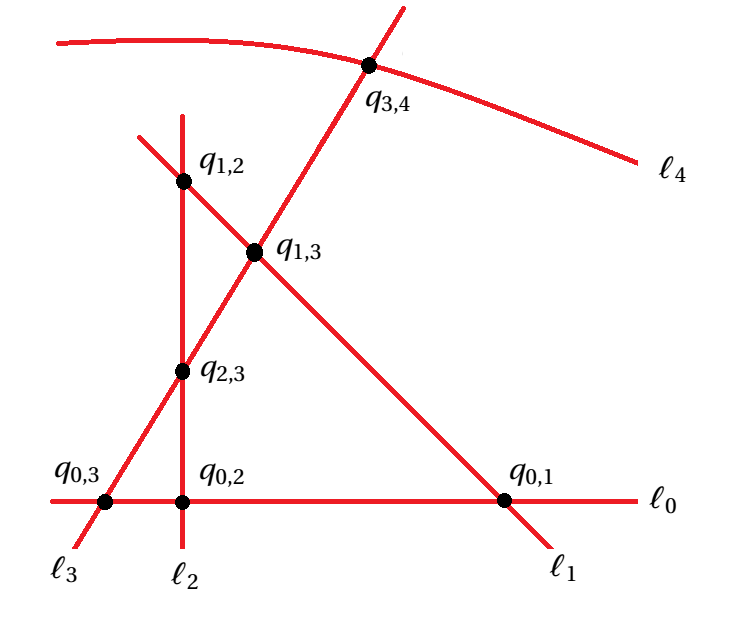}	
\caption{The lines $\ell_i$ and their intersections $q_{i,j}$.}
\label{fig_intersecciones_lineas}
\end{center}	
\end{figure}

Let $i,j,k,m\in\set{0,1,2,3,4}$ be pairwise different integers such that $i<j$, $k<m$. We define the complex projective lines
	$$H_{i,j,k,m}=\overleftrightarrow{q_{i,j},q_{k,m}}.$$
Observe that $H_{i,j,k,m}=H_{k,m,i,j}$. We denote by $\mathcal{H}$ the union of the lines $H_{i,j,k,m}$.\\ 

\begin{figure}[H]
\begin{center}	
\includegraphics[height=55mm]{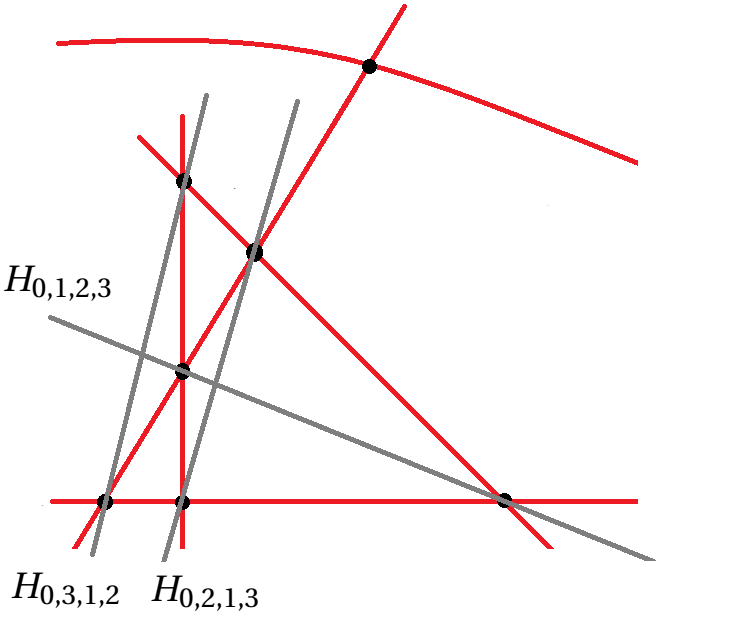}	
\caption{Some complex lines $H_{i,j,k,m}$.}
\label{fig_lineas_prohibidas_para_z}
\end{center}	
\end{figure}

If we only consider $i,j,k,m\in\set{0,1,2,3}$, we have 3 complex lines $H_{i,j,k,m}$, which are shown in figure (\ref{fig_lineas_prohibidas_para_z}).

For $z\nin \mathcal{H}$ we define $L_{i,j}=\overleftrightarrow{z,q_{i,j}}$, where $i<j$, with $i,j\in\set{0,...,4}$. Denote by $\mathcal{M}$ the union of these 10 lines $L_{i,j}$:

$$\mathcal{M}=\bigcup_{\begin{array}{c}
i,j\in\set{0,...,4}\\
i<j
\end{array}} L_{i,j}\subset \CP^2.$$

Denote by $L_{i,j}^{\ast}\in\parentesis{CP^2}^{\ast}$ to the dualization of the lines $L_{i,j}$. Let

$$\mathcal{M}^{\ast}=\bigcup_{\begin{array}{c}
i,j\in\set{0,...,4}\\
i<j
\end{array}} L^{\ast}_{i,j}\subset \parentesis{\CP^2}^{\ast}.$$
 
\begin{figure}[H]
\begin{center}	
\includegraphics[height=55mm]{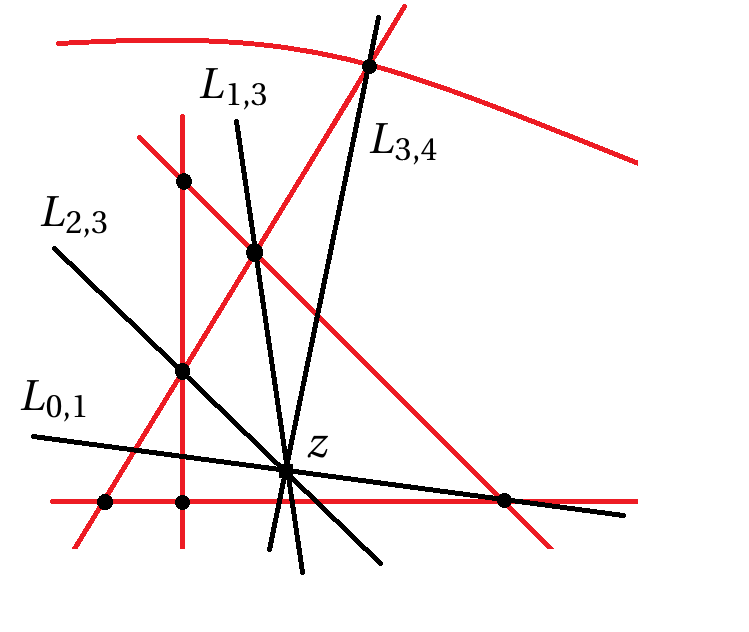}	
\caption{Some lines $L_{i,j}$.}
\label{fig_lineas_prohibidas}
\end{center}	
\end{figure}

Observe that the line $L_{i,4}$ has the same complex direction of $\ell_i$, for any $i=0,...,3$. This is because $L_{i,4}$, $\ell_i$ and $\ell_4$ are concurrent and $\ell_4$ is the line at infinity.

\begin{lem}\label{lem_casos_finitos_colapsos_intersecciones}
Let $z\in\Omega$ such that $z\nin \mathcal{H}$. Let $\eta\in\C$, if $L^{\ast}_\eta\nin\mathcal{M}^{\ast}$ then $L_\eta$ intersects $\mathcal{L}$ in exactly 5 points. If $L^{\ast}_\eta\in\mathcal{M}^{\ast}$ then $L_\eta$ intersects $\mathcal{L}$ in exactly 4 points.
\end{lem}

\begin{proof}
Let $z\in\Omega\setminus\mathcal{H}$ and $\xi\in\C$. We know from Propositions \ref{prop_interseccion_lineas_complejas} that $L_\eta$ intersects each line $\ell_i$ in exactly one point, therefore $\valorabs{L_\eta\cap \mathcal{L}}\leq 5$.\\

If $L^\ast_\eta\nin\mathcal{M}^\ast$, let us suppose that $\valorabs{L_\eta\cap \mathcal{L}}< 5$, this means that $L_\eta\cap \ell_i=L_\eta\cap \ell_j$ for some $i,j\in\set{0,...,4}$ and suppose, without loss of generality that $i<j$, denote $y=L_\eta\cap \ell_i$. Then $y\in\ell_i$ and $y\in\ell_j$, therefore $y=q_{i,j}$. Since $L_\eta$ passes through $z$ and $y\in L_\eta$ then $L_\eta=\overleftrightarrow{z,q_{i,j}}=L_{i,j}$, this contradicts that $L^\ast_\eta\nin\mathcal{M}^\ast$. Therefore, $L^\ast_\eta\nin\mathcal{M}^\ast$ implies $\valorabs{L_\eta\cap \mathcal{L}}= 5$.\\

If $L^\ast_\eta\in\mathcal{M}^\ast$ then $L_\eta=L_{i,j}$ for some $i,j\in\set{0,...,4}$ with $i<j$, then $L_\eta$ passes through $q_{i,j}$, in other words, $L_\eta\cap \ell_i=L_\eta\cap \ell_j$ and therefore $\valorabs{L_\eta\cap \mathcal{L}}\leq 4$. Suppose that $\valorabs{L_\eta\cap \mathcal{L}}<4$, this can happen in two ways:
\begin{enumerate}
\item $L_\eta\cap \ell_i=L_\eta\cap \ell_j=L_\eta\cap \ell_k$ for $i,j,k\in\set{0,...,4}$ pairwise different.
\item $L_\eta\cap \ell_i=L_\eta\cap \ell_j$ and $L_\eta\cap \ell_k=L_\eta\cap \ell_m$ for $i,j,k,m\in\set{0,...,4}$ pairwise different.
\end{enumerate}

In the first case, if we denote $y=L_\eta\cap \ell_i$, then $y\in \ell_i\cap\ell_j\cap\ell_k$, contradicting that the array of lines $\mathcal{L}$ is in general position.\\

In the second case, $q_{i,j}\in L_\eta$ and $q_{k,m}\in L_\eta$ then $L_\eta=H_{i,j,k,m}$ and therefore $z\in H_{i,j,k,m}$ contradicting that $z\nin\mathcal{H}$. Then, $L^\ast_\eta\in\mathcal{M}^\ast$ implies $\valorabs{L_\eta\cap \mathcal{L}}=4$.
\end{proof}

\begin{obs}
The previous lemma says that, except for 10 cases (whenever the lines $L_\eta$ coincide with some line $L_{i,j}$), any line $L_\eta$ intersects the array of lines $\mathcal{L}$ in exactly 5 points.
\end{obs}

The following lemma describes the complex parameters $\eta$ for which $L_\eta$ intersects the array of lines $\mathcal{L}$ in exactly 4 points.

\begin{lem}
Let $\Upsilon=(\zeta_1,\zeta_2)\in \mathcal{P}$ be a parameter determining an array of 4 complex lines in general position, $\mathcal{L}=\mathcal{L}_\Upsilon$. If $\eta\in\C_0$ then $L_\eta$ intersects the array $\mathcal{L}$ in exactly 4 points, where
	$$\C_0=\set{0,1,\zeta_1,\frac{z_1}{1-z_2},\frac{\zeta_2 z_1+\zeta_1 z_2 - \zeta_1\zeta_2(z_1+z_2)}{\zeta_2(1-\zeta_1)+z_2(\zeta_1-\zeta_2)},z_1+z_2,\frac{\zeta_2 z_1}{z_2-\zeta_2},z_1,\frac{\zeta_2 z_1+\zeta_1 z_2}{\zeta_2}}.$$
\end{lem}

\begin{proof}
A straight-forward calculation verifies that:
	\begin{center}
	\begin{tabular}{ll}
	If $\eta=1$, $L_\eta=L_{0,1}$. & If $\eta=0$, $L_\eta=L_{0,2}$. \\
	If $\eta=\zeta_1$, $L_\eta=L_{0,3}$. & If $\eta=\frac{z_1}{1-z_2}$, $L_\eta=L_{1,2}$.\\
	If $\eta=\frac{\zeta_2 z_1+\zeta_1 z_2 - \zeta_1\zeta_2(z_1+z_2)}{\zeta_2(1-\zeta_1)+z_2(\zeta_1-\zeta_2)}$, $L_\eta=L_{1,3}$. & If $\eta=z_1+z_2$, $L_\eta=L_{1,4}$.\\
	If $\eta=\frac{\zeta_2 z_1}{z_2-\zeta_2}$, $L_\eta=L_{2,3}$. & If $\eta=z_1$, $L_\eta=L_{2,4}$.\\
	If $L_\eta=\frac{\zeta_2 z_1+\zeta_1 z_2}{\zeta_2}$, $L_\eta=L_{3,4}$. & 
	\end{tabular}	
	\end{center}

As it has been proved in Lemma \ref{lem_casos_finitos_colapsos_intersecciones}, $L_\eta=L_{i,j}$ implies $\valorabs{L_\eta\cap\mathcal{L}}=4$. Observe that there is no $\eta\in\C$ such that $L_\eta=L_{0,4}$. However, for $\eta=\infty\in\hat{\C}$, $L_\eta=L_{0,4}$.
\end{proof}

From now on, consider the array of complex lines in general position $\mathcal{L}$ given by 
	$$\mathcal{L}=\ell_0\cup\ell_1\cup\ell_2\cup\ell_3.$$
As always, one has to keep in mind that the fifth complex projective line $\ell_4$ is the line at infinity.\\
	
As a consequence of Lemma \ref{lem_casos_finitos_colapsos_intersecciones}, let us denote 
	$$P_i^\eta= L_\eta \cap \ell_i$$ 
for $i=0,...,3$, the 4 intersection points $L_\eta\cap\ell_i$. Let 
	$$p_i^\eta=h_\eta^{-1}\parentesis{P_i^\eta}$$
for $i=0,...,3$. Since the complex lines $\ell_0,\ell_1,\ell_2$ are fixed, then

	$$\begin{array}{ll}
	P_0^\eta = (\eta,0) & p_0^\eta = 1\\
	P_1^\eta = (\frac{z_1+\eta(z_2-1)}{z_1+z_2-\eta},\frac{(\eta-1)z_2}{z_1+z_2-\eta}) & p_1^\eta = \frac{z_1+z_2-1}{z_1+z_2-\eta}\\
	P_2^\eta = (0,\frac{\eta z_2}{\eta-z_1}) & p_2^\eta = \frac{z_1}{z_1-\eta}.	
	\end{array}$$

Besides, these points do not depend on the parameter $\Upsilon$. On the other hand, $P_3$ and $p_3$ do depend on $\Upsilon$,

	$$\begin{array}{ll}
	P_3^{\Upsilon,\eta} = \parentesis{\frac{\zeta_1\parentesis{\eta(z_2-\zeta_2)+z_1\zeta_2}}{z_2\zeta_1+(z_1-\eta)\zeta_2},\frac{z_2\zeta_2(\zeta_1-\eta)}{z_2\zeta_1+(z_1-\eta)\zeta_2}} & p_3^{\Upsilon,\eta} = \frac{z_2\zeta_1+\zeta_2(z_1-\zeta_1)}{z_2\zeta_1+(z_1-\eta)\zeta_2}.
	\end{array}$$

\subsection{Estimate of the Kobayashi volume in domains of $\CP^1$}

In this subsection we will use indistinctly the terms \emph{volume} and \emph{area}. Let $S^{\Upsilon,\eta}=\C\setminus\set{p_0^\eta,p_1^\eta,p_2^\eta,p_3^{\Upsilon,\eta}}$, this region is homeomorphic to a 5-punctured sphere $S_{0,5}$. As long as there is no ambiguity, we will omit the dependence on $\Upsilon$ and $\eta$, and just write $S$ instead of $S^{\Upsilon,\eta}$; analogously, we will just write $L$, $h$, $P_i$ and $p_i$.\\ 

Observe that $S$ is Kobayashi hyperbolic (see \cite{kobayashi}) and, since 
	$$h:S\rightarrow h(S)=L\setminus\set{P_0,...,P_3}$$ 
is biholomorphic, it follows from Proposition \ref{prop_dist_contractante} that $\parentesis{S,d_S}$ and $\parentesis{L,d_{L}}$ are isometric. Furthermore, we have the following lemma.

\begin{lem}\label{lem_esfera_isometrica_linea}
$\parentesis{S,d_S}$ is isometric to $\parentesis{L,d_\Omega}$.
\end{lem}

\begin{proof}
Since $\parentesis{S,d_S}$ and $\parentesis{L\setminus\set{P_0,...,P_3},d_{L\setminus\set{P_0,...,P_3}}}$ are isometric, it will be enough to prove that the Kobayashi metrics $d_{L\setminus\set{P_0,...,P_3}}$ and $d_\Omega$ restricted to $L\setminus\set{P_0,...,P_3}$ are equal. To do this, we will prove that, for any chain of holomorphic disks in $\Omega$ there is a chain of holomorphic disks in $L\setminus\set{P_0,...,P_3}$ with the same length.\\

Let $p,q\in L\setminus\set{P_0,...,P_3}$ and let $\alpha$ be a chain of holomorphic disks in $\Omega$ given by
	$$\begin{array}{c}
	p=p_0,p_1,...,p_{k-1},p_k=q \in\Omega\\
	a_1,b_1,...,a_k,b_k\in \D\\
	f_1,...,f_k:\D\rightarrow \Omega\text{ holomorphic.}
	\end{array}$$

For any $n=0,...,k$ such that $p_n\in L\setminus\set{P_0,...,P_3}$, let $\tilde{p}_n=p_n$. For any $n=0,...,k$ such that $p_n\in\C^2\setminus\parentesis{L\setminus\set{P_0,...,P_3}}$, let $\tilde{p}_n\in L\setminus\set{P_0,...,P_3}$ such that $\tilde{p}_n\neq p_i$ for $i=0,...,k$.\\

For $i=1,...,k$, let $\tilde{f}_i:\D\rightarrow L\setminus\set{P_0,...,P_3}$ given by 
	$$\tilde{f}_i(z)=\tilde{p}_{i-1}+\frac{z-a_i}{b_i-a_i}\parentesis{\tilde{p}_i-\tilde{p}_{i-1}}.$$
Clearly $\tilde{f}_i$ is holomorphic and, furthermore
	\begin{align*}
	\tilde{f}_i(a_i)&=\tilde{p}_{i-1}\\
	\tilde{f}_i(b_i)&=\tilde{p}_{i}. 
	\end{align*}
Let $\tilde{a}$ be the chain of holomorphic disks formed with the points $\tilde{p}_i,a_i,b_i$ and the holomorphic mappings $\tilde{f}_i$. It is clear that $\ell(\alpha)=\ell(\tilde{\alpha})$. Then, for any chain $\alpha$ containing the points $p_i\nin L\setminus\set{P_0,...,P_3}$ we have a chain $\tilde{a}$ such that $\tilde{p}_i\in L\setminus\set{P_0,...,P_3}$ for $i=0,...,k$, with the same length. Therefore, if $p,q\in L\setminus\set{P_0,...,P_3}$, it holds
	$$d_\Omega(p,q)=d_{L\setminus\set{P_0,...,P_3}}(p,q).$$
Therefore, $\parentesis{S,d_S}$ and $\parentesis{L,d_\Omega}$ are isometric.	
\end{proof}

Since $S$ is hyperbolic, there is a covering map $\pi_{\Upsilon,\eta}:\D\rightarrow S$ and a Fuchsian model $\Gamma_{\Upsilon,\eta}$ of $S$, that is, a discrete subgroup of $\text{Aut}\parentesis{\D}$ such that
	$$S\cong\D/\Gamma_{\Upsilon,\eta}.$$ 
Once again, as long as there is no ambiguity, we will write $\pi$ and $\Gamma$ instead of $\pi_{\Upsilon,\eta}$ and $\Gamma_{\Upsilon,\eta}$ respectively.\\ 

We will need the following proposition (see Proposition 1.6 of \cite{kobayashi}).

\begin{prop}\label{prop_kobayashi_cubriente}
Let $M$ be a complex manifold and $\tilde{M}$ a covering manifold of $M$ with covering map $\pi:\tilde{M}\rightarrow M$. Let $p,q\in M$ and $\tilde{p},\tilde{q}\in \tilde{M}$ such that $\pi(\tilde{p})=p$ and $\pi(\tilde{q})=q$. Then
	$$d_M\parentesis{p,q}=\inf_{\tilde{q}}d_{\tilde{M}}\parentesis{\tilde{p},\tilde{q}},$$
where the infimum is taken over all points $\tilde{q}\in \tilde{M}$ such that $\pi(\tilde{q})=q$.
\end{prop}

Then it follows that, for any $w_1,w_2\in S$ and $\tilde{w}_1\in\pi^{-1}(w_1)$,
	
	\begin{equation}\label{eq_calculo_vol_entropia_2}
	d_{S}\parentesis{w_1,w_2}=\inf_{\tilde{w}_2\in\pi^{-1}(w_2) }\rho\parentesis{\tilde{w}_1,\tilde{w}_2},	
	\end{equation}

where $\rho$ is the Poincar\'e metric on $\D$ (which coincides with the Kobayashi metric on $\D$, see Proposition 1.3 of \cite{kobayashi}). Since $\Gamma$ is discrete, (\ref{eq_calculo_vol_entropia_2}) can be rewritten as

	\begin{equation}\label{eq_calculo_vol_entropia_3}
	d_{S}\parentesis{w_1,w_2}=\min_{\tilde{w}_2\in\pi^{-1}(w_2) }\rho\parentesis{\tilde{w}_1,\tilde{w}_2}.	
	\end{equation}

Let $p=h^{-1}(z)$ and let $\hat{p}\in\pi^{-1}(p)$ be a point in the pre-image $\pi^{-1}(p)$, suppose that $\hat{p}$ is not a fixed point of any $\gamma\in\Gamma\setminus\set{\id}$. Let $D=D_{\hat{p}}$ a Dirichlet region centered in $\hat{p}$. For any $q\in S$ we denote by $\tilde{q}\in\pi^{-1}(q)$ the element of the pre-image satisfying	
	\begin{equation}\label{eq_calculo_vol_entropia_4}
	d_{S}\parentesis{p,q}=\rho\parentesis{\hat{p},\tilde{q}}.	
	\end{equation}	 
As a consequence of this, an open ball centered in $p$, with respect to the Kobayashi metric looks like
	\begin{align*}
	\bola{S}{R}{p}&=\SET{q\in S}{d_S(p,q)<R}\\
		&=\SET{q\in S}{\rho\parentesis{\hat{p},\tilde{q}}<R}.
	\end{align*}		
In other words, the points in $\bola{S}{R}{p}$ are the points $q\in S$ such that $\tilde{q}\in\bola{\D}{R}{\hat{p}}$.\\

As a consequence of the previous observation and the definition of a fundamental domain, if $R<\rho\parentesis{\hat{p},\partial D}$ then $\bola{\D}{R}{\hat{p}}\subset D$ and then $\pi:\bola{\D}{R}{\hat{p}}\rightarrow \bola{S}{R}{p}$ is a biholomorphism and therefore, an isometry. Then
	$$\text{Vol}\parentesis{\bola{\D}{R}{\hat{p}}}=\text{Vol}\parentesis{\bola{S}{R}{p}}.$$
Let $B_{\gamma,R}=\bola{\D}{R}{\hat{p}}\cap \gamma(D)$. If $R>\rho\parentesis{\hat{p},\partial D}$ then $B_{\gamma,R}\neq \emptyset$ for some $\gamma\in\Gamma\setminus\set{\id}$. The restriction $\pi:\bola{\D}{R}{\hat{p}}\rightarrow S$ is no longer an isometry in general, since, for $w\in B_{\gamma,R}$ $\exists$ $w_0\in D$ such that $\pi(w)=\pi(w_0)$ and then this restriction of $\pi$ is not injective. However, the restriction 
	$$\pi:B_{\gamma,R}\rightarrow \pi\parentesis{B_{\gamma,R}}$$ 
is an isometry, since we can apply the same argument we used when $R<\rho\parentesis{\hat{p},\partial D}$ and as a consequence of
	$$d_S(p,q)=\rho\parentesis{\hat{p},\tilde{q}}=\rho\parentesis{\gamma(\hat{p}),\gamma(\tilde{p})}.$$
This means that we can estimate the volume of $\bola{S}{R}{p}$ by decomposing $\bola{\D}{R}{\hat{p}}$ into the intersections $B_{\gamma,R}\neq \emptyset$.\\

For $R>0$, we define 
	\begin{align*}
	\hat{N}(R)&=\valorabs{\SET{\gamma\in\Gamma}{\bola{\D}{R}{\hat{p}}\cap \gamma(D)\neq\emptyset}}\\
		&=\valorabs{\SET{\gamma\in\Gamma}{B_{\gamma,R}\neq\emptyset}}.
	\end{align*}		
This map counts the number of images of $D$ under elements of $\Gamma$ intersecting $\bola{\D}{R}{\hat{p}}$.
Observe that
$$\text{Vol}(D)=\text{Vol}\parentesis{\gamma(D)}$$
for any $\gamma\in\Gamma$. Therefore, we can estimate the area of $\bola{S}{R}{p}$ in the following way:
	\begin{equation}\label{eq_forma_cota_volumen_rebanada}
	\text{Vol}\parentesis{\bola{S}{R}{p}}\leq\hat{N}(R)\text{Vol}(D),
	\end{equation}		
whenever $\text{Vol}(D)<\infty$. The following proposition guarantees that this always happens.

\begin{figure}[H]
\begin{center}	
\includegraphics[height=35mm]{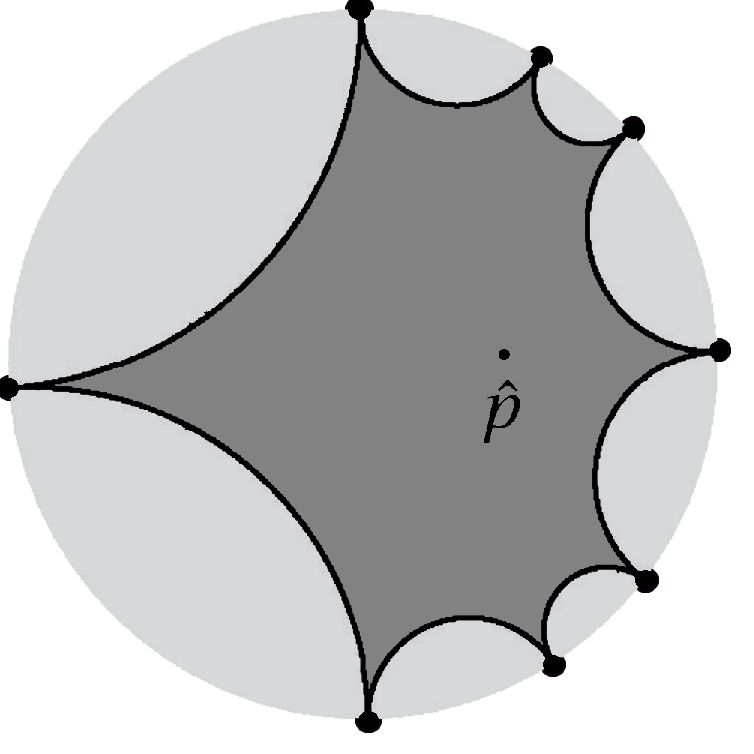}	
\caption{Dirichlet region for $\Gamma$.}
\label{fig_region_fundamental}
\end{center}	
\end{figure}

\begin{prop}\label{prop_volumen_finito}
For any $\eta\in\C$, other that the finite number of cases described in Lemma \ref{lem_casos_finitos_colapsos_intersecciones}, and for any $\Upsilon\in \mathcal{P}$, the Dirichlet region $D$ centered in $\hat{p}$ for $\Gamma=\Gamma_{\Upsilon,\eta}$ has finite volume. Furthermore,
	$$\text{Vol}(D)=6\pi.$$ 
Finally, $D$ is an ideal octagon like the one depicted in Figure (\ref{fig_region_fundamental}).
\end{prop}

\begin{figure}[H]
\begin{center}	
\includegraphics[height=45mm]{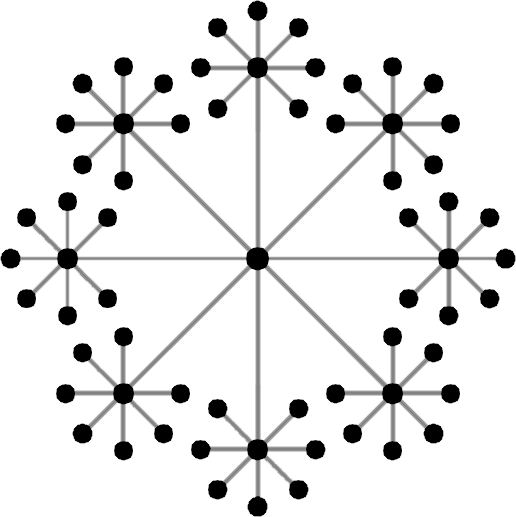}	
\caption{Part of the Cayley tree of $\Gamma$.}
\label{fig_arbol_cayley}
\end{center}	
\end{figure}

We know that $\Gamma$ is a free group with 8 generators. Consider the Cayley tree of $\Gamma$, as in Figure (\ref{fig_arbol_cayley}).

\begin{figure}[H]
\begin{center}	
\includegraphics[height=45mm]{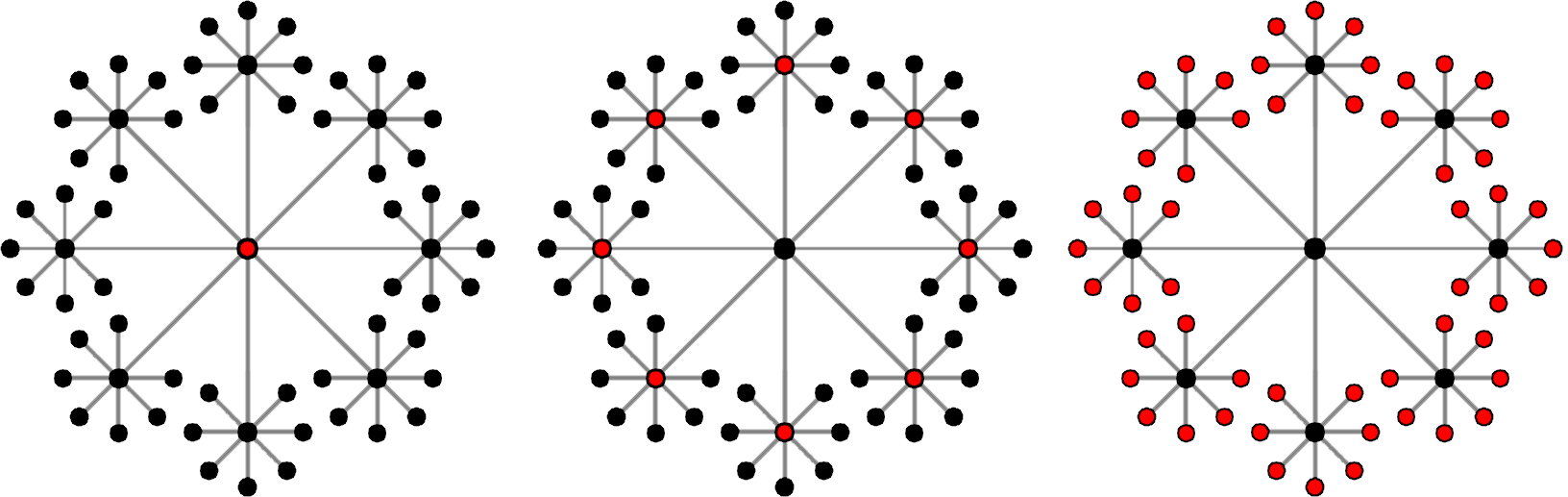}	
\caption{$S^\Gamma_{0}(\text{id})$, $S^\Gamma_{1}(\text{id})$ and $S^\Gamma_{2}(\text{id})$ respectively.}
\label{fig_niveles_arbol_cayley}
\end{center}	
\end{figure}

We denote by $\delta$ the word metric in $\Gamma$ and define, for $n\in\Z^+$ 
	\begin{align*}
	S_n^\Gamma(\id) &= \SET{\gamma\in \Gamma}{\delta(\id,\gamma)=n } \\ 
	B_n^\Gamma(\id) &= \SET{\gamma\in \Gamma}{\delta(\id,\gamma)\leq n}.	
	\end{align*}
In Figure (\ref{fig_niveles_arbol_cayley}), $S^\Gamma_{0}(\text{id})$, $S^\Gamma_{1}(\text{id})$ and $S^\Gamma_{2}(\text{id})$ are shown, in red, respectively.

The proof of the following lemma is straight-forward.

\begin{lem}\label{lem_conteo_gamma}
For $n\in \Z^+$ it holds
	$$\valorabs{S_n^\Gamma(\id)}=\begin{cases} 1 &, n=0.\\
	8\cdot 7^{n-1} &, n>0.\end{cases}$$
and, in consequence,
	$$\valorabs{B_n^\Gamma(\id)}=1+\frac{4}{3}(7^n-1).$$
\end{lem}

\begin{prop}
For any $R>0$ it holds $N(R,\hat{p},\hat{p})\leq \hat{N}(R)$
\end{prop}

\begin{proof}
Let
	\begin{align*}
	\Gamma_1 &= \SET{\gamma\in\Gamma}{\rho\parentesis{\hat{p},\gamma(\hat{p})}<R}\\
	\Gamma_2 &= \SET{\gamma\in\Gamma}{\bola{\D}{R}{\hat{p}}\cap \gamma(D)\neq\emptyset}
	\end{align*}
Let $\gamma\in \Gamma_1$, then $\gamma(\hat{p})\in\bola{\D}{R}{\hat{p}}$. Furthermore, since $\hat{p}\in D$, it follows that $\gamma(\hat{p})\in \gamma\parentesis{D}$ and therefore,
	$$\gamma(\hat{p})\in \bola{\D}{R}{\hat{p}}\cap \gamma(D),$$
then $\bola{\D}{R}{\hat{p}}\cap \gamma(D)\neq\emptyset$, and then $\gamma\in \Gamma_2$, this means that $\Gamma_1\subset \Gamma_2$. Hence, 
	$$N(R,\hat{p},\hat{p})=\valorabs{\Gamma_1}\leq\valorabs{\Gamma_2}= \hat{N}(R).$$
\end{proof}

The previous proposition gives a lower bound for $\hat{N}(R)$. Now, we need to find an upper bound for $\hat{N}(R)$. Let us define
	$$A(R)=\SET{\gamma\in\Gamma}{\gamma(D)\cap\bola{\D}{R}{\hat{p}}\neq\emptyset},$$
then $\hat{N}(R)=\valorabs{A(R)}$. We can decompose $A$ into the following disjoint subsets
	\begin{align*}
	A_1(R) &= \SET{\gamma\in\Gamma}{\gamma(\hat{p})\in\bola{\D}{R}{\hat{p}}} \\
	A_2(R) &= \SET{\gamma\in\Gamma}{\gamma(D)\cap\bola{\D}{R}{\hat{p}}\neq\emptyset\text{ y }\gamma(\hat{p})\nin\bola{\D}{R}{\hat{p}}}.
	\end{align*}
Observe that $N(R,\hat{p},\hat{p})=\valorabs{A_1(R)}$ and $A(R)=A_1(R)\cup A_2(R)$, then
	$$\hat{N}(R)=N(R,\hat{p},\hat{p})+\valorabs{A_2(R)}.$$
We will prove that 
	$$\hat{N}(R)\leq \valorabs{\bola{\Gamma}{K}{\id}}$$
for some $K\in\Z$, and therefore 
	$$A(R)\subset \bola{\Gamma}{K}{\id}.$$
Let $\gamma\in A(R)$, we have two possibilities:
	\begin{enumerate}[(i)]
	\item $\gamma\in A_1(R)$, i.e. $\rho(\hat{p},\gamma\hat{p})<R$. We have $N_0=N(R,\hat{p},\hat{p})$ elements of $\Gamma$ having this property. We seek the least $n\in\Z$ such that 
		$$N_0\leq\valorabs{B_n^\Gamma(\id)}.$$
	Namely, 
		$$N_0\leq 1+\frac{4}{3}(7^n-1).$$
	Manipulating this inequality we have 
		$$n\geq \log_7\parentesis{\frac{1}{4}(3N_0+1)}$$
	and then, the integer we are looking for is
		$$n=\big\lceil\log_7\parentesis{\frac{1}{4}(3N_0+1)}\big\rceil.$$
	Since $\Gamma$ is a free group and because it has a \emph{kissing-Schottky type} dynamics (see Chapter 6 of \cite{mumford2002indra}), $\bola{\D}{R}{\hat{p}}$ intersects the points $\gamma(\hat{p})$ in the order imposed by the \emph{layers} $S_n^\Gamma(\id)$. This is, it first intersects all the elements of the orbit $\Gamma\hat{p}$ given by $S^\Gamma_{0}(\id)$, then $S^\Gamma_{1}(\id)$, then $S^\Gamma_{2}(\id)$ y so on. This means that 
		$$\gamma\in \bola{\Gamma}{n}{\id},$$
	where $n=\big\lceil\log_7\parentesis{\frac{1}{4}(3N_0+1)}\big\rceil$.	
	\item $\gamma\in A_2(R)$, then $\gamma(\hat{p})\nin\bola{\D}{R}{\hat{p}}$. However, denoting again $N_0=N(R,\hat{p},\hat{p})$, the ball $\bola{\D}{R}{\hat{p}}$ already intersected $N_0$ traslations of $D$, all of which are determined by elements of $\Gamma$ of length at most $n=\big\lceil\log_7\parentesis{\frac{1}{4}(3N_0+1)}\big\rceil$. Then,
		$$\gamma\in \text{B}^{\Gamma}_{n+1}(\id).$$	 
	\end{enumerate}
	All of the above implies that 
		$$A(R)\subset \text{B}^{\Gamma}_{n+1}(\id).$$
	and then,
		\begin{equation}\label{eq_cota_superior_NR_1}
		\hat{N}(R)=\valorabs{A(R)}\leq \valorabs{\text{B}^{\Gamma}_{n+1}(\id)}=1+\frac{4}{3}(7^{\big\lceil\log_7\parentesis{\frac{1}{4}(3N_0+1)}\big\rceil+1}-1).
		\end{equation}
	This is the desired upper bound.\\
	
	Observe that $\lceil b \rceil < b+1$ for any $b\in\R$, using this inequality in (\ref{eq_cota_superior_NR_1}) we have 	
	\begin{equation}\label{eq_cota_superior_NR_2}
	\hat{N}(R)< 49N(R,\hat{p},\hat{p})+16.
	\end{equation}					
		
We will use the following theorem (see Theorem 1.5.1 of \cite{nicholls}).
		
\begin{thm}\label{thm_cota_sup_conteo_orbital}
Let $\Gamma$ be a discrete group acting on $\D$. There is a constant $A$, depending on $\Gamma$ and $y\in\D$, such that 
	$$N(r,x,y)< Ae^r$$
for any $x\in\D$ and $r>0$. 
\end{thm}		

Using the previous Theorem \ref{thm_cota_sup_conteo_orbital} in (\ref{eq_cota_superior_NR_2}) yields

	\begin{equation}\label{eq_cota_superior_NR_3}
	\hat{N}(R)\leq 49Ae^R +16. 
	\end{equation}

where $A$ is a constant depending on the group $\Gamma$ and $\hat{p}$. Specifically, $A$ has the form
	$$A=\frac{e^\varepsilon}{4Vol\parentesis{\bola{\D}{\varepsilon}{\hat{p}}}},$$
where $\varepsilon>0$ is such that the images of $\bola{\D}{\varepsilon}{\hat{p}}$ under $\Gamma$ are pairwise disjoint. Therefore, it is enough that 
	$$\varepsilon < \rho\parentesis{\hat{p},\partial D}.$$ 
Observe that $\varepsilon$ depends on the group $\Gamma$ and $\hat{p}$, namely, on $\eta$ and $\Upsilon$. As long as there is no ambiguity, we will omit this dependency and write $\varepsilon$ instead of $\varepsilon_{\eta,\Upsilon}$.\\

Finally, by (\ref{eq_cota_superior_NR_3}), we get
	\begin{equation}\label{eq_cota_superior_NR_4}
	\hat{N}(R)\leq \frac{49}{4} \frac{e^\varepsilon}{Vol\parentesis{\bola{\D}{\varepsilon}{\hat{p}}}} e^R +16. 
	\end{equation}

Combining (\ref{eq_cota_superior_NR_4}), (\ref{eq_forma_cota_volumen_rebanada}) and Proposition \ref{prop_volumen_finito} it follows
	
	\begin{equation}\label{forma_final_volumen_rebanada_1}
	\text{Vol}\parentesis{\bola{S}{R}{p}}\leq 6\pi \parentesis{\frac{49}{4} \frac{e^\varepsilon}{Vol\parentesis{\bola{\D}{\varepsilon}{\hat{p}}}} e^R +16}.	
	\end{equation}
	
Finally, as a immediate consequence of Lemma \ref{lem_esfera_isometrica_linea}, we have $\text{Vol}\parentesis{\hat{L}_{R,\eta}}=\text{Vol}\parentesis{\bola{S}{R}{p}}$. Therefore,

	\begin{equation}\label{forma_final_volumen_rebanada_2}
	\text{Vol}\parentesis{\hat{L}_{R,\eta}}\leq 6\pi \parentesis{\frac{49}{4} \frac{e^\varepsilon}{Vol\parentesis{\bola{\D}{\varepsilon}{\hat{p}}}} e^R +16}.	
	\end{equation}

This estimate bounds the Kobayashi volume of each \emph{diameter} of $\bola{\Omega}{R}{z}$.

\subsection*{Integrate the volume of $\hat{L}_{R,\eta}$ over $\eta$}
	
Now that we have an estimate of the volume of a \emph{diameter} of a Kobayashi ball, we need to find a way to integrate this estimate over all the \emph{diameters} passing through the center $z$. This task would be feasible if the Kobayashi volume in $\C^2$ could be decomposed into its two components. However, this is still unknown (see Problem B.2 of \cite{kobayashi1976intrinsic}), although there is evidence that it cannot be done. In Corollary 3.2 of \cite{graham1985some} it is proven that this does happen for the Eisenman instrinsic metric.

\section{Construction of the measure}\label{sec_construction}

In this final section we describe the construction of the desired measure up until the point where we need to guarantee the finiteness of the entropy volume of the the Kobayashi metric in complements of arrays of 5 complex projective lines in general position in $\CP^2$.\\

The next proposition says that there cannot be invariant probability measures and therefore, we must restrict our attention to quasi-invariant measures instead.\\

\begin{prop}\label{prop_quasi_invariant_measure}
Let $\Gamma\subset\PSL$ be a complex Kleinian group and suppose there exists a $\Gamma$-invariant probability measure supported on $\CG(\Gamma)$, then $\Gamma$ is elemental.
\end{prop}

\begin{proof}
Suppose that $\Gamma$ is not elemental, as a consequence of Observation \ref{obs_cg_noelemental} we can take an attractive fixed point $p\in\CG(\Gamma)$ of some loxodromic element $\gamma\in\Gamma$.\\
Since $\gamma$ is loxodromic, there is an open ball $R\subset\CP^2$ such that $p\in R$ and 
	$$\gamma(\overline{R})\subset R.$$
Let us denote $R_0=R\setminus \gamma(\overline{R})$ and $R_n=\gamma^{n}(R_0)$ for $n\in\Z$, then
	\begin{enumerate}[(i)]
	\item $R_n\cap R_m = \emptyset$ for any $n,m\in\Z$, $n\neq m$.
	\item $\underset{n\in\Z}{\bigcup}\overline{R}_n = \CP^2$.
	\end{enumerate}
Assume that there is a $\Gamma$-invariant probability measure $\mu$ supported on $\CG(\Gamma)$. Since $p$ is an attractive fixed point there are infinitely many integers $n$ such that $R_n\cap\CG(\Gamma)\neq \emptyset$ and therefore $\mu(R_n)>0$. Let $A\subset\Z$ be the infinite subset of integers such that $R_n\cap\CG(\Gamma)\neq\emptyset$. Then 
	$$\mu(R_n)=\mu(R_m)=k$$
for any $n,m\in A$. Using the countable additivity of $\mu$, it follows
	$$\sum_{n\in A}k=\sum_{n\in A}\mu(R_n)=\sum_{n\in\Z}\mu(R_n)=\mu\parentesis{\underset{n\in\Z}{\bigcup}R_n}\leq \mu(\CP^2)=\mu(\CG(\Gamma))=1.$$ 
This contradiction proves the proposition.
\end{proof}

Consider $\Gamma\subset\PSL$ a strongly irreducible complex Kleinian group. By Theorem \ref{thm_complemento_kob_hyp}, $\KulD(\Gamma)$ is complete Kobayashi hyperbolic. For any $z,w\in\KulD(\Gamma)$ and $s\in\R$ we define the series

	\begin{equation}\label{eq_defn_serie}
	f_s(z,w)=\sum_{\gamma\in\Gamma}e^{-s\rho(z,\gamma w)}
	\end{equation}		
	
where $\rho$ is the Kobayashi metric (or the Eisenman metric) on $\KulD(\Gamma)$. We now prove some necessary lemmas.

\begin{lem}\label{lem_orden_tiende_inf}
Let $\Gamma\subset\PSL$ be a complex Kleinian group and $z,w\in\KulD(\Gamma)$, then the elements of $\Gamma$ can be arranged as $\set{\gamma_1,\gamma_2,...}$ such that $\rho(z,\gamma_n w)\rightarrow \infty$ when $n\rightarrow\infty$.
\end{lem}

\begin{proof}
Let us suppose that there is no such ordering, then for any arbitrary enumeration $\Gamma=\set{\gamma_1,\gamma_2,...}$, it does not hold that 
	$$\lim_{n\rightarrow\infty}\rho(z,\gamma_n w)=\infty.$$
This means that there exists $R>0$ such that an infinite number of elements $\gamma$ of the $\Gamma$-orbit of $w$ satisfy $\gamma w\in \overline{B(z,R)}$. Let us write $K=\overline{B(z,R)}$, let 
	$$A=\SET{\gamma\in\Gamma}{\gamma w \in K},$$
since $A$ is infinite and $K$ is compact, there is a sequence $\set{\phi_n w}\subset K$ converging to a point $q\in K$, with $\phi_n\in A$. This means that $q$ is a point of accumulation of the orbit $\Gamma w$. Observe that $w\nin L_0(\Gamma)$ since $w\in\KulD(\Gamma)$, then $q\in L_1(\Gamma)$ but $q\in K\subset \KulD(\Gamma)$, this contradiction implies that there exists an enumeration $\Gamma=\set{\gamma_1,\gamma_2,...}$ such that 
	$$\rho(z,\gamma_n w)\rightarrow \infty$$ 
when $n\rightarrow\infty$.      
\end{proof}

As a consequence of this lemma we can arrange the terms
	$$\set{\rho(z,\gamma w)}_{\gamma\in\Gamma}$$ 
in an increasing order form such that they tend to $\infty$, thus the series (\ref{eq_defn_serie}) is a general Dirichlet series.\\

According to Chapter 8 \cite{apostol}, if there exists $s>0$ such that $f_s(z,w)$ converges and, for any $t<0$, $f_t(z,w)$ diverges then there is a convergence exponent $\delta_{z,w}\in\R$ such that $f_s(z,w)$ converges for $s> \delta_{z,w}$ and diverges for $s< \delta_{z,w}$. Furthermore, this exponent of convergence is given by
	
	\begin{equation}\label{eq_defn_deltazw}
	\delta_{z,w}=\limsup_{n\rightarrow\infty}\frac{\log n}{\rho(z,\gamma_n w)},
	\end{equation}

as a consequence of this, it holds that $\delta_{z,w}>0$. \\

Let $s<0$, since the terms $\rho(z,\gamma_n w)$ are arranged in increasing order, the sequence $e^{-s\rho(z,\gamma_n w)}$ does not converge to $0$ and therefore the series $f_s(z,w)$ diverges.

\begin{prop}
If it exists, the exponent of convergence $\delta_{z,w}$ depends only on $\Gamma$. 
\end{prop}

\begin{proof}
Let $z,w\in\KulD(\Gamma)$ and $\gamma\in\Gamma$. Using the triangle inequality it holds
	\begin{align*}
	\rho(z,\gamma(w)) &\leq \rho(z,w)+ \rho(w,\gamma(w)),\\
	\rho(w,\gamma(w)) &\leq \rho(z,w)+ \rho(z,\gamma(w)).
	\end{align*}
Then
	$$\rho(w,\gamma(w))-\rho(z,w)\leq \rho(z,\gamma(w)) \leq \rho(z,w)+ \rho(w,\gamma(w)) $$
and then, for $s>\delta_{z,w}$,

	\begin{align*}
	-\rho(z,w)-\rho(w,\gamma(w)) &\leq -\rho(z,\gamma(w)) \leq \rho(z,w)- \rho(w,\gamma(w))\\
	e^{-s\rho(z,w)}e^{-s\rho(w,\gamma(w))} &\leq e^{-s\rho(z,\gamma(w))} \leq e^{s\rho(z,w)}e^{- s\rho(w,\gamma(w))}\\
	e^{-s\rho(z,w)}\sum_{\gamma\in\Gamma}e^{-s\rho(w,\gamma(w))} &\leq \sum_{\gamma\in\Gamma}e^{-s\rho(z,\gamma(w))} \leq e^{s\rho(z,w)}\sum_{\gamma\in\Gamma}e^{- s\rho(w,\gamma(w))}\\
	e^{-s\rho(z,w)}f_s(w,w) &\leq f_s(z,w) \leq e^{s\rho(z,w)}f_s(w,w).	
	\end{align*}
This means that $f_s(z,w)$ converges if and only if $f_s(w,w)$ converges. Therefore $\delta_{z,w}$ does not depends of $z,w$; it only depends of $\Gamma$.
\end{proof}

As a consequence of the previous proposition, we can write $\delta_{\Gamma}$ instead of $\delta_{z,w}$. Equivalently, we can write $\delta_{\Gamma}$ in (\ref{eq_defn_deltazw}) as 
	\begin{equation}\label{eq_defn_alternativa_deltazw}
	\delta_{\Gamma}=\limsup_{n\rightarrow\infty}\frac{1}{r}\log N(r,z,w).
	\end{equation}

\begin{defn}
We say that $\Gamma$ is of \emph{divergence type} if $f_\delta(z,w)$ diverges for any $z,w\in\KulD(\Gamma)$.
\end{defn}

In the following theorem we give a generalization of Theorem \ref{thm_cota_sup_conteo_orbital}, assuming that the entropy volume for the Kobayashi metric is finite for some particular point $z\in\KulD(\Gamma)$.

\begin{thm}\label{thm_conteo_orbital_acotado_dim2}
Let $\Gamma\subset\PSL$ be a strongly irreducible complex Kleinian group acting on $\Omega=\KulD(\Gamma)$ and let $w\in\Omega$. If the entropy volume of $d_{\Omega}$ sa-tisfies $e(\Omega,z)<\infty$ for some $z\in\Omega$, then there is a constant $K$, depending on the metric $d_{\Omega}$ and $z$, and a constant $A$, depending on $\Gamma$, $K$ and $w$, such that
	$$N(r,z,w)<Ae^{rK},$$
for any $r>0$.  
\end{thm}

\begin{proof}
Let $\varepsilon>0$ small enough such that no two $\Gamma$-images of $B=\bola{\Omega}{\varepsilon}{w}$ overlap. Let $r>0$, then
	\begin{equation}\label{eq_thm_conteo_orbital_acotado_dim2_1}		
	\text{Vol}(B)N(r,w,z)<\text{Vol}\parentesis{\bola{\Omega}{r+\varepsilon}{z}}.
	\end{equation}

\begin{figure}[H]
\begin{center}	
\includegraphics[height=45mm]{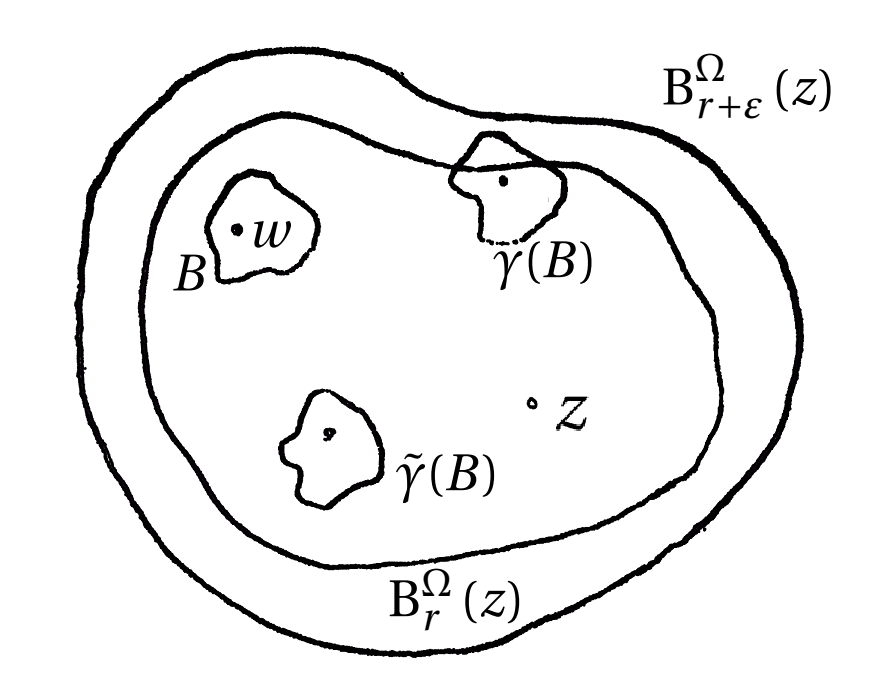}	
\end{center}	
\end{figure}

Since we are assuming that $e(\Omega,z)<\infty$, by definition there is a constant $K>0$ depending on the metric $d_{\Omega}$ such that 
	\begin{equation}\label{eq_thm_conteo_orbital_acotado_dim2_2}
	\text{Vol}\parentesis{\bola{\Omega}{R}{z}}<e^{KR}
	\end{equation}
for any $R>0$ and $z\in\Omega$. Then, using (\ref{eq_thm_conteo_orbital_acotado_dim2_2}) in (\ref{eq_thm_conteo_orbital_acotado_dim2_1}) yields
	$$N(r,w,z)<\frac{e^{K(r+\varepsilon)}}{\text{Vol}(B)}=\frac{e^{K\varepsilon}}{\text{Vol}(B)}e^{Kr}.$$
\end{proof}

Using Theorem \ref{thm_conteo_orbital_acotado_dim2} on (\ref{eq_defn_alternativa_deltazw}) we can conclude that, if the entropy volume of $d_{\Omega}$ is finite, then
	$$\delta_{\Gamma}<\infty$$
and then we can construct the Patterson-Sullivan measures in the same way as in Chapter 3 of \cite{nicholls}, we will do this in Theorem \ref{thm_existencia_PS}. Before we state and prove this theorem, we need some previous results, we will only give the proofs if they are essentially different from their respective counterparts given in \cite{nicholls}.\\

The proof of the following lemma is similar to the proof given in p. 46 of \cite{nicholls}.

\begin{lem}
Let $z,w\in\KulD(\Gamma)$. If $s>\delta_{\Gamma}$ then $f_s(z,w)$ converges, if $s>\delta_{\Gamma}$ then $f_s(z,w)$ diverges.
\end{lem}

Using Lemma \ref{lem_orden_tiende_inf}, the proof of the following lemma is similar to the proof of Lemma 3.1.1 of \cite{nicholls}.

\begin{lem}\label{lem_existencia_funcion_h}
Let $\Gamma\subset\PSL$ a complex Kleinian group with critical exponent $\delta=\delta_\Gamma$. There exists a continuous non-decreasing function $h:\R^+\rightarrow\R^+$ such that the series
		$$\sum_{\gamma\in\Gamma}e^{-s\rho(z,\gamma w)}h\parentesis{e^{\rho(z,\gamma w)}}$$
converges for $s>\delta$ and diverges for $s\leq\delta$.
\end{lem} 

The previous lemma gives a tool to construct a modified Poincar\'e series for the group $\Gamma$. This modified series has the same convergence exponent as the original series $f_s(z,w)$, however, if the group $\Gamma$ was of convergence type with $f_s(z,w)$, now it is of divergence type.\\

Ttheorem \ref{thm_helly} is known as the Helly's theorem, it is a compacteness theorem for the space of functions with bounded total variation (see p. 222 of \cite{natanson2016theory}, Theorem 1.126 of \cite{barbu2012convexity} or Theorem 19 of \cite{kreuzer2013bounded}). Before, we need to define the variation of a measure (see Section 1.6 of \cite{saks1937theory}).\\

Let $\mu$ be a measure on a measurable space $(X,\Sigma)$. For $E\in\Sigma$, we define
	\begin{align*}
	\overline{W}(\mu,E) &= \sup\SET{\mu(A)}{A\in\Sigma,\;A\subset E},\\
	\underline{W}(\mu,E) &= \inf\SET{\mu(A)}{A\in\Sigma,\;A\subset E} 
	\end{align*}
The \emph{variation of the measure} is given by
	$$v(\mu,E)=\overline{W}(\mu,E)+\valorabs{\underline{W}(\mu,E)},$$
and the \emph{total variation of measure} is given by
	$$V(\mu)=v(\mu,X).$$

\begin{thm}\label{thm_helly}
Let $U\subset \R$ be an open set and let $f_n:U\rightarrow\R$ be a sequence of functions. If $\set{f_n}$ has uniformly bounded variation on any compact set that is compactly embedded in $U$ and it is uniformly bounded at a point of $U$, then there exists a subsequence $\set{f_{n_k}}$ which converges pointwise to a function $f:U\rightarrow \R$ of bounded variation.
\end{thm}

In the next lemma we verify that the family of measures $\set{\mu_{z,w,s}}_{s>\delta_\Gamma}$ defined in (\ref{eq_defn_measure_s}) will satisfy the hypothesis of Theorem \ref{thm_helly}. We will use this lemma in the proof of Theorem \ref{thm_existencia_PS}, the proof is similar to the proof of Lemma 3.3.1 of \cite{nicholls}.

\begin{lem}\label{lem_measures_bounded}
Let $\Gamma\subset\PSL$ a complex Kleinian group, with critical exponent $\delta=\delta_\Gamma$. Denote $\Omega=\KulD(\Gamma)$, then $\mu_{x,s}(\Omega)$ is uniformly bounded for $s\in(\delta,\delta+1)$.
\end{lem} 

\begin{thm}\label{thm_existencia_PS}
Let $\Gamma\subset\PSL$ be a strongly irreducible complex Kleinian group acting on $\Omega=\KulD(\Gamma)$ and let $z\in\Omega$ such that $e(\Omega,z)<\infty$. Then, there exists a family of quasi-invariant measures supported on the Kulkarni limit set of $\Gamma$. 
\end{thm}

\begin{proof}
Let $z,w\in\KulD(\Gamma)$, $\delta=\delta_\Gamma$ and $s>\delta$. If $\Gamma$ is a group of convergence type, let $h$ be the function given by Lemma \ref{lem_existencia_funcion_h}, if $\Gamma$ is of divergence type, we take $h\equiv 1$. Let us consider,
	$$f^{\ast}_s(z,w)=\sum_{\gamma\in\Gamma}e^{-s\rho(z,\gamma w)}h\parentesis{e^{\rho(z,\gamma w)}}.$$
We define
	$$\mu_{z,w,s}=\frac{1}{f^{\ast}_s(w,w)}\sum_{\gamma\in\Gamma}e^{-s\rho(z,\gamma w)}h\parentesis{e^{\rho(z,\gamma w)}}D(\gamma w),$$
with $D(\gamma w)$ denoting the Dirac point mass of weight one at $\gamma w$. For a Borel set $E\subset\Omega$ we define the measure
	\begin{equation}\label{eq_defn_measure_s}
	\mu_{z,w,s}(E)=\frac{1}{f^{\ast}_s(w,w)}\sum_{\gamma\in\Gamma}e^{-s\rho(z,\gamma w)}h\parentesis{e^{\rho(z,\gamma w)}}1_E(\gamma w),
	\end{equation}		
where $1_E$ is the characteristic function of $E$. It is straight-forward to verify that $\mu_{z,w,s}$ is a measure of probability.
\begin{figure}[H]
\begin{center}	
\includegraphics[height=40mm]{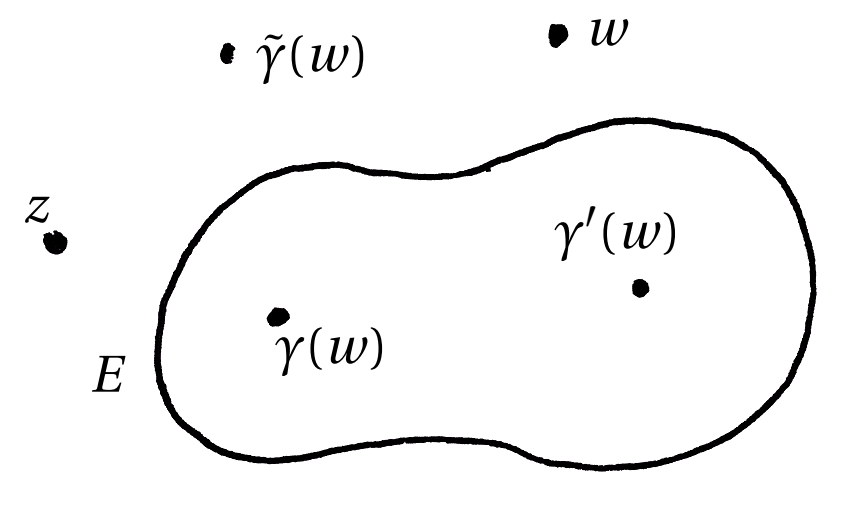}	
\end{center}	
\end{figure}

We will not be changing the point $w$, so we will just write $\mu_{z,s}$. Let $\set{s_n}\subset(\delta,\delta+1)$ be a decreasing sequence, using Lemma \ref{lem_measures_bounded} we conclude that the family of measures $\set{\mu_{z,s_n}}$ has uniformly bounded variation and then, by Helly's theorem (Theorem \ref{thm_helly}) there is a subsequence of measures $\set{\mu_{z,s_{n_k}}}$ converging to a measure $\mu_z$. This measure depends on the choice of the sequence $\set{s_n}$. We denote this family of measures by $M_z$.\\

Now we verify that this family of measures $M_z$ are supported on the Kulkarni limit set. In order to do this we choose a measure $\mu\in M_z$ and a ball $B$ around some $\gamma w$ such that $B$ does not contain any other element of the orbit $\Gamma w$, this is possible because the accumulation points of the orbit are in $\KulL(\Gamma)$. Then we need to prove that $\mu(B)=0$. This is done in the same way as in Theorem 3.3.2 of \cite{nicholls}.
\end{proof}

The family of measures $M_z$ obtained in the proof of the previous theorem, for a point $z\in\KulD(\Gamma)$ such that $e(\KulD(\Gamma),z)<\infty$, are measures with similar properties to the Patterson-Sullivan measures. However, it is important to note again that this generalization is possible only if one is able to show that the entropy volume of $d_{\Omega}$ is finite.

\bibliographystyle{alpha}
\bibliography{Bibliografia}

\end{document}